\documentclass[a4paper,12pt,amsart,frenchb]{article}

\usepackage{amsmath,amsbsy,amsfonts,amssymb}
\usepackage[french]{babel}
\newcommand{\ass}[2]{\vskip0.3cm\noindent
{\bf {#1}}. { \sl {#2}}\vskip0.3cm\noindent
}
  
\begin{document}

    \title{ Stabilisation de la formule des traces tordue III: int\'egrales orbitales et endoscopie sur un corps local non-archim\'edien; r\'eductions et preuves}
\author{ J.-L. Waldspurger}
\date{7 f\'evrier  2014}
\maketitle

\bigskip

{\bf Introduction}

Cet article est la suite de [I] et  [II] dont on adopte les notations. Soient $F$ un corps local non archim\'edien de caract\'eristique nulle, $G$ un groupe r\'eductif connexe d\'efini sur $F$, $\tilde{G}$ un espace tordu sous $G$ et ${\bf a}$ un \'el\'ement de $H^1(W_{F}; Z(\hat{G}))$, dont se d\'eduit un caract\`ere $\omega$ de $G(F)$. On impose les hypoth\`eses de [II] 1.1. Soit $\tilde{M}$ un espace de Levi de $\tilde{G}$. On a \'enonc\'e dans [II] diverses assertions concernant les int\'egrales pond\'er\'ees $\omega$-\'equivariantes $I_{\tilde{M}}^{\tilde{G}}(\boldsymbol{\gamma},{\bf f})$ et leurs variantes endoscopiques $I_{\tilde{M}}^{\tilde{G},{\cal E}}(\boldsymbol{\gamma},{\bf f})$. Le terme $\boldsymbol{\gamma}$ est un \'el\'ement de $D_{g\acute{e}om}(\tilde{M}(F),\omega)\otimes Mes(M(F))^*$, ce qui revient essentiellement \`a dire que c'est une combinaison lin\'eaire finie d'int\'egrales orbitales dans $\tilde{M}(F)$ (tordues par $\omega$). Le terme ${\bf f}$ est un \'el\'ement de $I(\tilde{G}(F),\omega)\otimes Mes(G(F))$, ce qui revient essentiellement \`a dire que c'est une fonction sur $\tilde{G}(F)$, localement constante et \`a support compact. Consid\'erons l'hypoth\`ese suivante:

{\bf (Hyp)}{\it soient $\boldsymbol{\gamma}\in D_{g\acute{e}om}(\tilde{M}(F),\omega)\otimes Mes(M(F))^*$ et ${\bf f}\in I(\tilde{G}(F),\omega)\otimes Mes(G(F))$; supposons que $\boldsymbol{\gamma}$ soit \`a support fortement r\'egulier dans $\tilde{G}(F)$; alors on a l'\'egalit\'e
$$I_{\tilde{M}}^{\tilde{G},{\cal E}}(\boldsymbol{\gamma},{\bf f})=I_{\tilde{M}}^{\tilde{G}}(\boldsymbol{\gamma},{\bf f}).$$}

Le but de l'article est de prouver toutes les assertions \'enonc\'ees en [II] sous cette hypoth\`ese. La preuve de celle-ci n\'ecessite un argument global et sera faite plus tard. Esquissons comment se d\'eduit le th\'eor\`eme [II] 1.16 de l'hypoth\`ese (Hyp). L'\'enonc\'e de ce th\'eor\`eme est le m\^eme que celui de (Hyp), sauf que l'on supprime l'hypoth\`ese sur le support de $\boldsymbol{\gamma}$. En utilisant la th\'eorie des germes de Shalika, on a prouv\'e en [II] 2.10 que (Hyp) entra\^{\i}nait une assertion plus forte, \`a savoir la m\^eme \'egalit\'e sous l'hypoth\`ese plus faible que le support de $\boldsymbol{\gamma}$ est form\'e d'\'el\'ements $\tilde{G}$-\'equisinguliers, c'est-\`a-dire d'\'el\'ements $\gamma\in \tilde{M}(F)$ tels que $M_{\gamma}=G_{\gamma}$. Passons au cas o\`u $\boldsymbol{\gamma}$ est quelconque. 
On peut fixer une classe de conjugaison stable ${\cal O}$ d'\'el\'ements semi-simples de $\tilde{M}(F)$ et supposer $\boldsymbol{\gamma}\in D_{g\acute{e}om}({\cal O},\omega)\otimes Mes(M(F))^*$, c'est-\`a-dire que le support de $\boldsymbol{\gamma}$ est form\'e d'\'el\'ements de parties semi-simples dans ${\cal O}$. Soit $a\in A_{\tilde{M}}(F)$ en position g\'en\'erale. On note $a\boldsymbol{\gamma}$ la translat\'ee de $\boldsymbol{\gamma}$ par $a$. Le support de cette distribution  est form\'e d'\'el\'ements $\tilde{G}$-\'equisinguliers. On a donc pour tout ${\bf f}$ l'\'egalit\'e
$$(1) \qquad I_{\tilde{M}}^{\tilde{G},{\cal E}}(a\boldsymbol{\gamma},{\bf f})=I_{\tilde{M}}^{\tilde{G}}(a\boldsymbol{\gamma},{\bf f}).$$
Faisons tendre $a$ vers $1$. On a \'etabli un d\'eveloppement de ces deux termes en [II] 3.2 et [II] 3.9. Pour les \'enoncer facilement, faisons l'hypoth\`ese simplificatrice que $\tilde{M}$ est un espace de Levi propre et maximal de $\tilde{G}$. Alors $I_{\tilde{M}}^{\tilde{G},{\cal E}}(a\boldsymbol{\gamma},{\bf f})$ est \'equivalent (en un sens d\'efini en [II] 3.1) \`a
$$(2) \qquad I_{\tilde{M}}^{\tilde{G}}(\boldsymbol{\gamma},{\bf f})+\sum_{J\in {\cal J}_{\tilde{M}}^{\tilde{G}}}I^{\tilde{G}}(\rho_{J}(\boldsymbol{\gamma},a)^{\tilde{G}},{\bf f}).$$
L'ensemble ${\cal J}_{\tilde{M}}^{\tilde{G}}$ est celui des racines de $A_{\tilde{M}}$ dans $\mathfrak{g}$, au signe pr\`es. Pour $J=\{\pm \alpha\}$, le terme $\rho_{J}(\boldsymbol{\gamma},a)$ est le produit d'un \'el\'ement de $D_{g\acute{e}om}({\cal O},\omega)\otimes Mes(M(F))^*$ ind\'ependant de $a$ et de 
$log(\vert \alpha(a)-\alpha(a)^{-1}\vert _{F})$. Enfin, l'exposant $\tilde{G}$ de $\rho_{J}(\boldsymbol{\gamma},a)^{\tilde{G}}$ indique que l'on induit cette distribution \`a $\tilde{G}(F)$. On a un d\'eveloppement parall\`ele
$$(3) \qquad I_{\tilde{M}}^{\tilde{G},{\cal E}}(\boldsymbol{\gamma},{\bf f})+\sum_{J\in {\cal J}_{\tilde{M}}^{\tilde{G}}}I^{\tilde{G}}(\rho_{J}^{{\cal E}}(\boldsymbol{\gamma},a)^{\tilde{G}},{\bf f}).$$
La notion d'\'equivalence \'evoqu\'ee plus haut est inoffensive: l'\'egalit\'e (1) entra\^{\i}ne que (2) et (3) sont \'egaux.  Pour $J\in {\cal J}_{\tilde{M}}^{\tilde{G}}$, consid\'erons l'assertion

$(4)_{J}$ on a $\rho_{J}^{{\cal E}}(\boldsymbol{\gamma},a)^{\tilde{G}}=\rho_{J}(\boldsymbol{\gamma},a)^{\tilde{G}}$ pour tout $\boldsymbol{\gamma}\in D_{g\acute{e}om}({\cal O},\omega)\otimes Mes(M(F))^*$ et tout $a\in A_{\tilde{M}}(F)$ en position g\'en\'erale et proche de $1$. 

L'ensemble ${\cal J}_{\tilde{M}}^{\tilde{G}}$ poss\`ede un \'el\'ement "maximal" $J_{max}=\{\pm \alpha\}$ form\'e des deux racines indivisibles. Supposons $(4)_{J}$ prouv\'e pour $J\not=J_{max}$. Alors l'\'egalit\'e de (2) et (3) entra\^{\i}ne
$$ I_{\tilde{M}}^{\tilde{G},{\cal E}}(\boldsymbol{\gamma},{\bf f})- I_{\tilde{M}}^{\tilde{G}}(\boldsymbol{\gamma},{\bf f})=I^{\tilde{G}}(\rho_{J_{max}}(\boldsymbol{\gamma},a)^{\tilde{G}} -\rho_{J_{max}}(\boldsymbol{\gamma},a)^{\tilde{G}},{\bf f}).$$
Or le membre de gauche est constant en $a$ tandis que celui de droite est proportionnel \`a $log(\vert \alpha(a)-\alpha(a)^{-1}\vert _{F})$, o\`u $\alpha$ est une racine indivisible. Il ne peuvent \^etre \'egaux que s'ils sont tous deux nuls. La nullit\'e du membre de gauche est l'assertion du th\'eor\`eme [II] 1.16. La nullit\'e du membre de droite est l'assertion $(4)_{J_{max}}$. Tout revient donc \`a prouver $(4)_{J}$ pour $J$ non maximal. Dans le cas o\`u $\tilde{G}=G$ et ${\bf a}=1$, on raisonne par r\'ecurrence sur la dimension de $G$. A un $J$ non maximal, on associe un certain sous-groupe $G_{J}\subsetneq G$, et on montre que l'assertion se d\'eduit de l'assertion analogue o\`u $G$ est remplac\'e par $G_{J}$. Le  cas g\'en\'eral utilise la descente. On fixe $\eta\in {\cal O}$ et on montre que les termes $\rho_{J}(\boldsymbol{\gamma},a)$ et $\rho_{J}^{{\cal E}}(\boldsymbol{\gamma},a)$ se d\'eduisent de termes analogues o\`u $\tilde{G}$ est remplac\'e par la composante neutre $G_{\eta}$ du commutant de $\eta$ dans $\tilde{G}$. Ce groupe $G_{\gamma}$ n'\'etant plus tordu, le r\'esultat pr\'ec\'edent s'applique et on peut conclure. La th\'eorie de la descente est facile pour le terme $\rho_{J}(\boldsymbol{\gamma},a)$ (du moins, elle est facile maintenant que Harish-Chandra et Arthur ont travaill\'e pour nous). C'est beaucoup plus d\'elicat pour le terme endoscopique $\rho_{J}^{{\cal E}}(\boldsymbol{\gamma},a)$ car, dans le cas tordu, m\'elanger descente et endoscopie fait appara\^{\i}tre des "triplets endoscopiques non standard", cf. 6.1. Il y a une analogue de l'assertion $(4)_{J}$ pour de tels triplets que nous d\'eduirons elle-aussi de l'hypoth\`ese (Hyp), modulo un raisonnement par r\'ecurrence assez sophistiqu\'e, cf. 6.4.

Il y a deux cas o\`u on obtient des r\'esultats non conditionnels, parce que l'on peut d\`es maintenant prouver la validit\'e de (Hyp). Le premier, d\'etaill\'e dans la section 1, est celui o\`u il n'y a pas de torsion: $\tilde{G}=G$ et ${\bf a}=1$. Dans ce cas, (Hyp) a \'et\'e prouv\'e par Arthur ([A1] local theorem 1). Le second, auquel est consacr\'e la section 2, est celui o\`u $G$ est quasi-d\'eploy\'e, $\tilde{G}$ est \`a torsion int\'erieure et ${\bf a}=1$. Dans ce cas, on montre que l'on peut plonger $\tilde{G}$ dans un groupe non tordu $H$ de sorte que l'hypoth\`ese (Hyp) pour $\tilde{G}$ se d\'eduise de la m\^eme hypoth\`ese pour ce groupe $H$. Comme on vient de le dire, cette derni\`ere a \'et\'e prouv\'ee par Arthur puisque $H$ n'est pas tordu.

Dans la  troisi\`eme section, limit\'ee au cas des groupes non tordus, on montre comment se comportent  nos objets (par exemple les termes $\rho_{J}(\boldsymbol{\gamma},a)$) par passage au rev\^etement simplement connexe du groupe d\'eriv\'e. Dans la quatri\`eme section, on \'enonce comment se comportent ces m\^emes objets par descente d'Harish-Chandra. Dans la  section 5, on reprend les constructions de [W1] qui permettent de relier descente et endoscopie. C'est l\`a qu'apparaissent les triplets endoscopiques non standard. La section 6 leur est consacr\'ee. Dans la section 7, on d\'eveloppe la d\'emonstration grossi\`erement \'evoqu\'ee ci-dessus, c'est-\`a-dire que l'on montre que la plupart des \'enonc\'es de [II] r\'esultent de l'hypoth\`ese (Hyp). La huiti\`eme et derni\`ere section concerne les \'enonc\'es restants, \`a savoir ceux concernant les germes de Shalika. Ces germes sont locaux, c'est-\`a-dire vivent au voisinage d'une classe de conjugaison stable semi-simple ${\cal O}$ fix\'ee dans $\tilde{M}(F)$. On prouve  ces r\'esultats sans recourir \`a l'hypoth\`ese (Hyp), pourvu que ${\cal O}$ n'appartienne pas \`a un ensemble au plus fini de telles classes. Ces r\'esultats n'ont pas de cons\'equence imm\'ediate mais le fait qu'ils  soient obtenus sans hypoth\`ese nous sera utile plus tard.

\bigskip

\section{Le cas des groupes non tordus}

\bigskip

\subsection{Rappel des r\'esultats d'Arthur}
Dans tout l'article, le corps de base $F$ est local non-archim\'edien de caract\'eristique nulle, sauf mention expresse du contraire.
On consid\`ere dans cette  section un triplet $(G,\tilde{G},{\bf a})$ non tordu, c'est-\`a-dire que $\tilde{G}=G$ et ${\bf a}=1$. Le triplet se r\'eduit donc \`a l'unique groupe $G$. On fixe une fonction $B$ comme en [II] 1.8. On peut affaiblir les hypoth\`eses de r\'ecurrence pos\'ees en [II] 1.1. En effet, en partant de notre groupe $G$, on ne peut faire appara\^{\i}tre par les constructions de [II] que des triplets non tordus, r\'eduits \`a un unique groupe. Les hypoth\`eses de r\'ecurrence suivantes sont donc suffisantes: si $G$ est quasi-d\'eploy\'e, on suppose connus tous les r\'esultats concernant des groupes  $G'$ quasi-d\'eploy\'es tels que $dim(G'_{SC})< dim(G_{SC})$; si $G$ n'est pas quasi-d\'eploy\'e, on suppose connus tous les r\'esultats concernant les groupes quasi-d\'eploy\'es $G'$ tels que $dim(G'_{SC})\leq dim(G_{SC})$ et tous les r\'esultats concernant les groupes quelconques tels que $dim(G'_{SC})< dim(G_{SC})$. Si une assertion est relative \`a un Levi $M$ de $G$, on suppose connues toutes les assertions concernant le m\^eme groupe $G$ et relatives \`a un Levi $L\in {\cal L}(M)$ tel que $L\not=M$.

Soit $M$ un Levi de $G$. Dans ce cas, les r\'esultats suivants ont \'et\'e prouv\'es par Arthur ([A1]  local theorem 1):

(1) soit $\boldsymbol{\gamma}\in D_{g\acute{e}om}(M(F))\otimes Mes(M(F))^*$ et ${\bf f}\in I(G(F))\otimes Mes(G(F))$; on suppose que le support de $\boldsymbol{\gamma}$ est form\'e d'\'el\'ements fortement $G$-r\'eguliers; alors on a l'\'egalit\'e
$$I^{G,{\cal E}}_{M}(\boldsymbol{\gamma},{\bf f})=I^G_{M}(\boldsymbol{\gamma},{\bf f});$$

(2) supposons $G$ quasi-d\'eploy\'e; soit $\boldsymbol{\delta}\in D_{g\acute{e}om}^{st}(M(F))\otimes Mes(M(F))^*$; on suppose que le support de $\boldsymbol{\delta}$ est form\'e d'\'el\'ements fortement $G$-r\'eguliers; alors la distribution ${\bf f}\mapsto S_{M}^G(\boldsymbol{\delta},{\bf f})$ est stable. 

En vertu de la proposition 2.10 de [II], les m\^emes r\'esultats valent sous des hypoth\`eses plus faibles concernant les supports des \'el\'ements $\boldsymbol{\gamma}$ ou $\boldsymbol{\delta}$: on peut y remplacer "$G$-r\'eguliers" par "$G$-\'equisinguliers".

\bigskip

\subsection{Int\'egrales orbitales pond\'er\'ees stables}
On suppose $G$ quasi-d\'eploy\'e. Soit $M$ un Levi de $G$. On a d\'efini en [II] 3.3 l'ensemble ${\cal J}_{M}^G(B)$. Pour $J\in {\cal J}_{M}^G(B)$, on peut appliquer la construction [II] 3.5 \`a la classe ${\cal O}=\{1\}$ et au syst\`eme de fonctions d\'eduit de notre fonction $B$. Cela d\'efinit une application lin\'eaire 
$$\sigma^G_{J}:D^{st}_{unip}(M(F))\otimes Mes(M(F))^*\to U_{J}\otimes (D_{unip}(M(F))\otimes Mes(M(F))^*)/Ann_{unip}^{G}.$$
On a aussi d\'efini en [II] 3.3 un groupe $G_{J}$ qui n'est pas, en g\'en\'eral, un sous-groupe de $G$. Le syst\`eme de racines de $G_{J}$ est un sous-syst\`eme de celui  de $G$, de m\^eme rang que celui-ci et l'inclusion est \'equivariante pour les actions galoisiennes. Il en r\'esulte que $Z(\hat{G})^{\Gamma_{F}}$ est un sous-groupe d'indice fini de $Z(\hat{G}_{J})^{\Gamma_{F}}$. On pose
$$i_{J}^G=[Z(\hat{G}_{J})^{\Gamma_{F}}:Z(\hat{G})^{\Gamma_{F}}]^{-1}.$$
Remarquons que, dans le cas o\`u $J$ est maximal, cf. [II] 3.1, on a $G_{J}=G$ et $i_{J}^G=1$. En particulier, si l'on remplace $G$ par $G_{J}$, $J$ devient l'\'el\'ement maximal de ${\cal J}_{M}^{G_{J}}(B)$ et $i_{J}^{G_{J}}=1$. 
On a prouv\'e en [II] 3.3 l'inclusion $Ann_{unip}^{G_{J}}\subset Ann_{unip}^{G}$.

\ass{Proposition}{(i) Pour tout $J\in {\cal J}_{M}^G(B)$, $ \sigma_{J}^G$ est la compos\'ee de $i_{J}^G\sigma_{J}^{G_{J}}$ et de la projection
$$U_{J}\otimes (D_{unip}(M(F))\otimes Mes(M(F))^*)/Ann_{unip}^{G_{J}}\to U_{J}\otimes (D_{unip}(M(F))\otimes Mes(M(F))^*)/Ann_{unip}^{G}.$$

(ii) Pour tout $J\in {\cal J}_{M}^G(B)$, $\sigma_{J}^G$ prend ses valeurs dans 
$$U_{J}\otimes (D^{st}_{unip}(M(F))\otimes Mes(M(F))^*)/Ann_{unip}^{G,st}.$$

(iii) Pour tout $\boldsymbol{\delta}\in D_{unip}^{st}(M(F))\otimes Mes(M(F))^*$, la distribution
$${\bf f}\mapsto S_{M}^G(\boldsymbol{\delta},B,{\bf f})$$
est stable.}

Preuve.  Rappelons que l'on note $\rho_{J,st}^G$ la restriction de $\rho^G_{J}$ \`a l'espace $D^{st}_{unip}(M(F))\otimes Mes(M(F))^*$. On peut reformuler la d\'efinition de $\sigma_{J}^G$ par l'\'egalit\'e
$$(1) \qquad \rho^G_{J,st}=\sum_{s\in Z(\hat{M})^{\Gamma_{F}}/Z(\hat{G})^{\Gamma_{F}}, J\in {\cal J}_{M}^{G'(s)}(B)}i_{M}(G,G'(s))\sigma_{J}^{G'(s)}.$$

Plus exactement, il s'agit des projections modulo $Ann_{unip}^G$ des $\sigma_{J}^{G'(s)}$. Pour simplifier, on oublie de telles projections dans la notation. Le membre de gauche est \'egal \`a $\rho^{G_{J}}_{J,st}$ d'apr\`es [II] 3.3(i) (en oubliant les projections). Pour $s\not=1$, on peut utiliser  par r\'ecurrence le (i) de l'\'enonc\'e: on a $\sigma_{J}^{G'(s)}=i_{J}^{G'(s)}\sigma_{J}^{G'(s)_{J}}$. Posons $x=\sigma_{J}^G-i_{J}^G\sigma_{J}^{G_{J}}$. Remarquons que le (i) de l'\'enonc\'e revient \`a prouver que $x=0$. L'\'egalit\'e (1) devient
$$(2) \qquad \rho_{J,st}^{G_{J}}=x+\sum_{s\in Z(\hat{M})^{\Gamma_{F}}/Z(\hat{G})^{\Gamma_{F}}, J\in {\cal J}_{M}^{G'(s)}(B)}i_{M}(G,G'(s))i_{J}^{G'(s)}\sigma_{J}^{G'(s)_{J}}.$$
On a une projection
$$Z(\hat{M})^{\Gamma_{F}}/Z(\hat{G})^{\Gamma_{F}}\to Z(\hat{M})^{\Gamma_{F}}/Z(\hat{G}_{J})^{\Gamma_{F}}.$$
Un \'el\'ement $s\in Z(\hat{M})^{\Gamma_{F}}/Z(\hat{G})^{\Gamma_{F}}$ d\'etermine donc \`a la fois un groupe endoscopique $G'(s)$ de $G$ et un groupe endoscopique $(G_{J})'(s)$ de $G_{J}$. 
Montrons que

(3) on a $J\in {\cal J}_{M}^{G'(s)}(B)$ si et seulement si $J\in {\cal J}_{M}^{(G_{J})'(s)}(B)$;

(4) si $J\in {\cal J}_{M}^{G'(s)}(B)$, on a l'\'egalit\'e $G'(s)_{J}=(G_{J})'(s)$.

 Rappelons que l'on a associ\'e \`a $J$ un r\'eseau  $R_{J}\subset \mathfrak{a}_{M}^*$ de rang $n=a_{M}-a_{G}$, cf. [II] 3.1. Fixons une paire de Borel $(B,T)$  de $G$ d\'efinie sur $F$ telle que $M$ soit standard. On note $\alpha_{M}$ la restriction \`a $\mathfrak{a}_{M}$ d'un \'el\'ement $\alpha\in \mathfrak{t}^*$. On a les \'egalites
  $$\Sigma^{G'(s)}(T)=\{\alpha\in \Sigma^{G}(T); \hat{\alpha}(s)=1\},$$
 $$\Sigma^{G_{J}}(T)=\{\alpha\in \Sigma^G(T); B(\alpha)^{-1}\alpha_{M}\in R_{J}\},$$
 $$\Sigma^{(G_{J})'(s)}(T)=\{\alpha\in \Sigma^G(T); B(\alpha)^{-1}\alpha_{M}\in R_{J},\hat{\alpha}(s)=1\}.$$
 On a $J\in {\cal J}_{M}^{G'(s)}(B)$ si et seulement si il existe $\alpha_{1},...,\alpha_{n}\in \Sigma^{G'(s)}(T)$ telles que les \'el\'ements $B(\alpha_{i})^{-1}\alpha_{i,M}$ pour $i=1,...,n$ engendrent $R_{J}$. Dans ce cas, ces \'el\'ements $\alpha_{i}$ appartiennent aussi \`a $\Sigma^{(G_{J})'(s)}(T)$, donc $J\in {\cal J}_{M}^{(G_{J})'(s)}(B)$. La r\'eciproque est \'evidente. Si ces conditions sont v\'erifi\'ees, les ensembles $\Sigma^{G'(s)_{J}}(T)$ et $\Sigma^{(G_{J})'(s)}(T)$ sont  \'egaux: conserver les racines $\alpha$ telles que $\hat{\alpha}(s)=1$ et conserver les racines $\alpha$ telles que $B(\alpha)^{-1}\alpha_{M}\in R_{J}$ sont des op\'erations qui commutent. Les actions galoisiennes sont aussi les m\^emes: ce sont les restrictions de l'action sur $\Sigma^G(T)$. Cela prouve (3) et (4). 
 
 On r\'ecrit (2) sous la forme
 $$\rho_{J,st}^{G_{J}}=x+\sum_{s\in Z(\hat{M})^{\Gamma_{F}}/Z(\hat{G})^{\Gamma_{F}}, J\in {\cal J}_{M}^{(G_{J})'(s)}(B)}i_{M}(G,G'(s))i_{J}^{G'(s)}\sigma_{J}^{(G_{J})'(s)}.$$
 Pour tout $s$ apparaissant, on a
 $$i_{M}(G,G'(s))i_{J}^{G'(s)}=[Z(\hat{G}'(s))^{\Gamma_{F}}:Z(\hat{G})^{\Gamma_{F}}]^{-1}[Z(\hat{G}'(s)_{J})^{\Gamma_{F}}:Z(\hat{G}'(s))^{\Gamma_{F}}]^{-1}.$$
 En utilisant (4), on obtient
 $$i_{M}(G,G'(s))i_{J}^{G'(s)}=[Z((\hat{G}_{J})'(s))^{\Gamma_{F}}:Z(\hat{G})^{\Gamma_{F}}]^{-1}$$
 $$=[Z((\hat{G}_{J})'(s))^{\Gamma_{F}}:Z(\hat{G}_{J})^{\Gamma_{F}}]^{-1}[Z(\hat{G}_{J})^{\Gamma_{F}}:Z(\hat{G})^{\Gamma_{F}}]^{-1}$$
 $$=i_{M}(G_{J},(G_{J})'(s))[Z(\hat{G}_{J})^{\Gamma_{F}}:Z(\hat{G})^{\Gamma_{F}}]^{-1}.$$
 Alors (2) se r\'ecrit
 $$\rho_{J,st}^{G_{J}}=x+ [Z(\hat{G}_{J})^{\Gamma_{F}}:Z(\hat{G})^{\Gamma_{F}}]^{-1}\sum_{s\in Z(\hat{M})^{\Gamma_{F}}/Z(\hat{G})^{\Gamma_{F}}, J\in {\cal J}_{M}^{(G_{J})'(s)}(B)}i_{M}(G_{J},(G_{J})'(s))\sigma_{J}^{(G_{J})'(s)}$$
 ou encore
$$\rho_{J,st}^{G_{J}}=x+ \sum_{s\in Z(\hat{M})^{\Gamma_{F}}/Z(\hat{G}_{J})^{\Gamma_{F}}, J\in {\cal J}_{M}^{(G_{J})'(s)}(B)}i_{M}(G_{J},(G_{J})'(s))\sigma_{J}^{(G_{J})'(s)}.$$
En comparant avec l'\'egalit\'e (1) appliqu\'ee au groupe $G_{J}$, on obtient $x=0$, ce qui prouve le (i) de l'\'enonc\'e. 
 
 En raisonnant par r\'ecurrence, ce r\'esultat implique (ii) pour tout $J$ tel que $G_{J}\not=G$. C'est-\`a-dire pour tout $J$ sauf l'unique \'el\'ement maximal $J_{max}$.
 
  Rappelons le d\'eveloppement [II] 3.7. Soient  $\boldsymbol{\delta}\in D^{st}_{unip}(M(F))\otimes Mes(M(F))^*$ et ${\bf f}\in I(G(F))\otimes Mes(G(F))$. Le germe en $1$ de la fonction $a\mapsto S_{M}^G(a\boldsymbol{\delta},{\bf f})$ est \'equivalent \`a
  $$\sum_{J\in {\cal J}_{M}^G(B)} I^G(\sigma_{J}^G(\boldsymbol{\delta},a)^{G},{\bf f})$$
  $$+\sum_{L\in {\cal L}(M),L\not=G}\sum_{J\in {\cal J}_{M}^L(B) } S_{L}^G(\sigma_{J}^L(\boldsymbol{\delta},a)^{L},B,{\bf f}).$$
   Pour $L\not=G$ et $J\in {\cal J}_{M}^L(B)$, on sait  par r\'ecurrence que $\sigma^L_{J}(\boldsymbol{\delta},a)^{L}$ est stable. Si de plus $L\not=M$, on sait  que la distribution ${\bf f}\mapsto S_{L}^G(\sigma_{J}^L(\boldsymbol{\delta},a)^{L},B,{\bf f})$ est stable.  Pour $J\in {\cal J}_{M}^G(B)$ non maximal, on sait que $\sigma_{J}^G(\boldsymbol{\delta},a)^G$ est stable.  Supposons que ${\bf f}$ soit "instable", c'est-\`a-dire que l'image de ${\bf f}$ dans $SI(G(F))\otimes Mes(G(F))$ soit nulle. Alors tous les termes du d\'eveloppement ci-dessus s'annulent, sauf deux: ceux pour l'\'el\'ement maximal $J_{max}$ de ${\cal J}_{M}^G(B)$ et pour  l'\'el\'ement $J=\emptyset$ de ${\cal J}_{M}^M(B)$ (on peut les supposer distincts, sinon $M=G$ et la proposition est tautologique). Pour $J=\emptyset$, on a  $\sigma^M_{\emptyset}(\boldsymbol{\delta},a)=\boldsymbol{\delta}$. Le d\'eveloppement ci-dessus se r\'eduit \`a
   $$I^G(\sigma^G_{J_{max}}(\boldsymbol{\delta},a)^G,{\bf f})+S_{M}^G(\boldsymbol{\delta},B,{\bf f}).$$
   Les r\'esultats d'Arthur impliquent que $S_{M}^G(a\boldsymbol{\delta},{\bf f})=0$, cf. 1.1. Donc la somme ci-dessus est \'equivalente \`a $0$. Comme fonction de $a$, le premier terme appartient \`a $U_{J_{max}}$ et le second est constant. La propri\'et\'e [II] 3.1(3) entra\^{\i}ne que les deux termes sont nuls. La nullit\'e du premier pour tout ${\bf f}$ instable signifie que $\sigma^G_{J_{max}}(\boldsymbol{\delta},a)^G$ est stable. En vertu du lemme [I] 5.13, cela  ach\`eve de prouver (ii). La nullit\'e du second terme implique le (iii) de l'\'enonc\'e. $\square$
   
   \bigskip
   
   \subsection{Germes stables}
   On suppose $G$ quasi-d\'eploy\'e. Soit $M$ un Levi de $G$.
   \ass{Corollaire}{Pour tout $\boldsymbol{\delta}\in D_{g\acute{e}om,G-\acute{e}qui}^{st}(M(F))\otimes Mes(M(F))^*$ assez proche de $1$, le terme $Sg_{M,unip}^G(\boldsymbol{\delta},B)$ appartient \`a $D_{unip}^{st}(G(F))\otimes Mes(G(F))^*$.}
   
   Preuve. On applique le lemme [II] 2.9, en prenant pour ${\cal D}^{st}$ l'espace $D_{g\acute{e}om,G-\acute{e}qui}^{st}(M(F))\otimes Mes(M(F))^*$ tout entier. L'hypoth\`ese de ce lemme est v\'erifi\'ee d'apr\`es 1.1. Le (i) du lemme l'est aussi d'apr\`es le (iii) de la proposition pr\'ec\'edente. Donc le (ii) aussi, ce qui est l'assertion de l'\'enonc\'e. $\square$

\bigskip

\subsection{Int\'egrales orbitales pond\'er\'ees endoscopiques}
Le groupe $G$ est quelconque. Soient $M$ un Levi de $G$ et ${\bf M}'=(M',{\cal M}',\zeta)$ une donn\'ee endoscopique de $M$, elliptique et relevante. Pour $J\in {\cal J}_{M}^G(B)$, on a d\'efini en [II] 3.8 un homomorphisme
$$\rho_{J}^{G,{\cal E}}({\bf M}'):D^{st}_{unip}({\bf M}')\otimes Mes(M'(F))^*\to U_{J}\otimes (D_{unip}(M(F))\otimes Mes(M(F))^*)/Ann_{unip}^{G}.$$

\ass{Proposition}{(i) Pour tout $J\in{\cal J}_{M}^G(B)$, $\rho_{J}^{G,{\cal E}}({\bf M}')$ est le compos\'e de $\rho_{J}^{G_{J},{\cal E}}({\bf M}')$ et de la projection
$$U_{J}\otimes (D_{unip}(M(F))\otimes Mes(M(F))^*)/Ann_{unip}^{G_{J}}\to U_{J}\otimes (D_{unip}(M(F))\otimes Mes(M(F))^*)/Ann_{unip}^{G}.$$

(ii) Pour tout $J\in{\cal J}_{M}^G(B)$, tout $\boldsymbol{\delta}\in D_{unip}^{st}({\bf M}')\otimes Mes(M'(F))^*$ et tout $a\in A_{M}(F)$ en position g\'en\'erale et assez proche de $1$, on a l'\'egalit\'e
$$\rho_{J}^{G,{\cal E}}({\bf M}',\boldsymbol{\delta},a)=\rho_{J}^G(transfert(\boldsymbol{\delta}),a).$$

(iii) Pour tout $\boldsymbol{\gamma}\in D_{unip}(M(F))\otimes Mes(M(F))^*$ et tout ${\bf f}\in I(G(F))\otimes Mes(G(F))$, on a l'\'egalit\'e
$$I_{M}^{G,{\cal E}}(\boldsymbol{\gamma},B,{\bf f})=I_{M}^{G}(\boldsymbol{\gamma},B,{\bf f}).$$}

Preuve.  Rappelons la d\'efinition, pour $\boldsymbol{\delta}\in D_{unip}^{st}({\bf M}')\otimes Mes(M'(F))^*$ et $a\in A_{M}(F)$ en position g\'en\'erale et proche de $1$:
$$\rho_{J}^{G,{\cal E}}({\bf M}',\boldsymbol{\delta},a)=\sum_{s\in\zeta Z(\hat{M})^{\Gamma_{F}}/Z(\hat{G})^{\Gamma_{F}}}i_{M'}(G,G'(s))\sum_{J'\in {\cal J}_{M'}^{G'(s)}(B);J'\mapsto J}transfert(\sigma_{J'}^{{\bf G}'(s)}(\boldsymbol{\delta},\xi(a))),$$
cf. [II] 3.8. Comme on l'a dit en [II] 3.6, il y a un isomorphisme
$$\iota:D_{unip}^{st}(M'(F))\otimes Mes(M'(F))^*\to D_{unip}^{st}({\bf M}')\otimes Mes(M'(F))^*$$
et on a l'\'egalit\'e
$$\sigma_{J'}^{{\bf G}'(s)}(\boldsymbol{\delta},\xi(a))=\iota(\sigma_{J'}^{G'(s)}(\iota^{-1}(\boldsymbol{\delta}),\xi(a)))$$
pour tout $s$ apparaissant ci-dessus. En posant $\boldsymbol{\delta}'=\iota^{-1}(\boldsymbol{\delta})$, la formule ci-dessus devient
$$\rho_{J}^{G,{\cal E}}({\bf M}',\boldsymbol{\delta},a)=transfert\circ \iota(X_{J}^G),$$
o\`u
$$(1) \qquad X_{J}^G=\sum_{s\in\zeta Z(\hat{M})^{\Gamma_{F}}/Z(\hat{G})^{\Gamma_{F}}}i_{M'}(G,G'(s))\sum_{J'\in {\cal J}_{M'}^{G'(s)}(B);J'\mapsto J}\sigma_{J'}^{ G'(s)}(\boldsymbol{\delta}',\xi(a)).$$
Plus exactement, ce qui nous importe est la projection de ce terme modulo le sous-espace $Ann'\subset D_{unip}^{st}(M'(F))\otimes Mes(M'(F))^*$ form\'e des \'el\'ements dont l'image par $transfert\circ \iota$ appartient \`a $Ann_{unip}^{G}$. En vertu des formules ci-dessus, il suffit pour prouver (i) de montrer que, modulo cet espace $Ann'$, $X_{J}^G$ co\"{\i}ncide avec $X_{J}^{G_{J}}$. La preuve est alors similaire \`a celle de 1.2. Un \'el\'ement $s\in\zeta Z(\hat{M})^{\Gamma_{F}}/Z(\hat{G})^{\Gamma_{F}}$ d\'etermine \`a la fois un groupe endoscopique $G'(s)$ de $G$ et un groupe endoscopique $(G_{J})'(s)$ de $G_{J}$. Les propri\'et\'es similaires \`a 1.2(3) et (4) sont v\'erifi\'ees:

(2) les ensembles (r\'eduits \`a au plus un \'el\'ement) $\{J'\in {\cal J}_{M'}^{G'(s)}(B);J'\mapsto J\}$ et $\{J'\in {\cal J}_{M'}^{(G_{J})'(s)}(B); J'\mapsto J\}$ co\"{\i}ncident;

(3) s'ils ne sont pas vides, notons $J'$ leur seul \'el\'ement; alors $G'(s)_{J'}=(G_{J})'(s)$.

En utilisant ces propri\'et\'es et en appliquant la proposition 1.2(i) aux termes du membre de droite de (1), on transforme (1) en
$$X_{J}^G=\sum_{s\in\zeta Z(\hat{M})^{\Gamma_{F}}/Z(\hat{G})^{\Gamma_{F}}}i_{M'}(G,G'(s))\sum_{J'\in {\cal J}_{M'}^{(G_{J})'(s)}(B);J'\mapsto J}i_{J'}^{G'(s)}\sigma_{J'}^{ (G_{J})'(s)}(\boldsymbol{\delta}',\xi(a))$$
du moins modulo $Ann'$. 
De nouveau, on calcule
$$i_{M'}(G,G'(s))i_{J'}^{G'(s)}=[Z(\hat{G}_{J})^{\Gamma_{F}}:Z(\hat{G})^{\Gamma_{F}}]^{-1}i_{M'}(G_{J},(G_{J})'(s)).$$
D'o\`u
$$X_{J}^G=[Z(\hat{G}_{J})^{\Gamma_{F}}:Z(\hat{G})^{\Gamma_{F}}]^{-1}\sum_{s\in\zeta Z(\hat{M})^{\Gamma_{F}}/Z(\hat{G})^{\Gamma_{F}}}i_{M'}(G_{J},(G_{J})'(s))$$
$$\sum_{J'\in {\cal J}_{M'}^{(G_{J})'(s)}(B);J'\mapsto J}\sigma_{J'}^{ (G_{J})'(s)}(\boldsymbol{\delta}',\xi(a))$$
$$=\sum_{s\in\zeta Z(\hat{M})^{\Gamma_{F}}/Z(\hat{G}_{J})^{\Gamma_{F}}}i_{M'}(G_{J},(G_{J})'(s))\sum_{J'\in {\cal J}_{M'}^{(G_{J})'(s)}(B);J'\mapsto J}\sigma_{J'}^{ (G_{J})'(s)}(\boldsymbol{\delta}',\xi(a))=X_{J}^{G_{J}}.$$
Cela prouve le (i) de l'\'enonc\'e. 

Posons $\boldsymbol{\gamma}=transfert(\boldsymbol{\delta})$. Pour tout $a\in A_{M}(F)$, on a $a\boldsymbol{\gamma}=transfert(\xi(a)\boldsymbol{\delta})$. Soit ${\bf f}\in I(G(F))\otimes Mes(G(F))$. Utilisons les d\'eveloppements [II] 3.3 de $I_{M}^G(a\boldsymbol{\gamma},{\bf f})$ et [II] 3.9 de $I_{M}^{G,{\cal E}}(a\boldsymbol{\gamma},{\bf f})=I_{M}^{G,{\cal E}}({\bf M}',\xi(a)\boldsymbol{\delta},{\bf f})$. On obtient que
$$(4) \qquad I_{M}^G(a\boldsymbol{\gamma},{\bf f})-I_{M}^{G,{\cal E}}(a\boldsymbol{\gamma},{\bf f})$$
est \'equivalent \`a
$$\sum_{L\in {\cal L}(M)}\sum_{J\in {\cal J}_{M}^L(B)}I_{L}^G(\rho_{J}^L(\boldsymbol{\gamma},a)^{L},B,{\bf f})-I_{L}^{G,{\cal E}}(\rho_{J}^{L,{\cal E}}({\bf M}',\boldsymbol{\delta},a)^{L},B,{\bf f}).$$
On sait d'apr\`es 1.1 que l'expression (4) est nulle. En raisonnant par r\'ecurrence et en utilisant le (i) de l'\'enonc\'e, on conna\^{\i}t l'assertion (ii) pour tout $J$ sauf  pour le terme maximal $J_{max}$ de ${\cal J}_{M}^G(B)$. On sait aussi que les int\'egrales orbitales $I_{L}^{G,{\cal E}}$ co\"{\i}ncident avec $I_{L}^{G}$ si $L\not=M$. L'expression ci-dessus se simplifie et on obtient  que l'expression 
$$I^G(\rho^G_{J_{max}}(\boldsymbol{\gamma},a)^G-\rho_{J_{max}}^{G,{\cal E}}({\bf M}',\boldsymbol{\delta},a)^G,{\bf f})$$
$$+I_{M}^G(\boldsymbol{\gamma},B,{\bf f})-I_{M}^{G,{\cal E}}(\boldsymbol{\gamma},B,{\bf f})$$
est nulle. La propri\'et\'e [II] 3.1(3) entra\^{\i}ne de nouveau la nullit\'e de ces deux termes. La nullit\'e du premier ach\`eve de prouver l'assertion (ii) de l'\'enonc\'e. La nullit\'e du second d\'emontre l'assertion (iii) pour notre distribution $\boldsymbol{\gamma}$. Mais tout \'el\'ement de $D_{unip}(M(F))\otimes Mes(M(F))^*$ est combinaison lin\'eaire de transferts de $\boldsymbol{\delta}$ comme ci-dessus, quand on fait varier la donn\'ee ${\bf M}'$. Il s'ensuit que l'assertion (iii) est v\'erifi\'ee pour tout $\boldsymbol{\gamma}\in D_{unip}(M(F))\otimes Mes(M(F))^*$. $\square$

\bigskip

\subsection{Germes endoscopiques}
Soit $M$ un Levi de $G$.

\ass{Corollaire}{Pour tout $\boldsymbol{\gamma}\in D_{g\acute{e}om}(M(F))\otimes Mes(M(F))^*$ assez proche de $1$, on a l'\'egalit\'e
$$g_{M,unip}^G(\boldsymbol{\gamma},B)=g_{M,unip}^{G,{\cal E}}(\boldsymbol{\gamma},B).$$}

Comme en 1.3, cela se d\'eduit de la proposition pr\'ec\'edente en utilisant le lemme [II] 2.8. $\square$

  \bigskip
  
  \section{Premiers r\'esultats dans le cas quasi-d\'eploy\'e et \`a torsion int\'erieure}

\subsection{Un lemme sur les groupes ab\'eliens finis}
Soient $X$ un groupe ab\'elien fini et $n$ un entier sup\'erieur ou \'egal \`a $1$. Tout \'el\'ement $\underline{m}=(m_{1},...,m_{n})\in {\mathbb Z}^n$ d\'etermine un homomorphisme 
$$\begin{array}{cccc}\varphi_{\underline{m}}:&X^n&\to& X\\ &\underline{x}=(x_{1},...,x_{n})&\mapsto&\sum_{i=1,...,n}m_{i}x_{i}\\ \end{array}$$
Evidemment, cet homomorphisme ne d\'epend que de l'image de $\underline{m}$ dans $({\mathbb Z}/N{\mathbb Z})^n$, o\`u $N$ est l'exposant de $X$ (c'est-\`a-dire le plus petit entier sup\'erieur ou \'egal \`a $1$ qui annule $X$).

\ass{Lemme}{Soit $\underline{x},\underline{y}\in X^n$. Alors $\underline{y}$ appartient au sous-groupe de $X^n$ engendr\'e par $\underline{x}$ si et seulement si, pour tout $\underline{m}\in {\mathbb Z}^n$, $\varphi_{\underline{m}}(\underline{y})$ appartient au sous-groupe de $X$ engendr\'e par $\varphi_{\underline{m}}(\underline{x})$. }

Preuve. Dans un sens, c'est \'evident: si $\underline{y}=r\underline{x}$, avec $r\in {\mathbb Z}$, alors $\varphi_{\underline{m}}(\underline{y})=r\varphi_{\underline{m}}(\underline{x})$ pour tout $\underline{m}$. Supposons inversement que, pour tout $\underline{m}\in {\mathbb Z}^n$, $\varphi_{\underline{m}}(\underline{y})$ appartient au sous-groupe de $X$ engendr\'e par $\varphi_{\underline{m}}(\underline{x})$. D\'ecomposons $X$ en somme directe  $\oplus_{p\in P}X_{p}$ o\`u $P$ est un ensemble fini de nombres premiers et $X_{p}$ est un $p$-groupe pour tout $p\in P$.   On d\'ecompose conform\'ement tout $z\in X$ en $z=\sum_{p\in P}z_{p}$. Le sous-groupe de $X$ engendr\'e par $z$ est l'ensemble des $z'=\sum_{p\in P}z'_{p}$ tels que, pour tout $p\in P$, $z'_{p}$ appartienne au sous-groupe de $X_{p}$ engendr\'e par $z_{p}$. La m\^eme propri\'et\'e s'applique \`a $X^n=\oplus_{p\in P}X_{p}^n$. Donc, pour tout $p\in P$ et pour tout $\underline{m}\in {\mathbb Z}^n$, $(\varphi_{\underline{m}}(\underline{y}))_{p}$ appartient au sous-groupe de $X_{p}$ engendr\'e par $(\varphi_{\underline{m}}(\underline{x}))_{p}$. Pour tout $\underline{z}\in X^n$, on a $(\varphi_{\underline{m}}(\underline{z}))_{p}=\varphi_{\underline{m}}(\underline{z}_{p})$.  Donc, pour tout $p\in P$ et pour tout $\underline{m}\in {\mathbb Z}^n$, $\varphi_{\underline{m}}(\underline{y}_{p})$ appartient au sous-groupe de $X_{p}$ engendr\'e par $\varphi_{\underline{m}}(\underline{x}_{p})$. Supposons le lemme d\'emontr\'e pour chaque $X_{p}$. Alors la propri\'et\'e pr\'ec\'edente entra\^{\i}ne que, pour tout $p\in P$, $\underline{y}_{p}$ appartient au sous-groupe de $X_{p}^n$ engendr\'e par $\underline{x}_{p}$. D'o\`u la conclusion. 

On est ainsi ramen\'e au cas o\`u $X$ est un $p$-groupe pour  un certain  nombre premier $p$. Ecrivons $\underline{x}=(x_{1},...,x_{n})$, $\underline{y}=(y_{1},...,y_{n})$.  Pour tout $i=1,...,n$, notons $a_{i}\in {\mathbb N}$ le plus petit entier   tel que $p^{a_{i}}x_{i}=0$. A permutation pr\`es, on peut supposer $a_{1}\geq...\geq a_{n}$. L'assertion \`a prouver \'etant \'evidente dans le cas $n=1$, on suppose $n\geq2$. On pose $\underline{x}'=(x_{1},...,x_{n-1})$, $\underline{y}'=(y_{1},...,y_{n-1})$. Ces elements v\'erifient la m\^eme hypoth\`ese que $\underline{x}$ et $\underline{y}$, mais pour $n-1$. En raisonnant par r\'ecurrence, on peut supposer qu'il existe $r\in {\mathbb Z} $ tel que $\underline{y}'=r\underline{x}'$. Alors l'\'el\'ement $(0,...,0,y_{n}-rx_{n})=\underline{y}-r\underline{x}$ v\'erifie la m\^eme hypoth\`ese que $\underline{y}$. On va montrer qu'il est nul. En oubliant cette construction, on suppose simplement que $\underline{y}=(0,...,0,y_{n})$ et on va prouver que $y_{n}=0$. En appliquant l'hypoth\`ese \`a $\underline{m}=(0,...,0,1)$, on voit qu'il existe $r\in {\mathbb Z}$ tel que $y_{n}=rx_{n}$. D'o\`u la conclusion si $x_{n}=0$. On suppose $x_{n}\not=0$. Soit $h\in {\mathbb N}$ le plus petit entier tel que $p^hx_{n}$ appartienne au sous-groupe de $X$ engendr\'e par $x_{1}$. On a $h\leq a_{n}$ et une \'egalit\'e $p^hx_{n}+p^{h'}ux_{1}=0$, o\`u $u\in {\mathbb Z}$ est premier \`a $p$ et $h'=h+a_{1}-a_{n}$. Posons $z=x_{n}+p^{a_{1}-a_{n}}ux_{1}$. On v\'erifie que l'application
$$\begin{array}{ccc}{\mathbb Z}/p^{a_{1}}{\mathbb Z}\oplus {\mathbb Z}/p^h{\mathbb Z}&\to&X\\ (e,f)&\mapsto& ex_{1}+fz\\ \end{array}$$
est injective. Appliquons l'hypoth\`ese \`a $\underline{m}=(m,0,...,0,1)$. On obtient que $rx_{n}$ appartient au groupe engendr\'e par $mx_{1}+x_{n}$. Autrement dit $r(z-p^{a_{1}-a_{n}}ux_{1})$ appartient au groupe engendr\'e par $z+(m-p^{a_{1}-a_{n}}u)x_{1}$. Posons $m=p^{a_{1}-a_{n}}u$. Alors $r(z-p^{a_{1}-a_{n}}ux_{1})$  appartient au groupe engendr\'e par $z$. D'apr\`es l'injectivit\'e pr\'ec\'edente, $p^{a_{1}}$ doit diviser $rup^{a_{1}-a_{n}}$. Donc $p^{a_{n}}$ divise $r$. D'o\`u $y_{n}=rx_{n}=0$. $\square$

\bigskip

\subsection{Un lemme sur les tores}
Dans ce paragraphe et les trois suivants, on l\`eve l'hypoth\`ese que $F$ est non-archim\'edien. Le corps $F$ est  un corps local de caract\'eristique nulle. 

\ass{Lemme}{Soient $T$ un tore d\'efini sur $F$ et $U\subset T(F)$ un sous-groupe ouvert d'indice fini. Alors il existe un tore $T'$ d\'efini sur $F$ et un homomorphisme $f:T'\to T$ d\'efini sur $F$ de sorte que $f(T'(F))=U$. }

Preuve dans le cas o\`u $F$ est archim\'edien.  Si $F={\mathbb C}$, $T({\mathbb C})$ est connexe. Donc  $U=T({\mathbb C})$. Le tore $T'=T$ et l'homomorphisme identit\'e conviennent. Supposons $F={\mathbb R}$. Introduisons les trois tores $T_{1}$, $T_{2}$ et $T_{3}$ sur ${\mathbb R}$ tels que $T_{1}({\mathbb R})={\mathbb R}^{\times}$, $T_{2}({\mathbb R})={\mathbb C}^{\times}$, $T_{3}({\mathbb R})= \{z\in {\mathbb C}; z\bar{z}=1\}$. On sait que $T$ est isomorphe \`a un produit de tels tores, disons $T=T_{1}^{a}\times T_{2}^b\times T_{3}^c$. Le sous-groupe $U$ est n\'ecessairement de la forme $U_{1}\times T_{2}({\mathbb R})^b\times T_{3}({\mathbb R})^c$, o\`u $U_{1}$ est un sous-groupe ouvert d'indice fini de $T_{1}({\mathbb R})^{a}$.   Si on trouve $T'_{1}$ et $f_{1}$ r\'esolvant le probl\`eme pour le tore $T_{1}^{a}$ et le sous-groupe $U_{1}$, on pose $T'=T'_{1}\times T_{2}^b\times T_{3}^c$, on \'etend $f_{1}$ en $f$ par l'identit\'e sur les autres composantes.  Cela r\'esout le probl\`eme initial. On est ainsi ramen\'e au cas o\`u $T=T_{1}^{a}$. Quitte \`a appliquer un automorphisme de $T$, on peut supposer qu'il existe un entier $e$ avec $0\leq e\leq a$ de sorte que $U=({\mathbb R}_{+}^{\times})^{e}\times ({\mathbb R}^{\times})^{a-e}$. On pose $T'=T_{2}^{e}\times T_{1}^{a-e}$, on d\'efinit $f$ comme \'etant la norme sur les $e$-premi\`eres composantes et l'identit\'e sur les $a-e$ derni\`eres. Cela r\'esout le probl\`eme.

Preuve dans le cas o\`u $F$ est non-archim\'edien. On fixe une extension finie $E$ de $F$ tel que $\Gamma_{E}$ agisse trivialement sur $X_{*}(T)$. On introduit le tore $S=Res_{E/F}(GL(1)_{E})$. Le groupe  $X_{*}(S)$ est le groupe des fonctions $\phi:\Gamma_{E}\backslash \Gamma_{F}\to {\mathbb Z}$, muni de l'action de $\Gamma_{F}$ par translations \`a droite. On a $S(F)=E^{\times}$. On introduit le tore $D$ tel que $X_{*}(D)$ soit $X_{*}(S)\otimes_{{\mathbb Z}} X_{*}(T)$, c'est-\`a-dire le groupe des fonctions $\phi:\Gamma_{E}\backslash \Gamma_{F}\to X_{*}(T)$, muni de l'action de $\Gamma_{F}$ par translations \`a droite. On a $D\simeq S^n$, o\`u $n$ est la dimension de $T$.  On a un plongement $\iota:T\to D$ ainsi d\'efini: pour $x_{*}\in X_{*}(T)$, $\iota\circ x_{*}$ est l'\'el\'ement $\phi$ de $X_{*}(D)$ tel que $\phi(\sigma)=\sigma(x_{*})$. Ce plongement est d\'efini sur $F$. Soit $N\geq1$ un entier. Pour tout groupe ab\'elien $Y$, notons $Y^{(N)}$ le groupe des puissances $N$-i\`emes dans $Y$. Montrons que

(1) il existe $N$ tel que $D(F)^{(N)}\cap \iota(T(F))\subset \iota(U)$.

Introduisons le sous-groupe compact maximal $D(F)_{c}$ de $D(F)$. Les sous-groupes $(D(F)_{c})^N$ forment un syst\`eme de voisinages ouverts de l'origine dans $D(F)$.  Le plongement $\iota:T(F)\to D(F)$ est une immersion ferm\'ee. Puisque $U$ est ouvert dans $T(F)$, il existe un  entier $N_{1}\geq1$ tel que $(D(F)_{c})^{N_{1}}\cap \iota(T(F))\subset \iota(U)$. 
On a une suite exacte
$$1\to D(F)_{c}\to D(F)\stackrel{\pi}{\to} {\mathbb Z}^n\to 0$$
Posons $L_{T}=\pi\circ\iota(T(F))$, $L_{U}=\pi\circ\iota(U)$ et $L_{0}={\mathbb Z}^n\cap (L_{T}\otimes_{{\mathbb Z}}{\mathbb Q})$. Le groupe $L_{U}$ est d'indice fini dans $L_{T}$ par hypoth\`ese et $L_{T}$ est d'indice fini dans $L_{0}$. Soit $N_{2}\geq1$ tel que $N_{2}L_{0}\subset L_{U}$. Soit $N=N_{1}N_{2}$ et soit $d\in D(F)$ tel que $d^N\in \iota(T(F))$. Alors $N\pi(d)\in L_{T}\subset L_{0}$. Le groupe ${\mathbb Z}^n/L_{0}$ est sans torsion. Donc $\pi(d)\in L_{0}$ puis $N_{2}\pi(d)\in L_{U}$. On peut donc trouver $u\in U$ et $d_{c}\in D(F)_{c}$ tels que $d^{N_{2}}= \iota(u)d_{c}$. On a $d_{c}^{N_{1}}=d^N\iota(u)^{-N_{1}}$. Ceci appartient \`a $\iota(T(F))$, donc \`a $(D(F)_{c})^{N_{1}}\cap \iota(T(F))$, donc \`a $\iota(U)$. Posons $d_{c}^{N_{1}}=\iota(v)$, avec $v\in U$.  Alors $d^N=\iota(u^{N_{1}}v)$ appartient \`a $\iota(U)$. Cela d\'emontre (1). 

Fixons $N$ v\'erifiant (1). Si $N=1$, on a $U=T(F)$ et le lemme est \'evident (on prend $T'=T$ et $f$ l'identit\'e). Supposons $N>1$. Soit $d\in D(F)$, notons $V$ le sous-groupe de $D(F)$ engendr\'e par $d$ et $D(F)^{(N)}$. Montrons que

(2) il existe un tore $D'$ d\'efini sur $F$ et un homomorphisme $g:D'\to D$ d\'efini sur $F$ de sorte que $g(D'(F))=V$ et que le noyau de $g$ soit connexe.

On identifie $D$ \`a $S^n$. Si $n=1$, $V$ est un sous-groupe ouvert d'indice fini de $S(F)=E^{\times}$. On sait qu'il existe une extension finie $E'$ de $E$ telle que, en notant $N_{E'/E}$ la norme, $V$ soit \'egal au groupe des normes $N_{E'/E}(E^{'\times})$ ([S] XIV.6 th\'eor\`eme 1). On pose $S'=Res_{E'/F}(GL(1)_{E'})$. On construit facilement un homomorphisme $g:S'\to S$ d\'efini sur $F$ dont l'homomorphisme d\'eduit de $S'(F)=E^{'\times}$ dans $S(F)=E^{\times}$ soit la norme $N_{E'/E}$. Avec la description donn\'ee plus haut de $X_{*}(S)$ et la description similaire de $X_{*}(S')$, pour $\phi'\in X_{*}(S')$, on a $g\circ\phi'(\sigma)=\sum_{\gamma\in \Gamma_{E'}\backslash \Gamma_{E}}\phi'(\gamma\sigma)$ pour tout $\sigma\in \Gamma_{F}$. Il en r\'esulte que l'on a une suite exacte
$$0\to Y\to X_{*}(S')\to X_{*}(S)\to 0,$$
o\`u $Y$ est un ${\mathbb Z}$-module libre. Le noyau de $g$ est donc connexe. Alors le tore $D'=S'$ et cet homomorphisme $g$ conviennent.

Supposons maintenant $n>1$. On choisit un sous-ensemble $\underline{M}\subset {\mathbb Z}^n$ qui s'envoie bijectivement sur ${\mathbb Z}^n/(N{\mathbb Z})^n$. On suppose que $\underline{M}$ contient les \'el\'ements de base de ${\mathbb Z}^n$, c'est-\`a-dire les \'el\'ements qui ont une coordonn\'ee \'egale \`a $1$ et dont les autres coordonn\'ees sont nulles. Pour tout $\underline{m}=(m_{1},...,m_{n})\in \underline{M}$, on d\'efinit l'homomorphisme
$$\begin{array}{cccc}\varphi_{\underline{m}}:&D=S^n&\to&S\\&(x_{1},...,x_{n})&\mapsto &\prod_{i=1,...,n}x_{i}^{m_{i}}\\ \end{array}$$
On note $V_{\underline{m}}$ le sous-groupe de $S(F)$ engendr\'e par $S(F)^{(N)}$ et par $\varphi_{\underline{m}}(d)$.  En appliquant le r\'esultat du cas $n=1$, on choisit un tore $S_{\underline{m}}$ d\'efini sur $F$ et un homomorphisme $g_{\underline{m}}:S_{\underline{m}}\to S$ d\'efini sur $F$ de sorte que $g_{\underline{m}}(S_{\underline{m}}(F))=V_{\underline{m}}$ et que le noyau de $g_{\underline{m}}$ soit connexe.  En posant $S_{\underline{M}}=\prod_{\underline{m}\in \underline{M}}S_{\underline{m}}$ et $S^{\underline{m}}=\prod_{\underline{m}\in \underline{M}}S$, les homomorphismes $g_{\underline{m}}$ se regroupent en un homomorphisme $g_{\underline{M}}:S_{\underline{M}}\to S^{\underline{M}}$. Son noyau est connexe et on a l'\'egalit\'e $g_{\underline{M}}(S_{\underline{M}}(F))=\prod_{\underline{m}\in \underline{M}}
V_{\underline{m}}$.  D'autre part, les homomorphismes $\varphi_{\underline{m}}$ se regroupent en un homomorphisme $\varphi_{\underline{M}}:D=S^n\to S^{\underline{M}}$. C'est une immersion ferm\'ee:   quand $\underline{m}$ d\'ecrit les \'el\'ements de base de ${\mathbb Z}^n$, les applications $\varphi_{\underline{m}}$ d\'ecrivent les applications coordonn\'ees naturelles sur $S^n$. Notons $D'$ le produit fibr\'e de $D$ et $S_{\underline{M}}$ au-dessus de $S^{\underline{M}}$. Autrement dit $D'(\bar{F})$ est le groupe des $(x,y)\in D(\bar{F})\times S_{\underline{M}}(\bar{F})$ tels que $\varphi_{\underline{M}}(x)=g_{\underline{M}}(y)$. Parce que $\varphi_{\underline{M}}$ est une immersion ferm\'ee et que le noyau de $g_{\underline{M}}$ est connexe, $D'$ est connexe. C'est donc un tore, qui est \'evidemment d\'efini sur $F$. On note $g:D'\to D$ la projection $(x,y)\mapsto x$. Cet homomorphisme est d\'efini sur $F$. Son noyau est celui de $g_{\underline{M}}$, donc est connexe.  Le groupe $g(D'(F))$ est celui des $x\in D(F)$ tels que, pour tout $\underline{m}\in \underline{M}$, $\varphi_{\underline{m}}(x)$ appartienne \`a $g_{\underline{m}}(S_{\underline{m}}(F))$, autrement dit \`a $V_{\underline{m}}$. En appliquant le lemme du paragraphe pr\'ec\'edent au groupe $X=S(F)/S(F)^{(N)}$, on obtient que $g(D'(F))=V$.  Cela prouve (2).

La propri\'et\'e (2) s'\'etend de la fa\c{c}on suivante. Soit $V$ un sous-groupe de $D(F)$ contenant $D(F)^{(N)}$. Alors 

(3) il existe un tore $D'$ d\'efini sur $F$ et un homomorphisme $g:D'\to D$ d\'efini sur $F$ de sorte que $g(D'(F))=V$ et que le noyau de $g$ soit connexe.

Le groupe $V$ est engendr\'e par $D(F)^{(N)}$ et un ensemble fini d'\'el\'ements $d_{1},...,d_{k}$. On peut supposer $k\geq1$, quitte \`a prendre $d_{1}=1$. Si $k=1$, on applique l'assertion (2). Si $k\geq2$, on note $V_{1}$, resp. $V_{2}$, le sous-groupe de $D(F)$ engendr\'e par $D(F)^{(N)}$ et les \'el\'ements $d_{1},...,d_{k-1}$, resp. $d_{k}$. En raisonnant par r\'ecurrence, on choisit $D'_{1}$ et $g_{1}$, resp. $D'_{2}$ et $g'_{2}$, v\'erifiant (3) pour le groupe $V_{1}$, resp. $V_{2}$. On pose $D'=D'_{1}\times D'_{2}$ et on prend pour $g$ le produit de $g_{1}$ et $g_{2}$. Il est clair que $D$ et $g$ sont d\'efinis sur $F$ et que $g(D(F))=V$. Le noyau de $g$ est fibr\'e au-dessus de $D'_{1}$, de fibres isomorphes au noyau de $g_{2}$. Donc ce noyau est connexe. Cela prouve (3).

Appliquons (3) au groupe $V =D(F)^{(N)}\iota(U)$. On en d\'eduit un tore $D'$ et un homomorphisme $g$. Soit $T'$ le produit fibr\'e de $T$ et $D'$ au-dessus de $D$. C'est-\`a-dire que $T'(\bar{F})$ est le groupe des $(t,d')\in T(\bar{F})\times D'(\bar{F})$ tels que $\iota(t)=g(d')$. Notons $f:T'\to T$ la projection $(t,d')\mapsto t$. Cette projection est surjective et son noyau est isomorphe \`a celui de $g$, donc est connexe. Donc le groupe $T'$ est lui-m\^eme connexe et c'est un tore. Il est clair que $T$ et $f$ sont d\'efinis sur $F$. L'image $f(T'(F))$ est le sous-groupe des $t\in T(F)$ tels que $\iota(t)$ appartienne \`a $g(D'(F))$, autrement dit \`a $D(F)^{(N)}\iota(U)$. En appliquant (1), on obtient que $f(T'(F))=U$. Cela ach\`eve la d\'emonstration. $\square$

\bigskip

\subsection{D\'etordre un triplet $(G,\tilde{G},{\bf a})$ quasi-d\'eploy\'e et \`a torsion int\'erieure}

Soit $(G,\tilde{G},{\bf a})$ un triplet quasi-d\'eploy\'e et \`a torsion int\'erieure. On suppose comme toujours $\tilde{G}(F)\not=\emptyset$.

\ass{Proposition}{Il existe des objets $H$, $D$, $d$, $\iota$, $\tilde{\iota}$, $q$ v\'erifiant les conditions suivantes:

(i) $H$ est un groupe r\'eductif connexe d\'efini et quasi-d\'eploy\'e sur $F$;

(ii) $D$ est un tore d\'efini sur $F$;

(iii) $d\in D(F)$;

(iv) $\iota:G\to H$ est un plongement d\'efini sur $F$ dont l'image est un sous-groupe distingu\'e de $H$;

(v) $q:H\to D$ est un homomorphisme;

(vi) la suite 
$$1\to G\stackrel{\iota}{\to}H\stackrel{q}{\to}D\to 1$$
est exacte;

(vii) en notant $H_{d}=\{h\in H; q(h)=d\}$, $\tilde{\iota}:\tilde{G}\to H_{d}$ est un isomorphisme de vari\'et\'es alg\'ebriques d\'efinies sur $F$ tel que $\tilde{\iota}(g\gamma g')=\iota(g)\tilde{\iota}(\gamma)\iota(g')$ pour tous $g,g'\in G$, $\gamma\in \tilde{G}$;

(viii)  $q(H(F))=D(F)$ et ce groupe est engendr\'e par $q(Z(H;F))$ et par $d$.}

Preuve. On a construit en [W1] 1.3(6) et (7) des objets $H'$, $D'$, $\iota'$, $q'$ v\'erifiant les analogues de (i), (ii), (iv), (v), et tels que $Z(H')$ soit connexe et soit un tore induit.  Comme on l'a dit en [I] 1.9, l'ensemble ${\cal Z}(\tilde{G})$ s'identifie \`a celui des $e\in \tilde{G}$ tels que $ad_{e}$ soit l'identit\'e. Fixons un \'el\'ement $e\in {\cal Z}(\tilde{G})$. Il y a un cocycle $z:\Gamma_{F}\to Z(G)$ tel que $\sigma(e)=z(\sigma)^{-1}e$ pour tout $\sigma\in \Gamma_{F}$. Puisque $Z(H')$ est induit, le cocycle $\iota\circ z$ est un bord. On peut fixer $e_{H'}\in Z(H')$ tel que $\sigma(e_{H'})=\iota\circ z(\sigma)^{-1}e_{H'}$ pour tout $\sigma$. On pose $d'=q'(e_{H'})$. On d\'efinit $\tilde{\iota}':\tilde{G}\to H'_{d'}$ par $\tilde{\iota}'(ge)=\iota'(g)e_{H'}$. On voit que toutes nos conditions sont v\'erifi\'ees, sauf \'eventuellement la huiti\`eme. Appliquons le lemme du paragraphe pr\'ec\'edent au tore $D'$ et au groupe  $U$ engendr\'e par $q'(Z(H';F))$ et $d'$. On obtient un tore que nous notons $D$ et un homomorphisme $f:D\to D'$. Notons $H$ le produit fibr\'e de $H'$ et $D$ au-dessus de $D'$. C'est-\`a-dire que $H(\bar{F})$ est le groupe des $(x,y)\in H'(\bar{F})\times D(\bar{F})$ tels que $q'(x)=f(y)$.  On note $\iota:G\to H$ le plongement  $g\mapsto (\iota'(g),1)$ et $q:H\to D$ la projection $(x,y)\mapsto y$. La suite 
$$1\to G\stackrel{\iota}{\to}H\stackrel{q}{\to}D\to 1$$
est exacte. Cela prouve que $H$ est connexe. On sait qu'il existe $d\in D(F)$ tel que $f(d)=d'$. On fixe un tel $d$ et on d\'efinit $\tilde{\iota}:\tilde{G}\to H_{d}$ par $\tilde{\iota}(\gamma)=(\tilde{\iota}'(\gamma),d)$. Les sept premi\`eres propri\'et\'es de l'\'enonc\'e sont v\'erifi\'ees. Le groupe $Z(H)$ est le produit fibr\'e de $Z(H')$ et de $D$ au-dessus de $D'$. Remarquons que $d$ appartient \`a $q(H(F))$: on a $d=q\circ\tilde{\iota}(\gamma)$ pour tout $\gamma\in \tilde{G}(F)$.   Soit $y\in D(F)$. On a $f(y)\in U$. On peut donc \'ecrire $f(y)=q'(z')(d')^n$, avec $z'\in Z(H';F)$ et $n\in {\mathbb Z}$. Posons $y_{1}=yd^{-n}$. Alors $y_{1}\in D(F)$ et $f(y_{1})=q'(z')$. L'\'el\'ement $z=(z',y_{1})$ appartient \`a $Z(H;F)$.  Alors $y=q(z)d^n$. Cela prouve que $D(F)$ est engendr\'e par $q(Z(H;F))$ et par $d$.  Puisque $d$ appartient \`a $q(H(F))$ et que $q(Z(H;F))$ est inclus dans ce groupe, on a aussi $D(F)=q(H(F))$. $\square$

Pour la suite de la section, les hypoth\`eses sont celles de ce paragraphe et on fixe des objets v\'erifiant la proposition. Pour simplifier les notations, on oublie $\iota$ et $\tilde{\iota}$ en identifiant $G$ et $\tilde{G}$ \`a des sous-ensembles de $H$ via ces plongements. 

Soit $\gamma\in \tilde{G}(F)$. La propri\'et\'e (viii) entra\^{\i}ne que

(1)  tout \'el\'ement de  $H(F)$ peut s'\'ecrire $\gamma^nzg$, avec $n\in {\mathbb Z}$, $z\in Z(H;F)$ et $g\in G(F)$. 

Puisque $\gamma$ appartient \`a $Z_{H}(\gamma;F)$, il en r\'esulte que 

(2)   $Z_{H}(\gamma;F)$ est le sous-groupe de $H(F)$ engendr\'e par $\gamma$, $Z(H;F)$ et $Z_{G}(\gamma;F)$.

On a aussi

(3) l'application naturelle $Z_{G}(\gamma;F)\backslash G(F)\to  Z_{H}(\gamma;F)\backslash H(F)$ est bijective. 

En effet, elle est \'evidemment injective. Les propri\'et\'es (1) et (2) entra\^{\i}nent sa surjectivit\'e. 
 
Il y a une bijection $M\mapsto M^{H}$ entre Levi de $G$ et Levi de $H$: $M^H$ est engendr\'e par $M$ et $Z(H)$; inversement, $M=G\cap M^H$. Puisque $\tilde{G}$ est \`a torsion int\'erieure, il y a aussi une bijection $M\mapsto \tilde{M}$ entre Levi de $G$ et espaces de Levi de $\tilde{G}$. On a simplement $\tilde{M}=M^H\cap \tilde{G}$. Il est clair que, pour tout espace de Levi $\tilde{M}$, le groupe $M^H$ et les m\^emes objets $D$, $q$, $d$ v\'erifient la proposition du paragraphe pr\'ec\'edent relativement \`a $\tilde{M}$. 

\bigskip

\subsection{Fonctions, int\'egrales orbitales, repr\'esentations}
De l'inclusion $\tilde{G}(F)\subset H(F)$ se d\'eduit un homomorphisme de restriction $res_{\tilde{G}}^H:C_{c}^{\infty}(H(F))\to C_{c}^{\infty}(\tilde{G}(F))$. Dans le cas o\`u $F$ est archim\'edien, ces espaces sont munis d'une topologie et cet homomorphisme est continu. On a donc un homomorphisme dual  qui, \`a une distribution sur $\tilde{G}(F)$, associe une distribution sur $H(F)$. Soit $\gamma\in \tilde{G}(F)$, munissons $Z_{G}(\gamma;F)\backslash G(F)$ d'une mesure. A ces donn\'ees est associ\'ee une distribution sur $\tilde{G}(F)$, qui \`a $f\in C_{c}^{\infty}(\tilde{G}(F))$ associe l'int\'egrale orbitale
$$\int_{Z_{G}(\gamma;F)\backslash G(F)}f(g^{-1}\gamma g)\,dg.$$

D'apr\`es 2.3(2), son image dans l'espace des distributions sur $H(F)$ est l'int\'egrale orbitale sur $H(F)$ associ\'ee \`a $\gamma$ et la mesure sur $Z_{H}(\gamma;F)\backslash H(F)$ transport\'ee de celle fix\'ee sur $Z_{G}(\gamma;F)\backslash G(F)$ par l'isomorphisme entre ces deux quotients. Il en r\'esulte que l'homomorphisme $res_{\tilde{G}}^H$ se quotiente en un homomorphisme $res_{\tilde{G}}^H:I(H(F))\to I(\tilde{G}(F))$. En sens inverse, une distribution invariante sur $\tilde{G}(F)$ \`a support dans un nombre fini de classes de conjugaison par $G(F)$ s'envoie sur une distribution sur $H(F)$ \`a support dans un nombre fini de classes de conjugaison par $H(F)$. Autrement dit, on obtient un homomorphisme $res_{\tilde{G}}^{H,*}:D_{g\acute{e}om}(\tilde{G}(F))\to D_{g\acute{e}om}(H(F))$. 

On a vu qu'il \'etait plus canonique de consid\'erer les espaces $I(\tilde{G}(F))\otimes Mes(G(F))$ et $D_{g\acute{e}om}(\tilde{G}(F))\otimes Mes(G(F))^*$. La suite exacte
$$1\to G(F)\to H(F)\to D(F)\to 1$$
induit un isomorphisme $Mes(H(F))\simeq Mes(G(F))\otimes Mes(D(F))$.  On choisit une fois pour toutes une mesure de Haar sur $D(F)$. L'isomorphisme ci-dessus devient simplement un isomorphisme $Mes(H(F))\simeq Mes(G(F))$. On a aussi un isomorphisme dual $Mes(H(F))^*\simeq Mes(G(F))^*$. On peut voir les homomorphismes ci-dessus sous la forme
$$res_{\tilde{G}}^H:I(H(F))\otimes Mes(H(F))\to I(\tilde{G}(F))\otimes Mes(G(F)) ,$$
$$res_{\tilde{G}}^{H,*}:D_{g\acute{e}om}(\tilde{G}(F))\otimes Mes(G(F))^* \to D_{g\acute{e}om}(H(F))\otimes Mes(H(F))^*.$$

Pour un \'el\'ement $\gamma\in \tilde{G}_{reg}(F)$, la classe de conjugaison stable de $\gamma$ dans $\tilde{G}(F)$ est \'egale \`a la classe de conjugaison stable de $\gamma$ dans $H(F)$. La propri\'et\'e 2.3(3) implique d'ailleurs que, si $\dot{{\cal X}}(\gamma)$ est un ensemble de repr\'esentants des classes de conjugaison par $G(F)$ dans cette classe de conjugaison  stable, c'est aussi un ensemble de repr\'esentants des classes de conjugaison par $H(F)$. En choisissant des mesures comme ci-dessus, on voit que
$$S^{\tilde{G}}(\gamma,res_{\tilde{G}}^H(f))=S^H(\gamma,f)$$
pour tout $f\in C_{c}^{\infty}(H(F))$. Il en r\'esulte que l'homomorphisme de restriction se quotiente en un homomorphisme
$$res_{\tilde{G}}^H:SI(H(F))\otimes Mes(H(F))\to SI(\tilde{G}(F))\otimes Mes(G(F)).$$
On a un homomorphisme dual
$$res_{\tilde{G}}^{H,*}:D_{g\acute{e}om}^{st}(\tilde{G}(F))\otimes Mes(G(F))^*\to D_{g\acute{e}om}^{st}(H(F))\otimes Mes(H(F))^*.$$
\bigskip

Les choses sont moins simples du c\^ot\'e spectral. Soit $(\pi,\tilde{\pi})$ une repr\'esentation $G(F)$-irr\'eductible de $\tilde{G}(F)$. C'est-\`a-dire que $\pi$ est une repr\'esentation admissible irr\'eductible de $G(F)$ dans un espace complexe $V_{\pi}$ et $\tilde{\pi}$ est une application de $\tilde{G}(F)$ dans le groupe des automorphismes lin\'eaires de $V_{\pi}$ telle que $\tilde{\pi}(g\gamma g')=\pi(g)\tilde{\pi}(\gamma)\pi(g')$ pour tous $g,g'\in G(F)$ et $\gamma\in \tilde{G}(F)$. Notons $\chi_{\pi}$ le caract\`ere central de $\pi$ et prolongeons-le en un caract\`ere $\chi_{\pi}^H$ de $Z(H;F)$. Fixons $\gamma_{0}\in \tilde{G}(F)$, notons $N$ l'ordre de $d$ dans le groupe fini $D(F)/q(Z(H;F))$. Alors $\gamma_{0}^N$ appartient \`a $Z(H;F)G(F)$ et on peut l'\'ecrire conform\'ement $\gamma_{0}^N=z_{0}g_{0}$. Le lemme de Schur implique qu'il existe $c_{0}\in {\mathbb C}^{\times}$ tel que $\tilde{\pi}(\gamma_{0})^N =c_{0}\chi_{\pi}^H(z_{0})\pi(g_{0})$. Fixons une racine $N$-i\`eme  $c$ de $c_{0}$. Pour $h\in H(F)$, \'ecrivons $h=zg\gamma^n$, avec $z\in Z(H;F)$, $g\in G(F)$ et $n\in {\mathbb Z}$. Posons $\pi^H(h)=\chi_{\pi}^H(z)\pi(g)(c^{-1}\tilde{\pi}(\gamma))^n$. On v\'erifie que cela ne d\'epend pas de la d\'ecomposition choisie de $h$ et que l'application $\pi^H$ ainsi d\'efinie est une repr\'esentation admissible de $H(F)$ dans $V_{\pi}$. Elle est irr\'eductible puisque sa restriction $\pi$ \`a $G(F)$ l'est. Introduisons le groupe  localement compact $D(F)^{\vee}$ des caract\`eres unitaires de $D(F)$.    La th\'eorie de la dualit\'e pour les groupes ab\'eliens localement compacts nous dit que $Mes(D(F)^{\vee})$ est isomorphe \`a $Mes(D(F))^*$. Autrement dit, de  la mesure que l'on a fix\'ee sur $D(F)$ se d\'eduit une mesure duale $d\kappa$ sur $D(F)^{\vee}$. Fixons une mesure de Haar $dh$ sur $H(F)$, qui d\'etermine une telle  mesure $dg$ sur $G(F)$. De $\tilde{\pi}$, resp. $\pi^H$, se d\'eduit un caract\`ere-distribution $I^{\tilde{G}}(\tilde{\pi},.)$, resp. $I^H(\pi^H,.)$, sur $C_{c}^{\infty}(\tilde{G}(F))\otimes Mes(G(F))$, resp. sur $C_{c}^{\infty}(H(F))\otimes Mes(H(F))$. Pour $f\in C_{c}^{\infty}(H(F))$, on v\'erifie que l'int\'egrale de gauche ci-dessous est absolument convergente et que l'on a l'\'egalit\'e
$$(1) \qquad \int_{D(F)^{\vee}}I^H(\pi^H,f(\kappa\circ q)\otimes dh)\kappa(d)^{-1}\,d\kappa=I^{\tilde{G}}(res_{\tilde{G}}^H(f)\otimes dg).$$
   Comme on le sait, les caract\`eres-distributions sont associ\'es \`a des fonctions localement int\'egrables $trace\, \tilde{\pi}$ sur $\tilde{G}(F)$ et $trace\,\pi^H$ sur $H(F)$. Alors $trace\, \tilde{\pi}$ n'est autre que la restriction de $trace\,\pi^H$ \`a $\tilde{G}(F)$. 

Inversement, soit $\pi^H$ une repr\'esentation admissible irr\'eductible de $H(F)$. Utilisons la th\'eorie de Mackey appliqu\'ee \`a $H(F)$ et \`a son sous-groupe distingu\'e $G(F)Z(H;F)$, dont le quotient $D(F)/q(Z(H;F))$ est engendr\'e par l'image de $d$. Cette th\'eorie nous dit que, si la restriction de  $trace\,\pi^H$ n'est  pas identiquement nulle sur $\tilde{G}(F)$, alors la restriction $\pi$ de $\pi^H$ \`a $G(F)$ est irr\'eductible. Notons dans ce cas $\tilde{\pi}$ la restriction de $\pi^H$ \`a $\tilde{G}(F)$. Le couple $(\pi,\tilde{\pi})$ est une repr\'esentation $G(F)$-irr\'eductible de $\tilde{G}(F)$. Le proc\'ed\'e ci-dessus appliqu\'e \`a ce couple, en prenant pour  caract\`ere $\chi_{\pi}^H$ le caract\`ere central de $\pi^H$, reconstruit $\pi^H$. On obtient que l'application qui, \`a $trace\,\pi^H$, associe sa restriction \`a $\tilde{G}(F)$, est un homomorphisme surjectif
$$(2) \qquad D_{spec}(H(F))\to D_{spec}(\tilde{G}(F)).$$
Il est clair qu'il se restreint en un homomorphisme surjectif
$$D_{temp}(H(F))\to D_{temp}(\tilde{G}(F)),$$
les indices $temp$ signifiant que l'on se limite aux repr\'esentations temp\'er\'ees.  On a introduit en [W2] 2.12 le sous-espace $D_{ell}(\tilde{G}(F))\subset D_{temp}(\tilde{G}(F))$ engendr\'e par les caract\`eres de repr\'esentations elliptiques au sens d'Arthur. On a le sous-espace analogue $D_{ell}(H(F))\subset D_{temp}(H(F))$.

\ass{Lemme}{L'homomorphisme pr\'ec\'edent se restreint en un homomorphisme surjectif
$$D_{ell}(H(F))\to D_{ell}(\tilde{G}(F)).$$}

Preuve.  Soient $M$ un Levi semi-standard de $G$, $\sigma$ une repr\'esentation irr\'eductible et de la s\'erie discr\`ete de $M(F)$ et $(A,\gamma)\in {\cal N}^{\tilde{G}}(\sigma)$. C'est-\`a-dire que $A$ est un automorphisme unitaire de l'espace $V_{\sigma}$ de $\sigma$, $\gamma\in \tilde{G}(F)$ normalise $M$ et on a la relation
$$\sigma(ad_{\gamma}(x))\circ A=A\circ \sigma(x)$$
pour tout $x\in M(F)$. On note $W_{0}(\sigma)$ le groupe habituel de la th\'eorie des $R$-groupes (cf. [W2] 1.11) et on suppose $W_{0}(\sigma)=\{1\}$. On suppose aussi que l'automorphisme de ${\cal A}_{M}/{\cal A}_{G}$ d\'efini par $\gamma$ n'a pas de point fixe non nul. Fixons $P\in {\cal P}(M)$.  A l'aide de $(A,\gamma)$, on a d\'efini en [W2] 2.9 une repr\'esentation $(\pi,\tilde{\pi})$ de $\tilde{G}(F)$. La repr\'esentation $\pi$ n'est autre que l'induite $Ind_{P}^G(\sigma)$. Elle n'est pas irr\'eductible en g\'en\'eral. En la r\'ealisant dans son mod\`ele $V_{\pi}$ habituel, l'op\'erateur $\tilde{\pi}(\gamma)$ est le compos\'e des trois op\'erateurs

- $e\mapsto A\circ e$ de $V_{\pi}$ dans $V_{\pi'}$, o\`u $\pi'=Ind_{P}^G(\sigma\circ ad_{\gamma})$;

- $e\mapsto \partial_{P}(\gamma)^{1/2}e\circ ad_{\gamma}^{-1}$ de $V_{\pi'}$ dans $V_{\pi''}$, o\`u $\pi''=Ind_{ad_{\gamma}(P)}^G(\sigma)$; $\partial_{P}(\gamma)^{1/2}$ est un facteur de normalisation sans importance pour nous;

- l'op\'erateur d'entrelacement normalis\'e $R_{P\vert  ad_{\gamma}(P)}(\sigma):V_{\pi''}\to V_{\pi}$.

Le caract\`ere de $\tilde{\pi}$ appartient \`a $D_{ell}(\tilde{G}(F))$ (il peut \'eventuellement \^etre nul) et cet espace est engendr\'e par de tels caract\`eres. 

Remarquons que cette construction s'applique pour construire $D_{ell}(H(F))$, en posant simplement $\tilde{H}=H$. 

Notons $M^H$ le Levi de $H$ associ\'e \`a $M$. Fixons un caract\`ere unitaire $\chi$ de $Z(H;F)$ qui co\"{\i}ncide sur $Z(M;F)\cap Z(H;F)$ avec le caract\`ere central de $\sigma$. On prolonge $\sigma$ en une repr\'esentation encore not\'ee $\sigma$ de $Z(H;F)M(F)$ par $\sigma(zx)=\chi(z)\sigma(x)$ pour tous $z\in Z(H;F)$ et $x\in M(F)$. Posons $\sigma^H=Ind_{Z(H;F)M(F)}^{M^H(F)}(\sigma)$ que l'on r\'ealise dans son espace habituel $V_{\sigma^H}$. On d\'efinit un op\'erateur $A^H$ de $V_{\sigma^H}$ par
$$(A^Hf)(x)=Af(\gamma^{-1}x\gamma).$$
Il v\'erifie la relation
$$\sigma^H(ad_{\gamma}(x))\circ A^H=A^H\circ \sigma^H(x)$$
pour tout $x\in M^H(F)$. Fixons $\delta\in \tilde{M}(F)$, notons $l$ le plus petit entier strictement positif tel que $\sigma\circ (ad_{\delta})^l\simeq \sigma$. D'apr\`es la th\'eorie de Mackey, la repr\'esentation  $\sigma^H$ se d\'ecompose en une somme
 $$\Sigma\oplus (\Sigma\otimes (\kappa\circ q))\oplus...\otimes (\Sigma\otimes( \kappa\circ q)^{\frac{N}{l}-1})$$
 de repr\'esentations irr\'eductibles et deux \`a deux in\'equivalentes, o\`u $N$ est l'ordre de $d$ dans $D(F)/q(Z(H;F))$ et $\kappa$ est un  caract\`ere primitif de $D(F)/q(Z(H;F))$. Notons $P^H$ le sous-groupe parabolique de $H$  d\'eduit de $P$. On peut utiliser pour chaque composante $\Sigma\otimes (\kappa\circ q)^n$ les m\^emes facteurs de normalisation que pour $\sigma$ et d\'efinir ainsi l'op\'erateur $R_{P^H\vert  ad_{\gamma}(P^H)}(\Sigma\otimes (\kappa\circ q)^n)$. Ces op\'erateurs se regroupent en un op\'erateur $R_{P^H\vert  ad_{\gamma}(P^H)}(\sigma^H)$.  Posons $\pi^H=Ind_{P^H}^H(\sigma^H)$. On copie la d\'efinition ci-dessus pour d\'efinir un op\'erateur $\tilde{\pi}^H(\gamma)$ de l'espace $V_{\pi^H}$, puis une repr\'esentation $(\pi^H,\tilde{\pi}^H)$ de $\tilde{H}(F)=H(F)$ (on a par d\'efinition $\tilde{\pi}^H(x\gamma)=\pi^H(x)\tilde{\pi}^H(\gamma)$ pour tout $x\in H(F)$). Montrons que
 
 (3) l'image par l'homomorphisme (2) du caract\`ere de $\tilde{\pi}^H$ est $N$ fois le caract\`ere de $\tilde{\pi}$.
 
Pour la simplicit\'e de l'\'ecriture, on ne distingue pas les repr\'esentations de leurs espaces naturels.  Notons $\epsilon:\sigma^H\to \sigma$ l'\'evaluation $f\mapsto f(1)$.  Pour $\varphi\in \sigma^H$ et $n\in {\mathbb Z}$, d\'efinissons une fonction $\varphi_{n}$ sur $G(F)$ par $\varphi_{n}(g)=\epsilon\circ \varphi(\gamma^ng)$. Elle appartient \`a $Ind_{ad_{\gamma}^{-n}(P)}^G(\sigma\circ ad_{\gamma}^n)$. Puisque $H(F)$ est r\'eunion disjointe des $\gamma^nZ(H;F)G(F)$ pour $n=0,...,N-1$ et $M^H(F)$ est r\'eunion disjointe des $\gamma^nZ(H;F)M(F)$ pour les m\^emes $n$, on v\'erifie que l'application
 $$(4)\qquad\begin{array}{ccc}\pi^H=Ind_{P^H}^H(\sigma^H)&\to &\oplus_{n=0,...,N-1}  Ind_{ad_{\gamma}^{-n}(P)}^G(\sigma\circ ad_{\gamma}^n)\\ \varphi&\mapsto& (\varphi_{n})_{n=0,...N-1}\\ \end{array}$$
 est un isomorphisme. Il est \'equivariant pour les actions de $G(F)$. Pour tout $n$, d\'efinissons une application
 $$ \begin{array}{ccc}Ind_{ad_{\gamma}^{-n}(P)}^G(\sigma\circ ad_{\gamma}^n)&\to& Ind_{P}^G(\sigma)\\\varphi_{n}&\mapsto&\psi_{n}= R_{P\vert ad_{\gamma}^{-n}(P)}(\sigma)(A^{-n}\circ \varphi_{n})\\ \end{array}$$
 C'est un isomorphisme. En composant (4) avec ces isomorphismes, on obtient un isomorphisme $G(F)$-\'equivariant
$$ \begin{array}{ccc}\pi^H=Ind_{P^H}^H(\sigma^H)&\to &\oplus_{n=0,...,N-1}  Ind_{P}^G(\sigma)\\ \varphi&\mapsto& (\psi_{n})_{n=0,...N-1}\\ \end{array}$$ 
Un simple calcul montre que cet isomorphisme transporte l'op\'erateur $\tilde{\pi}^H(\gamma)$ de l'espace de gauche sur l'op\'erateur diagonal de l'espace de droite, dont chaque composante est $\tilde{\pi}(\gamma)$.  L'assertion (3) en r\'esulte. 

Montrons que:

(5) ou bien la repr\'esentation $(\pi^H,\tilde{\pi}^H)$ de $\tilde{H}(F)=H(F)$ est somme de repr\'esentations elliptiques,  ou bien son caract\`ere est nul.

Remarquons d'abord que, de la d\'ecomposition de $\sigma^H$ en composantes irr\'eductibles r\'esulte une d\'ecomposition en composantes pas forc\'ement irr\'eductibles
$$(6) \qquad \pi^H=Ind_{P^H}^H(\Sigma)\oplus...\oplus Ind_{P^H}^H(\Sigma\otimes (\kappa\circ q)^{\frac{N}{l}-1})).$$
Puisque $\sigma^H\circ ad_{\gamma}$ est isomorphe $\sigma^H$, sa composante irr\'eductible $\Sigma\circ ad_{\gamma}$ est isomorphe \`a une autre composante irr\'eductible $\Sigma\otimes (\kappa\circ q)^j$, avec $j\in \{0,...,\frac{N}{l}-1\}$. Par tensorisation, $(\Sigma\otimes (\kappa\circ q)^n)\circ ad_{\gamma}$ est isomorphe \`a $\Sigma\otimes (\kappa\circ q)^{j+n}$ pour tout $n\in \{0,...,\frac{N}{l}-1\}$. Supposons d'abord $j\not=0$. Puisque l'op\'erateur $A^H$ r\'ealise les isomorphismes ci-dessus, il permute sans point fixe l'ensemble des composantes irr\'eductibles de $\sigma^H$. Il r\'esulte de sa construction que l'op\'erateur $\tilde{\pi}^H(\gamma)$ permute sans point fixe les diff\'erentes composantes  du membre de droite de (6). Il en est de m\^eme de $\tilde{\pi}^H(x\gamma)$ pour tout $x\in H(F)$. Il est alors  clair que le caract\`ere de $\tilde{\pi}^H$ est nul. Supposons maintenant $j=0$. Pour la m\^eme raison, l'op\'erateur $\tilde{\pi}^H(\gamma)$ conserve chaque composante du membre de droite de (6) et $\tilde{\pi}$ se d\'ecompose en repr\'esentations agissant dans chaque composante. Il suffit de voir que chacune de ces sous-repr\'esentations est elliptique ou de caract\`ere nul. On ne perd rien \`a se limiter \`a la premi\`ere composante $Ind_{P^H}^H(\Sigma)$. La repr\'esentation $\Sigma$ est de la s\'erie discr\`ete. L'op\'erateur $A^H$ se restreint \`a $\Sigma$ en un op\'erateur $B$ qui v\'erifie 
$$\Sigma(ad_{\gamma}(x))\circ B= B\circ \Sigma(x)$$
pour tout $x\in M^H(F)$. Donc $(B,\gamma)\in {\cal N}^H(\Sigma)$. On voit que la restriction de $\tilde{\pi}^H(\gamma)$ \`a $Ind_{P^H}^H(\Sigma)$ s'obtient \`a partir du couple $(B,\gamma)$ par le m\^eme proc\'ed\'e rappel\'e plus haut qui construit les repr\'esentations elliptiques. Il suffit de montrer que le  $\Sigma$ et $\gamma$  v\'erifie les conditions requises plus haut, \`a savoir que l'action d\'eduite de $\gamma$ dans ${\cal A}_{M^H}/{\cal A}_{H}$ est sans point fixe non nul et que $W_{0}(\Sigma)=\{1\}$. La premi\`ere condition r\'esulte de l'hypoth\`ese sur $\gamma$ et de l'isomorphisme
${\cal A}_{M}/{\cal A}_{G}\simeq {\cal A}_{M^H}/ {\cal A}_{H}$. Supposons $W_{0}(\Sigma)\not=\{1\}$. Par d\'efinition de ce groupe, on peut trouver un Levi $L^H$ de $H$ contenant strictement $M^H$, tel que ${\cal A}_{M^H}/{\cal A}_{L^H}$ soit de dimension $1$ et tel que  la condition suivante soit v\'erifi\'ee. L'ensemble ${\cal P}^{L^H}(M^H)$ a deux \'el\'ements, disons $Q^H$ et $\bar{Q}^H$. Pour $\lambda\in {\cal A}_{M^H,{\mathbb C}}$, on d\'efinit $\Sigma_{\lambda}$ par tensorisation de $\Sigma$ avec le caract\`ere $x\mapsto e^{<H_{M^H}(x),\lambda>}$ de $M^H(F)$. On d\'efinit l'op\'erateur d'entrelacement usuel $J_{\bar{Q}^H\vert Q^H}(\Sigma _{\lambda})$, qui est m\'eromorphe en $\lambda$. Alors cet op\'erateur a un p\^ole en $\lambda=0$. Posons $L=L^H\cap G$, $Q=Q^H\cap L$, $\bar{Q}=\bar{Q}^H\cap L$.  Les op\'erateurs d'entrelacement  vivent dans le groupe d\'eriv\'e de $L^H$, a fortiori dans $L$. Puisque la restriction de $\Sigma$ \`a $M(F)$ se d\'ecompose en $\sigma\oplus...\oplus \sigma\circ ad_{\delta^{l-1}}$, il existe $j\in \{0,...,l-1\}$ tel que l'op\'erateur $J_{\bar{Q}\vert Q}(\sigma_{\lambda}\circ ad_{\delta^j})$ ait un p\^ole en $\lambda=0$.  Conjuguer par $\delta^j$ ne change pas cette propri\'et\'e. Donc $J_{\bar{Q}\vert  Q}(\sigma_{\lambda})$ a un p\^ole en $\lambda=0$. Mais alors $W_{0}(\sigma)\not=\{1\}$ contrairement \`a l'hypoth\`ese. Cela ach\`eve la preuve de (5).

Il r\'esulte de (3) et (5) que, si le caract\`ere de $\tilde{\pi}$ n'est pas nul, c'est l'image par l'homomorphisme (1) d'un \'el\'ement de $D_{ell}(H(F))$.

La r\'eciproque est similaire. On part cette fois d'un Levi semi-standard $M^H$ de $H$, d'une repr\'esentation $\Sigma$ de $M^H(F)$ irr\'eductible et de la s\'erie discr\`ete et d'un couple $(B,\gamma)\in {\cal N}^H(\Sigma)$. On suppose que l'automorphisme de ${\cal A}_{M^H}/{\cal A}_{H}$ d\'eduit de $\gamma$ n'a pas de point fixe non nul et que $W_{0}(\Sigma)=\{1\}$. On d\'eduit de ces donn\'ees une repr\'esentation $\tilde{\pi}^H$ de $\tilde{H}(F)=H(F)$ dans $\pi^H=Ind_{M^H}^H(\Sigma)$. L'espace $D_{ell}(H(F))$ est engendr\'e par les caract\`eres de telles repr\'esentations.   On ne change pas la repr\'esentation $\tilde{\pi}^H$ si l'on remplace le couple $(B,\gamma)$ par $(\Sigma(x)B,x\gamma)$ pour un $x\in M^H(F)$. Notons que cette op\'eration ne change pas l'automorphisme de ${\cal A}_{M^H}/{\cal A}_{H}$ d\'eduit de $\gamma$. Puisque $M^H(F)$ s'envoie surjectivement sur $D(F)$, on peut par un tel changement supposer $q(\gamma)=d$, autrement dit $\gamma\in \tilde{G}(F)$. On pose $M=M^H\cap G$ et on fixe $\delta\in M^H(F)\cap \tilde{G}(F)$. La th\'eorie de Mackey nous dit que la restriction de $\Sigma$ \`a $M(F)$ se d\'ecompose en une somme $\sigma_{1}\oplus...\oplus \sigma_{l}$ de repr\'esentations irr\'eductibles deux-\`a-deux non \'equivalentes.
Pour $n=1,...,l$, on note $\epsilon_{n}$ la projection sur la composante $\sigma _{n}$. Pour $\varphi\in \pi^H$ et $n\in\{1,...,l\}$, d\'efinissons une fonction $\varphi_{n}$ sur $G(F)$ par $\varphi_{n}(g)=\epsilon_{n}\circ \varphi(g)$. On v\'erifie que l'application
$$(7) \qquad \begin{array}{ccc}\pi^H&\to&\oplus_{n=0,...,l-1}Ind_{M}^G(\sigma_{n})\\\varphi&\mapsto&(\varphi_{n})_{n=1,...,l}\\ \end{array}$$
est un isomorphisme \'equivariant pour les actions de $G(F)$. L'op\'erateur $B$ permute les composantes $\sigma_{n}$. Par construction, l'op\'erateur $\tilde{\pi}^H(\gamma)$ permute conform\'ement les composantes de la d\'ecomposition (7). Il en est de m\^eme de $\tilde{\pi}^H(x\gamma)$ pour tout $x\in G(F)$, autrement dit de $\tilde{\pi}^H(\gamma')$ pour tout $\gamma'\in \tilde{G}(F)$. On veut montrer que l'image du caract\`ere de $\tilde{\pi}^H$ par l'homomorphisme (1) appartient \`a $D_{ell}(\tilde{G}(F))$. Cette image ne d\'epend que la restriction de ce caract\`ere \`a $\tilde{G}(F)$. Les composantes de (7) permut\'ees non trivialement ne contribuent pas \`a cette restriction. On peut donc se limiter aux  $n$ tels que $\sigma_{n}$  est conserv\'e par $B$. Pour un tel $n$, la restriction de $\tilde{\pi}^H$ \`a $\tilde{G}(F)$ conserve la composante $Ind_{M}^G(\sigma_{n})$ du membre de droite de (7). Il suffit de prouver que cette action de $\tilde{G}(F)$ dans cette composante est elliptique ou de trace nulle. Fixons un tel $n$, notons simplement $\sigma=\sigma_{n}$, $\pi=Ind_{M}^G(\sigma)$ et $A$ la restriction de $B$ \`a $\sigma$. La repr\'esentation $\sigma$ est de la s\'erie discr\`ete et on a $(A,\gamma)\in {\cal N}^{\tilde{G}}(\sigma)$. On voit que l'action de $\tilde{G}(F)$ dans $\pi$ est d\'eduite de $(A,\gamma)$ par le m\^eme proc\'ed\'e d\'ecrit au d\'ebut de la preuve. Il suffit de prouver que $(\sigma,\gamma)$ v\'erifie les conditions requises, \`a savoir que $W_{0}(\sigma)=\{1\}$ et que l'action sur ${\cal A}_{M}/{\cal A}_{G}$ d\'eduite de $\gamma$ est sans point fixe non nul. Ces deux propri\'et\'es r\'esultent comme dans la preuve de (5) des propri\'et\'es analogues de $(\Sigma,\gamma)$. Cela ach\`eve la preuve. $\square$

\bigskip

\subsection{Endoscopie}
Du c\^ot\'e dual, on a une suite exacte
$$1\to \hat{D}\to \hat{H}\stackrel{\hat{\iota}}{\to} \hat{G}\to 1$$
qui est \'equivariante pour les actions galoisiennes.  Cette suite s'\'etend en une suite exacte
$$1\to \hat{D}\to {^LH}\stackrel{^L\iota}{\to} {^LG}\to 1.$$
Soit ${\bf G}'=(G',{\cal G}',s)$ une donn\'ee endoscopique pour $(G,\tilde{G})$ (on oublie ${\bf a}$ qui est trivial). Fixons $s^{H}\in \hat{H}$ qui s'envoie sur $s$. Cet \'el\'ement est uniquement d\'etermin\'e modulo $\hat{D}$, a fortiori modulo $Z(\hat{H})$. On a une suite exacte
$$1\to \hat{D}\to \hat{H}_{s^H}\to \hat{G}_{s}\to 1.$$
Notons ${\cal H}'$ l'image inverse de ${\cal G}'$ dans $^LH$. Pour tout $w\in W_{F}$, fixons $g_{w}=(g(w),w)\in {\cal G}'$ tel que $ad_{g_{w}}$ agisse par $w_{G'}$ sur $\hat{G}'=\hat{G}_{s}$. Relevons $g_{w}$ en $h_{w}=(h(w),w)\in {\cal H}'$. Alors l'action $w\mapsto ad_{h_{w}}$ munit $\hat{H}_{s^H}$ d'une action de $W_{F}$, qui se quotiente en une action de $\Gamma_{F}$ qui pr\'eserve une paire de Borel \'epingl\'ee. On peut introduire un groupe reductif $H'$ d\'efini et quasi-d\'eploy\'e sur $F$ de sorte que $\hat{H}_{s^H}$, muni de l'action pr\'ec\'edente, s'identifie au groupe dual $\hat{H}'$ de $H'$. On a une suite exacte
$$1\to G'\to H'\stackrel{q'}{\to} D\to 1.$$ 
Pour $(g,w)\in {\cal G}'$, on a une \'egalit\'e $sgw(s)^{-1} =a(w)g$, o\`u $a:W_{F}\to Z(\hat{G})$ est un cocycle qui est un cobord. Il en r\'esulte qu'il existe un cocycle $a^H::W_{F}\to Z(\hat{H})$ de sorte que $s_{H}hw(s_{H})^{-1}=a^H(w)h$ pour tout $(h,w)\in {\cal H}'$. En notant ${\bf a}^H$ la classe de $a^H$, on voit que ${\bf H}'=(H',{\cal H}',s_{H})$ est une donn\'ee endoscopique pour le triplet $(H,\tilde{H}=H,{\bf a}^H)$ (nous noterons simplement ce triplet comme un couple $(H,{\bf a}^H)$). Evidemment, ${\bf a}^H$ appartient au noyau $Ker$ de l'homomorphisme
$$H^1(W_{F};Z(\hat{H}))\to H^1(W_{F};Z(\hat{G})),$$
ou encore \`a l'image de l'homomorphisme
$$H^1(W_{F};\hat{D})\to H^1(W_{F};Z(\hat{H})).$$
\bigskip
Inversement, soit ${\bf a}^H$ un \'el\'ement de $Ker$ et soit ${\bf H}'=(H',{\cal H}',s_{H})$ une donn\'ee endoscopique pour le couple $(H,{\bf a}^H)$. On a une injection $\hat{D}\subset Z(\hat{H})\subset Z(\hat{H}')$, d'o\`u une surjection $H'\to D$ dont le noyau est connexe. On note $G'$ ce noyau. On note ${\cal G}'$ la projection de ${\cal H}'$ dans $^LG$ et $s$ la projection de $s^H$ dans $\hat{G}$. Alors $(G',{\cal G}',s)$ est une donn\'ee endoscopique de $(G,\tilde{G})$. 

Il est assez clair que les correspondances ci-dessus se quotientent en des bijections entre l'ensemble des classes d'\'equivalence de donn\'ees endoscopiques pour $(G,\tilde{G})$ et la r\'eunion sur les \'el\'ements ${\bf a}^H\in Ker$ des ensembles de classes d'\'equivalence de donn\'ees endoscopiques pour $(H,{\bf a}^H)$. Cette bijection pr\'eserve l'ellipticit\'e: ${\bf G}'$ est elliptique si et seulement si ${\bf H}'$ l'est.

Soient ${\bf G}'$ et ${\bf H}'$ comme ci-dessus. On a construit l'espace endoscopique $\tilde{G}'=G'\times_{Z(G)}{\cal Z}(\tilde{G})$. Or ${\cal Z}(\tilde{G})\subset \tilde{G}$ s'envoie dans $H$ par $\tilde{\iota}$, plus pr\'ecis\'ement dans $Z(H)$. On a une injection $Z(H)\to Z(H')$. Puisque $G'$ s'envoie lui-aussi dans $H'$, on en d\'eduit une application naturelle $\tilde{\iota}':\tilde{G}'\to H'$. Elle est d\'efinie sur $F$. Par construction, son image est contenue dans l'image r\'eciproque de $d\in D(F)$ dans $H'$. Puisque cette image r\'eciproque est une unique classe \`a gauche modulo $G'$, l'image de $\tilde{\iota}'$ est exactement cette image r\'eciproque. Soient $\gamma\in \tilde{G}_{ss}(F)$ et $\delta\in \tilde{G}'_{ss}(F)$. On a d\'efini en [I] 1.10 la propri\'et\'e: $\gamma$ et $\delta$ se correspondent. Cette notion est relative aux donn\'ees ambiantes $\tilde{G}$ et $\tilde{G}'$. C'est-\`a-dire que, si on consid\`ere maintenant $\gamma$ comme un \'el\'ement de $H_{ss}(F)$ et $\delta$ comme un \'el\'ement de $H'_{ss}(F)$, on a une autre notion de correspondance relative aux donn\'ees ambiantes $H$ et $H'$. On v\'erifie qu'en fait, ces deux notions co\"{\i}ncident. Cela r\'esulte du fait qu'il y a une bijection \'evidente entre paires de Borel pour $G$, resp. $G'$, et paires de Borel pour $H$, resp. $H'$. 

\ass{Lemme}{Les conditions suivantes sont \'equivalentes:

(i) la donn\'ee ${\bf G}'$ est relevante;

(ii) $d$ appartient \`a $q'(H'(F))$;

(iii) $d$ appartient \`a $q'(H'(F))$ et   ${\bf a}^H=1$.}

Preuve. La donn\'ee ${\bf G}'$ est relevante si et seulement si $\tilde{G}'(F)$ n'est pas vide ([I] lemme 1.9). Puisque $\tilde{G}'(F)$ est l'ensemble des $h\in H'(F)$ tels que $q'(h)=d$, on obtient l'\'equivalence entre (i) et (ii). Evidemment, (iii) entra\^{\i}ne (ii). Supposons ${\bf G}'$ relevante.
 Soient $\gamma\in \tilde{G}_{ss}(F)$ et $\delta\in \tilde{G}_{ss}'(F)$ deux \'el\'ements qui se correspondent, avec $\gamma\in \tilde{G}_{reg}(F)$. Comme on vient de le dire, il se correspondent aussi pour les donn\'ees ambiantes $H$ et $H'$, et on a encore $\gamma\in \tilde{H}_{reg}(F)$. L'\'el\'ement ${\bf a}^H$ provient d'un \'el\'ement ${\bf a}^D\in H^1(W_{F};\hat{D})$. Ces \'el\'ements d\'eterminent des caract\`eres $\omega^H$ de $H(F)$ et $\omega^D$ de $D(F)$. On a $\omega^H=\omega^D\circ q$.  D'apr\`es [KS]
lemme 4.4.C, $\omega^H$ est trivial sur $H_{\gamma}(F)$. D'apr\`es 2.3(2), $\omega^H$ est trivial sur $\gamma$ et sur $Z(H;F)$. Donc $\omega^D$ est trivial sur $q(\gamma)=d$ et sur $q(Z(H;F))$. Il l'est aussi sur le groupe engendr\'e par $d$ et $q(Z(H;F))$, c'est-\`a-dire $D(F)$ tout entier. Donc $\omega^D=1$. Puisque $D$ est un tore, cela entra\^{\i}ne que ${\bf a}^D=1$, donc aussi ${\bf a}^H=1$. $\square$

Supposons ${\bf G}'$ relevante. Alors les objets $H'$, $D$, $d$, $q'$ et les plongements que l'on a d\'efinis de $G'$  et  $\tilde{G}'$ dans $H'$ v\'erifient la proposition 2.3 relativement \`a $(G',\tilde{G}')$. En effet, puisque $D(F)$ est engendr\'e par $q(Z(H;F))$ et $d$ et puisque $Z(H)\subset Z(H')$, $D(F)$ est a fortiori engendr\'e par $q'(Z(H';F))$ et $d$. 
Fixons des donn\'ees auxiliaires $H'_{1}$, $\tilde{H}'_{1}=H'_{1}$, $C_{1}$, $\hat{\xi}^H_{1}$, $\Delta_{1}^H$ pour ${\bf H}'$. On note $G'_{1}$ et $\tilde{G}'_{1}$ les images r\'eciproques de $G'$ et $\tilde{G}'$ dans $H'_{1}$. On a une suite exacte
$$1\to G'_{1}\to H'_{1}\to D\to 1$$
d'o\`u dualement
$$1\to \hat{D}\to {^LH}'_{1}\to {^LG}'_{1}\to 1.$$
Le plongement $\hat{\xi}^H_{1}:{\cal H}'\to {^LH}'_{1}$ se quotiente en un plongement
$$\hat{\xi}_{1}:{\cal G}'={\cal H}'/\hat{D}\to {^LG}'_{1}={^LH}'_{1}/\hat{D}.$$
Les donn\'ees $G'_{1}$, $\tilde{G}'_{1}$, $C_{1}$, $\hat{\xi}_{1}$ sont des donn\'ees auxiliaires pour ${\bf G}'$. Notons $\Delta_{1}$ la restriction de $\Delta_{1}^H$ aux couples $(\delta,\gamma)$ d'\'el\'ements qui se correspondent tels que $\delta\in \tilde{G}'_{1}(F)$, $\gamma\in \tilde{G}(F)$. Il est facile quoique fastidieux de v\'erifier que $\Delta_{1}$ est un facteur de transfert compl\'etant nos donn\'ees auxiliaires. On a un homomorphisme de restriction
$$res_{\tilde{G}'_{1}}^{H'_{1}}:C_{c,\lambda_{1}}^{\infty}(H'_{1}(F))\to C_{c,\lambda_{1}}^{\infty}(\tilde{G}'_{1}(F)),$$
ou encore
$$res_{\tilde{G}'_{1}}^{H'_{1}}:C_{c,\lambda_{1}}^{\infty}(H'_{1}(F))\otimes Mes(H'(F))\to C_{c,\lambda_{1}}^{\infty}(\tilde{G}'_{1}(F))\otimes Mes(G'(F)).$$
Comme en 2.4, il se quotiente en un homomorphisme
$$res_{\tilde{G}'_{1}}^{H'_{1}}:SI_{\lambda_{1}}(H'_{1}(F))\otimes Mes(H'(F))\to SI_{\lambda_{1}}(\tilde{G}'_{1}(F))\otimes Mes(G'(F)).$$
Du fait que $\Delta_{1}$ est la restriction de $\Delta_{1}^H$ r\'esulte que le diagramme suivant est commutatif
$$\begin{array}{ccc}C_{c}^{\infty}(H(F))\otimes Mes(H(F))&\stackrel{res_{\tilde{G}}^H}{\to}&C_{c}^{\infty}(\tilde{G}(F))\otimes Mes(G(F))\\ \,\,\downarrow transfert&&\downarrow transfert\\ SI_{\lambda_{1}}(H'_{1}(F))\otimes Mes(H'(F))&\stackrel{res_{\tilde{G}'_{1}}^{H'_{1}}}{\to}& SI_{\lambda_{1}}(\tilde{G}'_{1}(F))\otimes Mes(G'(F))\\ \end{array}$$
Si on fait varier les donn\'ees auxiliaires pour ${\bf H}'$, on voit que les applications $res_{\tilde{G}'_{1}}^{H'_{1}}$ se recollent en un homomorphisme
$$res_{{\bf G}'}^{{\bf H}'}:SI({\bf H}')\otimes Mes(H'(F))\to SI({\bf G}')\otimes Mes(G'(F)).$$
Le diagramme ci-dessus devient un diagramme commutatif
$$\begin{array}{ccc}C_{c}^{\infty}(H(F))\otimes Mes(H(F))&\stackrel{res_{\tilde{G}}^H}{\to}&C_{c}^{\infty}(\tilde{G}(F))\otimes Mes(G(F))\\ \,\,\downarrow transfert&&\downarrow transfert\\ SI({\bf H}')\otimes Mes(H'(F))&\stackrel{res_{{\bf G}'}^{{\bf H}'}}{\to}& SI({\bf G}')\otimes Mes(G'(F))\\ \end{array}$$

\bigskip

\subsection{L'application $\phi_{\tilde{M}}$}
On suppose de nouveau, et jusqu'\`a la fin de l'article, que $F$ est non-archim\'edien.

{\bf Remarque.} On reprend cette hypoth\`ese parce que nous allons travailler avec des objets que nous n'avons d\'efini que dans le cas non-archim\'edien. Mais s'il anticipe les d\'efinitions n\'ecessaires dans le cas archim\'edien, le lecteur verra que la suite de cette section vaut aussi dans ce cas. 
\bigskip

 On a la suite exacte
$$0\to {\cal A}_{G}\to {\cal A}_{H}\to {\cal A}_{D}\to 0.$$
On a introduit en  [II] 1.6 l'ensemble $\tilde{{\cal A}}_{\tilde{G},F}$ et l'application $\tilde{H}_{\tilde{G}}:\tilde{G}(F)\to \tilde{{\cal A}}_{\tilde{G},F}$. On peut identifier $\tilde{{\cal A}}_{\tilde{G},F}$ \`a l'image de $\tilde{G}(F)$ dans ${\cal A}_{H}$ par l'application $H_{H}$ (avec une double signification de la lettre $H$) et $\tilde{H}_{\tilde{G},F}$ \`a la restriction de cette application $H_{H}$ \`a $\tilde{G}(F)$. On a introduit en  [II] 1.6 les espaces $C_{ac}^{\infty}(\tilde{G}(F))$ et $I_{ac}(\tilde{G}(F))$. On voit que l'application lin\'eaire $res_{\tilde{G}}^H:C_{c}^{\infty}(H(F))\to C_{c}^{\infty}(\tilde{G}(F))$ se prolonge en une application lin\'eaire $res_{\tilde{G}}^H:C_{ac}^{\infty}(H(F))\to C_{ac}^{\infty}(\tilde{G}(F))$. Celle-ci se quotiente en une application lin\'eaire $res_{\tilde{G}}^H:I_{ac}(H(F))\to I_{ac}(\tilde{G}(F))$.

Soit $\tilde{M}$ un espace de Levi de $\tilde{G}$. On fixe un sous-groupe compact maximal sp\'ecial $K^H$ de $H(F)$ en bonne position relativement \`a $M^H$. Posons $K=G(F)\cap K^H$. C'est un sous-groupe compact maximal de $G(F)$ en bonne position relativement \`a $M$. On a d\'efini en  [W2] 6.4 une application lin\'eaire
$$\phi_{\tilde{M}}:C_{ac}^{\infty}(\tilde{G}(F))\otimes Mes(G(F))\to I_{ac}(\tilde{M}(F))\otimes Mes(M(F)).$$
Remarquons que la mesure fix\'ee sur $D(F)$ d\'etermine encore un isomorphisme $Mes(M^H(F))\simeq Mes(M(F))$. 
On a un diagramme
$$\begin{array}{ccc}C_{ac}^{\infty}(H(F))\otimes Mes(H(F))&\stackrel{res_{\tilde{G}}^H}{\to}& C_{ac}^{\infty}(\tilde{G}(F))\otimes Mes(G(F))\\ \phi_{M^H}\downarrow\,\,&&\,\,\downarrow \phi_{\tilde{M}}\\ I_{ac}(M^H(F))\otimes Mes(M^H(F))&\stackrel{res_{\tilde{M}}^{M^H}}{\to}&I_{ac}(\tilde{M}(F))\otimes Mes(M(F))\\ \end{array}$$

\ass{Lemme}{Le diagramme ci-dessus est commutatif.}

Preuve. Pour simplifier, on fixe sur  les groupes $G(F)$ et $M(F)$  des mesures de Haar compatibles. Il s'en d\'eduit des mesures de Haar sur $H(F)$ et $M^H(F)$. Cela nous d\'ebarrasse des espaces de mesures.   Soit $\pi^H$ une repr\'esentation temp\'er\'ee irr\'eductible de $M^H(F)$. Notons $\tilde{\pi}$ sa restriction \`a $\tilde{M}(F)$ et supposons que le caract\`ere de cette restriction n'est pas nul. Soit $X\in \tilde{{\cal A}}_{\tilde{M},F}\subset {\cal A}_{M^H}$. Soit $f\in C_{c}^{\infty}(H(F))$, posons $\varphi=res_{\tilde{G}}^H(f)$.  On a d\'efini en [W2] 6.4 les termes $J_{M^H}^H(\pi^H,X,f)$ et $J_{\tilde{M}}^{\tilde{G}}(\tilde{\pi},X,\varphi)$. On va les comparer. Notons $D(F)_{c}$ le plus grand sous-groupe compact de $D(F)$ et $D(F)_{c}^{\vee}$ son groupe des caract\`eres. On prolonge chaque \'el\'ement de $D(F)_{c}^{\vee}$ en un caract\`ere unitaire de $D(F)$. Pour chaque $\kappa\in D(F)_{c}^{\vee}$, on dispose de la repr\'esentation $\pi^H\otimes (\kappa\circ q)$ de $H(F)$. On note $mes(D(F)_{c})$ la mesure de $D(F)_{c}$, vu comme sous-groupe ouvert de $D(F)$. Montrons que

$$(1) \qquad J_{\tilde{M}}^{\tilde{G}}(\tilde{\pi},X,\varphi)=mes(D(F)_{c})^{-1}\sum_{\kappa\in D(F)_{c}^{\vee}}J_{M^H}^H(\pi^H\otimes(\kappa\circ q),X,f )\kappa(d)^{-1}.$$

Montrons d'abord que la somme est finie. Soit $P^H\in {\cal P}(M^H)$. Pour $\lambda\in i{\cal A}_{M^H}^*$, posons  $\Pi_{\lambda}^H=Ind_{P^H}^H(\pi^H_{\lambda})$. La fonction $f$ n'intervient dans la d\'efinition de 
 $J_{M^H}^H(\pi^H\otimes(\kappa\circ q),X,f )$ que via des op\'erateurs $(\Pi^H_{\lambda}\otimes (\kappa\circ q))(f)$. Notons $U$ le plus grand sous-groupe compact de $Z(H;F)$ et 
 $\chi$ la restriction \`a $U$ du caract\`ere central de $\pi^H$.  Les op\'erateurs ci-dessus ne d\'ependent que de la fonction $f_{*}$ sur $H(F)$ 
  d\'efinie par
$$f_{*}(h)=\int_{U}f(zh)\chi(z)\kappa\circ q(z)\,dz.$$
Fixons un sous-groupe   $U'\subset U$ ouvert et d'indice fini tel que $f$ et $\chi$ soient invariants par $U'$. Alors $f_{*}$ est nulle si $\kappa$ n'est pas trivial sur $q(U')$. Puisque $q(Z(H;F))$ est d'indice fini dans $D(F)$, $q(U')$ est d'indice fini dans $D(F)_{c}$. Il n'y a qu'un nombre fini de $\kappa$ triviaux sur $q(U')$, d'o\`u l'assertion de finitude.  

Par d\'efinition,
$$J_{\tilde{M}}^{\tilde{G}}(\tilde{\pi},X,\varphi)=\int_{i{\cal A}^*_{M,F}}J_{\tilde{M}}^{\tilde{G}}(\tilde{\pi}_{\tilde{\lambda}},\varphi)e^{-<\tilde{\lambda},X>}\,d\lambda.$$
Expliquons cette formule. Ici ${\cal A}^*_{M,F}=i{\cal A}_{M}^*/i{\cal A}_{M,F}^{\vee}$, o\`u $i{\cal A}_{M,F}^{\vee}$ est le sous-groupe des $\lambda\in i{\cal A}_{M}^*$ tels que $<\lambda,H_{M}(x)>\in 2i\pi{\mathbb Z}$ pour tout $x\in M(F)$. On a relev\'e tout \'el\'ement $\lambda\in i{\cal A}^*_{M,F}$ en un \'el\'ement $\tilde{\lambda}\in i{\cal A}^*_{M^H}$. L'expression ci-dessus ne d\'epend pas du rel\`evement choisi. Enfin, la mesure sur $i{\cal A}^*_{\tilde{M},F}$ est de masse totale $1$.  On a une suite exacte
$$0\to i{\cal A}_{D}^*\to i{\cal A}_{M^H}^*\to i{\cal A}_{M}^*\to 0$$
dont on d\'eduit une suite exacte
$$0\to i{\cal A}_{D,F}^*\to i{\cal A}_{M^H,F}^*\to i{\cal A}_{M,F}^*\to 0,$$
avec des notations imit\'ees des pr\'ec\'edentes. La formule ci-dessus se r\'ecrit
$$J_{\tilde{M}}^{\tilde{G}}(\tilde{\pi},X,\varphi)=\int_{i{\cal A}^*_{M^H,F}/i{\cal A}^*_{D,F}}J_{\tilde{M}}^{\tilde{G}}(\tilde{\pi}_{\lambda},\varphi)e^{-<\lambda,X>}\,d\lambda.$$

L'expression $J_{\tilde{M}}^{\tilde{G}}(\tilde{\pi}_{\lambda},\varphi)$ est construite \`a l'aide d'op\'erateurs d'entrelacement et de l'op\'erateur $\tilde{\Pi}_{\lambda}(\varphi)$, o\`u $\tilde{\Pi}_{\lambda}=Ind_{\tilde{P}}^{\tilde{G}}(\tilde{\pi}_{\lambda})$ pour un \'el\'ement fix\'e $\tilde{P}\in {\cal P}(\tilde{M})$. Les op\'erateurs pour les repr\'esentations induites de $\tilde{\pi}_{\lambda}$ \'etant les restrictions des m\^emes op\'erateurs pour les repr\'esentations induites de $\pi^H_{\lambda}$, la formule 2.4(1) se g\'en\'eralise en
$$J_{\tilde{M}}^{\tilde{G}}(\tilde{\pi}_{\lambda},\varphi)=\int_{D(F)^{\vee}}J_{M^H}^H(\pi^H_{\lambda},f(\kappa\circ q)) \kappa(d)^{-1}\,d\kappa.$$
Le groupe $D(F)^{\vee}$ s'identifie au produit  $D(F)_{c}^{\vee}\times i{\cal A}_{D,F}^*$. Si on munit  le groupe discret $D(F)_{c}^{\vee}$ de la mesure de comptage et le groupe $i{\cal A}_{D,F}^*$   de la mesure de masse totale $1$, l'identification ci-dessus envoie la mesure sur $D(F)^{\vee}$ sur $mes(D(F)_{c})^{-1}$ fois celle sur $D(F)_{c}^{\vee}\times i{\cal A}_{D,F}^*$. D'o\`u
$$J_{\tilde{M}}^{\tilde{G}}(\tilde{\pi}_{\lambda},\varphi)=mes(D(F)_{c})^{-1}\sum_{\kappa\in D(F)_{c}^{\vee}}\int_{i{\cal A}_{D,F}^*}J_{M^H}^H(\pi^H_{\lambda},f_{\mu}(\kappa\circ q)) \kappa(d)^{-1}e^{-<\mu,H_{D}(d)>}\,d\mu,$$
o\`u on a not\'e $f_{\mu}$ la fonction $x\mapsto f(x)e^{<\mu,H_{H}(x)>}$ sur $H(F)$.   La somme en $\kappa$ est finie pour la m\^eme raison que ci-dessus.  Le terme $J_{M^H}^H(\pi^H_{\lambda},f_{\mu}(\kappa\circ q))$ est construit \`a l'aide d'op\'erateurs d'entrelacement et de l'op\'erateur $\Pi^H_{\lambda}(f_{\mu}(\kappa\circ q))$, o\`u $\Pi^H_{\lambda}=Ind_{P^H}^H(\pi_{\lambda})$ pour un \'el\'ement fix\'e $P^H\in {\cal P}(M^H)$. Les op\'erateurs d'entrelacement vivent dans $G(F)$ et sont insensibles \`a la torsion par un caract\`ere se factorisant par $q$. On a aussi $\Pi^H_{\lambda}(f_{\mu}(\kappa\circ q))=(\Pi^H_{\lambda+\mu}\otimes(\kappa\circ q))(f)$ et $\Pi^H_{\lambda+\mu}\otimes(\kappa\circ q)$ n'est autre que $Ind_{P^H}^H((\pi^H\otimes(\kappa\circ q))_{\lambda+\mu})$. On obtient
$$J_{M^H}^H(\pi^H_{\lambda},f_{\mu}(\kappa\circ q)) =J_{M^H}^H((\pi^H\otimes(\kappa\circ q))_{\lambda+\mu},f) .$$
D'o\`u
$$J_{\tilde{M}}^{\tilde{G}}(\tilde{\pi},X,\varphi)=mes(D(F)_{c})^{-1}\sum_{\kappa\in D(F)_{c}^{\vee}}\int_{i{\cal A}^*_{M^H,F}/i{\cal A}^*_{D,F}}$$
$$\int_{i{\cal A}_{D,F}^*}J_{M^H}^H((\pi^H\otimes(\kappa\circ q))_{\lambda+\mu},f) \kappa(d)^{-1}e^{-<\mu,H_{D}(d)>}\,d\mu\, e^{-<\lambda,X>}\,d\lambda.$$
Puisque $X\in \tilde{{\cal A}}_{\tilde{M}}$, sa projection dans ${\cal A}_{D}$ est $H_{D}(d)$, donc $<\mu,H_{D}(d)>=<\mu,X>$ pour tout $\mu\in i{\cal A}_{D,F}^*$. La double int\'egrale ci-dessus se recompose en une int\'egrale unique
$$J_{\tilde{M}}^{\tilde{G}}(\tilde{\pi},X,\varphi)=mes(D(F)_{c})^{-1}\sum_{\kappa\in D(F)_{c}^{\vee}}\kappa(d)^{-1}\int_{i{\cal A}^*_{M^H,F}}J_{M^H}^H((\pi^H\otimes(\kappa\circ q))_{\lambda},f) e^{-<\lambda,X>}\,d\lambda.$$
cette derni\`ere int\'egrale n'est autre que $J_{M^H}^H(\pi^H\otimes (\kappa\circ q),X,f)$, ce qui prouve (1). 

Un raisonnement facile, similaire \`a celui fait en [W2] 6.4, permet d'\'etendre la relation (1) \`a une fonction $f\in C_{ac}^{\infty}(H(F))$. Posons $\varphi_{\tilde{M}}=res_{\tilde{M}}^{M^H}(\phi_{M^H}(f))$. En appliquant (1) au cas $G=M$ et \`a la fonction $\phi_{M^H}(f))$, on obtient
$$(2) \qquad I^{\tilde{M}}(\tilde{\pi},X,\varphi_{\tilde{M}})=mes(D(F)_{c})^{-1}\sum_{\kappa\in D(F)_{c}^{\vee}}\kappa(d)^{-1}I^{M^H}(\pi^H\otimes(\kappa\circ q),X,\phi_{M^H}(f)).$$
Par d\'efinition de $\phi_{M^H}$, les membres de droite de (1) et (2) sont \'egaux. Donc aussi les membres de gauche. Par d\'efinition de $\phi_{\tilde{M}}$, cela signifie que $\varphi_{\tilde{M}}=\phi_{\tilde{M}}(\varphi)$, autrement dit $res_{\tilde{M}}^{M^H}(\phi_{M^H}(f))=\phi_{\tilde{M}}(res_{\tilde{G}}^H(f))$. $\square$

\bigskip

\subsection{Int\'egrales orbitales pond\'er\'ees \'equivariantes}

Soit $\tilde{M}$ un espace de Levi de $\tilde{G}$. On fixe un sous-groupe compact maximal sp\'ecial $K^H$ de $H(F)$ en bonne position relativement \`a $M^H$. Posons $K=G(F)\cap K^H$. C'est un sous-groupe compact maximal de $G(F)$ en bonne position relativement \`a $M$. Soient $f\in C_{c}^{\infty}(H(F))$ et $\gamma\in \tilde{M}(F)\cap \tilde{G}_{reg}(F)$. On d\'efinit comme en [II] 1.2 les int\'egrales orbitales pond\'er\'ees $J_{M^H}^H(\gamma,f)$ et
$J_{\tilde{M}}^{\tilde{G}}(\gamma,res_{\tilde{G}}^H(f))$. On suppose que les mesures sur $M_{\gamma}(F)\backslash G(F)$ et $M_{\gamma}^H(F)\backslash H(F)$ se correspondent par la bijection 2.3(3). On suppose aussi que les mesures sur ${\cal A}_{M}^G$ et ${\cal A}_{M^H}^H$ n\'ecessaires pour d\'efinir des int\'egrales orbitales pond\'er\'ees se correspondent via l'isomorphisme naturel entre ces deux espaces. On v\'erifie que les fonctions poids $v_{\tilde{M}}$ et $v_{M^H}$, qui sont d\'efinies sur ces deux quotients, se correspondent par la bijection. Il en r\'esulte que les deux int\'egrales ci-dessus sont \'egales.

 En appliquant la d\'efinition de [II] 1.6, le lemme du paragraphe pr\'ec\'edent permet d'en d\'eduire par r\'ecurrence la m\^eme \'egalit\'e des int\'egrales pond\'er\'ees invariantes
$$I_{\tilde{M}}^{\tilde{G}}(\gamma,res_{\tilde{G}}^H(f))=I_{M^H}^H(\gamma,f).$$
A partir de cette \'egalit\'e, les choix effectu\'es ci-dessus de sous-groupes compacts maximaux n'ont plus d'importance. Les choix de mesures disparaissent aussi:   on a  l'\'egalit\'e
$$I_{\tilde{M}}^{\tilde{G}}(\boldsymbol{\gamma},res_{\tilde{G}}^H({\bf f}))=I_{M^H}^H(res_{\tilde{M}}^{M^H,*}(\boldsymbol{\gamma}),{\bf f}),$$
pour tout $\boldsymbol{\gamma}\in  D_{g\acute{e}om}(\tilde{M}(F))\otimes Mes(M(F))^*\otimes Mes(D(F))^*$ \`a support form\'e d'\'el\'ements $\tilde{G}$-fortement r\'eguliers et tout ${\bf f}\in C_{c}^{\infty}(H(F))\otimes Mes(H(F))$. 

\bigskip

\subsection{Int\'egrales orbitales pond\'er\'ees stables}

On utilisera plus loin la propri\'et\'e suivante:

(1) soit $\varphi\in C_{c}^{\infty}(\tilde{G}(F))$ dont l'image dans $SI(\tilde{G}(F))$ est nulle; alors il existe $f\in C_{c}^{\infty}(H(F))$ dont l'image dans $SI(H(F))$ est nulle et telle que $\varphi=res_{\tilde{G}}^H(f)$.

Fixons un sous-espace $\mathfrak{s}$ d\'efini sur $F$ de $\mathfrak{z}(H)$  suppl\'ementaire de $\mathfrak{z}(G)$. Fixons un voisinage ouvert $\mathfrak{u}$ de $0$ dans $\mathfrak{s}(F)$. Si $\mathfrak{u}$ est assez petit, l'application
$$\begin{array}{ccc}\tilde{G}(F)\times \mathfrak{u}&\to&H(F)\\ (\gamma,X)&\mapsto& exp(X)\gamma\\ \end{array}$$
est un isomorphisme de $\tilde{G}(F)\times \mathfrak{u}$ sur un voisinage ouvert  ${\cal U}$ de $\tilde{G}(F)$ dans $H(F)$ invariant par conjugaison  et par conjugaison stable (si $\gamma\in {\cal U}$  est fortement $H$-r\'egulier, sa classe de conjugaison stable dans $H(F)$ est contenue dans ${\cal U}$). On fixe une fonction $\psi\in C_{c}^{\infty}(\mathfrak{s}(F))$ \`a support dans $\mathfrak{u}$ et telle que $\psi$ est constante de valeur $1$ dans un voisinage de $0$. On d\'efinit $f_{1}$ sur $\tilde{G}(F)\times \mathfrak{u}$ par $f_{1}(\gamma,X)=\psi(X)\varphi(\gamma)$. On transporte $f_{1}$ par l'isomorphisme ci-dessus en une fonction sur ${\cal U}$, que l'on prolonge par $0$ hors de ${\cal U}$ en une fonction $f$ sur $H(F)$. Cette fonction r\'epond \`a la question.

Soit $\tilde{M}$ un espace de Levi de $\tilde{G}$.  On a d\'efini $S_{\tilde{M}}^{\tilde{G}}(\boldsymbol{\delta},{\bf f})$ pour tout $\boldsymbol{\delta}\in D_{g\acute{e}om}^{st}(\tilde{M}(F))\otimes Mes(M(F))^*$ et tout ${\bf f}\in C_{c}^{\infty}(\tilde{G}(F))\otimes Mes(G(F))$. On sait aussi d\'efinir $S_{M^H}^H(\boldsymbol{\delta},{\bf f})$ pour tout $\boldsymbol{\delta}\in D_{g\acute{e}om}^{st}(M^H(F))\otimes Mes(M^H(F))^*$ et tout ${\bf f}\in C_{c}^{\infty}(H(F))\otimes Mes(H(F))$.  

\ass{Proposition}{Soit $\boldsymbol{\delta}\in D_{g\acute{e}om}^{st}(\tilde{M}(F))\otimes Mes(M(F))^*$. On suppose que le support de $\boldsymbol{\delta}$ est form\'e d'\'el\'ements fortement r\'eguliers dans $\tilde{G}(F)$.

(i)   Pour tout  ${\bf f}\in  C_{c}^{\infty}(H(F))\otimes Mes(H(F))$, on a l'\'egalit\'e
$$ S_{\tilde{M}}^{\tilde{G}}(\boldsymbol{\delta},res_{\tilde{G}}^H({\bf f}))=S_{M^H}^H(res_{\tilde{M}}^{M^H}(\boldsymbol{\delta}),{\bf f}).$$

(ii) La distribution $\boldsymbol{\varphi}\mapsto S_{\tilde{M}}^{\tilde{G}}(\boldsymbol{\delta},\boldsymbol{\varphi})$ est stable.}

Preuve. Posons $\boldsymbol{\varphi}=res_{\tilde{G}}^H({\bf f})$ et recopions la d\'efinition [II] 1.10(8):
$$(2) \qquad S_{\tilde{M}}^{\tilde{G}}(\boldsymbol{\delta},\boldsymbol{\varphi})=I_{\tilde{M}}^{\tilde{G}}(\boldsymbol{\delta},\boldsymbol{\varphi})-\sum_{s\in Z(\hat{M})^{\Gamma_{F}}/Z(\hat{G})^{\Gamma_{F}}, s\not=1}i_{\tilde{M}}(\tilde{G},\tilde{G}'(s))S_{{\bf M}}^{{\bf G}'(s)}(\boldsymbol{\delta},\boldsymbol{\varphi}^{{\bf G}'(s)}).$$
D'apr\`es 2.7, le premier terme est \'egal \`a $I_{M^H}^{H}(res_{\tilde{M}}^{M^H}(\boldsymbol{\delta}),{\bf f})$. Parce que $Z(\hat{M}_{ad})$ est un tore induit, $Z(\hat{M}_{ad})^{\Gamma_{F}}$ est connexe et l'homomorphisme naturel
$$Z(\hat{M})^{\Gamma_{F}}/Z(\hat{G})^{\Gamma_{F}}\to Z(\hat{M}_{ad})^{\Gamma_{F}}$$
est bijectif. Mais $\hat{M}_{ad}=\hat{M}^H_{ad}$. D'o\`u un isomorphisme
$$Z(\hat{M}^H)^{\Gamma_{F}}/Z(\hat{H})^{\Gamma_{F}}=Z(\hat{M})^{\Gamma_{F}}/Z(\hat{G})^{\Gamma_{F}}.$$
Pour un \'el\'ement $s$ de cet ensemble, on a une donn\'ee endoscopique ${\bf G}'(s)$ de $(G,\tilde{G})$ d\'eduite de ${\bf M}$ et de $s$ et on a une donn\'ee endoscopique ${\bf H}'(s)$ de $H$ d\'eduite de ${\bf M}^H$ et de $s$. 
On voit que la donn\'ee ${\bf H}'(s)$ est d\'eduite de ${\bf G}'(s)$ par la correspondance d\'efinie en 2.5. Que l'une des donn\'ees soit elliptique \'equivaut \`a ce que l'autre le soit. Montrons que

(3) $i_{\tilde{M}}(\tilde{G},\tilde{G}'(s))=i_{M^H}(H,H'(s))$.

On peut supposer les donn\'ees elliptiques, sinon les deux membres sont nuls. 
 On a la suite exacte
$$1\to \hat{D}\to \hat{H}'(s)\to \hat{G}'(s)\to 1.$$
Les groupes adjoints $\hat{H}'(s)_{AD}$ et $\hat{G}'(s)_{AD}$ sont \'egaux. L'image de $\hat{M}^H$ dans le premier est \'egale \`a celle de $\hat{M}$ dans le second. On a donc comme plus haut l'\'egalit\'e
$$Z(\hat{M}^H)^{\Gamma_{F}}/Z(\hat{H}'(s))^{\Gamma_{F}}=Z(\hat{M})^{\Gamma_{F}}/Z(\hat{G}'(s))^{\Gamma_{F}}.$$
Les homomorphismes
$$Z(\hat{M})^{\Gamma_{F}}/Z(\hat{G})^{\Gamma_{F}}\to Z(\hat{M})^{\Gamma_{F}}/Z(\hat{G}'(s))^{\Gamma_{F}}$$
et
$$Z(\hat{M}^H)^{\Gamma_{F}}/Z(\hat{H})^{\Gamma_{F}}\to Z(\hat{M}^H)^{\Gamma_{F}}/Z(\hat{H}'(s))^{\Gamma_{F}}$$
s'identifient. Puisque $i_{\tilde{M}}(\tilde{G},\tilde{G}'(s))$, resp. $i_{M^H}(H,H'(s))$, est l'inverse du nombre d'\'el\'ement du noyau du premier homomorphisme, resp. du second, (3) s'ensuit.

 Pour $s\not=1$, on peut admettre par r\'ecurrence la proposition que l'on cherche \`a prouver. Modulo quelques formalit\'es, elle affirme que
$$S_{{\bf M}}^{{\bf G}'(s)}(\boldsymbol{\delta},res_{{\bf G}'(s)}^{{\bf H}'(s)}({\bf f}^{{\bf H}'(s)}))=S_{{\bf M}^H}^{{\bf H}'(s)}(res_{\tilde{M}}^{M^H,*}(\boldsymbol{\delta}),{\bf f}^{{\bf H}'(s)}).$$
Comme on l'a dit en 2.5, on a l'\'egalit\'e $\boldsymbol{\varphi}^{{\bf G}'(s)}=res_{{\bf G}'(s)}^{{\bf H}'(s)}({\bf f}^{{\bf H}'(s)})$. L'\'egalit\'e pr\'ec\'edente devient
$$S_{{\bf M}}^{{\bf G}'(s)}(\boldsymbol{\delta},\boldsymbol{\varphi}^{{\bf G}'(s)})=S_{{\bf M}^H}^{{\bf H}'(s)}(res_{\tilde{M}}^{M^H,*}(\boldsymbol{\delta}),{\bf f}^{{\bf H}'(s)}).$$
Le membre de droite de (2) devient
$$I_{M^H}^{H}(res_{\tilde{M}}^{M^H}(\boldsymbol{\delta}),{\bf f})-\sum_{s\in Z(\hat{M}^H)^{\Gamma_{F}}/Z(\hat{H})^{\Gamma_{F}}, s\not=1}i_{M^H}(H,H'(s))S_{{\bf M}^H}^{{\bf H}'(s)}(res_{\tilde{M}}^{M^H,*}(\boldsymbol{\delta}),{\bf f}^{{\bf H}'(s)}),$$
ce qui n'est autre que $S_{M^H}^H(res_{\tilde{M}}^{M^H}(\boldsymbol{\delta}),{\bf f})$. Cela prouve le (i) de l'\'enonc\'e.

Soit $\boldsymbol{\varphi}\in C_{c}^{\infty}(\tilde{G}(F))\otimes Mes(G(F))$ dont l'image dans $SI(\tilde{G}(F))\otimes Mes(G(F))$ est nulle. D'apr\`es (1), on peut choisir ${\bf f}\in C_{c}^{\infty}(H(F))\otimes Mes(H(F))$ dont l'image dans $SI(H(F))\otimes Mes(H(F))$ est nulle et telle que $\boldsymbol{\varphi}=res_{\tilde{G}}^H({\bf f})$. On veut prouver que $S_{\tilde{M}}^{\tilde{G}}(\boldsymbol{\delta},\boldsymbol{\varphi})=0$. D'apr\`es (i), il suffit de prouver que $S_{M^H}^H(res_{\tilde{M}}^{M^H,*}(\boldsymbol{\delta}),{\bf f})=0$. Mais on est maintenant dans la situation d'un groupe non tordu et l'assertion a \'et\'e prouv\'ee par Arthur, cf. 1.1.$\square$

\bigskip

\subsection{Int\'egrales orbitales pond\'er\'ees endoscopiques}
Soit $\tilde{M}$ un espace de Levi de $\tilde{G}$.  

\ass{Proposition}{Soit $\boldsymbol{\gamma}\in D_{g\acute{e}om}(\tilde{M}(F))\otimes Mes(M(F))$. On suppose que le support de $\boldsymbol{\gamma}$ est form\'e d'\'el\'ements fortement r\'eguliers dans $\tilde{G}(F)$.

(i) Pour tout ${\bf f}\in C_{c}^{\infty}(H(F))\otimes Mes(H(F))$, on a l'\'egalit\'e
$$I_{\tilde{M}}^{\tilde{G},{\cal E}}(\boldsymbol{\gamma},res_{\tilde{G}}^H({\bf f}))=I_{M^H}^{H,{\cal E}}(res_{\tilde{M}}^{M^H,*}(\boldsymbol{\gamma}),{\bf f}).$$

(ii) Pour tout $\boldsymbol{\varphi}\in C_{c}^{\infty}(\tilde{G}(F))\otimes Mes(G(F))$, on a l'\'egalit\'e
$$I_{\tilde{M}}^{\tilde{G},{\cal E}}(\boldsymbol{\gamma},\boldsymbol{\varphi})=I_{\tilde{M}}^{\tilde{G}}(\boldsymbol{\gamma},\boldsymbol{ \varphi}).$$}

Preuve. On peut fixer une donn\'ee endoscopique ${\bf M}'=(M',{\cal M}',\zeta)$ de $(M,\tilde{M})$, qui est elliptique et relevante, et un \'el\'ement $\boldsymbol{\delta}\in D_{g\acute{e}om}^{st}({\bf M}')$ et supposer que $\boldsymbol{\gamma}$ est le transfert de $\boldsymbol{\delta}$. On a alors
$$I_{\tilde{M}}^{\tilde{G},{\cal E}}(\boldsymbol{\gamma},res_{\tilde{G}}^H({\bf f}))=I_{\tilde{M}}^{\tilde{G},{\cal E}}({\bf M}',\boldsymbol{\delta},res_{\tilde{G}}^H({\bf f})).$$
Il y a un cocycle $a:W_{F}\to Z(\hat{M})$ d\'efini par $\zeta m w_{M'}(\zeta)^{-1}=a(w)m$ pour tout $(m,w)\in {\cal M}'$. Comme toujours, on suppose qu'il prend ses valeurs dans $Z(\hat{G})$, ce qui est possible quitte \`a multiplier $\zeta$ par un \'el\'ement de $Z(\hat{M})$.  On a alors
$$(1) \qquad I_{\tilde{M}}^{\tilde{G},{\cal E}}(\boldsymbol{\gamma},res_{\tilde{G}}^H({\bf f}))=\sum_{s\in \zeta Z(\hat{M})^{\Gamma_{F}}/Z(\hat{G})^{\Gamma_{F}}}i_{\tilde{M}'}(\tilde{G},\tilde{G}'(s))S_{{\bf M}'}^{{\bf G}'(s)}(\boldsymbol{\delta},(res_{\tilde{G}}^H({\bf f}))^{{\bf G}'(s)}).$$
A partir de ${\bf M}'$, on construit une donn\'ee endoscopique ${\bf M}^{'H}=(M^{'H},{\cal M}^{'H},\zeta^H)$ de $M^H$. Le lemme 2.5 et l'hypoth\`ese de relevance de ${\bf M}'$ assurent que c'est bien une donn\'ee endoscopique pour le caract\`ere trivial de $M^H$.  De plus, le cocycle $a^H$ associ\'e \`a ces donn\'ees v\'erifie automatiquement la condition analogue \`a celle v\'erifi\'ee par $a$. Comme dans le paragraphe pr\'ec\'edent, la projection naturelle induit un isomorphisme
$$ \zeta^HZ(\hat{M}^H)^{\Gamma_{F}}/Z(\hat{H})^{\Gamma_{F}}\simeq \zeta Z(\hat{M})^{\Gamma_{F}}/Z(\hat{G})^{\Gamma_{F}}.$$
Identifions ces deux ensembles. Pour un \'el\'ement $s$ dans cet ensemble commun, la donn\'ee ${\bf H}'(s)$ se d\'eduit de ${\bf G}'(s)$ par le proc\'ed\'e de 2.5. L'une de ces donn\'ees est elliptique si et seulement si l'autre l'est. Montrons que

(2) $i_{\tilde{M}'}(\tilde{G},\tilde{G}'(s))=i_{M^{'H}}(H',H'(s))$.

L'argument est le m\^eme qu'en 2.8(3). On peut supposer les donn\'ees elliptiques sinon les deux membres sont nuls.  Le premier terme est le nombre d'\'el\'ements du noyau de l'homomorphisme
$$Z(\hat{M})^{\Gamma_{F}}/Z(\hat{G})^{\Gamma_{F}}\to Z(\hat{M}')^{\Gamma_{F}}/ Z(\hat{G}'(s))^{\Gamma_{F}}.$$
Le deuxi\`eme terme est le nombre d'\'el\'ements du noyau de l'homomorphisme
$$Z(\hat{M}^H)^{\Gamma_{F}}/Z(\hat{H})^{\Gamma_{F}}\to Z(\hat{M}^{_{'}H})^{\Gamma_{F}}/ Z(\hat{H}'(s))^{\Gamma_{F}}.$$
Mais ces homomorphismes s'identifient et (2) en r\'esulte.

Comme on l'a dit en 2.5, on a l'\'egalit\'e $(res_{\tilde{G}}^H({\bf f}))^{{\bf G}'(s)}=res_{{\bf G}'(s)}^{{\bf H}'(s)}({\bf f}^{{\bf H}'(s)})$. 
Remarquons que l'on peut supposer que $\boldsymbol{\delta}$ est \`a support dans l'ensemble des \'el\'ements semi-simples de $\tilde{M}'(F)$ qui correspondent \`a un \'el\'ement du support de $\boldsymbol{\gamma}$.  Alors le support de $\boldsymbol{\delta}$ est form\'e d'\'el\'ements qui sont fortement r\'eguliers dans $H'(s)$. Modulo quelques formalit\'es, la proposition 2.8(i) nous dit que
$$S_{{\bf M}'}^{{\bf G}'(s)}(\boldsymbol{\delta},(res_{\tilde{G}}^H({\bf f}))^{{\bf G}'(s)})=S_{{\bf M}^{'H}}^{{\bf H}'(s)}(res_{{\bf M}'}^{{\bf M}^{'H},*}(\boldsymbol{\delta}),{\bf f}^{{\bf H}'(s)}).$$
Le membre de droite de (1) devient
$$\sum_{s\in \zeta^HZ(\hat{M}^H)^{\Gamma_{F}}/Z(\hat{H})^{\Gamma_{F}}}i_{M^{'H}}(H,H'(s))S_{{\bf M}^{'H}}^{{\bf H}'(s)}(res_{{\bf M}'}^{{\bf M}^{'H},*}(\boldsymbol{\delta}),{\bf f}^{{\bf H}'(s)}).$$
Ceci n'est autre que $I_{M^H}^{H,{\cal E}}({\bf M}^{'H},res_{{\bf M}'}^{{\bf M}^{'H},*}(\boldsymbol{\delta}),{\bf f})$, ou encore $I_{M^H}^{H,{\cal E}}(transfert(res_{{\bf M}'}^{{\bf M}^{'H},*}(\boldsymbol{\delta})),{\bf f})$ . Mais, par dualit\'e \`a partir du dernier diagramme de 2.5, on a l'\'egalit\'e
$$transfert(res_{{\bf M}'}^{{\bf M}^{'H},*}(\boldsymbol{\delta}))=res_{\tilde{M}}^{M^H,*}(transfert(\boldsymbol{\delta}))=res_{\tilde{M}}^{M^H,*}(\boldsymbol{\gamma}).$$
Le membre de droite de (1) est donc \'egal \`a $I_{M^H}^{H,{\cal E}}(res_{\tilde{M}}^{M^H,*}(\boldsymbol{\gamma}),{\bf f})$ et cela d\'emontre le (i) de l'\'enonc\'e.

Pour $\boldsymbol{\varphi}\in C_{c}^{\infty}(\tilde{G}(F))\otimes Mes(G(F))$, on choisit ${\bf f}\in C_{c}^{\infty}(H(F))\otimes Mes(H(F))$ tel que $\boldsymbol{\varphi}=res_{\tilde{G}}^H({\bf f})$.  Dans la situation d'un groupe non tordu, on peut appliquer 1.1: on a l'\'egalit\'e
$$I_{M^H}^{H,{\cal E}}(res_{\tilde{M}}^{M^H,*}(\boldsymbol{\gamma}),{\bf f})=I_{M^H}^{H}(res_{\tilde{M}}^{M^H,*}(\boldsymbol{\gamma}),{\bf f}).$$
On a vu en 2.7 que le membre de droite \'etait \'egal \`a $I_{\tilde{M}}^{\tilde{G}}(\boldsymbol{\gamma},\boldsymbol{\varphi})$. Le (i) de l'\'enonc\'e nous dit que le membre de gauche est \'egal \`a $I_{\tilde{M}}^{\tilde{G},{\cal E}}(\boldsymbol{\gamma},\boldsymbol{\varphi})$. D'o\`u le (ii) de l'\'enonc\'e. $\square$

\section{Passage \`a un rev\^etement}

\bigskip

\subsection{ D\'efinition des homomorphismes de passage}
On fixe pour toute la section un triplet $(G,\tilde{G},{\bf a})$ tel que $\tilde{G}=G$. Mais ${\bf a}$ est quelconque. On consid\`ere un sous-tore $Z\subset Z(G)$  et un groupe r\'eductif connexe $G_{\sharp}$. On suppose donn\'e un homomorphisme $q:G_{\sharp}\to G$. Ces trois donn\'ees sont d\'efinies sur $F$.  On pose $G_{\flat}=Z\times G_{\sharp}$ et on prolonge $q$ par l'identit\'e sur $Z$. On obtient ainsi un homomorphisme encore not\'e $q:G_{\flat}\to G$. On suppose qu'il sinscrit dans une suite exacte
$$1\to \Xi_{\flat}\to G_{\flat}\stackrel{q}{\to} G\to 1$$
o\`u $\Xi_{\flat}$ est un sous-groupe fini central.  On note $\Xi$ la projection de $\Xi_{\flat}$ dans $G_{\sharp}$. Notons que $\Xi_{\flat}\to \Xi$ est bijective puisque $Z$ est inclus dans $G$.

{\bf Exemples.} On peut prendre $Z=Z(G)^0$ et $G_{\sharp}=G_{SC}$. Ou bien, supposons que $G$ soit un Levi d'un groupe $H$. On note $G_{sc}$ son image r\'eciproque dans $H_{SC}$. On peut prendre $Z=Z(H)^0$ et $G_{\sharp}=G_{sc}$.
\bigskip

On suppose que $\omega$ est trivial sur $q(G_{\flat}(F))$.  On fixe une mesure de Haar sur $Z(F)$. Il s'en d\'eduit une identification $Mes(G_{\sharp}(F))\simeq Mes(G_{\flat}(F))$. Puisque l'homomorphisme $q:G_{\flat}(F)\to G(F)$ est un isomorphisme local, on a aussi un isomorphisme $Mes(G_{\flat}(F))\simeq Mes(G(F))$: deux mesures se correspondent si elles se correspondent localement. D'o\`u aussi $Mes(G_{\sharp}(F))\simeq Mes(G(F))$.

L'action adjointe de $G(F)$ sur lui-m\^eme se remonte en une action de $G(F) $ sur $G_{\sharp}(F)$. Fixons un voisinage ouvert  $V_{\sharp}$ de $1$ dans $G_{\sharp}(F)$ invariant par cette action de $G(F)$ et tel que   $x\in V_{\sharp}$ si et seulement si la partie semi-simple de $x$ appartient \`a $V_{\sharp}$. On suppose $V_{\sharp}$ assez petit pour que $V_{\sharp}\cap \xi V_{\sharp}=\emptyset$ pour tout $\xi\in \Xi(F)-\{1\}$. On pose $V=q(Z(F)\times V_{\sharp})$.
   Alors $q$ se restreint en un isomorphisme de $Z(F)\times V_{\sharp}$ sur $V$. 
   
   Rappelons que l'on a d\'efini en [II] 1.6 les espaces $C_{ac}^{\infty}(G(F))$ et $I_{ac}(G(F),\omega)$. On note $C_{c}^{\infty}(V)$, resp $C_{ac}^{\infty}(V)$, $I(V,\omega)$, $I_{ac}(V,\omega)$, $D_{g\acute{e}om}(V,\omega)$ , l'espace des \'el\'ements de $C_{c}^{\infty}(G(F))$, resp $C_{ac}^{\infty}(G(F))$, $I(G(F),\omega)$, $I_{ac}(G(F),\omega)$,  $D_{g\acute{e}om}(G(F),\omega)$, \`a support dans $V$.

Le groupe $q(G_{\flat}(F))$ est un sous-groupe distingu\'e de $G(F)$, qui est ouvert  et d'indice fini. Fixons un ensemble ${\cal U}$ de repr\'esentants du quotient
$$q(G_{\sharp}(F))\backslash G(F).$$
Pour $f\in C_{c}^{\infty}(V)$ et $u\in {\cal U}$, on d\'efinit une fonction $(^{u}f)_{G_{\sharp}}$ sur $V_{\sharp}$ par $(^{u}f)_{G_{\sharp}}(x)=f(u^{-1}q(x)u)$ pour tout $x\in V_{\sharp}$. On d\'efinit une application lin\'eaire
$$\iota_{G_{\sharp},G}:C_{c}^{\infty}(V)\to C_{c}^{\infty}(V_{\sharp})$$
par $\iota_{G_{\sharp},G}(f)=\vert {\cal U}\vert ^{-1}\sum_{u\in {\cal U}}\omega(u)(^{u}f)_{G_{\sharp}}$. Elle d\'epend du choix de ${\cal U}$.

En sens inverse,  on  d\'efinit une application
$$\iota_{G,G_{\sharp}}:C_{c}^{\infty}(V_{\sharp})\to C_{ac}^{\infty}(G(F))$$
de la fa\c{c}on suivante. Pour $\varphi\in C_{c}^{\infty}(G_{\sharp}(F))$, $f=\iota_{G,G_{\sharp}}(\varphi)$ est la fonction sur $G(F)$ qui est nulle hors de $V$ et qui v\'erifie $f(zq(x_{\sharp}))=\varphi(x_{\sharp})$ pour tout $z\in Z(F)$ et $x_{\sharp}\in V_{\sharp}$. 

 Fixons un sous-tore maximal $T$ de $G$ et notons $T_{\sharp}$ son image r\'eciproque dans $G_{\sharp}$. De la m\^eme fa\c{c}on que ci-dessus, on a un isomorphisme $Mes(T(F))\simeq Mes(T_{\sharp}(F))$. Fixons des mesures de Haar sur $T_{\sharp}(F)$ et $G_{\sharp}(F)$, donc aussi sur $T(F)$ et $G(F)$. Alors, pour $t_{\sharp}\in T_{\sharp}(F)\cap G_{\sharp,reg}(F)$, l'int\'egrale orbitale $I^{G_{\sharp}}(t_{\sharp},.)$ est bien d\'efinie. C'est un \'el\'ement de $D_{g\acute{e}om}(G_{\sharp}(F))$. De m\^eme, pour $t\in T(F)\cap G_{reg}(F)$, l'int\'egrale orbitale $I^G(t,\omega,.)$ est bien d\'efinie. C'est un \'el\'ement de $D_{g\acute{e}om}(G(F),\omega)$.  Notons $c_{T}$ le nombre d'\'el\'ements de l'ensemble de doubles classes
$$q(G_{\flat}(F))\backslash G(F)/ T(F).$$

\ass{Lemme}{ Soient  $t_{\sharp}\in T_{\sharp}(F)\cap V_{\sharp} $  et $z\in Z(F)$. Posons $t=q(t_{\sharp})$ et supposons $t\in G_{reg}(F)$.  

(i)  Soit $f\in C_{c}^{\infty}(V)$, posons $\varphi=\iota_{G_{\sharp},G}(\varphi)$. On a l'\'egalit\'e 
$$c_{T}^{-1}I^G(t,\omega,f)=I^{G_{\sharp}}(t_{\sharp},\varphi).$$

(ii) Soit $\varphi\in C_{c}^{\infty}(V_{\sharp})$, posons $f=\iota_{G,G_{\sharp}}(\varphi)$. On a l'\'egalit\'e    
$$c_{T} \vert {\cal U}\vert ^{-1} \sum_{u\in {\cal U}}\omega(u)I^{G_{\sharp}}( ad_{u^{-1}}(t_{\sharp}),\varphi)=I^G(zt,\omega,f).$$}

Preuve de (i). On a par d\'efinition
$$I^{G_{\sharp}}(t_{\sharp},\varphi)\vert {\cal U}\vert ^{-1}\sum_{u\in {\cal U}}\omega(u)D^{G_{\sharp}}(t_{\sharp})^{1/2}\int_{T_{\sharp}(F)\backslash G_{\sharp}(F)}f(u^{-1}q(x^{-1}t_{\sharp}x)u)\,dx.$$
On voit que cette expression ne d\'epend pas du choix de ${\cal U}$. On peut fixer des ensembles de repr\'esentants ${\cal U}'$ du quotient $q(G_{\flat}(F))\backslash q(G_{\flat}(F))T(F)$ et ${\cal U}''$ du quotient $q(G_{\flat}(F))T(F)\backslash G(F)$ et supposer que ${\cal U}$ est l'ensemble des produits $u'u''$ avec $u'\in {\cal U}'$ et $u''\in {\cal U}''$. On peut de plus supposer que ${\cal U}'\subset T(F)$. On obtient
$$I^{G_{\sharp}}(t_{\sharp},\varphi)=\vert {\cal U}''\vert ^{-1}\sum_{u''\in {\cal U}''}\omega(u'')\vert {\cal U}'\vert ^{-1}\sum_{u'\in {\cal U}'}\omega(u')D^{G_{\sharp}}(t_{\sharp})^{1/2}$$
$$\int_{T_{\sharp}(F)\backslash G_{\sharp}(F)}f((u'')^{-1}(u')^{-1}q(x^{-1}t_{\sharp}x)u'u'')\,dx.$$
Pour $u'\in {\cal U}'\subset T(F)$, l'action $ad_{u'}^{-1}$ sur $G_{\sharp}(F)$ normalise $T_{\sharp}(F)$ et d\'efinit un automorphisme de  $T_{\sharp}(F)\backslash G_{\sharp}(F)$ qui pr\'eserve la mesure. D'autre part, cette action fixe $t_{\sharp}$. Par changement de variables, on voit que le terme $u'$ dispara\^{\i}t de l'int\'egrale int\'erieure. L'expression ci-dessus devient 
$$I^{G_{\sharp}}(t_{\sharp},\varphi)=d\vert {\cal U}''\vert ^{-1}\sum_{u''\in {\cal U}''}\omega(u'') D^{G_{\sharp}}(t_{\sharp})^{1/2}\int_{T_{\sharp}(F)\backslash G_{\sharp}(F)}f((u'')^{-1}q(x^{-1}t_{\sharp}x)u'')\,dx,$$
o\`u
$$d=\vert {\cal U}'\vert ^{-1}\sum_{u'\in {\cal U}'}\omega(u').$$
Si $\omega$ est non trivial sur $T(F)$, $d$ est nul et $I^{G_{\sharp}}(t_{\sharp},\varphi)=0$. Mais l'int\'egrale orbitale $I^G(t,\omega,f)$ est nulle elle aussi, d'o\`u l'\'egalit\'e voulue dans ce cas. Supposons que  $\omega$ est trivial sur $T(F)$. Alors $d=1$. Pour tout $u''\in {\cal U}''$, l'application
$$\begin{array}{ccc}T_{\sharp}(F)\backslash G_{\sharp}(F)&\to &T(F)\backslash G(F)\\ x&\mapsto& q(x)u''\\ \end{array}$$
est un isomorphisme de l'espace de d\'epart sur un ouvert ferm\'e de l'espace d'arriv\'ee. Il respecte les mesures par d\'efinition de celles-ci. Par d\'efinition de ${\cal U}''$ et parce que $Z\subset T$, $T(F)\backslash G(F)$ est r\'eunion disjointe des images de ces applications quand $u''$ d\'ecrit ${\cal U}''$. On obtient
$$I^{G_{\sharp}}(t_{\sharp},\varphi)=\vert {\cal U}''\vert ^{-1} D^{G_{\sharp}}(t_{\sharp})^{1/2}\int_{T(F)\backslash G(F)}f(y^{-1}ty)\omega(y)\,dy.$$
Il est clair que $D^{G_{\sharp}}(t_{\sharp})=D^G(t)$. Par d\'efinition, on a $c_{T}=\vert {\cal U}''\vert $. Alors la formule ci-dessus \'equivaut \`a
$$I^{G_{\sharp}}(t_{\sharp},\varphi)=c_{T}^{-1}I^G(t,\omega,f),$$
d'o\`u le (i) de l'\'enonc\'e.

Le (ii) se d\'emontre de fa\c{c}on analogue. $\square$ 

Ce lemme entra\^{\i}ne que les applications lin\'eaires d\'efinies ci-dessus se quotientent en des applications lin\'eaires
$$\iota_{G_{\sharp},G}:I(V,\omega)\to I(V_{\sharp})$$
et
$$\iota_{G,G_{\sharp}}:I(V_{\sharp})\to I_{ac}(V,\omega).$$
Ces applications ne d\'ependent pas du choix de l'ensemble ${\cal U}$.
Dualement, on a des applications lin\'eaires
$$\iota^*_{G_{\sharp},G}:D_{g\acute{e}om}(V_{\sharp})\to D_{g\acute{e}om}(V,\omega),$$
$$\iota^*_{G,G_{\sharp}}:D_{g\acute{e}om}(V,\omega)\to D_{g\acute{e}om}(V_{\sharp}).$$
D\'ecrivons plus compl\`etement ces applications. De l'action adjointe de $G(F)$ sur $G_{\sharp}(F)$ se d\'eduit une action de $G(F)$ sur $D_{g\acute{e}om}(V_{\sharp})$. Notons $D_{g\acute{e}om}(V_{\sharp})^{G(F),\omega}$ le sous-espace des $\boldsymbol{\gamma}\in D_{g\acute{e}om}(V_{\sharp})$ tels que $ad(g)(\boldsymbol{\gamma})=\omega(g)\boldsymbol{\gamma}$ pour tout $g\in G(F)$. Puisque cette action se quotiente en l'action du groupe fini $q(G_{\flat}(F))\backslash G(F)$, on a une projection naturelle
$$p:D_{g\acute{e}om}(V_{\sharp})\to D_{g\acute{e}om}(V_{\sharp})^{G(F),\omega}.$$ 
Notons d'autre part $D_{g\acute{e}om}(V,\omega)_{\sharp}$ le sous-espace des \'el\'ements de $D_{g\acute{e}om}(V,\omega)$ \`a support dans $q(V_{\sharp})$. Alors $\iota^*_{G_{\sharp},G}$ se factorise en
$$D_{g\acute{e}om}(V_{\sharp})\stackrel{p}{\to} D_{g\acute{e}om}(V_{\sharp})^{G(F),\omega}\stackrel{\iota^*_{G_{\sharp},G}}{\simeq}D_{g\acute{e}om}(V,\omega)_{\sharp}\subset D_{g\acute{e}om}(V,\omega).$$
En sens inverse, tout \'el\'ement $\boldsymbol{\gamma}\in D_{g\acute{e}om}(V,\omega) $ s'\'ecrit de fa\c{c}on unique $\boldsymbol{\gamma}=\sum_{z\in Z(F)}z\boldsymbol{\gamma}_{z}$, o\`u $\boldsymbol{\gamma}_{z}\in D_{g\acute{e}om}(V,\omega)_{\sharp}$ et $\boldsymbol{\gamma}_{z}=0$ pour presque tout $z$. Notons $\boldsymbol{\gamma}_{z,\sharp}$ l'\'el\'ement de $D_{g\acute{e}om}(V_{\sharp})^{G(F),\omega}$ tel que $\boldsymbol{\gamma}_{z}=\iota^*_{G_{\sharp},G}(\boldsymbol{\gamma}_{z,\sharp})$. Alors
$$\iota^*_{G,G_{\sharp}}(\boldsymbol{\gamma})=\sum_{z\in Z(F)}\boldsymbol{\gamma}_{z,\sharp}.$$

Remarquons que nos applications d\'ependent du choix de $V_{\sharp}$. Mais, pour deux tels voisinages, les applications relatives \`a chacun de ces voisinages co\"{\i}ncident sur l'intersection de leurs domaines de d\'efinition.
En particulier, $\iota_{G_{\sharp},G}^*$ se restreint en une application surjective
$$D_{unip}(G_{\sharp}(F))\to D_{unip}(G(F),\omega).$$
On note encore $\iota_{G_{\sharp},G}^*$ l'application obtenue en  tensorisant ces espaces  par les espaces de mesures ad\'equats. 

\bigskip

\subsection{Les termes $\rho_{J}$}

Soit $M$ un espace de Levi de $G$. On note $M_{\sharp}$, resp. $M_{\flat}$, son image r\'eciproque dans $G_{\sharp}$, resp. $G_{\flat}$. Remarquons que   l'application naturelle
$$q(M_{\flat}(F))\backslash M(F)\simeq q(G_{\flat}(F))\backslash G(F)$$
est bijective. Son injectivit\'e est imm\'ediate. Pour la surjectivit\'e, il suffit de traiter le cas o\`u $M$ est un Levi minimal de $G$. On fixe $P\in {\cal P}(M)$ et on note $P_{\flat}$ son image r\'eciproque dans $G_{\flat}$. Alors $(P_{\flat},M_{\flat})$ est aussi une paire parabolique minimale de $G_{\flat}$. 
Pour $g\in G(F)$, $ad_{g}(P_{\flat},M_{\flat})$ est encore une paire parabolique minimale de $G_{\sharp}$. Deux telles paires \'etant toujours conjugu\'ees par un \'el\'ement de $G_{\flat}(F)$, 
on peut multiplier $g$ \`a gauche par un \'el\'ement de $q(G_{\flat}(F))$ de sorte que $ad_{g} $ conserve $(P_{\flat},M_{\flat})$. Mais alors $ad_{g}$ conserve $(P,M)$ donc $g\in M(F)$. D'o\`u la surjectivit\'e requise.

En cons\'equence, on peut supposer que l'ensemble ${\cal U}$ de 3.1 est contenu dans $M(F)$. 
On voit que les applications lin\'eaires d\'efinies en 3.1 commutent au passage au terme constant. C'est-\`a-dire que, pour ${\bf f}\in I(V,\omega)\otimes Mes(G(F))$, on a l'\'egalit\'e
$$\iota_{M_{\sharp},M}({\bf f}_{M,\omega})=(\iota_{G_{\sharp},G}({\bf f})_{M_{\sharp}})$$
et, pour $\boldsymbol{\varphi}\in I(V_{\sharp}(F))\otimes Mes(G_{\sharp}(F))$, on a l'\'egalit\'e
$$\iota_{M,M_{\sharp}}(\boldsymbol{\varphi}_{M_{\sharp}})=(\iota_{G,G_{\sharp}}(\boldsymbol{\varphi}))_{M,\omega}.$$
On a des formules duales de commutation \`a l'induction.

L'application $q$ se restreint en une bijection entre les ensembles d'\'el\'ements unipotents de $M_{\sharp}(F)$ et de $M(F)$. On a aussi une bijection entre les ensembles de racines $\Sigma(A_{M_{\sharp}})$ et $\Sigma(A_{M})$. Enfin, de $q$ se d\'eduit un plongement $q:{\cal A}_{M_{\sharp}}\to {\cal A}_{M}$. Soit $u$ un \'el\'ement unipotent de $M_{\sharp}(F)$ et $\alpha\in \Sigma(A_{M_{\sharp}})$. On a d\'efini un \'el\'ement $\rho^{G_{\sharp}}(\alpha,u)\in {\cal A}_{M_{\sharp}}$ en [II] 1.4. Modulo l'identification ci-dessus, on a aussi un \'el\'ement $\rho^G(\alpha,q(u))\in {\cal A}_{M}$. On a

(1) $\rho^G(\alpha,q(u))=q(\rho^{G_{\sharp}}(\alpha,u))$.

Preuve. On n'en donne qu'une esquisse. On montre d'abord la m\^eme \'egalit\'e pour les termes primitifs d\'efinis par Arthur, c'est-\`a-dire l'\'egalit\'e $\rho^{G,Art}(\alpha,q(u))=q(\rho^{G_{\sharp},Art}(\alpha,u))$. Pour cela, on applique comme dans les preuves de [II] 1.4 la caract\'erisation de ces termes par les fonctions $W_{\omega}(a,\pi)$ de [A2]   (3.8). On compare ais\'ement ces fonctions pour $G$ et $G_{\sharp}$ et l'assertion en r\'esulte. Ensuite, \`a tout \'el\'ement $\alpha\in \Sigma(A_{M_{\sharp}})$ sont associ\'es deux groupes $G_{\alpha}$ et $G_{\sharp,\alpha}$. On voit que ces groupes sont reli\'es de la m\^eme fa\c{c}on que $G$ et $G_{\sharp}$, c'est-\`a-dire que l'on a une suite exacte
$$1\to \Xi_{\flat}\to Z\times G_{\alpha,\sharp}\stackrel{q}{\to}G_{\alpha}\to 1.$$
La d\'efinition inductive des termes $\rho^G(\alpha,q(u))$ et $\rho^{G_{\sharp}}(\alpha,u)$ conduit alors au r\'esultat. $\square$

De l'identification  ci-dessus entre ensembles de racines se d\'eduit une identification entre ${\cal J}_{M}^G$ et ${\cal J}_{M_{\sharp}}^{G_{\sharp}}$. Pour un \'el\'ement $J$ de cet ensemble, on a d\'efini un espace $U_{J}$ de germes de fonctions sur $A_{M}(F)$ et, de m\^eme, un espace $U_{J,\sharp}$ de germes de fonctions sur $A_{M_{\sharp}}(F)$. L'espace $U_{J,\sharp}$ est celui des fonctions $u\circ q$ pour $u\in U_{J}$. On identifie ainsi ces deux espaces.  On a d\'efini en [II] 3.2 des applications $\rho_{J}^G$ et $\rho_{J}^{G_{\sharp}}$. Dans leur d\'efinition interviennent des mesures sur ${\cal A}_{M}^G$ et ${\cal A}_{M_{\sharp}}^{G_{\sharp}}$. On suppose que ces mesures se correspondent par la bijection d\'eduite de $q$ entre ces espaces. On a

(2) le diagramme suivant est commutatif
$$\begin{array}{ccc}D_{unip}(M(F),\omega)\otimes Mes(M(F))^*&\stackrel{\rho_{J}^G}{\to}&U_{J}\otimes(D_{unip}(M(F),\omega)\otimes Mes(M(F))^*)/Ann^{G}\\ \quad \uparrow \iota^*_{M_{\sharp},M}&&\quad \uparrow \iota^*_{M_{\sharp},M}\\ D_{unip}(M_{\sharp}(F))\otimes Mes(M_{\sharp}(F))^*&\stackrel{\rho_{J}^{G_{\sharp}}}{\to}&U_{J}\otimes(D_{unip}(M_{\sharp}(F))\otimes Mes(M_{\sharp}(F))^*)/Ann^{G_{\sharp}}.\\ \end{array}$$

Cela r\'esulte des d\'efinitions des applications et de (1). 

{\bf Variante.} Supposons ${\bf a}=1$ et supposons donn\'ee une fonction $B$ comme en [II] 1.8. Cette fonction se remonte \`a $G_{\sharp}$. On sait alors d\'efinir les variantes $\rho^G(\alpha,u,B)$,  $\Sigma(A_{M},B)$ etc... des termes consid\'er\'es ci-dessus. Ces variantes v\'erifient les m\^emes propri\'et\'es. 

 \bigskip

\subsection{Int\'egrales orbitales pond\'er\'ees et rev\^etement}

Soit $M$ un Levi de $G$.
Rappelons que les int\'egrales orbitales pond\'er\'ees d\'ependent du choix d'un sous-groupe compact sp\'ecial $K$ de $G(F)$ en bonne position relativement \`a $M$. Si besoin est, on introduit ce compact dans la notation pour la pr\'eciser: $J_{M}^G(\boldsymbol{\gamma},\omega,{\bf f},K)$. Ce compact \'etant fix\'e, on note $K_{\sharp}$ son image r\'eciproque dans $G_{\sharp}(F)$. C'est  un sous-groupe compact sp\'ecial de $G_{\sharp}(F)$   (il est associ\'e au m\^eme point sp\'ecial de l'immeuble de $G_{\sharp,AD}=G_{AD}$). 
 Pour $g\in G(F)$, on pose $^{g}K_{\sharp}=ad_{g}(K_{\sharp})$. Les int\'egrales orbitales pond\'er\'ees d\'ependent aussi de mesures sur ${\cal A}_{M}^G$ et ${\cal A}_{M_{\sharp}}^{G_{\sharp}}$. Comme dans le paragraphe pr\'ec\'edent, on suppose que ces mesures se correspondent par l'isomorphisme d\'eduit de $q$ entre ces espaces.  

Fixons des mesures  sur nos diff\'erents groupes $G(F)$, $G_{\sharp}(F)$, etc... qui se correspondent comme en 3.1.

\ass{Lemme}{Pour tout $\boldsymbol{\gamma}_{\sharp}\in D_{g\acute{e}om}(V_{\sharp})$ \`a support $G_{\sharp}$-\'equisingulier ou unipotent et tout $f\in C_{c}^{\infty}(G(F))$, on a l'\'egalit\'e
$$J_{M}^G(\iota_{M_{\sharp},M}^*(\boldsymbol{\gamma}_{\sharp}),\omega,f)=\vert {\cal U}\vert ^{-1}\sum_{u\in {\cal U}}\omega(u)J_{M_{\sharp}}^{G_{\sharp}}(\boldsymbol{\gamma}_{\sharp},(^{u}f)_{G_{\sharp}},^{u}K_{\sharp}).$$}

  Preuve. On peut supposer que $\boldsymbol{\gamma}_{\sharp}$ est l'int\'egrale orbitale associ\'ee \`a un \'el\'ement $\gamma_{\sharp}\in M_{\sharp}(F)$  et \`a une certaine mesure sur $M_{\sharp,\gamma_{\sharp}}(F)$. Par un calcul similaire \`a celui de la preuve du lemme 3.1,  $\iota_{M_{\sharp},M}(\boldsymbol{\gamma}_{\sharp})$ est l'int\'egrale orbitale associ\'ee \`a l'\'el\'ement $\gamma=q(\gamma_{\sharp})$,  et \`a une certaine mesure sur $M_{\gamma}(F)$. Remarquons en passant que, si on remplace $\gamma_{\sharp}$ par $a\gamma_{\sharp}$, avec $a\in A_{M_{\sharp}}(F)$ et si on conserve la m\^eme mesure sur $M_{\sharp,a\gamma_{\sharp}}(F)=M_{\sharp,\gamma_{\sharp}}(F)$, la mesure d\'eduite sur $M_{q(a)\gamma}(F)=M_{\gamma}(F)$ ne change pas. Supposons d'abord que $\gamma_{\sharp}$ soit $G_{\sharp}$-\'equisingulier. Si $\omega$ n'est pas trivial sur $M_{\gamma}(F)$, on v\'erifie facilement que les deux membres de l'\'egalit\'e de l'\'enonc\'e sont nuls. Supposons que $\omega$ soit trivial sur $M_{\gamma}(F)$. Alors
$$J_{M}^G(\iota_{M_{\sharp},M}^*(\boldsymbol{\gamma}_{\sharp}),\omega, f)=J_{M}^G(\gamma,\omega,f)=D^G(\gamma)^{1/2}\int_{M_{\gamma}(F)\backslash G(F)}f(x^{-1}\gamma x)\omega(x)v_{M}^G(x)\,dx$$
$$=D^G(\gamma)^{1/2}D^M(\gamma)^{-1/2}\int_{M(F)\backslash G(F)}I^{M}(  \gamma,\omega,(^xf)_{M})\omega(x)v_{M}^G(x)\,dx,$$
o\`u $(^xf)_{M}$ est la fonction $y\mapsto f(x^{-1}yx)$ sur $M(F)$. Puisque l'int\'egrale orbitale $I^M(\gamma,\omega,.)$ est par d\'efinition l'image par $\iota_{M_{\sharp},M}^*$ de l'int\'egrale orbitale $I^{M_{\sharp}}(\gamma_{\sharp},.)$, la d\'efinition de $\iota_{M_{\sharp},M}^*$ entra\^{\i}ne
$$J_{M}^G(\iota_{M_{\sharp},M}^*(\boldsymbol{\gamma}_{\sharp}),\omega, f)=D^G(\gamma)^{1/2}D^M(\gamma)^{-1/2}\int_{M(F)\backslash G(F)}\vert {\cal U}\vert ^{-1}$$
$$\sum_{u\in {\cal U}}I^{M_{\sharp}}(\gamma_{\sharp},(^{ux}f)_{M_{\sharp}})\omega(ux)v_{M}^G(ux)\,dx,$$
avec une d\'efinition \'evidente de $(^{ux}f)_{M_{\sharp}}$.
Par d\'efinition de ${\cal U}$, cette expression se r\'ecrit
$$J_{M}^G(\iota_{M_{\sharp},M}^*(\boldsymbol{\gamma}_{\sharp}),\omega, f)=D^G(\gamma)^{1/2}D^M(\gamma)^{-1/2}\vert {\cal U}\vert ^{-1}$$
$$\int_{q(M_{\flat}(F))\backslash G(F)}I^{M_{\sharp}}(\gamma_{\sharp},(^{x}f)_{M_{\sharp}})\omega(x)v_{M}^G(x)\,dx.$$ 
L'application
$$\begin{array}{ccc}(M_{\sharp}(F)\backslash G_{\sharp}(F))\times {\cal U}&\to &q(M_{\flat}(F))\backslash G(F)\\ (x,u)&\mapsto &q(x)u\\ \end{array}$$
est bijective. On v\'erifie qu'elle pr\'eserve les mesures, si on met sur ${\cal U}$ la mesure de comptage. L'expression pr\'ec\'edente se r\'ecrit
$$J_{M}^G(\iota_{M_{\sharp},M}^*(\boldsymbol{\gamma}_{\sharp}),\omega, f)= D^G(\gamma)^{1/2}D^M(\gamma)^{-1/2}\vert {\cal U}\vert ^{-1}\sum_{u\in {\cal U}}$$
$$\int_{M_{\sharp}(F)\backslash G_{\sharp}(F)}I^{M_{\sharp}}(\gamma_{\sharp},(^{xu}f)_{M_{\sharp}})\omega(xu)v_{M}^G(xu)\,dx.$$
Pour $u\in {\cal U}$ et $x\in G_{\sharp}(F)$, on v\'erifie que $v_{M}^G(xu)$ est \'egal au poids $v_{M_{\sharp}}^{G_{\sharp}}(x,{^{u}K}_{\sharp})$ calcul\'e relativement au compact ${^{u}K}_{\sharp}$ (on utilise ici la compatibilit\'e entre les mesures sur ${\cal A}_{M_{\sharp}}^{G_{\sharp}} $ et sur ${\cal A}_{M}^G $).  On a de plus l'\'egalit\'e
$$D^G(\gamma)^{1/2}D^M(\gamma)^{-1/2}=D^{G_{\sharp}}(\gamma_{\sharp})^{1/2}D^{M_{\sharp}}(\gamma_{\sharp})^{-1/2}.$$
L'int\'egrale int\'erieure ci-dessus multipli\'ee par ce facteur devient $J_{M_{\sharp}}^{G_{\sharp}}(\boldsymbol{\gamma}_{\sharp},(^{u}f)_{G_{\sharp}},{^{u}K}_{\sharp})$ et on obtient la formule de l'\'enonc\'e.

Si $\gamma_{\sharp}$ est unipotent, on a une \'egalit\'e
$$J_{M}^G(\iota_{M_{\sharp},M}^*(\boldsymbol{\gamma}_{\sharp}),\omega,f)=lim_{a\to 1}\sum_{L\in {\cal L}(M)}r_{M}^L(\gamma,q(a))J_{L}^G(\iota_{M_{\sharp},M}^*(q(a)\boldsymbol{\gamma}_{\sharp}),\omega,f),$$
o\`u $a\in A_{M_{\sharp}}(F)$ est en position g\'en\'erale,
cf. [A1] 6.5 ou [II] 1.5. On a des formules similaires pour les termes du membre de droite de l'\'enonc\'e. Remarquons que, d'apr\`es une remarque faite plus haut, les mesures qui interviennent implicitement dans ces formules ne d\'ependent pas de $a$. En utilisant 3.2(1), on voit que $r_{M}^L(\gamma,q(a))=r_{M_{\sharp}}^{L_{\sharp}}(\gamma_{\sharp},a)$. L'assertion de l'\'enonc\'e  pour $\gamma_{\sharp}$ unipotent se d\'eduit alors par passage \`a la limite du cas particulier trait\'e pr\'ec\'edemment. $\square$

{\bf Remarque.} La proposition vaut en fait pour tout $\boldsymbol{\gamma}_{\sharp}\in D_{g\acute{e}om}(V_{\sharp}\cap M_{\sharp}(F))$ mais nous ne nous en servirons que pour les $\boldsymbol{\gamma}_{\sharp}$ indiqu\'es.

\bigskip

\subsection{Germes de Shalika et  rev\^etement }
On conserve les donn\'ees du paragraphe pr\'ec\'edent.   
\ass{Proposition}{(i) Pour $\boldsymbol{\gamma}_{\sharp}\in D_{g\acute{e}om,G_{\sharp}-\acute{e}qui}(M_{\sharp}(F))\otimes Mes(M_{\sharp}(F))^*$ assez voisin de l'origine, on a l'\'egalit\'e
$$g_{M,unip}^G\circ \iota^*_{M_{\sharp},M}(\boldsymbol{\gamma}_{\sharp})=\iota^*_{G_{\sharp},G}\circ g_{M_{\sharp},unip}^{G_{\sharp}}(\boldsymbol{\gamma}_{\sharp}).$$

(ii) Pour $\boldsymbol{\gamma}\in D_{g\acute{e}om,G-\acute{e}qui}(M(F))\otimes Mes(M(F))^*$ assez voisin de l'origine, on a l'\'egalit\'e
$$g_{M,unip}^G(\boldsymbol{\gamma})=g_{M,unip}^G\circ \iota_{M_{\sharp},M}^*\circ \iota^*_{M,M_{\sharp}}(\boldsymbol{\gamma}).$$}

Preuve. On fixe des mesures de Haar sur les groupes $G(F)$, $G_{\sharp}(F)$ etc... se correspondant comme en 3.1. Soient $\boldsymbol{\gamma}_{\sharp}\in D_{g\acute{e}om,G_{\sharp}-\acute{e}qui}(M_{\sharp}(F))$ et $f\in C_{c}^{\infty}(G(F))$. Si $\boldsymbol{\gamma}_{\sharp}$ est assez voisin de l'origine, on a l'\'egalit\'e
$$J_{M}^G( \iota^*_{M_{\sharp},M}(\boldsymbol{\gamma}_{\sharp}),\omega,f)=\sum_{L\in {\cal L}(M)}J_{L}^G(g_{M,unip}^L\circ \iota^*_{M_{\sharp},M}(\boldsymbol{\gamma}_{\sharp}),\omega,f).$$
Pour $L\not=G$, nos hypoth\`eses de r\'ecurrence nous permettent d'appliquer la proposition \`a prouver: on a
$$g_{M,unip}^L\circ \iota^*_{M_{\sharp},M}(\boldsymbol{\gamma}_{\sharp})=\iota^*_{L_{\sharp},L}\circ g_{M_{\sharp},unip}^{L_{\sharp}}(\boldsymbol{\gamma}_{\sharp}).$$
Pour $L=G$, notons $X$ la diff\'erence entre le membre de gauche de l'\'enonc\'e et celui de droite.   On obtient
$$J_{M}^G( \iota^*_{M_{\sharp},M}(\boldsymbol{\gamma}_{\sharp}),\omega,f)=I^G(X,\omega,f)+\sum_{L\in {\cal L}(M)}J_{L}^G( \iota^*_{L_{\sharp},L}\circ g_{M_{\sharp},unip}^{L_{\sharp}}(\boldsymbol{\gamma}_{\sharp}),\omega,f).$$
On applique la proposition 3.3 au membre de gauche et aux termes de la somme de droite. On obtient
$$I^G(X,\omega,f)=\vert {\cal U}\vert ^{-1}\sum_{u\in {\cal U}}\left(J_{M_{\sharp}}^{G_{\sharp}}(\boldsymbol{\gamma}_{\sharp},(^{u}f)_{\sharp},{^{u}K_{\sharp}})-\sum_{L_{\sharp}\in {\cal L}(M_{\sharp})}J_{L_{\sharp}}^{G_{\sharp}}(g_{M_{\sharp},unip}^{L_{\sharp}}(\boldsymbol{\gamma}_{\sharp}),{^{u}f}_{G_{\sharp}},{^{u}K}_{\sharp})\right).$$
Bien que les int\'egrales orbitales pond\'er\'ees d\'ependent du choix d'un sous-groupe compact sp\'ecial, les germes n'en d\'ependent pas. Cela r\'esulte ais\'ement des formules de passage entre int\'egrales orbitales pond\'er\'ees relatives \`a diff\'erents choix de sous-groupes compacts (et c'est ce qui permet aux germes pour les int\'egrales orbitales pond\'er\'ees non $\omega$-\'equivariantes d'\^etre aussi les germes pour leurs avatars $\omega$-\'equivariants).  Donc tous les termes de la somme en $u$ sont nuls. Donc aussi le membre de gauche, c'est-\`a-dire $X=0$.  Cela prouve  le (i) de la proposition. 

Il est clair que les germes sont insensibles aux centres, c'est-\`a-dire que l'on a l'\'egalit\'e $g_{M,unip}^G(z\boldsymbol{\gamma})=g_{M,unip}^G(\boldsymbol{\gamma})$ pour tout $\boldsymbol{\gamma}\in D_{g\acute{e}om,G-\acute{e}qui}(M(F))\otimes Mes(M(F))^*$ et $z\in Z(F)$, pourvu que $z$ et le support de $\boldsymbol{\gamma}$ soient assez proches de $1$. L'assertion (ii) r\'esulte alors de (i) et de la description explicite donn\'ee en 3.1 de l'application $\iota_{M,M_{\sharp}}^*$. $\square$
 
{\bf Variante.} Supposons $\omega=1$ et soit $B$ une fonction comme en [II].1.8. La m\^eme proposition vaut pour les germes $g_{M,unip}^{G}(.,B)$ et $g_{M_{\sharp},unip}^{G_{\sharp}}(.,B)$.  

\bigskip

\subsection{ Rev\^etement  et stabilit\'e}
On suppose que $G$ est quasi-d\'eploy\'e et que ${\bf a}=1$. On suppose que le voisinage $V_{\sharp}$ utilis\'e en 3.1 est invariant par conjugaison stable (si $x\in V_{\sharp}$ est fortement r\'egulier et si $y$ est stablement conjugu\'e \`a $x$, alors $y\in V_{\sharp}$). Notons qu'il existe de tels voisinages v\'erifiant de plus notre condition $V_{\sharp}\cap \xi V_{\sharp}=\emptyset$ pour tout $\xi\in \Xi(F)-\{1\}$ puisqu'un tel $\xi$ n'est pas stablement conjugu\'e \`a $1$. On a un analogue du lemme 3.1 pour les int\'egrales orbitales stables. Reprenons les hypoth\`eses de ce lemme.

\ass{Lemme}{Soient $t_{\sharp}\in T_{\sharp}(F)\cap V_{\sharp}$ et $z\in Z(F)$. Posons $t=q(t_{\sharp})$ et supposons   $t\in G_{reg}(F)$. Alors

(i) l'image de $S^{G_{\sharp}}(t_{\sharp},.)$ par l'application $\iota^*_{G_{\sharp},G}$ est $S^G(t,.)$;

(ii) l'image de $S^G( zt,.)$ par l'application $\iota^*_{G,G_{\sharp}}$ est
$ S^{G_{\sharp}}( t_{\sharp},.)$.}

Preuve.     Notons ${\cal O}_{\sharp}$ la classe de conjugaison stable de $t_{\sharp}$ dans $G_{\sharp}(F)$ et ${\cal O}$ celle de $t$ dans $G(F)$. Il est clair que $q({\cal O}_{\sharp})\subset {\cal O}$. Montrons que

(1)   la restriction de $q$ \`a ${\cal O}_{\sharp}$ est une bijection de ${\cal O}_{\sharp}$ sur ${\cal O}$.

 Par hypoth\`ese sur $V_{\sharp}$, on a ${\cal O}_{\sharp}\subset V_{\sharp}$ et on sait que $q$ est injective sur $V_{\sharp}$, a fortiori sur ${\cal O}_{\sharp}$. 
Soit   $t'\in {\cal O}$. On peut fixer $x\in G$ tel que $x^{-1}tx=t'$. Soit $x_{\sharp}\in G_{\sharp}$ ayant m\^eme image que $x$ dans $G_{AD}=G_{\sharp,AD}$. Pour tout $\sigma\in \Gamma_{F}$, on a la relation $x\sigma(x)^{-1}\in T$, donc $x_{\sharp}\sigma(x_{\sharp})^{-1}\in T_{\sharp}$. Il en r\'esulte que l'\'el\'ement $t'_{\sharp}=x_{\sharp}^{-1}t_{\sharp}x_{\sharp}$ appartient \`a  ${\cal O}_{\sharp}$. Et on a $t'=q(t'_{\sharp})$. D'o\`u la surjectivit\'e, ce qui prouve (1). 

  Fixons un ensemble de repr\'esentants $\dot{{\cal X}}$ des classes de conjugaison par $G(F)$ dans ${\cal O}$. Pour tout $x\in \dot{{\cal X}}$, fixons un ensemble de repr\'esentants $\dot{{\cal X}}_{x}$ des classes de conjugaison par $G_{\sharp}(F)$ dans l'image r\'eciproque par $q$ de la classe de conjugaison par $G(F)$ de $x$. Posons $\dot{{\cal X}}_{\sharp}=\sqcup_{x\in \dot{{\cal X}}}\dot{{\cal X}}_{x}$. C'est un ensemble de repr\'esentants des classes de conjugaison par $G_{\sharp}(F)$ dans ${\cal O}_{\sharp}$. Par d\'efinition, on a l'\'egalit\'e
$$S^{G_{\sharp}}(t_{\sharp},.)=\sum_{x_{\sharp}\in \dot{{\cal X}}_{\sharp}}I^{G_{\sharp}}(x_{\sharp},.)=\sum_{x\in \dot{{\cal X}}}\sum_{x_{\sharp}\in \dot{{\cal X}}_{x}}I^{G_{\sharp}}(x_{\sharp},.).$$
On applique le lemme 3.1(i). Par l'application $\iota^*_{G_{\sharp},G}$, la derni\`ere int\'egrale $I^{G_{\sharp}}(x_{\sharp},.)$ s'envoie sur $c_{T_{x}}^{-1}I^G(x,.)$, o\`u $T_{x}$ est le commutant de $x$ dans $G$. Donc
$$\iota^*_{G_{\sharp},G}(S^{G_{\sharp}}(t_{\sharp},.))=\sum_{x\in \dot{{\cal X}}}c_{T_{x}}^{-1}\vert \dot{{\cal X}}_{x}\vert I^G(x,.).$$
La classe de conjugaison de $x$ par $G(F)$ s'identifie \`a $T_{x}(F)\backslash G(F)$. Donc $\dot{{\cal X}}_{x}$ est un ensemble de repr\'esentants du quotient $T_{x}(F)\backslash G(F)/q(G_{\sharp}(F))=T_{x}(F)\backslash G(F)/q(G_{\flat}(F))$. D'o\`u $\vert \dot{{\cal X}}_{x}\vert =c_{T_{x}}$ et le membre de droite de l'\'egalit\'e ci-dessus devient
$$\sum_{x\in \dot{{\cal X}}}I^G(x,.),$$
ce qui n'est autre que $S^G(t,.)$. Cela prouve le (i) de l'\'enonc\'e. La preuve de (ii) est similaire. $\square$

 Il  r\'esulte de ce lemme que les applications de 3.1 se quotientent en des applications lin\'eaires
$$\iota_{G_{\sharp},G}:SI(V)\to SI(V_{\sharp})$$
et
$$\iota_{G,G_{\sharp}}:SI(V_{\sharp})\to SI(V),$$
avec des d\'efinitions \'evidentes de ces espaces.
On a dualement des applications lin\'eaires
$$\iota_{G_{\sharp},G}^*:D_{g\acute{e}om}^{st}(V_{\sharp})\to D_{g\acute{e}om}^{st}(V)$$
et
$$\iota_{G,G_{\sharp}}^*:D_{g\acute{e}om}^{st}(V)\to D_{g\acute{e}om}^{st}(V_{\sharp}).$$
  La description de 3.1 se simplifie pour les distributions stables: une distribution stable sur $V_{\sharp}$ est forc\'ement invariante par l'action du groupe adjoint $G_{\sharp,AD}(F)$, a fortiori par celle de $G(F)$. Il en r\'esulte que $\iota_{G_{\sharp},G}^*$ est injective et a pour image le sous-espace des \'el\'ements de $D^{st}_{g\acute{e}om}(V)$ \`a support dans $q(V_{\sharp})$. En particulier, cette application se restreint en un isomorphisme
$$\iota^*_{G_{\sharp},G}:D_{unip}^{st}(G_{\sharp}(F))\simeq D_{unip}^{st}(G(F)).$$
On peut aussi tensoriser les applications ci-dessus par des espaces de mesures. 

\bigskip

\subsection{Les termes $\sigma_{J}$}
On suppose $G$ quasi-d\'eploy\'e et ${\bf a}=1$. On suppose fix\'ee une fonction $B$ comme en [II] 1.8. Soit $M$ un Levi de $G$.  L'assertion suivante est similaire \`a 3.2(2).

\ass{Lemme}{Pour tout $J\in {\cal J}_{M}^G(B)$, le diagramme suivant est commutatif
$$\begin{array}{ccc}D_{unip}^{st}(M(F),\omega)\otimes Mes(M(F))^*&\stackrel{\sigma_{J}^G}{\to}&U_{J}\otimes(D_{unip}(M(F),\omega)\otimes Mes(M(F))^*)/Ann^{G}\\ \quad \uparrow \iota^*_{M_{\sharp},M}&&\quad \uparrow \iota^*_{M_{\sharp},M}\\ D_{unip}^{st}(M_{\sharp}(F))\otimes Mes(M_{\sharp}(F))^*&\stackrel{\sigma_{J}^{G_{\sharp}}}{\to}&U_{J}\otimes(D_{unip}(M_{\sharp}(F))\otimes Mes(M_{\sharp}(F))^*)/Ann^{G_{\sharp}}.\\ \end{array}$$}

{\bf Remarque.} Ainsi qu'on l'a dit ci-dessus, les applications verticales sont des isomorphismes.
\bigskip

Preuve. Dualement \`a la suite exacte
$$1\to \Xi_{\flat}\to Z\times G_{\sharp}\to G\to 1$$
on a une suite exacte
$$1\to \hat{\Xi}_{\flat}\to \hat{G}\to \hat{Z}\times \hat{G}_{\sharp}\to 1,$$
o\`u $\hat{\Xi}_{\flat}$ est un certain sous-groupe fini central. On a les isomorphismes
$$Z(\hat{M})^{\Gamma_{F}}/Z(\hat{G})^{\Gamma_{F}}\simeq Z(\hat{M}_{ad})^{\Gamma_{F}}\simeq Z(\hat{M}_{\sharp,ad})^{\Gamma_{F}}\simeq Z(\hat{M}_{\sharp})^{\Gamma_{F}}/Z(\hat{G}_{\sharp})^{\Gamma_{F}}.$$
En effet, les deux fl\`eches extr\^emes sont bijectives car les groupes  $Z(\hat{M}_{ad})^{\Gamma_{F}}$ et $Z(\hat{M}_{\sharp,ad})^{\Gamma_{F}}$ sont connexes. Celle du milieu est bijective car l'\'egalit\'e $\hat{G}_{AD}=\hat{G}_{\sharp,AD}$ entra\^{\i}ne $\hat{M}_{ad}=\hat{M}_{\sharp,ad}$. Ainsi, pour $s\in Z(\hat{M})^{\Gamma_{F}}/Z(\hat{G})^{\Gamma_{F}}$, on a \`a la fois une donn\'ee endoscopique ${\bf G}'(s)$ de $G$ et une donn\'ee endoscopique ${\bf G}'_{\sharp}(s)$ de $G_{\sharp}$. On a une suite exacte
$$1\to \hat{\Xi}_{\flat}\to \hat{G}'(s)\to \hat{Z}\times \hat{G}'_{\sharp}(s)\to 1.$$
Dualement, on a une suite exacte
$$1\to \Xi_{\flat}\to Z\times G'_{\sharp}(s)\stackrel{q_{s}}{\to} G'(s)\to 1.$$
C'est bien le groupe $\Xi_{\flat}$ qui intervient ici. En effet,  $Z\times G'_{\sharp}(s)$ a pour Levi $M_{\flat}=Z\times M_{\sharp}$, tandis que $G'(s)$ a pour Levi $M$. L'homomorphisme $ q_{s}$ se restreint en l'homomorphisme de d\'epart $M_{\flat}\to M$, dont le noyau est $\Xi_{\flat}$. Montrons que l'on a l'\'egalit\'e
$$(1) \qquad i_{M}(G,G'(s))=i_{M_{\sharp}}(G_{\sharp},G'_{\sharp}(s)).$$
Puisque ${\cal A}_{Z}\oplus{\cal A}_{G'_{\sharp}(s)}={\cal A}_{G'(s)}$ et ${\cal A}_{Z}\oplus{\cal A}_{G_{\sharp}}={\cal A}_{G}$, les donn\'ees ${\bf G}'(s)$ et ${\bf G}'_{\sharp}(s)$ sont simultan\'ement elliptiques ou non. Si elles ne le sont pas, les deux termes ci-dessus sont nuls. Supposons qu'elles soient elliptiques.  
Le m\^eme argument que ci-dessus fournit l'isomorphisme
$$Z(\hat{M})^{\Gamma_{F}}/Z(\hat{G}'(s))^{\Gamma_{F}}=Z(\hat{M}_{\sharp})^{\Gamma_{F}}/Z(\hat{G}'_{\sharp}(s))^{\Gamma_{F}}.$$
Le terme  $i_{M}(G,G'(s))$ est l'inverse du nombre d'\'el\'ements du noyau de l'homomorphisme
$$Z(\hat{M})^{\Gamma_{F}}/Z(\hat{G})^{\Gamma_{F}}\to Z(\hat{M})^{\Gamma_{F}}/Z(\hat{G}'(s))^{\Gamma_{F}}.$$
Le terme  $i_{M_{\sharp}}(G_{\sharp},G_{\sharp}'(s))$ est l'inverse du nombre d'\'el\'ements du noyau de l'homomorphisme
$$Z(\hat{M}_{\sharp})^{\Gamma_{F}}/Z(\hat{G}_{\sharp})^{\Gamma_{F}}\to Z(\hat{M}_{\sharp})^{\Gamma_{F}}/Z(\hat{G}'_{\sharp}(s))^{\Gamma_{F}}.$$
Ces homomorpismes s'identifient et (1) en r\'esulte.

Soient $\boldsymbol{\delta}_{\sharp}\in D_{unip}^{st}(M_{\sharp}(F))\otimes Mes(M_{\sharp}(F))^*$ et $a_{\sharp}\in A_{M_{\sharp}}(F)$ en position g\'en\'erale et proche de $1$. Posons $\boldsymbol{\delta}=\iota_{M_{\sharp},M}^*(\boldsymbol{\delta}_{\sharp})$ et $a=q(a_{\sharp})$. Posons 
$$X=\sigma_{J}^G(\boldsymbol{\delta},a)-\iota_{M_{\sharp},M}^*(\sigma_{J}^{G_{\sharp}}(\boldsymbol{\delta}_{\sharp},a_{\sharp})).$$
L'assertion de l'\'enonc\'e est que $X=0$.  On a 
$$\sigma_{J}^G(\boldsymbol{\delta},a)=\rho_{J}^G(\boldsymbol{\delta},a)-\sum_{s\in Z(\hat{M})^{\Gamma_{F}}/Z(\hat{G})^{\Gamma_{F}},s\not=1, J\in {\cal J}_{M}^{G'(s)}(B)}i_{M}(G,G'(s))\sigma_{J}^{G'(s)}(\boldsymbol{\delta},a).$$
On peut appliquer au premier terme la propri\'et\'e 3.2(2). Par r\'ecurrence, on peut appliquer le pr\'esent lemme aux autres termes. On obtient que $\sigma_{J}^G(\boldsymbol{\delta},a)$ est l'image par $\iota_{M_{\sharp},M}^*$ de
$$\rho_{J}^{G_{\sharp}}(\boldsymbol{\delta}_{\sharp},a_{\sharp})-\sum_{s\in Z(\hat{M})^{\Gamma_{F}}/Z(\hat{G})^{\Gamma_{F}},s\not=1, J\in {\cal J}_{M}^{G'(s)}(B)}i_{M}(G,G'(s))\sigma_{J}^{G_{\sharp}'(s)}(\boldsymbol{\delta}_{\sharp},a_{\sharp}).$$
Il est clair que $J\in {\cal J}_{M}^{G'(s)}(B)$ si et seulement si $J\in {\cal J}_{M_{\sharp}}^{G_{\sharp}'(s)}(B)$. Les consid\'erations qui pr\'ec\`edent transforment l'expression ci-dessus en
$$\rho_{J}^{G_{\sharp}}(\boldsymbol{\delta}_{\sharp},a_{\sharp})-\sum_{s\in Z(\hat{M}_{\sharp})^{\Gamma_{F}}/Z(\hat{G}_{\sharp})^{\Gamma_{F}},s\not=1, J\in {\cal J}_{M_{\sharp}}^{G_{\sharp}'(s)}(B)}i_{M_{\sharp}}(G_{\sharp},G_{\sharp}'(s))\sigma_{J}^{G_{\sharp}'(s)}(\boldsymbol{\delta}_{\sharp},a_{\sharp})$$
ce qui n'est autre que $\sigma_{J}^{G_{\sharp}}(\boldsymbol{\delta}_{\sharp},a_{\sharp})$. Cela prouve $X=0$ et le lemme. $\square$

\bigskip
\subsection{Rev\^etement  et  germes stables}

On  conserve les m\^emes hypoth\`eses.  Soit $M$ un Levi de $G$.

\ass{Proposition}{(i) Pour $\boldsymbol{\delta}_{\sharp}\in D_{g\acute{e}om,G_{\sharp}-\acute{e}qui}^{st}(M_{\sharp}(F))\otimes Mes(M_{\sharp}(F))^*$ assez voisin de l'origine, on a l'\'egalit\'e
$$Sg_{M,unip}^G( \iota^*_{M_{\sharp},M}(\boldsymbol{\delta}_{\sharp}),B)=\iota^*_{G_{\sharp},G}( Sg_{M_{\sharp},unip}^{G_{\sharp}}(\boldsymbol{\delta}_{\sharp},B)).$$

(ii) Pour $\boldsymbol{\delta}\in D_{g\acute{e}om,G-\acute{e}qui}^{st}(M(F))\otimes Mes(M(F))^*$ assez voisin de l'origine, on a l'\'egalit\'e
$$Sg_{M,unip}^G(\boldsymbol{\delta},B)=Sg_{M,unip}^G(\boldsymbol{\delta}',B),$$
o\`u $\boldsymbol{\delta}'=\iota_{M_{\sharp},M}^*\circ \iota_{M,M_{\sharp}}^*(\boldsymbol{\delta})$.}

Preuve. Posons $\boldsymbol{\delta}= \iota^*_{M_{\sharp},M}(\boldsymbol{\delta}_{\sharp})$. On utilise la d\'efinition de [II] 2.4(1):
$$(1) \qquad  Sg_{M,unip}^{G}(\boldsymbol{\delta},B)=g_{M,unip}^{G}(\boldsymbol{\delta},B)-\sum_{s\in Z(\hat{M})^{\Gamma_{F}}/Z(\hat{G})^{\Gamma_{F}}; s\not=1}i_{M}(G,G'(s))$$
$$transfert(Sg_{{\bf M},unip}^{{\bf G}'(s)}(\boldsymbol{\delta},B)).$$
  Fixons $s\in Z(\hat{M})^{\Gamma_{F}}/Z(\hat{G})^{\Gamma_{F}}$ avec $s\not=1$. Comme on l'a dit dans la preuve pr\'ec\'edente, on a \`a la fois une donn\'ee endoscopique ${\bf G}'(s)$ de $G$ et une donn\'ee endoscopique ${\bf G}'_{\sharp}(s)$ de $G_{\sharp}$. On a vu en [II] 1.10 que l'on pouvait choisir comme donn\'ees auxiliaires pour ${\bf G}'(s)$ le groupe $G'(s)_{1}=G'(s)$, le tore  $C(s)_{1}=\{1\}$ et un plongement $\hat{\xi}(s)_{1}:{\cal G}'(s)\to{^LG}'(s)$ de la forme $(x,w)\mapsto (x\chi(w),w)$ o\`u $\chi$ est un cocycle de $W_{F}$ dans $Z(\hat{M})$. On doit fixer de plus un facteur de transfert $\Delta(s)$. On a une projection $\hat{q}_{s}:\hat{G}'(s)\to \hat{G}'_{\sharp}(s)$ qui est \'equivariante pour les actions galoisiennes. D\'efinissons $\hat{\xi}_{\sharp}(s)_{1}:{\cal G}'_{\sharp}(s)\to {^LG}'_{\sharp}(s)$ par  $\hat{\xi}_{\sharp}(s)_{1}(x,w)=(x\hat{q}_{s}(\chi(w)),w)$. Alors les donn\'ees $G'_{\sharp}(s)_{1}=G'_{\sharp}(s)$, $C_{\sharp}(s)_{1}=\{1\}$ et $\hat{\xi}_{\sharp}(s)_{1}$ sont des donn\'ees auxiliaires pour ${\bf G}'_{\sharp}(s)$. Notons $q_{s}:G'_{\sharp}(s)\to G'(s)$ un homomorphisme dual de  $\hat{q}_{s}$. Pour deux \'el\'ements assez r\'eguliers $\delta\in G'_{\sharp}(s)(F)$ et $\gamma\in G_{\sharp}(F)$ qui se correspondent, posons $\Delta_{\sharp}(s)(\delta,\gamma)=\Delta(s)(q_{s}(\delta),q(\gamma))$. On v\'erifie que $\Delta_{\sharp}(s)$ est un facteur de transfert.
  
  {\bf Remarque.} Ce facteur n'est pas tout-\`a-fait l'image r\'eciproque de $\Delta(s)$. A cause du noyau $\Xi$, il y a des couples $(\delta,\gamma)$ d'\'el\'ements qui ne se correspondent pas mais pour lesquels $q_{s}(\delta)$ et $q(\gamma)$ se correspondent. Pour un tel couple, on a $\Delta(s)(q_{s}(\delta),q(\gamma))\not=0$ mais $\Delta_{\sharp}(s)(\delta,\gamma)=0$.
  
  \bigskip
  
  L'\'el\'ement $\boldsymbol{\delta}$ appartient au d\'epart \`a l'espace $D_{g\acute{e}om}^{st}(M(F))\otimes Mes(M(F))^*$. Dans la formule (1), on l'a identifi\'e \`a un \'el\'ement de $D_{g\acute{e}om}^{st}({\bf M})\otimes Mes(M(F))^*$ en utilisant pour la donn\'ee ${\bf M}$ les donn\'ees auxiliaires "triviales". Une  fois fix\'e $s$, on  utilise les donn\'ees auxiliaires pour ${\bf M}$ d\'eduites par restriction des donn\'ees auxiliaires pour ${\bf G}'(s)$. Alors $\boldsymbol{\delta}$ s'identifie de nouveau \`a un \'el\'ement de $D_{g\acute{e}om}^{st}(M(F))\otimes Mes(M(F))^*$, qui n'est pas en g\'en\'eral l'\'el\'ement de d\'epart et que l'on note $\boldsymbol{\delta}(s)$. De la m\^eme fa\c{c}on, l'\'el\'ement $\boldsymbol{\delta}_{\sharp}$ d\'etermine un \'el\'ement $\boldsymbol{\delta}_{\sharp}(s)\in D_{g\acute{e}om}^{st}(M_{\sharp}(F))\otimes Mes(M_{\sharp}(F))^*$. 
 Montrons que
  
  (2) $\iota_{M_{\sharp},M}^*(\boldsymbol{\delta}_{\sharp}(s))=\boldsymbol{\delta}(s)$.

   Pour $m\in M(F)$ assez r\'egulier, on a une \'egalit\'e $\Delta(s)(m,m)=c\underline{\chi}(m)$, o\`u $c$ est une constante non nulle et $\underline{\chi}$ est le caract\`ere de $M(F)$ d\'eduit de $\chi$ ( \`a moins que ce ne soit l'inverse, peu importe). Par construction, $\boldsymbol{\delta}(s)$ se d\'eduit de $\boldsymbol{\delta}$ par un automorphisme de $D_{g\acute{e}om}^{st}(M(F))\otimes Mes(M(F))^*$.  Celui-ci  envoie l'int\'egrale orbitale stable associ\'ee \`a un \'el\'ement $m\in M(F)$ assez r\'egulier sur la m\^eme int\'egrale orbitale stable, multipli\'ee par $c\underline{\chi}(m)$. Le caract\`ere $\underline{\chi}$ est localement constant. Il vaut $1$ sur les \'el\'ements assez voisins de l'origine. Restreint aux \'el\'ements \`a support proche de l'origine, l'automorphisme est donc l'homoth\'etie de rapport $c$. D'o\`u $\boldsymbol{\delta}(s)=c\boldsymbol{\delta}$. D'apr\`es la d\'efinition de $\Delta_{\sharp}(s)$, on a de m\^eme $\boldsymbol{\delta}_{\sharp}(s)=c\boldsymbol{\delta}_{\sharp}$ avec la m\^eme constante $c$. L'assertion (2) s'ensuit.   
  
  Par d\'efinition, on a l'\'egalit\'e
  $$Sg_{{\bf M},unip}^{{\bf G}'(s)}(\boldsymbol{\delta},B)=Sg_{M,unip}^{G'(s)}(\boldsymbol{\delta}(s),B).$$
 Puisque $s\not=1$, on peut appliquer le pr\'esent \'enonc\'e par r\'ecurrence. Gr\^ace \`a (2), on obtient
$$(3) \qquad Sg_{{\bf M},unip}^{{\bf G}'(s)}(\boldsymbol{\delta},B)=\iota_{G_{\sharp}'(s),G'(s)}^*(Sg_{M_{\sharp},unip}^{G'_{\sharp}(s)}(\boldsymbol{\delta}_{\sharp}(s),B)).$$

Montrons que 
 
 (4) on a l'\'egalit\'e $transfert\circ \iota _{G_{\sharp}'(s),G'(s)}^*=\iota_{G_{\sharp},G}^*\circ transfert$. 
 
Rappelons que nos applications  $\iota _{G_{\sharp}'(s),G'(s)}^*$ et $\iota_{G_{\sharp},G}^*$ sont d\'efinies au voisinage de l'unit\'e.
 Consid\'erons d'abord un \'el\'ement $t'_{\sharp}\in G'_{\sharp}(s;F)$ semi-simple assez r\'egulier et proche de $1$. Soit $t_{\sharp}\in G_{\sharp}(F)$ un \'el\'ement qui lui correspond. On note $T'_{\sharp}$ et $T_{\sharp}$ leurs commutants dans $G'_{\sharp}(s)$ et $G_{\sharp}$. On note $t'=q_{s}(t'_{\sharp})$ et $t=q(t_{\sharp})$ les images de $t'_{\sharp}$ et $t_{\sharp}$ dans $G'(s;F)$ et $G(F)$. On note  $T'$ et $T$ les commutants de  $t'$ et $t$ dans  $G'(s)$ et $G$. Les tores $T'_{\sharp}$ et $T_{\sharp}$ sont isomorphes, ainsi que les tores $T'$ et $T$. Fixons une mesure de Haar sur $T(F)$ que l'on transporte en une mesure sur $T'(F)$.  On fixe aussi des mesures de Haar sur $G(F)$ et $G'(s;F)$. On en d\'eduit comme en 3.1 des mesures sur $T_{\sharp}(F)$, $T'_{\sharp}(F)$, $G_{\sharp}(F)$ et $G'_{\sharp}(s)(F)$. 
Toutes ces mesures permettent de d\'efinir les int\'egrales orbitales qui interviennnent ci-dessous. En particulier, on a une int\'egrale orbitale stable $S^{G'_{\sharp}(s)}(t'_{\sharp},.)$.  Le lemme  3.5 montre que son image par $\iota^*_{G'_{\sharp}(s),G'(s)}$ est $S^{G'(s)}(t',.)$. Utilisons les notations introduites dans la preuve du lemme 3.5 pour les \'el\'ements $t_{\sharp}$ et $t$. L'image par transfert de $S^{G'(s)}(t',.)$ est
$$(5) \qquad \sum_{x\in \dot{{\cal X}}}\Delta(s)(t',x)I^G(x,\omega,.).$$
D'autre part, l'image par transfert de $S^{G'_{\sharp}(s)}(t'_{\sharp},.)$ est
$$\sum_{x_{\sharp}\in \dot{{\cal X}}_{\sharp}}\Delta_{\sharp}(s)(t'_{\sharp},x_{\sharp})I^{G_{\sharp}}(x_{\sharp},.), $$
ou encore
$$\sum_{x\in \dot{{\cal X}}}\sum_{x_{\sharp}\in \dot{{\cal X}}_{x}}\Delta_{\sharp}(s)(t'_{\sharp},x_{\sharp})I^{G_{\sharp}}(x_{\sharp},.).$$
En utilisant le lemme 3.1, l'image de cette expression par $\iota^*_{G_{\sharp},G}$ est
$$\sum_{x\in \dot{{\cal X}}} \sum_{x_{\sharp}\in \dot{{\cal X}}_{x}}\Delta_{\sharp}(s)(t'_{\sharp},x_{\sharp})c_{T_{q(x_{\sharp})}}^{-1}I^G(q(x_{\sharp}),\omega,.).$$
Pour $x_{\sharp}$ intervenant ci-dessus, on a $\Delta_{\sharp}(s)(t'_{\sharp},x_{\sharp})=\Delta(s)(t,q(x_{\sharp}))$ par d\'efinition. Puisque la fonction 
$$y\mapsto \Delta(s)(t,y)c_{T_{y}}^{-1}I^G(y,\omega,.)$$
 est invariante par conjugaison et puisque $q(x_{\sharp})$ est conjugu\'e \`a $x$, on peut remplacer $q(x_{\sharp})$ par $x$ dans l'expression ci-dessus. On obtient
 $$\sum_{x\in \dot{{\cal X}}} c_{T_{x}}^{-1}\vert \dot{{\cal X}}_{x}\vert \Delta(s)(t,x)I^G(x,\omega,.).$$
 On a vu dans la preuve du lemme 3.5 que $c_{T_{x}}=\vert \dot{{\cal X}}_{x}\vert $. Alors l'expression ci-dessus devient (5). 
Cela prouve l'\'egalit\'e (4) sur les distributions stables \`a support assez r\'egulier. Elle se g\'en\'eralise sans hypoth\`ese de support par bidualit\'e. En effet, puisque le transfert d'une fonction est d\'etermin\'e par  ses int\'egrales orbitales stables assez r\'eguli\`eres, cela implique  que la relation duale \`a (4) vaut pour les fonctions. Par dualit\'e, cela entra\^{\i}ne l'assertion (4) sans restriction sur le support des distributions. Cela prouve (4). 
  
  En utilisant (3), (4) et la proposition 3.4, on transforme l'expression (1) sous la forme suivante.  Le terme $Sg_{M,unip}(\boldsymbol{\delta},B)$ est l'image par $\iota^*_{G_{\sharp},G}$ de
  $$(5) \qquad g_{M_{\sharp},unip}^{G_{\sharp}}(\boldsymbol{\delta}_{\sharp},B)-\sum_{s\in Z(\hat{M})^{\Gamma_{F}}/Z(\hat{G})^{\Gamma_{F}},s\not=1}i_{M}(G,G'(s))transfert(Sg_{M_{\sharp},unip}^{G'_{\sharp}(s)}(\boldsymbol{\delta}_{\sharp}(s),B)).$$
  Comme ci-dessus, on a
  $$transfert(Sg_{M_{\sharp},unip}^{G'_{\sharp}(s)}(\boldsymbol{\delta}_{\sharp}(s),B))=transfert(Sg_{{\bf M}_{\sharp},unip}^{{\bf G}'_{\sharp}(s)}(\boldsymbol{\delta}_{\sharp},B)).$$
  Les m\^emes calculs que dans la preuve pr\'ec\'edente transforment  (5) en
  $$g_{M_{\sharp},unip}^{G_{\sharp}}(\boldsymbol{\delta}_{\sharp},B)-\sum_{s\in Z(\hat{M}_{\sharp})^{\Gamma_{F}}/Z(\hat{G}_{\sharp})^{\Gamma_{F}},s\not=1}i_{M_{\sharp}}(G_{\sharp},G_{\sharp}'(s))transfert(Sg_{{\bf M}_{\sharp},unip}^{{\bf G}'_{\sharp}(s)}(\boldsymbol{\delta}_{\sharp},B)).$$
  Ou encore, par une formule similaire \`a (1), en $Sg_{M_{\sharp},unip}^{G_{\sharp}}(\boldsymbol{\delta}_{\sharp},B)$. On a ainsi obtenu l'\'egalit\'e
  $$Sg^G_{M,unip}(\boldsymbol{\delta},B)=\iota^*_{G_{\sharp},G}(Sg_{M_{\sharp},unip}^{G_{\sharp}}(\boldsymbol{\delta}_{\sharp},B)).$$
  C'est ce qu'affirme le (i) de l'\'enonc\'e. L'assertion (ii) s'en d\'eduit comme en 3.4. $\square$

  \bigskip

\section{Germes et descente d'Harish-Chandra }

\bigskip

\subsection{Formule de descente pour les termes $\rho^{\tilde{G}}_{J}$}
On consid\`ere un triplet $(G,\tilde{G},{\bf a})$ quelconque. Soient $\tilde{M}$ un espace de Levi de $\tilde{G}$ et $\eta$ un \'el\'ement semi-simple de $\tilde{M}(F)$. 
Il y a une application naturelle $Z(\hat{G})\to Z(\hat{G}_{\eta})$, \'equivariante pour les actions galoisiennes. La classe de cocycle ${\bf a}$ d\'etermine par composition avec cette application une classe dans $H^1(\Gamma_{F};Z(\hat{G}_{\eta}))$, que nous noterons encore ${\bf a}$ pour simplifier. Le caract\`ere de $G_{\eta}(F)$ qui s'en d\'eduit est la restriction \`a ce groupe du caract\`ere $\omega$ de $G(F)$ d\'eduit du ${\bf a}$ initial. En consid\'erant $G_{\eta}$ comme un espace tordu sur lui-m\^eme, le triplet $(G_{\eta},G_{\eta},{\bf a})$ v\'erifie les m\^emes hypoth\`eses que notre triplet initial mais est "sans torsion".

Consid\'erons un voisinage ouvert $U_{\eta}$ de l'origine dans $G_{\eta}(F)$ qui est invariant par l'action de $Z_{G}(\eta;F)$, qui est tel que, pour $x\in G_{\eta}(F)$, $x$ appartient \`a $U_{\eta}$ si et seulement si la partie semi-simple de $x$ appartient \`a $U_{\eta}$ et qui est "assez petit". 
La descente d'Harish-Chandra fournit des applications transpos\'ees
$$desc_{\eta}^{\tilde{G}}:I(\tilde{G}(F),\omega)\otimes Mes(G(F))\to I(U_{\eta},\omega)\otimes Mes(G_{\eta}(F)),$$
$$desc_{\eta}^{\tilde{G},*}:D_{g\acute{e}om}(U_{\eta},\omega)\otimes Mes(G_{\eta}(F))^*\to D_{g\acute{e}om}(\tilde{G}(F),\omega)\otimes Mes(G(F))^*.$$
Rappelons que la donn\'ee d'un \'el\'ement $x\in  U_{\eta}$ et d'une mesure $dh$  sur $(G_{\eta})_{x}(F)$ d\'efinit un \'el\'ement  ${\bf x}\in D_{g\acute{e}om}(U_{\eta},\omega)\otimes Mes(G_{\eta}(F))^*$. C'est  l'int\'egrale orbitale qui, \`a $f\in C^{\infty}_{c}(U_{\eta})$ et \`a une mesure $dg$ sur $G_{\eta}(F)$, associe l'int\'egrale
$$I^{G_{\eta}}({\bf x},\omega,f\otimes dg)=D^{G_{\eta}}(x)^{1/2}\int_{(G_{\eta})_{x}(F)\backslash G_{\eta}(F)}f(y^{-1}xy)\omega(y)\,dy,$$
o\`u $dy$ est d\'eduite de $dg$ et $dh$. Si $x$ est assez proche de $1$, on a $(G_{\eta})_{x}=G_{x\eta}$ et le couple $(x\eta,dh)$ d\'efinit de m\^eme un \'el\'ement de $D_{g\acute{e}om}(\tilde{G}(F),\omega)\otimes Mes(G(F))^*$. Alors $desc_{\eta}^{\tilde{G},*}$ envoie l'\'el\'ement de $D_{g\acute{e}om}( U_{\eta},\omega)\otimes Mes(G_{\eta}(F))^*$ d\'efini par $(x,dh)$ sur l'\'el\'ement de $D_{g\acute{e}om}(\tilde{G}(F),\omega)\otimes Mes(G(F))^*$ d\'efini par $(x\eta,dh)$. Pour simplifier les notations, on oubliera le voisinage $U_{\eta}$ et on notera 
$$desc_{\eta}^{\tilde{G}}:I(\tilde{G}(F),\omega)\otimes Mes(G(F))\to I(G_{\eta}(F),\omega)\otimes Mes(G_{\eta}(F))$$
et
$$desc_{\eta}^{\tilde{G},*}:D_{g\acute{e}om}(G_{\eta}(F),\omega)\otimes Mes(G_{\eta}(F))^*\to D_{g\acute{e}om}(\tilde{G}(F),\omega)\otimes Mes(G(F))^*$$
les applications ci-dessus. On consid\'erera toutefois que, pour ${\bf f}\in I(\tilde{G}(F),\omega)\otimes Mes(G(F))$, les int\'egrales orbitales de $desc_{\eta}^{\tilde{G}}({\bf f})$ n'ont un sens que dans un voisinage de l'origine et que, de m\^eme, pour $\boldsymbol{\gamma}\in D_{g\acute{e}om}(G_{\eta}(F),\omega)$, $desc_{\eta}^{\tilde{G},*}(\boldsymbol{\gamma})$ n'est d\'efini que si le support de $\boldsymbol{\gamma}$ est assez voisin de l'origine. 

 On note ${\cal O}$ la classe de conjugaison de $\eta$ dans $\tilde{M}(F)$ et ${\cal O}^{\tilde{G}}$ sa classe de conjugaison dans $\tilde{G}(F)$. Rappelons que l'on note $D_{g\acute{e}om}({\cal O}^{\tilde{G}},\omega)$ le sous-espace de $D_{g\acute{e}om}(\tilde{G}(F),\omega)$ engendr\'e par les int\'egrales orbitales associ\'ees \`a des \'el\'ements $\gamma\in \tilde{G}(F)$ dont les parties semi-simples appartiennent \`a ${\cal O}^{\tilde{G}}$.   Il est clair que $desc_{\eta}^{\tilde{G},*}$ envoie $D_{unip}(G_{\eta}(F),\omega)\otimes Mes(G_{\eta}(F))^*$ dans (et m\^eme sur) $D_{g\acute{e}om}({\cal O}^{\tilde{G}},\omega)\otimes Mes(G(F))^*$.

 {\bf  On suppose $\eta$ elliptique dans $\tilde{M}(F)$}, c'est-\`a-dire $A_{\tilde{M}}=A_{M_{\eta}}$.  De cette \'egalit\'e se  d\'eduit une  injection $\Sigma^{G_{\eta}}(A_{M_{\eta}})\to \Sigma^G(A_{\tilde{M}})$. Rappelons que les ensembles ${\cal J}_{M_{\eta}}^{G_{\eta}}$, resp. ${\cal J}_{\tilde{M}}^{\tilde{G}}$, sont des classes d'\'equivalence d'ensembles $\{\alpha_{1},...,\alpha_{n}\}$ form\'es d'\'el\'ements lin\'eairement ind\'ependants de $\Sigma^{G_{\eta}}(A_{M_{\eta}})$, resp. $ \Sigma^G(A_{\tilde{M}})$, tels que $n=a_{M_{\eta}}-a_{G_{\eta}}$, resp. $n=a_{\tilde{M}}-a_{\tilde{G}}$. Deux ensembles sont \'equivalents s'ils engendrent le m\^eme ${\mathbb Z}$-module. L'\'egalit\'e $A_{\tilde{M}}= A_{M_{\eta}}$ \'equivaut \`a $a_{M_{\eta}}=a_{\tilde{M}}$. Si $A_{\tilde{G}}\not=A_{G_{\eta}}$, il n'y a pas de correspondance entre les ensembles ${\cal J}_{M_{\eta}}^{G_{\eta}}$, resp. ${\cal J}_{\tilde{M}}^{\tilde{G}}$ et ${\cal J}_{\tilde{M}}^{\tilde{G}}$. Mais, si l'on suppose $A_{\tilde{G}}=A_{G_{\eta}}$, 
 de l'injection pr\'ec\'edente se d\'eduit une injection ${\cal J}_{M_{\eta}}^{G_{\eta}}\to {\cal J}_{\tilde{M}}^{\tilde{G}}$ que l'on note $J'\mapsto J$. Pour de tels $J'\mapsto J$, l'espace $U_{J'}$ associ\'e \`a $J'$ est \'egal \`a l'espace $U_{J}$ associ\'e \`a $J$.

\ass{Lemme}{Soient $J\in {\cal J}_{\tilde{M}}^{\tilde{G}}$, $\boldsymbol{\gamma}'\in D_{unip}(M_{\eta}(F),\omega)\otimes Mes(M_{\eta}(F))^*$ et $a\in A_{\tilde{M}}(F)$ en position g\'en\'erale et assez proche de $1$. Posons $\boldsymbol{\gamma}=desc_{\eta}^{\tilde{M},*}(\boldsymbol{\gamma}')$. On a l'\'egalit\'e
$$\rho_{J}^{\tilde{G}}(\boldsymbol{\gamma},a)=\left\lbrace\begin{array}{cc}desc_{\eta}^{\tilde{M},*}(\rho_{J'}^{G_{\eta}}(\boldsymbol{\gamma}',a)),&\text{ si } A_{\tilde{G}}=A_{G_{\eta}} \text{ et }J{\text \,\,provient\,\, de\,\,  } J'\in {\cal J}_{M_{\eta}}^{G_{\eta}},\\ 0,&\text{ sinon. }\\ \end{array}\right.$$}

Preuve. Par lin\'earit\'e, on peut supposer que $\boldsymbol{\gamma}'$ est une int\'egrale orbitale associ\'ee \`a un \'el\'ement unipotent $u\in M_{\eta}(F)$. Alors $\boldsymbol{\gamma}$ est une int\'egrale orbitale associ\'e \`a l'\'el\'ement $u\eta\in \tilde{M}(F)$. On applique la formule de d\'efinition [II] 3.2(5):
$$(1) \qquad \rho_{J}^{\tilde{G}}(\boldsymbol{\gamma},a)=\sum_{\underline{\alpha}\in J}m(\underline{\alpha},u\eta)sgn(\underline{\alpha},u\eta)u_{\underline{\alpha}}(a)\boldsymbol{\gamma}.$$
Consid\'erons $\underline{\alpha}=\{\alpha_{1},...,\alpha_{n}\}\in J$. Pour tout $i=1,...,n$, fixons une "coracine" $\check{\alpha}_{i}$  que nous normalisons par la condition $<\alpha_{i},\check{\alpha}_{i}>=1$ (sic!).  Notons $m$ le volume du quotient de ${\cal A}_{\tilde{M}}^{\tilde{G}[J]}$ par le ${\mathbb Z}$-module engendr\'e par ces $\check{\alpha}_{i}$, pour $i=1,...,n$. Le terme $\rho^{\tilde{G}}(\alpha_{i},u\eta)$ d\'efini en [II] 1.5 est proportionnel \`a $\check{\alpha}_{i}$. Il r\'esulte des d\'efinitions que
$$ m(\underline{\alpha},u\eta)sgn(\underline{\alpha},u\eta)=m\prod_{i=1,...,n}<\alpha_{i},\rho^{\tilde{G}}(\alpha_{i},u\eta)>.$$
Soit $i\in \{1,...,n\}$. Parce que $A_{\tilde{M}}=A_{M_{\eta}}$, la d\'efinition [II] 1.5 entra\^{\i}ne que $\rho^{\tilde{G}}(\alpha_{i},u\eta)=0$ si $\alpha_{i}$ ne  provient pas de $ \Sigma_{M_{\eta}}^{G_{\eta}}(A_{M_{\eta}})$, tandis que 
$$\rho^{\tilde{G}}(\alpha_{i},u\eta)= \rho^{G_{\eta}}(\alpha'_{i},u)$$
si $\alpha_{i}$ est l'image de $\alpha'_{i}\in  \Sigma_{M_{\eta}}^{G_{\eta}}(A_{M_{\eta}})$. Si $A_{\tilde{G}}\not=A_{G_{\eta}}$ ou si $A_{\tilde{G}}=A_{G_{\eta}}$ mais $J$ ne provient pas de  $ {\cal J}_{M_{\eta}}^{G_{\eta}}$, il n'y a aucun \'el\'ement $\underline{\alpha}\in J$ qui v\'erifie cette condition pour tout $i$ et on obtient $ \rho_{J}^{\tilde{G}}(\boldsymbol{\gamma},a)=0$.  Supposons d\'esormais que  $A_{\tilde{G}}=A_{G_{\eta}}$ et que $J$ est l'image de $J'\in {\cal J}_{M_{\eta}}^{G_{\eta}}$. Les $\underline{\alpha}$ qui contribuent \`a (1) sont exactement les images d'\'el\'ements $\underline{\alpha}'\in J'$.  Si $\underline{\alpha}$ provient de $\underline{\alpha}'$, les formules ci-dessus montrent que 
$$m(\underline{\alpha},u\eta)sgn(\underline{\alpha},u\eta)=m(\underline{\alpha}',u)sgn(\underline{\alpha}',u).$$
On a aussi $u_{\underline{\alpha}}(a)=u_{\underline{\alpha}'}(a)$. 
On obtient 
$$\rho_{J}^{\tilde{G}}(\boldsymbol{\gamma},a)=\sum_{\underline{\alpha}'\in J'}m(\underline{\alpha}',u)sgn(\underline{\alpha}',u)u_{\underline{\alpha}'}(a)\boldsymbol{\gamma}.$$
 On a une formule analogue \`a (1):
 $$\rho_{J'}^{G_{\eta}}(\boldsymbol{\gamma}',a)=\sum_{\underline{\alpha}'\in J'}m(\underline{\alpha}',u)sgn(\underline{\alpha}',u)u_{\underline{\alpha}'}(a)\boldsymbol{\gamma}'.$$
 Puisque $\boldsymbol{\gamma}=desc_{\eta}^{\tilde{M},*}(\boldsymbol{\gamma}')$, on en d\'eduit
$$\rho_{J}^{\tilde{G}}(\boldsymbol{\gamma},a)=desc_{\eta}^{\tilde{M},*}(\rho_{J'}^{G_{\eta}}(\boldsymbol{\gamma}',a)),$$
ce qui prouve le lemme. $\square$

 {\bf Variante.} Supposons $(G,\tilde{G},{\bf a})$ quasi-d\'eploy\'e et \`a torsion int\'erieure. Fixons un syst\`eme de fonctions $B$ comme en [II] 1.9. Ce syst\`eme d\'etermine une fonction $B_{\eta}$ sur le syst\`eme de racines de $G_{\eta}$, que l'on notera $B_{{\cal O}}$ pour simplifier. Ainsi on dispose des ensembles ${\cal J}_{\tilde{M}}^{\tilde{G}}(B_{{\cal O}})$ et ${\cal J}_{M_{\eta}}^{G_{\eta}}(B_{\eta})$. On a une proposition similaire relative \`a ces ensembles. En fait, sa conclusion se simplifie car, avec la d\'efinition que l'on a donn\'ee de l'ensemble ${\cal J}_{\tilde{M}}^{\tilde{G}}(B_{{\cal O}})$, l'injection ${\cal J}_{M_{\eta}}^{G_{\eta}}(B_{{\cal O}})\to {\cal J}_{\tilde{M}}^{\tilde{G}}(B_{{\cal O}})$ est bijective. L'existence d'un $J\in {\cal J}_{\tilde{M}}^{\tilde{G}}(B_{{\cal O}})$ entra\^{\i}ne donc que $A_{\tilde{G}}=A_{G_{\eta}}$ et que $J$ provient d'un $J'\in {\cal J}_{M_{\eta}}^{G_{\eta}}(B_{{\cal O}})$. 
 
 \bigskip
 
 \subsection{Descente des germes d'int\'egrales orbitales pond\'er\'ees}
 On conserve les m\^emes donn\'ees. On suppose toujours  que $\eta$ est elliptique dans $\tilde{M}$.

\ass{Proposition}{On a l'\'egalit\'e 
$$g_{\tilde{M},{\cal O}}^{\tilde{G}}\circ desc_{\eta}^{\tilde{M},*}=\left\lbrace\begin{array}{cc}desc_{\eta}^{\tilde{G},*}\circ g_{M_{\eta},unip}^{G_{\eta}},&\text{ si }A_{\tilde{G}}=A_{G_{\eta}},\\ 0,&\text{ sinon. }\\ \end{array}\right. .$$}

Les deux termes sont des germes d'applications lin\'eaires d\'efinies sur $D_{g\acute{e}om,G_{\eta}-\acute{e}qui}(M_{\eta}(F),\omega)$  au voisinage de l'\'el\'ement neutre de $M_{\eta}(F)$. L'assertion n'est qu'une reformulation du lemme 9.2 de [A2].  Un examen de la preuve de ce lemme montre que l'hypoth\`ese de r\'egularit\'e figurant dans l'\'enonc\'e d'Arthur ne sert pas. D'autre part, on a modifi\'e les d\'efinitions d'Arthur en [II] 1.5 et on doit montrer que cela n'influe pas sur la preuve.  On se rend compte qu'il suffit de prouver que le lemme 8.2 de [A2] reste vrai avec notre d\'efinition des fonctions $r_{\tilde{M}}^{\tilde{G}}(\gamma,a)$. Pr\'ecis\'ement, soit $x\in M_{\eta}(F)$ assez proche de $1$. Soit $\tilde{P}\in {\cal P}(\tilde{M})$. Posons $P_{\eta}=P\cap G_{\eta}$. On doit montrer que, pour $\lambda\in i{\cal A}_{\tilde{M}}^*$ et pour $a\in A_{\tilde{M}}(F)\simeq A_{M_{\eta}}(F)$ en position g\'en\'erale et assez proche de $1$, on a l'\'egalit\'e
$$(1) \qquad r_{\tilde{P}}(x\eta,a;\lambda)=r_{P_{\eta}}(x,a;\lambda).$$
Cela r\'esulte de notre d\'efinition: si on note $t$ la partie semi-simple de $x$ et $u$ sa partie unipotente,  les deux fonctions se d\'efinissent \`a l'aide des m\^emes termes $\rho^{G_{t\eta}}(\beta,u)$, pour $\beta\in \Sigma^{G_{t\eta}}(A_{M_{\eta}})$. $\square$

{\bf Variante.} Supposons $(G,\tilde{G},{\bf a})$ quasi-d\'eploy\'e et \`a torsion int\'erieure. Fixons un syst\`eme de fonctions $B$ comme en [II] 1.9. On a

(2) soit $\boldsymbol{\gamma}\in D_{g\acute{e}om,G_{\eta}-\acute{e}qui}(M_{\eta}(F),\omega)\otimes Mes(M_{\eta}(F))^*$  assez proche de l'origine; alors on a l'\'egalit\'e
$$g_{\tilde{M},{\cal O}}^{\tilde{G}}(desc_{\eta}^{\tilde{M},*}(\boldsymbol{\gamma}),B)=\left\lbrace\begin{array}{cc}desc_{\eta}^{\tilde{G},*}( g_{M_{\eta},unip}^{G_{\eta}}(\boldsymbol{\gamma},B_{\eta})),&\text{ si }A_{\tilde{G}}=A_{G_{\eta}},\\ 0,&\text{ sinon. }\\ \end{array}\right.$$
 
La preuve est la m\^eme. Il y a toutefois une subtilit\'e. Il ne doit intervenir dans la preuve que des \'el\'ements $x\eta$, avec $x\in M_{\eta}(F)$, pour lesquels on a l'analogue de (1), \`a savoir
$$(2) \qquad  r_{\tilde{P}}(x\eta,a,B;\lambda)=r_{P_{\eta}}(x,a,B_{\eta};\lambda).$$
Avec les notations plus haut, la fonction de gauche est d\'efinie \`a l'aide des termes $\rho^{G_{t\eta}}(\beta,u)$ pour $\beta\in \Sigma^{G_{t\eta}}(A_{M_{\eta}},B_{t\eta})$ tandis que celle de droite est d\'efinie \`a l'aide des termes $\rho^{G_{t\eta}}(\beta,u)$ pour $\beta\in \Sigma^{G_{t\eta}}(A_{M_{\eta}},B_{\eta})$. Pour $x$ quelconque,  il n'y a gu\`ere de raison pour que ces termes soient \'egaux. Mais il n'intervient que des $x$ unipotents, pour lesquels $t=1$ et les termes co\"{\i}ncident, et des $x$ qui sont $G_{\eta}$-\'equisinguliers pour lesquels les termes co\"{\i}ncident aussi car les deux ensembles de racines sont vides.

\bigskip

\subsection{Formule de descente pour les termes $\sigma_{J}$}
Soient $(G,\tilde{G},{\bf a})$ un triplet quasi-d\'eploy\'e et \`a torsion int\'erieure, $B$ un syst\`eme de fonctions comme en [II] 1.9, $\tilde{M}$ un espace de Levi de $\tilde{G}$ et $\eta$ un \'el\'ement semi-simple de $\tilde{M}(F)$.   On note ${\cal O}$ la classe de conjugaison stable de $\eta$ dans $\tilde{M}(F)$ et ${\cal O}^{\tilde{G}}$ sa classe de conjugaison stable dans $\tilde{G}(F)$. 

La descente d'Harish-Chandra s'adapte aux distributions stables, cf. [I] 4.8.  Pour d\'efinir correctement les applications de descente, il faudrait introduire  comme en 4.1 un voisinage convenable $U_{\eta}$ de l'origine dans $G_{\eta}(F)$. Pour simplifier la notation, on consid\`ere comme dans ce paragraphe que l'on a des applications transpos\'ees
$$desc_{\eta}^{st,\tilde{G}}:SI(\tilde{G}(F))\otimes Mes(G(F))\to SI(G_{\eta}(F))\otimes Mes(G_{\eta}(F)),$$
$$desc_{\eta}^{st,\tilde{G},*}:D_{g\acute{e}om}^{st}(G_{\eta}(F))\otimes Mes(G_{\eta}(F))^*\to D_{g\acute{e}om}^{st}(G(F))\otimes Mes(G(F))^*.$$
Mais, pour ${\bf f}\in SI(\tilde{G}(F))\otimes Mes(G(F))$, les int\'egrales orbitales stables de $desc_{\eta}^{st,\tilde{G}}({\bf f})$ n'ont de sens que que dans un voisinage de l'origine. De m\^eme, pour $\boldsymbol{\delta}\in D_{g\acute{e}om}^{st}(G_{\eta}(F))\otimes Mes(G_{\eta}(F))^*$, $desc_{\eta}^{st,\tilde{G},*}(\boldsymbol{\delta})$ n'est d\'efini que si le support de $\boldsymbol{\delta}$ est assez voisin de l'origine. 
En particulier, $desc_{\eta}^{st,\tilde{G},*}$ se restreint en une surjection de $D_{unip}^{st}(G_{\eta}(F))\otimes Mes(G_{\eta}(F))^*$ sur $D_{g\acute{e}om}^{st}({\cal O}^{\tilde{G}})\otimes Mes(G(F))^*$.

{\bf Attention.} Le diagramme 
$$\begin{array}{ccc}I(\tilde{G}(F))\otimes Mes(G(F))&\stackrel{desc_{\eta}^{\tilde{G}}}{\to}& I(G_{\eta}(F))\otimes Mes(G_{\eta}(F))\\ \downarrow&&\downarrow\\ SI(\tilde{G}(F))\otimes Mes(G(F))&\stackrel{desc_{\eta}^{st,\tilde{G}}}{\to}& SI(G_{\eta}(F))\otimes Mes(G_{\eta}(F))\\ \end{array}$$
n'est pas commutatif, cf. [I] 5.10. 

\bigskip

{\bf On suppose $\eta$ elliptique dans $\tilde{M}(F)$}. Le groupe dual $\hat{M}_{\eta}$ s'identifie \`a un sous-groupe de $\hat{M}$. Dualement, l'ellipticit\'e de $\eta$  entra\^{\i}ne que $Z(\hat{M})^{\Gamma_{F}}$ est un sous-groupe d'indice fini de $Z(\hat{M}_{\eta})^{\Gamma_{F}}$. Dans le cas o\`u $A_{G}=A_{G_{\eta}}$, le groupe  $Z(\hat{G})^{\Gamma_{F}}$ est aussi  un sous-groupe d'indice fini de $Z(\hat{G}_{\eta})^{\Gamma_{F}}$. L'homomorphisme
 $$Z(\hat{M})^{\Gamma_{F}}/Z(\hat{G})^{\Gamma_{F}}\to Z(\hat{M}_{\eta})^{\Gamma_{F}}/Z(\hat{G}_{\eta})^{\Gamma_{F}}$$
 est surjectif et son noyau est fini. On note $e_{\tilde{M}}^{\tilde{G}}(\eta) $ l'inverse du nombre d'\'el\'ements de ce noyau.
 
 Rappelons que l'on a d\'efini en 4.1 une bijection    de ${\cal J}_{M_{\eta}}^{G_{\eta}}(B_{{\cal O}})$ dans ${\cal J}_{\tilde{M}}^{\tilde{G}}(B_{{\cal O}})$. Pour simplifier, on identifie ces deux ensembles. L'hypoth\`ese que ${\cal J}_{\tilde{M}}^{\tilde{G}}(B_{{\cal O}})$ est non vide implique automatiquement que $A_{G}=A_{G_{\eta}}$.

\ass{Proposition (\`a prouver)}{Soient $J\in {\cal J}_{\tilde{M}}^{\tilde{G}}(B_{{\cal O}})$, $\boldsymbol{\delta}'\in D_{unip}^{st}(M_{\eta}(F),\omega)\otimes Mes(M_{\eta}(F))^*$ et $a\in A_{\tilde{M}}(F)$ en position g\'en\'erale et assez proche de $1$. Posons $\boldsymbol{\delta}=desc_{\eta}^{st,\tilde{M},*}(\boldsymbol{\delta}')$.  On a l'\'egalit\'e
$$\sigma_{J}^{\tilde{G}}(\boldsymbol{\delta},a)= e_{\tilde{M}}^{\tilde{G}}(\eta)desc_{\eta}^{st,\tilde{M},*}(\sigma_{J}^{G_{\eta}}(\boldsymbol{\delta}',a)) .$$}

\bigskip

\subsection{Formule de descente pour les germes des int\'egrales orbitales pond\'er\'ees stables}
On conserve les m\^emes hypoth\`eses. On suppose $\eta$ elliptique dans $\tilde{M}(F)$.

\ass{Proposition (\`a prouver)}{Soit $\boldsymbol{\delta}\in D_{g\acute{e}om,G_{\eta}-\acute{e}qui}^{st}(M_{\eta}(F))\otimes Mes(M_{\eta}(F))^*$  \`a support assez proche de l'origine. Alors on a l'\'egalit\'e
$$Sg_{\tilde{M},{\cal O}}^{\tilde{G}}( desc_{\eta}^{st,\tilde{M},*}(\boldsymbol{\delta}),B)=\left\lbrace\begin{array}{cc} e_{\tilde{M}}^{\tilde{G}}(\eta) desc_{\eta}^{st,\tilde{G},*}( Sg_{M_{\eta},unip}^{G_{\eta}}(\boldsymbol{\delta},B_{\eta})),&\text{ si }A_{G}=A_{G_{\eta}},\\ 0,&\text{ sinon.}\\ \end{array}\right.$$}

\bigskip

\section{Descente et endoscopie}

\bigskip

\subsection{Descente de donn\'ees endoscopiques}
Soient  $(G,\tilde{G},{\bf a})$ un triplet quelconque et ${\bf G}'=(G',{\cal G}',\tilde{s})$ une donn\'ee endoscopique elliptique et relevante. 

D'apr\`es [I] paragraphe 3, on peut  r\'ealiser les objets duaux de la fa\c{c}on suivante. On fixe une paire de Borel \'epingl\'ee  $\hat{{\cal E}}=(\hat{B},\hat{T},(\hat{E}_{\alpha})_{\alpha\in \Delta})$ de $\hat{G}$. On note $\hat{\theta}$ l'automorphisme dual de $\theta^*$ qui conserve cette paire. On modifie l'action galoisienne de sorte que $\hat{{\cal E}}$ soit stable par cette action. On suppose $\tilde{s}=s\hat{\theta}$ avec $s\in \hat{T}$. On munit $\hat{G}'$ d'une paire de Borel \'epingl\'ee  conserv\'ee  par l'action galoisienne sur ce groupe, dont la paire de Borel sous-jacente est $(\hat{B}\cap \hat{G}',\hat{T}^{\hat{\theta},0})$. 

Fixons un diagramme
 $$(\epsilon,B,T',B,T,\eta)$$
 cf. [I] 1.10.  Le terme $\eta$ est un \'el\'ement semi-simple de $\tilde{G}(F)$ et
$\epsilon$ est un \'el\'ement semi-simple de $\tilde{G}'(F)$. On suppose que $G'_{\epsilon}$ est  quasi-d\'eploy\'e. Rappelons que l'homomorphisme $\xi:T\to T'$ d\'eduit du diagramme fix\'e se d\'eduit de l'homomorphisme compos\'e
$$X_{*}(T)\simeq X^*(\hat{T})\to X^*(\hat{T}^{\hat{\theta},0})\simeq X_{*}(T').$$
  Posons $(\bar{B},\bar{T})=(B\cap G_{\eta},T\cap G_{\eta})$, compl\'etons cette paire  en une paire de Borel \'epingl\'ee de $G_{\eta}$ et introduisons l'action galoisienne quasi-d\'eploy\'ee $\sigma\mapsto \sigma_{\bar{G}}$ sur $G_{\eta}$ qui la conserve, cf. [I] 1.2. Elle est de la forme $\sigma\mapsto \sigma_{\bar{G}}=ad_{\bar{u}(\sigma)}\circ \sigma_{G}$, o\`u $\bar{u}(\sigma)\in G_{\eta,SC}$.On note $\bar{G}$ le groupe $G_{\eta}$ muni de cette action galoisienne quasi-d\'eploy\'ee. Compl\'etons de m\^eme la paire de Borel $(B,T)$ de $G$ en une paire de Borel \'epingl\'ee et introduisons l'action galoisienne quasi-d\'eploy\'ee sur $G$ qui conserve celle-ci. Elle est de la forme $\sigma\mapsto  \sigma_{G^*}=ad_{u(\sigma)}\circ \sigma_{G}$, o\`u $u(\sigma)\in G_{SC}$. 
On a $\sigma_{\bar{G}}=ad_{\bar{u}(\sigma)u(\sigma)^{-1}}\circ \sigma_{G^*}$. Puisque les deux actions intervenant ici conservent $T$, l'\'el\'ement $\bar{u}(\sigma)u(\sigma)^{-1}$ normalise $T$. Il s'envoie sur un \'el\'ement du groupe de Weyl $W$ que l'on note $\omega_{\bar{G}}(\sigma)$. En utilisant le fait que $ad_{\eta}$ conserve $(B,T)$, on v\'erifie que $\omega_{\bar{G}}(\sigma)$ est fixe par $\theta$, cf. [I] 1.3(3).   On peut identifier un tore maximal $\hat{\bar{T}}$ du groupe $\hat{\bar{G}}$ a $\hat{T}/(1-\hat{\theta})(\hat{T})$, muni de l'action $\sigma\mapsto \omega_{\bar{G}}(\sigma)\circ \sigma_{G^*}$ (en notant encore $\sigma_{G^*}$ l'action sur le groupe $\hat{G}$).

Dans [W1] 3.5, on a d\'efini une donn\'ee endoscopique de $\bar{G}_{SC}$ associ\'e \`a la donn\'ee ${\bf G}'$ et au diagramme $(\epsilon,B,T',B,T,\eta)$. Notons-la $\bar{{\bf G}}'=(\bar{G}',\bar{{\cal G}}',\bar{s})$. Rappelons le point cl\'e de sa d\'efinition.    Par la suite d'homomorphismes
 $$\hat{T}\to \hat{T}/(1-\hat{\theta})(\hat{T})\simeq \hat{\bar{T}}\to \hat{\bar{T}}_{ad}=\hat{\bar{T}}/Z(\hat{\bar{G}}),$$
 l'\'el\'ement $s$ s'envoie sur un \'el\'ement $\bar{s}\in \hat{\bar{T}}_{ad}$.  Le groupe $\hat{\bar{G}}'$ est la composante neutre du commutant de $\bar{s}$ dans $\hat{\bar{G}}_{AD}$. L'action galoisienne sur $\hat{\bar{G}}'$  est   compos\'ee de celle sur $\hat{\bar{G}}_{AD}$ et d'un cocycle \`a valeurs dans ce groupe. On renvoie \`a [W1] 3.5 pour une description plus compl\`ete.    Fixons une paire de Borel  $(B^{\flat},T^{\flat})$ de $G'_{\epsilon}$ d\'efinie sur $F$ et un \'el\'ement $g'\in G'_{\epsilon,SC}$ tel que $ad_{g'}$ envoie cette paire sur $(B'\cap M'_{\epsilon},T')$.  Fixons une paire de Borel $(\bar{B}',\bar{T}')$ de $\bar{G}'$ d\'efinie sur $F$.     On note  $T^{\flat}_{sc}$, resp.  $\bar{T}'_{sc}$,  l'image r\'eciproque de $T^{\flat}$ dans $ G'_{\epsilon,SC}$, resp.    l'image r\'eciproque de $ \bar{T}'$ dans $\bar{G}'_{SC}$.  On a la suite d'homomorphismes
  $$ \bar{T}'_{sc}\times Z(\bar{G}')^0\times Z(\bar{G})^0\to  \bar{T}'\times Z(\bar{G})^0,$$
 $$ \bar{T}'\times Z(\bar{G})^0\to  \bar{T}_{sc}\times Z(\bar{G})^0,$$
 celui-ci provenant du fait que $\bar{{\bf  G}}'$ est une donn\'ee endoscopique de $\bar{G}_{SC}$,
 $$\bar{T}_{sc}\times Z(\bar{G})^0\to \bar{T},$$
 $$\bar{T}\stackrel{\xi}{\to}T'\stackrel{ad_{g'}^{-1}}{\to}T^{\flat},$$
 et l'homomorphisme en sens inverse
 $$T^{\flat}\leftarrow  T^{\flat}_{sc}\times Z(G'_{\epsilon})^0.$$
 Par composition, on en d\'eduit un isomorphisme d'alg\`ebres de Lie 
 $$j:\bar{\mathfrak{t}}'_{sc}\oplus \mathfrak{z}(\bar{G}')\oplus \mathfrak{z}(\bar{G})\simeq \mathfrak{t}^{\flat}_{sc}\oplus \mathfrak{z}(G'_{\epsilon}).$$
 On a vu en [W1] 3.5 que cet isomorphisme \'etait \'equivariant pour les actions galoisiennes et se restreignait en un isomorphisme de $\bar{\mathfrak{t}}'_{sc}$ sur $\mathfrak{t}^{\flat}_{sc}$ et un isomorphisme
  $$(1) \qquad  \mathfrak{z}(\bar{G})\oplus  \mathfrak{z}(\bar{G}') \simeq\mathfrak{z}(G'_{\epsilon}).$$
  On note $j_{*}:X_{*}(\bar{T}'_{sc})\otimes {\mathbb Q}\to X_{*}(T^{\flat}_{sc})\otimes {\mathbb Q}$ l'isomorphisme sous-jacent \`a la restriction de l'homomorphisme $j$ \`a $\bar{\mathfrak{t}}'_{sc}$. Le triplet $(\bar{G}'_{SC},G'_{\epsilon,SC},j_{*})$  est un triplet endoscopique non standard, cf. [W1] 1.7.

{\bf Cas particulier.} Supposons $(G,\tilde{G},{\bf a})$ quasi-d\'eploy\'e et \`a torsion int\'erieure. Alors les triplets endoscopiques non standard sont triviaux, c'est-\`a-dire que $ \bar{G}'_{SC}=G'_{\epsilon,SC}$. En fait, dans ce cas, on n'a pas besoin de passer au rev\^etements simplement connexes des groupes d\'eriv\'es (ce passage sert en g\'en\'eral \`a trivialiser le caract\`ere $\omega$).  Le groupe $G'_{\epsilon}$ est celui d'une donn\'ee endoscopique de $\bar{G}$. 

\bigskip
Fixons des donn\'ees auxiliaires $G'_{1},...,\Delta_{1}$ pour la donn\'ee ${\bf G}'$. On fixe une image r\'eciproque $\epsilon_{1}$ de $\epsilon$ dans $\tilde{G}'_{1}(F)$ ainsi qu'une d\'ecomposition d'alg\`ebres de Lie
$$(2) \qquad \mathfrak{g}'_{1}=\mathfrak{c}_{1}\oplus \mathfrak{g}'.$$
Notons ${\cal Y}$ l'ensemble des $y\in G$ tels que $y\sigma(y)^{-1}\in I_{\eta}$ pour tout $\sigma\in \Gamma_{F}$, o\`u $I_{\eta}=Z(G)^{\theta}G_{\eta}$.  Pour tout $y\in {\cal Y}$, on pose $\eta[y]=y^{-1}\eta y$. Notons pour plus de clart\'e $\psi:G_{\eta}\to \bar{G}$ l'identit\'e, qui est un torseur pour les actions galoisiennes d\'efinies sur ces deux groupes.   On a d\'efini la donn\'ee endoscopique $\bar{{\bf G}}'$ de $\bar{G}_{SC}$. L'application $\psi\circ ad_{y}:G_{\eta[y]}\to \bar{G}$ est un torseur int\'erieur gr\^ace auquel $\bar{{\bf G}}'$ appara\^{\i}t aussi comme une donn\'ee endoscopique pour $G_{\eta[y],SC}$.  Supposons que cette donn\'ee soit  relevante. Il s'agit ici d'endoscopie non tordue. De plus, le groupe $G_{\eta[y],SC}$ est simplement connexe. On sait qu'alors on n'a pas besoin de donn\'ees auxiliaires, on peut d\'efinir un facteur de transfert sur $\bar{G}'(F)\times G_{\eta[y],SC}(F)$. 
On fixe un tel facteur de transfert $\Delta(y)$. Soit $Y\in \mathfrak{g}'_{\epsilon}(F)$ en position g\'en\'erale et assez proche de $0$, que l'on identifie par (2) \`a un \'el\'ement de $\mathfrak{g}'_{1,\epsilon_{1}}(F)$.  On le d\'ecompose en $Y_{sc}+Z$ avec $Y_{sc}\in \mathfrak{g}'_{\epsilon,SC}(F)$ et $Z\in \mathfrak{z}(G'_{\epsilon};F)$. Via l'isomorphisme (1), on 
    d\'ecompose $Z$ en $Z_{1}+Z_{2}$, avec $Z_{1}\in \mathfrak{z}(\bar{G};F)$ et $Z_{2}\in \mathfrak{z}(\bar{G}';F)$. Par endoscopie non standard, $Y_{sc}$ d\'etermine une classe de conjugaison stable dans $\mathfrak{\bar{g}}'_{SC}(F)$. Fixons $\bar{Y}_{sc}$ dans cette classe. Posons $\bar{Y}=\bar{Y}_{sc}+Z_{2}$. Par endoscopie ordinaire, il peut correspondre ou pas \`a $\bar{Y}$ une classe de conjugaison stable dans $\mathfrak{g}_{\eta[y],SC}(F)$. Supposons que cette classe existe et fixons un \'el\'ement $X[y]_{sc}$ de cette classe.  Puisque $\bar{G}$ est une forme int\'erieure de $G_{\eta[y]}$, les espaces $\mathfrak{z}(\bar{G})$ et $\mathfrak{z}(G_{\eta[y]})$ s'identifient et on peut consid\'erer $Z_{1}$ comme un \'el\'ement de $\mathfrak{z}(G_{\eta[y]};F)$. On pose $X[y]=X[y]_{sc}+Z_{1}$. 
  Rappelons le th\'eor\`eme 3.9 de [W1]:  il existe $d(y)\in {\mathbb C}^{\times}$ tel que,  pour des donn\'ees comme ci-dessus, on ait l'\'egalit\'e
 $$(3) \qquad d(y)\Delta(y)(exp(\bar{Y}),exp(X[y]_{sc}))=\Delta_{1}(exp(Y)\epsilon_{1},exp(X[y])\eta[y])$$
 pourvu que $X$ et $Y$ soient assez proches de $0$.

\bigskip
 
 \subsection{Transfert des fonctions et des distributions}
 
Fixons des mesures de Haar sur  tous les groupes  intervenant, auxquelles on imposera quelques conditions de compatibilit\'e, cf. remarque (5) ci-dessous.  Soient  $f\in I(\tilde{G}(F),\omega)$ et $Y\in \mathfrak{g}'_{\epsilon}(F)$ comme en 5.1.  On utilise pour  $Y$ les d\'efinitions de ce paragraphe. On d\'efinit

-   le transfert $f'$ de $f$, qui est une fonction sur  $\tilde{G}'_{1}(F)$;

-   l'image $f'_{\epsilon_{1}}$ de $f^{\tilde{G}'_{1}}$ par l'application de descente $desc_{\epsilon_{1}}^{st,\tilde{G}'_{1}}$; c'est une fonction sur $G'_{1,\epsilon_{1}}(F)$;

  -   l'image $f'(Z)_{sc}$ de la fonction $x\mapsto f'_{\epsilon_{1}}(exp(Z)x)$ par l'application $\iota_{G'_{\epsilon,SC},G'_{1,\epsilon_{1}}}$; c'est une fonction sur $G'_{\epsilon,SC}(F)$ (on remarque que $G'_{\epsilon,SC}$ est aussi le rev\^etement simplement connexe du groupe d\'eriv\'e de $G'_{1,\epsilon_{1}}$).
  
  Soit $y\in {\cal Y}$. D\'efinissons 
  
  - l'image $f[y]$ de $f$ par l'application de descente $desc^{\tilde{G}}_{\eta[y]}$; c'est une fonction sur $G_{\eta[y]}(F)$;
  
  - l'image $f[y,Z_{1}]_{sc}$ de la fonction $x\mapsto f[y](exp(Z_{1})x)$ par l'application $\iota_{G_{\eta[y],SC},G_{\eta[y]}}$; c'est une fonction sur $G_{\eta[y],SC}(F)$;
  
  - le transfert $\bar{f}[y,Z_{1}]$ de $f[y,Z_{1}]_{sc}$; c'est une fonction sur $\bar{G}'(F)$ qui est nulle si $\bar{G}'$ n'est pas relevante pour $G_{\eta[y],SC}$;
  
  - l'image $\bar{f}[y,Z]_{sc}$ de la fonction $x\mapsto \bar{f}[y,Z_{1}](exp(Z_{2})x)$ par l'application $\iota_{\bar{G}'_{SC},\bar{G}'}$; c'est une fonction sur $\bar{G}'_{SC}(F)$. 
  
  On pose
  $$c[y]=[I_{\eta[y]}(F):G_{\eta[y]}(F)]^{-1}.$$
   Fixons un ensemble de repr\'esentants $\dot{{\cal Y}}$ de l'ensemble de doubles classes $I_{\eta}\backslash {\cal Y}/G(F)$. 
 Avec les notations ci-dessus, le lemme 3.11 de [W1] affirme l'\'egalit\'e
 $$   S_{\lambda_{1}}^{\tilde{G}'_{1}}(exp(Y)\epsilon_{1},f^{\tilde{G}'_{1}}) =\sum_{y\in \dot{{\cal Y}}}c[y]d(y)S^{\bar{G}'}(\bar{Y},\bar{f}[y,Z_{1}]) $$
 pourvu que $Y$ soit assez proche de $0$, cette notion \'etant ind\'ependante de $f$. Le terme $d(y)$ n'a \'et\'e d\'efini que dans le cas o\`u $\bar{G}'$ est relevante pour $G_{\eta[y],SC}$. Par convention, il est nul dans le cas contraire (la fonction $\bar{f}[y,Z_{1}]$ est nulle dans ce cas). 
 D'apr\`es les d\'efinitions, on a
$$S_{\lambda_{1}}^{\tilde{G}'_{1}}(exp(Y)\epsilon_{1},f^{\tilde{G}'_{1}})= S^{G'_{\epsilon,SC}}(exp(Y_{sc}),f'(Z)_{sc}),$$
$$S^{\bar{G}'}(\bar{Y},\bar{f}[y,Z_{1}])=S^{\bar{G}'_{SC}}(exp(\bar{Y}_{sc}),\bar{f}[y,Z]_{sc}).$$
L'\'egalit\'e pr\'ec\'edente devient
$$(1)  \qquad S^{G'_{\epsilon,SC}}(exp(Y_{sc}),f'(Z)_{sc})=\sum_{y\in \dot{{\cal Y}}}c[y]d(y)S^{\bar{G}'_{SC}}(exp(\bar{Y}_{sc}),\bar{f}[y,Z]_{sc}).$$
Remarquons que les termes $Y_{sc}$ et $Z$ jouent le r\^ole de variables ind\'ependantes. Les d\'efinitions des fonctions $f'(Z)_{sc}$ et $\bar{f}[y,Z]_{sc}$ conservent un sens pour $Z=0$. On note $f'_{sc}$ et $\bar{f}[y]_{sc}$ leurs valeurs en $Z=0$. On obtient l'\'egalit\'e
$$(2) \qquad  S^{G'_{\epsilon,SC}}(exp(Y_{sc}),f'_{sc})=\sum_{y\in \dot{{\cal Y}}}c[y]d(y)S^{\bar{G}'_{SC}}(exp(\bar{Y}_{sc}),\bar{f}[y]_{sc}).$$

{\bf Remarques.} (3) Dans [W1], on n'avait pas de constantes $d(y)$ simplement parce qu'on avait normalis\'e les facteurs de transfert de sorte que ces constantes vaillent $1$.
 
  (4)  La d\'efinition de la fonction $f[y]_{sc}$ dans  [W1]  est plus directe que celle ci-dessus.  Elle consiste \`a remplacer l'application $\iota_{G_{\eta[y],SC},G_{\eta[y]}}$ utilis\'ee ici par la simple image r\'eciproque. C'est-\`a-dire que l'on ne moyenne pas cette image r\'eciproque par l'action   adjointe de $G_{\eta[y]}(F)$. Mais la formule 5.1(3) implique que le facteur de transfert $\Delta(y)$   est $\omega$-\'equivariant par  cette action. Il en r\'esulte ais\'ement que  les deux d\'efinitions possibles de $f[y]_{sc}$ donnent le m\^eme transfert  $\bar{f}[y]$.
  
  (5) La formule (2) n\'ecessite des compatibilit\'es entre les mesures choisies. Notons $T_{Y_{sc}}$ le commmutant de $Y_{sc}$ dans $G'_{\epsilon,SC}$ et $T_{\bar{Y}_{sc}}$ celui de $\bar{Y}_{sc}$ dans $\bar{G}'_{SC}$. 
  Pour d\'efinir les int\'egrales orbitales stables, on doit fixer des mesures sur  $T_{Y_{sc}}(F)$ et 
  $T_{\bar{Y}_{sc}}$. Via l'exponentielle, il revient au m\^eme de fixer des mesures sur $\mathfrak{t}_{Y_{sc}}(F)$ et $\mathfrak{t}_{\bar{Y}_{sc}}(F)$. Or on est dans une situation d'endoscopie non standard et ces alg\`ebres de Lie sont naturellement isomorphes. On choisit des mesures qui se correspondent par cet isomorphisme. D'autre part, dans les passages entre groupes et rev\^etements simplement connexes, il est n\'ecessaire de fixer des mesures sur les centres. Pr\'ecis\'ement, on doit fixer des mesures sur $Z(G'_{\epsilon})^0(F)$, $Z(G_{\eta[y]})^0(F)$ et $Z(\bar{G}')^0(F)$. De nouveau, on peut remplacer ces groupes par les alg\`ebres de Lie correspondantes. Mais alors la premi\`ere est naturellement isomorphe \`a la somme directe des deux autres, cf. 5.1(1). On choisit des mesures qui se correspondent par cet isomorphisme.
  
 \bigskip

  Seul compte pour nous le comportement des fonctions $f'_{sc}$  et $\bar{f}[y]_{sc}$   dans des  voisinages de l'unit\'e invariants par conjugaison stable. On peut donc aussi bien descendre ces fonctions par l'exponentielle et  consid\'erer qu'elles sont d\'efinies sur des alg\`ebres de Lie. Pour simplifier, on ne change pas leur notation.
  Alors l'\'egalit\'e (2) \'equivaut \`a l'\'egalit\'e
$$(6) \qquad f'_{sc}=transfert(\sum_{y\in \dot{{\cal Y}}}c[y]d(y)\bar{f}[y]_{sc}),$$
o\`u $transfert$ est ici le transfert non standard de $\bar{G}'_{SC}$ \`a $G'_{\epsilon,SC}$. 

  Fixons un \'el\'ement $Z\in \mathfrak{z}(G'_{\epsilon},F)$ assez proche de $0$, que l'on \'ecrit $Z=Z_{1}+Z_{2}$ comme en 5.1. Soit $\boldsymbol{\delta}_{SC} \in D^{st}_{g\acute{e}om}(G'_{\epsilon,SC}(F))$ \`a support assez proche de $1$. On d\'efinit

- l'image $\boldsymbol{\delta}(Z)$ de $\boldsymbol{\delta}_{SC} $ par la compos\'ee de l'application $\iota_{G'_{\epsilon,SC},G'_{1,\epsilon_{1}}}^*$ et de la translation par $exp(Z)$;

- l'image $\boldsymbol{\delta}(Z)^{\tilde{G}'_{1}}$ de $\boldsymbol{\delta}(Z)$  par l'application $desc_{\epsilon_{1}}^{st,\tilde{G}'_{1}}$. 

C'est une distribution stable sur $\tilde{G}'_{1}(F)$, \`a support proche de la classe de conjugaison stable ${\cal O}^{\tilde{G}_{1}}$ de $\epsilon_{1}$.

On d\'efinit:

- l'image  $\bar{\boldsymbol{\delta}}_{SC}$ de $\boldsymbol{\delta}_{SC}$ par transfert endoscopique non standard (on doit descendre les distributions aux alg\`ebres de Lie pour d\'efinir ce transfert); c'est une distribution stable sur $\bar{G}'_{SC}(F)$;

- l'image $\bar{\boldsymbol{\delta}}(Z_{2})$ de $\bar{\boldsymbol{\delta}}_{SC}$ par la compos\'ee de l'application $\iota^*_{\bar{G}'_{SC},\bar{G}'}$ et de la translation par $exp(Z_{2})$;

 - pour $y\in {\cal Y}$, l'image $\boldsymbol{\delta}[y,Z_{2}]$ de $\bar{\boldsymbol{\delta}}(Z_{2})$ par transfert \`a $G_{\eta[y],SC}(F)$, avec la convention que ce transfert est nul si $\bar{{\bf G}}'$ n'est pas relevante pour $G_{\eta[y],SC}$;

- l'image $\boldsymbol{\delta}[y,Z]$ de $\boldsymbol{\delta}[y,Z_{2}]$ par la compos\'ee de l'application $\iota^*_{G_{\eta[y],SC},G_{\eta[y]}}$ et de la translation par $exp(Z_{1})$; c'est une distribution sur $G_{\eta[y]}(F)$;

- l'image $\boldsymbol{\delta}[y,Z]^{\tilde{G}}$ de  $\boldsymbol{\delta}[y,Z]$  par l'application $desc_{\eta[y]}^{\tilde{G},*}$.

C'est une distribution sur $\tilde{G}(F)$, \`a support proche de la classe de conjugaison stable ${\cal O}^{\tilde{G}}$ de $\eta$. 

A partir de l'\'egalit\'e (1), un calcul formel conduit \`a l'\'egalit\'e
$$(7) \qquad transfert(\boldsymbol{\delta}(Z)^{\tilde{G}'_{1}})=\sum_{y\in \dot{{\cal Y}}}c[y]d(y)\boldsymbol{\delta}[y,Z]^{\tilde{G}}.$$
On supprime de la notation les termes $Z$ dans le cas o\`u $Z=0$.

 \bigskip
 
 \subsection{Levi et descente de donn\'ees endoscopiques}

Soient $(G,\tilde{G},{\bf a})$ un triplet quelconque, $\tilde{M}$ un espace de Levi de $\tilde{G}$ et ${\bf M}'=(M',{\cal M}',\tilde{\zeta})$ une donn\'ee endoscopique elliptique et relevante de $(M,\tilde{M},{\bf a})$.  

D'apr\`es [I] paragraphe 3, on peut  r\'ealiser les objets duaux de la fa\c{c}on suivante. On fixe une paire de Borel \'epingl\'ee  $\hat{{\cal E}}=(\hat{B},\hat{T},(\hat{E}_{\alpha})_{\alpha\in \Delta})$ de $\hat{G}$. On note $\hat{\theta}$ l'automorphisme dual de $\theta^*$ qui conserve cette paire. On modifie l'action galoisienne de sorte que $\hat{{\cal E}}$ soit stable par cette action. On suppose que $\hat{M}$ est standard et on munit ce groupe de la paire de Borel \'epingl\'ee  $\hat{{\cal E}}^{M}=(\hat{B}\cap \hat{M},\hat{T},(\hat{E}_{\alpha})_{\alpha\in \Delta}^M)$. On suppose $\tilde{\zeta}=\zeta\hat{\theta}$ avec $\zeta\in \hat{T}$. On suppose $\zeta$ tel que le cocycle qui intervient dans la d\'efinition d'une donn\'ee endoscopique prenne ses valeurs dans $Z(\hat{G})$. On munit $\hat{M}'$ d'une paire de Borel \'epingl\'ee $\hat{{\cal E}}^{M'}$ dont la paire de Borel sous-jacente soit $(\hat{B}\cap \hat{M}',\hat{T}^{\hat{\theta},0})$.
 
Fixons un diagramme
 $$(\epsilon,B^{M'},T',B^M,T,\eta)$$
 cf. [I] 1.10. Les espaces ambiants sont $\tilde{M}'$ et $\tilde{M}$.
 Le terme $\eta$ est un \'el\'ement semi-simple de $\tilde{M}(F)$ et
$\epsilon$ est un \'el\'ement semi-simple de $\tilde{M}'(F)$. On suppose que $M'_{\epsilon}$ est  quasi-d\'eploy\'e.

Les deux paires de Borel $(B^M,T)$ de $M$ et $(\hat{B}\cap \hat{M},\hat{T})$ de $\hat{M}$ fournissent une identification $X_{*}(T)\simeq X^*(\hat{T})$, qui transporte l'action de $\theta=ad_{\eta}$ en $\hat{\theta}$. On introduit le groupe de Borel $B$ de $G$ contenant $T$  dont l'ensemble de coracines simples s'identifie \`a l'ensemble des racines simples de $\hat{B}$.  Il est clair que les deux paires de Borel $(B,T)$ de $G$ et $(\hat{B},\hat{T})$ de $\hat{G}$ fournissent la m\^eme identification  $X_{*}(T)\simeq X^*(\hat{T})$ que ci-dessus. On v\'erifie d'ailleurs que le groupe $P$ engendr\'e par $B$ et $M$ est un \'el\'ement de ${\cal P}(M)$ et que l'ensemble $\tilde{P}=P\eta$ est un \'el\'ement de ${\cal P}(\tilde{M})$.  De la paire de Borel $(\hat{B}\cap \hat{M}',\hat{T}^{\hat{\theta},0})$ de $\hat{M}'$ et de la paire $(B^{M'},T')$ de $M'$ se d\'eduit une identification $X_{*}(T')\simeq X^*(\hat{T}^{\hat{\theta},0})$.

  Posons $(\bar{B},\bar{T})=(B\cap G_{\eta},T\cap G_{\eta})$, compl\'etons cette paire  en une paire de Borel \'epingl\'ee de $G_{\eta}$ et introduisons l'action galoisienne quasi-d\'eploy\'ee sur $G_{\eta}$ qui la conserve. Elle est de la forme $\sigma\mapsto \sigma_{\bar{G}}=ad_{\bar{u}(\sigma)}\circ \sigma_{G}$, cf. 5.1. Parce que $P\cap G_{\eta}$ est d\'efini sur $F$ pour l'action galoisienne naturelle, on a n\'ecessairement $\bar{u}(\sigma)\in M_{\eta,sc}$. Donc $M_{\eta}$ est stable par  l'action $\sigma\mapsto \sigma_{\bar{G}}$ et les deux actions galoisiennes co\"{\i}ncident sur $Z(M_{\eta})$.  On note  $\bar{G}$, resp. $\bar{M}$, le groupe $G_{\eta}$, resp.  $M_{\eta}$, muni de  l'action $\sigma\mapsto \sigma_{\bar{G}}$. Compl\'etons de m\^eme la paire de Borel $(B,T)$ de $G$ en une paire de Borel \'epingl\'ee et introduisons l'action galoisienne quasi-d\'eploy\'ee sur $G$ qui conserve celle-ci. Elle est de la forme $\sigma\mapsto  \sigma_{G^*}=ad_{u(\sigma)}\circ \sigma_{G}$, o\`u $u(\sigma)\in G_{SC}$. Toujours parce que $P$ est d\'efini sur $F$ pour l'action galoisienne naturelle, on a $u(\sigma)\in M_{sc}$. Donc $M$ est stable par l'action $\sigma\mapsto \sigma_{G^*}$ et les deux actions co\"{\i}ncident sur $Z(M)$. On a $\sigma_{\bar{G}}=ad_{\bar{u}(\sigma)u(\sigma)^{-1}}\circ \sigma_{G^*}$.  L'\'el\'ement $\bar{u}(\sigma)u(\sigma)^{-1}$  appartient \`a $M_{sc}$. Il s'envoie donc sur un \'el\'ement du groupe de Weyl $W^M$ que l'on note $\omega_{\bar{M}}(\sigma)$. Il est fixe par $\theta$.  
    Le groupe $\hat{\bar{M}}$ est le commutant dans $\hat{\bar{G}}$ de l'image dans $\hat{\bar{T}}$ de $Z(\hat{M})$.

  On peut d\'efinir comme en 5.1 une donn\'ee endoscopique   de $\bar{M}_{SC}$  associ\'ee \`a la donn\'ee ${\bf M}'$ et au diagramme   $(\epsilon,B^{M'},T',B^M,T, \eta)$.    Mais on peut aussi  refaire la m\^eme construction en rempla\c{c}ant le groupe $\bar{M}_{SC}$ par l'image r\'eciproque $\bar{M}_{sc}$ de $\bar{M}$ dans $\bar{G}_{SC}$. On obtient une donn\'ee endoscopique de $\bar{M}_{sc}$. C'est celle-l\`a que nous noterons  $\bar{{\bf M}}'=(\bar{M}',\bar{{\cal M}}',\bar{\zeta})$. 
  
   Fixons une paire de Borel  $(B^{M',\flat},T^{\flat})$ de $M'_{\epsilon}$ d\'efinie sur $F$  et un \'el\'ement  $m'\in M'_{\epsilon,SC}$ tel que $ad_{m'}$ envoie cette paire sur $(B^{M'}\cap M'_{\epsilon},T')$. Fixons une paire de Borel $(\bar{B}^{\bar{M}'},\bar{T}')$ de $\bar{M}'$ d\'efinie sur $F$.  Pour une raison qui appara\^{\i}tra plus loin,  on note  $T^{\flat}_{M-sc}$, resp.  $\bar{T}'_{M'-sc}$,  l'image r\'eciproque de $T^{\flat}$ dans $ M'_{\epsilon,SC}$, resp.    l'image r\'eciproque de $ \bar{T}'$ dans $\bar{M}'_{SC}$.  Comme en 5.1, on a un isomorphisme d'alg\`ebres de Lie 
 $$j:\bar{\mathfrak{t}}'_{M'-sc}\oplus \mathfrak{z}(\bar{M}')\oplus \mathfrak{z}(\bar{G})\simeq \mathfrak{t}^{\flat}_{M-sc}\oplus \mathfrak{z}(M'_{\epsilon})$$
  qui est \'equivariant pour les actions galoisiennes. On note $j_{*}:X_{*}(\bar{T}'_{M'-sc})\otimes {\mathbb Q}\to X_{*}(T^{\flat}_{M-sc})\otimes {\mathbb Q}$ l'isomorphisme sous-jacent \`a la restriction de l'homomorphisme $j$ \`a $\bar{\mathfrak{t}}'_{M'-sc}$. Le triplet $(\bar{M}'_{SC},M'_{\epsilon,SC},j_{*})$  est un triplet endoscopique non standard.
On a  le diagramme suivant
 $$\begin{array}{ccc}\mathfrak{z}(M)^{\Gamma_{F},\theta}&=&\mathfrak{z}(M)^{\Gamma_{F},\theta}\\ \downarrow&&\\ \mathfrak{z}(\bar{M})^{\Gamma_{F}}&&\downarrow\\ &&\\ \parallel&&\mathfrak{z}(M')^{\Gamma_{F}}\\ &&\\ \mathfrak{z}(\bar{M}_{sc})^{\Gamma_{F}}\oplus\mathfrak{z}(\bar{G})^{\Gamma_{F}}&& \downarrow\\ \downarrow&&\\ \mathfrak{z}(\bar{M}')^{\Gamma_{F}}\oplus \mathfrak{z}(\bar{G})^{\Gamma_{F}}&\simeq&\mathfrak{z}(M'_{\epsilon})^{\Gamma_{F}}\\ \end{array}$$
 On suppose d\'esormais
 
 (1) $\epsilon$ est elliptique dans $\tilde{M}'(F)$.
 
 Les fl\`eches de droite du diagramme sont des isomorphismes parce que ${\bf M}'$ est une donn\'ee elliptique de $(M,\tilde{M},{\bf a})$ et que $\epsilon$ est elliptique dans $\tilde{M}'(F)$. Les fl\`eches de gauche sont donc aussi des isomorphismes. Pour la fl\`eche du haut, cela entra\^{\i}ne que $\eta$ est elliptique dans $\tilde{M}(F)$. Pour celle du bas, cela entra\^{\i}ne que ${\bf M}'$ est une donn\'ee endoscopique elliptique de $\bar{M}_{sc}$. Rappelons que, dans cette situation non tordue et sur un corps local non-archim\'edien, une telle donn\'ee est automatiquement relevante.  
  
 L'ellipticit\'e de $\eta$ entra\^{\i}ne dualement 
    que
  
  (2)  l'homomorphisme naturel
  $$Z(\hat{M})^{\Gamma_{F},\hat{\theta},0}\to Z(\hat{\bar{M}})^{\Gamma_{F},0}$$
  est surjectif et de noyau fini. 
  
  Puisque $Z(\hat{G})^{\Gamma_{F},\hat{\theta}}$ s'envoie \'evidemment dans $Z(\bar{G})^{\Gamma_{F}}$ et puisque les quotients ci-dessous sont connexes, on a aussi
  
  (3) l'homomorphisme
  $$Z(\hat{M})^{\Gamma_{F},\hat{\theta}}/Z(\hat{G})^{\Gamma_{F},\hat{\theta}}\to Z(\hat{\bar{M}})^{\Gamma_{F}}/Z(\hat{\bar{G}})^{\Gamma_{F}}$$
  est surjectif.

   Soit $\tilde{s}\in \tilde{\zeta}Z(\hat{M})^{\Gamma_{F},\hat{\theta}}/Z(\hat{G})^{\Gamma_{F},\hat{\theta}}$. Le groupe $\hat{G'}(\tilde{s})$ est composante neutre du centralisateur de $\tilde{s}$ dans $\hat{G}$. Il est muni de la paire de Borel $(\hat{B}\cap \hat{G}'(\tilde{s}),\hat{T}^{\hat{\theta},0})$. De m\^eme que l'on a construit le groupe $B$, on construit un groupe de Borel $B'$ de $G'(\tilde{s})$ de sorte que l'identification   $X_{*}(T')\simeq X^*(\hat{T}^{\hat{\theta},0})$ d\'eduite des  paires $(B',T')$ et  $(\hat{B}\cap \hat{G}'(\tilde{s}),\hat{T}^{\hat{\theta},0})$ soit la m\^eme que celle ci-dessus. De nouveau, le groupe $P'$ engendr\'e par $M'$ et $B'$ est un \'el\'ement de ${\cal P}^{G'(\tilde{s})}(M')$.  
On voit alors que le sextuplet $(\epsilon,B',T',B,T,\eta)$ est  encore un diagramme, avec cette fois pour groupes ambiants $G'(\tilde{s})$ et $G$, donnant naissance au m\^eme homomorphisme $\xi$.  A partir du diagramme $(\epsilon,B',T',B,T,\eta)$ et de ${\bf G}'(\tilde{s})$, on d\'efinit  une donn\'ee endoscopique de $\bar{G}_{SC}$. L'\'el\'ement $s$ s'envoie sur un \'el\'ement  $\bar{s}\in\hat{\bar{T}}_{ad}$, plus pr\'ecis\'ement $\bar{s}\in \bar{\zeta}Z(\hat{\bar{M}})^{\Gamma_{F}}$.  En consid\'erant les d\'efinitions de [W1] 3.5, on s'aper\c{c}oit que la donn\'ee endoscopique de $\bar{G}_{SC}$ que l'on construit ainsi n'est autre que celle associ\'ee \`a la donn\'ee $\bar{{\bf M}}'$ du Levi $\bar{M}_{sc}$ et \`a cet \'el\'ement $\bar{s}$.  Nous la noterons $\bar{{\bf G}}'(\bar{s})=(\bar{G}'(\bar{s}),\bar{{\cal G}}'(\bar{s}),\bar{s})$.  Fixons un sous-groupe parabolique de $\bar{G}'(\bar{s})$ d\'efini sur $F$ et de composante de Levi $\bar{M}'$. On note $\bar{B}'(\bar{s})$ le sous-groupe de Borel de $\bar{G}'(\bar{s})$ qui est inclus dans ce parabolique et v\'erifie $\bar{B}'(\bar{s})\cap \bar{M}'=\bar{B}^{\bar{M}'}$. Le couple $(\bar{B}'(\bar{s}),\bar{T}')$ est une paire de Borel de $\bar{G}'(\bar{s})$ qui est d\'efinie sur $F$.     On pose $B^{\flat}=ad_{m'}^{-1}(B')$. Alors $(B^{\flat},T^{\flat})$ est une paire de Borel de $G'(\tilde{s})_{\epsilon}$ d\'efinie sur $F$.  On note $T^{\flat}_{sc}$, resp.  $\bar{T}'_{sc}$,  l'image r\'eciproque de $T^{\flat}$ dans $ G'(\tilde{s})_{\epsilon,SC}$, resp.    l'image r\'eciproque de $ \bar{T}'$ dans $\bar{G}'(\bar{s})_{SC}$. On a un isomorphisme
$$j:\bar{\mathfrak{t}}'_{sc}\oplus \mathfrak{z}(\bar{G}'(\bar{s}))\oplus \mathfrak{z}(\bar{G})\simeq \mathfrak{t}^{\flat}_{sc}\times \mathfrak{z}(G'(\tilde{s})_{\epsilon})$$
qui v\'erifie des propri\'et\'es analogues \`a celles ci-dessus. C'est le m\^eme que plus haut, modulo les identifications
$$\bar{\mathfrak{t}}'_{sc}\simeq \bar{\mathfrak{t}}'_{M-sc}\oplus \mathfrak{z}(\bar{M}'(\bar{s})_{sc}),$$
o\`u $\bar{M}'(\bar{s})_{sc}$ est l'image r\'eciproque de $\bar{M}'$ dans $\bar{G}'(\bar{s})_{SC}$ et
$$\mathfrak{t}^{\flat}_{sc}\simeq \mathfrak{t}^{\flat}_{M'-sc}\oplus \mathfrak{z}(M'(\tilde{s})_{\epsilon,sc}),$$
o\`u $M'(\tilde{s})_{\epsilon,sc}$ est l'image r\'eciproque de $M'_{\epsilon}$ dans $G'(\tilde{s})_{\epsilon,SC}$. De nouveau (avec un l\'eger abus de notation),  on note $j_{*}: X_{*}(\bar{T}'_{sc})\otimes {\mathbb Q}\to X_{*}(T^{\flat}_{sc})\otimes {\mathbb Q}$ l'isomorphisme sous-jacent \`a la restriction de l'homomorphisme $j$ \`a $\bar{\mathfrak{t}}'_{sc}$. Le triplet $(\bar{G}'(\bar{s})_{SC},G'(\tilde{s})_{\epsilon,SC},j_{*})$  est un triplet endoscopique non standard. R\'esumons les isomorphismes obtenus:
$$(4) \quad \begin{array}{ccccccccc}&&\mathfrak{z}(G'(\tilde{s})_{\epsilon})&\to&\mathfrak{z}(\bar{G})\oplus \mathfrak{z}(\bar{G}'(\bar{s}))&&\mathfrak{z}(\bar{G})&\to &\mathfrak{z}(\bar{G})\\\mathfrak{z}(M'_{\epsilon})&\to&\oplus&&\oplus&\to&\oplus&&\oplus\\&&\mathfrak{z}(M'(\tilde{s})_{\epsilon,sc})&\to&\mathfrak{z}(\bar{M}'(\bar{s})_{sc})&&\mathfrak{z}(\bar{G}'(\bar{s}))\oplus\mathfrak{z}(\bar{M}'(\bar{s})_{sc})&\to&\mathfrak{z}(\bar{M}')\\ \end{array}$$

Notons ${\cal S}$ l'ensemble des $\tilde{s}\in \tilde{\zeta}Z(\hat{M})^{\Gamma_{F},\hat{\theta}}/Z(\hat{G})^{\Gamma_{F},\hat{\theta}}$ tels que ${\bf G}'(\tilde{s})$ soit elliptique et $A_{G'(\tilde{s})}=A_{G'(\tilde{s})_{\epsilon}}$. Montrons que

(5) si $A_{G_{\eta}}\not=A_{\tilde{G}}$, l'ensemble ${\cal S}$ est vide;

(6) supposons $A_{G_{\eta}}=A_{\tilde{G}}$; alors ${\cal S}$ est \'egal \`a l'ensemble des $\tilde{s}\in \tilde{\zeta}Z(\hat{M})^{\Gamma_{F},\hat{\theta}}/Z(\hat{G})^{\Gamma_{F},\hat{\theta}}$ tels que la donn\'ee d\'eduite ${\bf G}'(\bar{s})$ de $\bar{G}_{SC}$ soit elliptique.

Pour $\tilde{s}\in \tilde{\zeta}Z(\hat{M})^{\Gamma_{F},\hat{\theta}}/Z(\hat{G})^{\Gamma_{F},\hat{\theta}}$, on a un diagramme similaire \`a celui \'ecrit plus haut
 $$\begin{array}{ccc}\mathfrak{z}(G)^{\Gamma_{F},\theta}&=&\mathfrak{z}(G)^{\Gamma_{F},\theta}\\ \downarrow&& \downarrow\\ \mathfrak{z}(\bar{G})^{\Gamma_{F}}&&\mathfrak{z}(G'(\tilde{s}))^{\Gamma_{F}}\\ \downarrow&&\downarrow \\ &&\\ \mathfrak{z}(\bar{G}'(\bar{s}))^{\Gamma_{F}}\oplus\mathfrak{z}(\bar{G})^{\Gamma_{F}}&\simeq&\mathfrak{z}(G'(\tilde{s})_{\epsilon})^{\Gamma_{F}}\\ \end{array}$$
 L'\'el\'ement $\tilde{s}$ appartient \`a ${\cal S}$ si et seulement si les deux fl\`eches de droite sont des isomorphismes. C'est \'equivalent \`a ce que celles de gauche  soient elles-aussi des isomorphismes. Que la fl\`eche du haut \`a gauche soit un isomorphisme signifie que $A_{G_{\eta}}=A_{\tilde{G}}$. Cette condition est ind\'ependante de $\tilde{s}$. Si elle n'est pas v\'erifi\'ee, l'ensemble ${\cal S}$ est donc vide. Si elle est v\'erifi\'ee, la seule condition est que la fl\`eche du bas \`a gauche soit un isomorphisme, ce qui \'equivaut \`a l'ellipticit\'e de ${\bf G}'(\bar{s})$.
 \bigskip

\bigskip

\subsection{Facteurs de transfert et transfert des distributions}
Pour fixer les notations, reprenons bri\`evement les constructions de 5.1 et 5.2.
Soit $\tilde{s}\in \tilde{\zeta}Z(\hat{M})^{\Gamma_{F},\hat{\theta}}/Z(\hat{G})^{\Gamma_{F},\hat{\theta}}$. On fixe des donn\'ees auxiliaires $G_{1}'(\tilde{s}),...,\Delta_{1}(\tilde{s})$ pour la donn\'ee ${\bf G}'(\tilde{s})$. On note $M_{1}'(\tilde{s})$, resp. $\tilde{M}_{1}'(\tilde{s})$, l'image r\'eciproque de $M'$ dans $G_{1}'(\tilde{s})$, resp. de $\tilde{M}'$ dans $\tilde{G}'_{1}(\tilde{s})$. On fixe une image r\'eciproque $\epsilon_{1}(\tilde{s})$ de $\epsilon$ dans $\tilde{M}_{1}'(\tilde{s};F)$, ainsi qu'une d\'ecomposition d'alg\`ebres de Lie
 $$  \mathfrak{g}'_{1}(\tilde{s})=\mathfrak{c}_{1}(\tilde{s})\oplus \mathfrak{g}'(\tilde{s}).$$
 On en d\'eduit des d\'ecompositions
 $$\mathfrak{g}_{1}'(\tilde{s})_{\epsilon_{1}(\tilde{s})}=\mathfrak{c}_{1}(\tilde{s})\oplus \mathfrak{g}'(\tilde{s})_{\epsilon},$$
$$\mathfrak{m}'_{1}(\tilde{s})=\mathfrak{c}_{1}(\tilde{s})\oplus \mathfrak{m}',$$
 $$\mathfrak{m}_{1}'(\tilde{s})_{\epsilon_{1}(\tilde{s})}=\mathfrak{c}_{1}(\tilde{s})\oplus \mathfrak{m}'_{\epsilon}.$$

 Pour $y\in {\cal Y}$, on fixe  un facteur de transfert sur $G'(\bar{s},F)\times G_{\eta[y],SC}(F)$ que l'on note $\Delta(\bar{s},y)$. Alors il existe $d(\tilde{s},y)\in {\mathbb C}^{\times}$ tel que,  pour des donn\'ees comme en 5.1, on ait l'\'egalit\'e
 $$(1) \qquad d(\tilde{s},y)\Delta(\bar{s},y)(exp(\bar{Y}),exp(X[y]_{sc}))=\Delta(\tilde{s})_{1}(exp(Y)\epsilon_{1}(\tilde{s}),exp(X[y])\eta[y])$$
 pourvu que $X$ et $Y$ soient assez proches de $0$. 
 
 Pour appliquer la formule 5.2(7) \`a la donn\'ee ${\bf G}'(\tilde{s})$, on glisse des termes $\tilde{s}$ ou $\bar{s}$ dans les notations.  Fixons un \'el\'ement $Z\in \mathfrak{z}(G'(\tilde{s})_{\epsilon},F)$ assez proche de $0$, que l'on \'ecrit $Z=Z_{1}+Z_{2}$  avec $Z_{1}\in \mathfrak{z}(\bar{G};F)$ et $Z_{2}\in \mathfrak{z}(\bar{G}'(\bar{s});F)$. Soit $\boldsymbol{\delta}(\tilde{s})_{SC} \in D^{st}_{g\acute{e}om}(G'(\tilde{s})_{\epsilon,SC}(F))$ \`a support assez proche de $1$. On d\'efinit

- l'image $\boldsymbol{\delta}(\tilde{s},Z)$ de $\boldsymbol{\delta}(\tilde{s})_{SC} $ par la compos\'ee de l'application $\iota_{G'(\tilde{s})_{\epsilon,SC},G'_{1}(\tilde{s})_{\epsilon_{1}(\tilde{s})}}^*$ et de la translation par $exp(Z)$;

- l'image $\boldsymbol{\delta}(\tilde{s},Z)^{\tilde{G}'_{1}(\tilde{s})}$ de $\boldsymbol{\delta}(\tilde{s},Z)$  par l'application $desc_{\epsilon_{1}(\tilde{s})}^{st,\tilde{G}'_{1}(\tilde{s})}$;

- l'image  $\bar{\boldsymbol{\delta}}(\bar{s})_{SC}$ de $\boldsymbol{\delta}(\tilde{s})_{SC}$ par transfert endoscopique non standard;

- l'image $\bar{\boldsymbol{\delta}}(\bar{s},Z_{2})$ de $\bar{\boldsymbol{\delta}}(\bar{s})_{SC}$ par la compos\'ee de l'application $\iota^*_{\bar{G}'(\bar{s})_{SC},\bar{G}'(\bar{s})}$ et de la translation par $exp(Z_{2})$;

 - pour $y\in {\cal Y}$, l'image $\boldsymbol{\delta}[y,Z_{2}]$ de $\bar{\boldsymbol{\delta}}(\bar{s},Z_{2})$ par transfert \`a $G_{\eta[y],SC}(F)$;

- l'image $\boldsymbol{\delta}[y,Z]$ de $\boldsymbol{\delta}[y,Z_{2}]$ par la compos\'ee de l'application $\iota^*_{G_{\eta[y],SC},G_{\eta[y]}}$ et de la translation par $exp(Z_{1})$;

- l'image $\boldsymbol{\delta}[y,Z]^{\tilde{G}}$ de  $\boldsymbol{\delta}[y,Z]$  par l'application $desc_{\eta[y]}^{\tilde{G},*}$.

 La formule 5.2(7) devient 
$$(2) \qquad transfert(\boldsymbol{\delta}(\tilde{s},Z)^{\tilde{G}'_{1}(\tilde{s})})=\sum_{y\in \dot{{\cal Y}}}c[y]d(\tilde{s},y)\boldsymbol{\delta}[y,Z]^{\tilde{G}}.$$
On supprime de la notation les termes $Z$ dans le cas o\`u $Z=0$.

On peut remplacer $\tilde{G}$ par $\tilde{M}$ dans les constructions ci-dessus. On  ajoute des exposants $M$ pour indiquer les analogues pour $\tilde{M}$ des termes pr\'ec\'edemment d\'efinis pour $\tilde{G}$.
On obtient une assertion o\`u la donn\'ee de d\'epart est un \'el\'ement $\boldsymbol{\delta}_{SC}\in D^{st}_{g\acute{e}om}(M'_{\epsilon,SC}(F))$. On aura besoin d'une variante o\`u le groupe $M'_{\epsilon,SC}$ est remplac\'ee par le groupe interm\'ediaire $M'(\tilde{s})_{\epsilon,sc}$. Partons ainsi d'un \'el\'ement  $\boldsymbol{\delta}(\tilde{s})_{sc} \in D^{st}_{unip}(M'(\tilde{s})_{\epsilon,sc}(F))$ (le cas unipotent nous suffira). On d\'efinit

-  $\boldsymbol{\delta}(\tilde{s})=\iota_{M'(\tilde{s})_{\epsilon,sc},M'_{1}(\tilde{s}){\epsilon_{1}(\tilde{s})}}^*(\boldsymbol{\delta}(\tilde{s})_{sc} )$;

- l'image $\boldsymbol{\delta}(\tilde{s})^{\tilde{M}'_{1}(\tilde{s})}$ de $\boldsymbol{\delta}(\tilde{s})$  par l'application $desc_{\epsilon_{1}(\tilde{s})}^{st,\tilde{M}'_{1}(\tilde{s})}$. 

C'est une distribution stable sur $\tilde{M}'_{1}(\tilde{s};F)$,  dont les \'el\'ements du support ont des parties semi-simples dans la classe de conjugaison stable ${\cal O}_{\tilde{M}'_{1}(\tilde{s})}$ de $\epsilon_{1}(\tilde{s})$.

On d\'efinit:

- l'image  $\bar{\boldsymbol{\delta}}(\bar{s})_{sc}$ de $\boldsymbol{\delta}(\tilde{s})_{sc}$ par transfert endoscopique non standard; c'est une distribution stable sur $\bar{M}'(\bar{s})_{sc}(F)$;

-   $\bar{\boldsymbol{\delta}}(\bar{s}) =\iota^*_{\bar{M}'(\bar{s})_{sc},\bar{M}'}(\bar{\boldsymbol{\delta}}(\bar{s})_{sc})$;

 - pour $y\in {\cal Y}^M$, l'image $\boldsymbol{\delta}[y]$ de $\bar{\boldsymbol{\delta}}(\bar{s})$ par la compos\'ee du transfert \`a $M_{\eta[y],sc}(F)$ et de l'application $\iota^*_{M_{\eta[y],sc},M_{\eta[y]}}$; c'est une distribution sur $M_{\eta[y]}(F)$;

- l'image $\boldsymbol{\delta}[y]^{\tilde{M}}$ de  $\boldsymbol{\delta}[y,]$  par l'application $desc_{\eta[y]}^{\tilde{M},*}$.

C'est une distribution sur $\tilde{M}(F)$,  dont les \'el\'ements du support ont des parties semi-simples dans  la classe de conjugaison stable ${\cal O}$ de $\eta$. 

 On a alors l'\'egalit\'e
$$(3) \qquad transfert(\boldsymbol{\delta}(\tilde{s})^{\tilde{M}'_{1}(\tilde{s})})=\sum_{y\in \dot{{\cal Y}}^M}c^M[y]d(\tilde{s},y)\boldsymbol{\delta}[y]^{\tilde{M}}.$$

Cette variante se d\'eduit facilement de la formule pr\'ec\'edente. Il suffit d'utiliser la formule 3.7(4) qui permet de permuter un transfert endoscopique avec une application telle que $\iota^*_{M'_{\epsilon,SC},M'(\tilde{s})_{\epsilon,sc}}$, ainsi qu'une formule analogue concernant le transfert non standard, laquelle se prouve de la m\^eme fa\c{c}on. On laisse les d\'etails au lecteur.

\bigskip

\subsection{Applications de transition}
Notons ${\cal O}_{M'}$ la classe de conjugaison stable de $\epsilon$. Soit $\boldsymbol{\delta}\in D_{g\acute{e}om}^{st}({\bf M}')$ \`a support proche de ${\cal O}_{M'}$. Fixons des donn\'ees auxiliaires $M'_{1},...,\Delta_{1}$ pour ${\bf M}'$ et un \'el\'ement $\epsilon_{1}\in \tilde{M}'_{1}(F)$ se projetant sur $\epsilon$. Alors $\boldsymbol{\delta}$ s'identifie \`a un \'el\'ement de $D_{g\acute{e}om,\lambda_{1}}^{st}(\tilde{M}'_{1}(F))$. Comme on l'a dit en [II] 1.10,  on a une surjection $D_{g\acute{e}om}^{st}(\tilde{M}'_{1}(F))\to D_{g\acute{e}om,\lambda_{1}}^{st}(\tilde{M}'_{1}(F))$. L'action  par translations sur l'espace de d\'epart d'un \'el\'ement $c\in C_{1}(F)$ se traduit sur l'espace d'arriv\'ee par la multiplication par $\lambda_{1}(c)$. Puisque le support de $\boldsymbol{\delta}$ est voisin de ${\cal O}_{M'}$, on peut donc relever cette distribution en un \'el\'ement $\boldsymbol{\delta}_{1}\in D_{g\acute{e}om}^{st}(\tilde{M}'_{1}(F))$ \`a support voisin de la classe de conjugaison stable de $\epsilon_{1}$.  C'est l'image par l'application $desc^{st,\tilde{M}'_{1}}_{\epsilon_{1}}$ 
 d'un \'el\'ement de $D_{g\acute{e}om}^{st}(M'_{1,\epsilon_{1}}(F))$ \`a support proche de $1$. 
 Un tel \'el\'ement est combinaison lin\'eaire de termes $exp(Z)\iota^*_{M'_{\epsilon,SC},M'_{1,\epsilon_{1}}}(\boldsymbol{\delta}_{\epsilon,SC})$, avec $Z\in \mathfrak{z}(M'_{1,\epsilon_{1}};F)$ et $\boldsymbol{\delta}_{\epsilon,SC}\in D^{st}_{g\acute{e}om}(M'_{\epsilon,SC}(F))$, \`a support proche de $1$ (le produit par $exp(Z)$ signifie la translation par cet \'el\'ement). Pour simplifier, on suppose que la combinaison lin\'eaire est r\'eduite \`a un seul terme, c'est-\`a-dire que 
 $$\boldsymbol{\delta}_{1}=desc^{st,\tilde{M}'_{1}}_{\epsilon_{1}}(exp(Z)\iota^*_{M'_{\epsilon,SC},M'_{1,\epsilon_{1}}}(\boldsymbol{\delta}_{\epsilon,SC})).$$
 Fixons une d\'ecomposition
$$\mathfrak{m}'_{1,\epsilon_{1}}=\mathfrak{c}_{1}\oplus \mathfrak{m}'_{\epsilon}.$$ 
Les translations par $\mathfrak{c}_{1}(F)$ ne comptent pas par le m\^eme argument que ci-dessus (d'ailleurs, le compos\'e de $\lambda_{1}$ et de l'exponentielle est \'egal \`a $1$ au voisinage de $0$). On peut donc supposer que l'\'el\'ement $Z$ appartient \`a $\mathfrak{z}(M'_{\epsilon};F)$. Posons alors $\boldsymbol{\delta}_{\epsilon}=exp(Z)\iota^*_{M'_{\epsilon,SC},M'_{1,\epsilon_{1}}}(\boldsymbol{\delta}_{\epsilon,SC})$. On peut consid\'erer que c'est un \'el\'ement de $D_{g\acute{e}om}^{st}(M'_{\epsilon}(F))$ \`a support proche de $1$.

Soit $\tilde{s}\in \tilde{\zeta}Z(\hat{M})^{\Gamma_{F},\hat{\theta}}/Z(\hat{G})^{\Gamma_{F},\hat{\theta}}$. Utilisons les donn\'ees auxiliaires introduites au paragraphe pr\'ec\'edent, que l'on restreint \`a $M'$. En rempla\c{c}ant $\tilde{M}'_{1}$ par $\tilde{M}'_{1}(\tilde{s})$, on construit un autre \'el\'ement $\boldsymbol{\delta}(\tilde{s})_{\epsilon}$.  Les distributions $\boldsymbol{\delta}_{\epsilon}$ et $\boldsymbol{\delta}(\tilde{s})_{\epsilon}$ appartiennent au m\^eme espace $D_{g\acute{e}om}^{st}(M'_{\epsilon}(F))$ et proviennent par descente d'un m\^eme \'el\'ement $\boldsymbol{\delta}\in D_{g\acute{e}om}^{st}({\bf M}')$. Mais elles ne sont pas forc\'ement \'egales.  Rappelons en effet la construction de l'isomorphisme compos\'e
$$SI_{\lambda_{1}}(\tilde{M}'_{1}(F))\simeq SI({\bf M}')\simeq SI_{\lambda_{1}(\tilde{s})}(\tilde{M}_{1}'(\tilde{s};F)).$$
Il provient d'un isomorphisme
 $$\begin{array}{ccc}C_{c,\lambda_{1}}^{\infty}(\tilde{M}'_{1}(F))&\simeq &C_{c,\lambda_{1}(\tilde{s})}^{\infty}(\tilde{M}_{1}'(\tilde{s};F))\\ \phi&\mapsto&\phi(\tilde{s}).\\ \end{array}$$
 Soient $\delta_{1}\in \tilde{M}'_{1}(F)$ et $\delta_{1}(\tilde{s})\in \tilde{M}'_{1}(\tilde{s};F)$ deux \'el\'ements au-dessus d'un m\^eme point de $\tilde{M}'(F)$. Alors on a une \'egalit\'e
 $$\phi(\tilde{s})(\delta_{1}(\tilde{s}))=\tilde{\lambda}(\tilde{s})(\delta_{1},\delta_{1}(\tilde{s}))\phi(\delta_{1}),$$
 o\`u $\tilde{\lambda}(\tilde{s})$ est une fonction de recollement d\'efinie en [I] 2.5. Cette fonction de recollement est localement constante. Parce que l'on travaille dualement avec des distributions, un d\'evissage des d\'efinitions conduit \`a l'\'egalit\'e
 $$(1) \qquad \boldsymbol{\delta}(\tilde{s})_{\epsilon}=d(\tilde{s})\boldsymbol{\delta}_{\epsilon},$$
 o\`u
 $$d(\tilde{s})=\tilde{\lambda}(\tilde{s})(\epsilon_{1},\epsilon_{1}(\tilde{s}))^{-1}.$$
 On peut calculer ce terme en fixant un \'el\'ement $Y\in \mathfrak{m}'_{\epsilon}(F)$ en position g\'en\'erale et elliptique. Il lui correspond un \'el\'ement $X\in \mathfrak{m}_{\eta}(F)$ par la construction de 5.1 appliqu\'ee en rempla\c{c}ant $\tilde{G}$ par $\tilde{M}$ et $y$ par $1$. 
  A l'aide des d\'ecompositions fix\'ees, on peut identifier $Y$ soit \`a un \'el\'ement de $\mathfrak{m}'_{1,\epsilon_{1}}(F)$, soit \`a un \'el\'ement de $\mathfrak{m}_{1}'(\tilde{s})_{\epsilon_{1}(\tilde{s})}(F)$. D'apr\`es [I] lemme 2.5, on a l'\'egalit\'e
 $$(2) \qquad d(\tilde{s})=\frac{\Delta_{1}(exp(Y)\epsilon_{1},exp(X)\eta)}{\Delta_{1}(\tilde{s})(exp(Y)\epsilon_{1}(\tilde{s}),exp(X)\eta)}.$$

 \section{Triplets endoscopiques non standard}

\bigskip
\subsection{Apparition des triplets endoscopiques non standard}
   Rappelons la notion de triplet endoscopique non standard introduite dans [W1] 1.7. On consid\`ere deux groupes r\'eductifs connexes $G_{1}$ et $G_{2}$ d\'efinis sur $F$, simplement connexes et quasi-d\'eploy\'es. Pour $i=1,2$, soit $(B_{i},T_{i})$ une paire de Borel de $G_{i}$ d\'efinie sur $F$. On introduit l'ensemble $\Sigma(T_{i})$ des racines de $T_{i}$ dans $G_{i}$ et l'ensemble $\check{\Sigma}(T_{i})$ des coracines. On note $\alpha\mapsto \check{\alpha}$ la bijection usuelle entre ces ensembles. On suppose donn\'ees une bijection $\tau:\Sigma(T_{2})\to \Sigma(T_{1})$, une fonction $b:\Sigma(T_{2})\to {\mathbb Q}_{>0}$ et  un isomorphisme $j_{*}:X_{*}(T_{1})\otimes_{{\mathbb Z}}{\mathbb Q}\to X_{*}(T_{2})\otimes_{{\mathbb Z}}{\mathbb Q}$. On note $j^{*}:X^{*}(T_{2})\otimes_{{\mathbb Z}}{\mathbb Q}\to X^{*}(T_{1})\otimes_{{\mathbb Z}}{\mathbb Q}$ l'isomorphisme dual. On impose les conditions suivantes:

- pour $\alpha_{2}\in \Sigma(T_{2})$, $\alpha_{2}$ est positif pour $B_{2}$ si et seulement si $\tau(\alpha_{2})$ est positif pour $B_{1}$ (cette condition ne figure pas dans [W1] 1.7 mais peut \^etre ajout\'ee d'apr\`es le lemme de cette r\'ef\'erence);

- $j_{*}$ (donc aussi $j^*$) est \'equivariant pour les actions galoisiennes;

- pour tout $\alpha_{2}\in \Sigma(T_{2})$, $j^*(\alpha_{2})=b(\alpha_{2})\tau(\alpha_{2})$;

- pour tout $\alpha_{1}\in \Sigma(T_{1})$, $j_{*}(\check{\alpha}_{1})=b(\alpha_{2})\check{\alpha}_{2}$, o\`u $\alpha_{2}=\tau^{-1}(\alpha_{1})$.

A ces conditions, on dit que $(G_{1},G_{2},j_{*})$ est un triplet endoscopique non standard.

 De tels triplets  interviennent dans nos constructions comme on l'a vu en 5.1. Soit $(G,\tilde{G},{\bf a})$ un de nos triplets  comme en [II] 1.1 et  ${\bf G}'=(G',{\cal G}',\tilde{s})$ une donn\'ee endoscopique elliptique de $(G,\tilde{G},{\bf a})$. Soient $(\epsilon,B',T',B,T,\eta)$ un diagramme reliant  un \'el\'ement semi-simple $\epsilon\in \tilde{G}'(F)$ \`a un \'el\'ement   semi-simple  $\eta\in \tilde{G}'(F)$. Supposons  $G'_{\epsilon}$ quasi-d\'eploy\'e. A l'aide de ces donn\'ees, on construit une donn\'ee endoscopique $\bar{{\bf G}}' =(\bar{G}',\bar{{\cal G}}',\bar{s})$ de $G_{\eta,SC}$. Posons   $G_{1}=\bar{G}'_{SC}$ et $G_{2}=G'_{\epsilon,SC}$. Alors le couple $(G_{1},G_{2})$ se compl\`ete naturellement en un triplet endoscopique non standard $(G_{1},G_{2},j_{*})$.  Dans [W1] 1.7, on a classifi\'e les triplets endoscopiques non standard. Consid\'erons les triplets \'el\'ementaires suivants:

(1) $G_{1}=G_{2}$ et $j_{*}$ est l'identit\'e;

(2)   $G_{1}$ est de type $B_{n}$ avec $n\geq2$, $G_{2}$ est de type $C_{n}$ et $j_{*}$ envoie une coracine  courte sur une coracine longue et envoie la coracine longue sur $2$ fois la coracine courte;

(3)   $G_{1}$ est de type $C_{n}$ avec $n\geq2$, $G_{1}$ est de type $B_{n}$ et $j_{*}$ envoie la coracine  courte sur la coracine longue et envoie une coracine longue sur $2$ fois une coracine courte;

(4) $G_{1}$ et $G_{2}$ sont de type $F_{4}$ et $j_{*}$ envoie une coracine courte sur  une coracine longue et une coracine longue sur $2$ fois une coracine courte;

(5) $G_{1}$ et $G_{2}$ sont de type $G_{2}$ (sic!) et $j_{*}$ envoie une coracine courte sur  une coracine longue et une coracine longue sur $3$ fois une coracine courte.

Disons qu'un triplet est quasi-\'el\'ementaire s'il se d\'eduit par restriction des scalaires d'un triplet \'el\'ementaire. Disons que deux triplets $(G_{1},G_{2},j_{*})$ et $(G'_{1},G'_{2},j'_{*})$ sont \'equivalents si, \`a isomorphismes pr\`es, on a $G_{1}=G'_{1}$, $G_{2}=G'_{2}$ et $j_{*}=cj'_{*}$ o\`u $c$ est un rationnel strictement positif. Alors tout triplet est isomorphe \`a un produit de triplets \'equivalents \`a un triplet quasi-\'el\'ementaire. 

Pour chacun des triplets \'el\'ementaires ci-dessus, on pose $N(G_{1},G_{2},j_{*})=0$ dans le cas (1), $(n+1)(2n+1)$ dans le cas (2),  $4n^2-1$ dans le cas (3), $78$ dans le cas (4), $28$ dans le cas (5).  Pour un triplet quasi-\'el\'ementaire $(G_{1},G_{2},j_{*})$, d\'eduit par restriction des scalaires disons de $F'$ \`a $F$ d'un triplet \'el\'ementaire $(G'_{1},G'_{2},j'_{*})$, on pose $N(G_{1},G_{2},j_{*})=[F':F]N(G'_{1},G'_{2},j'_{*})$. Pour un triplet g\'en\'eral, produit sur $i=1,...,n$ de triplets \'equivalents \`a des triplets quasi-\'el\'ementaires $(G_{1,i},G_{2,i},j_{*,i})$, on pose 
$$N(G_{1},G_{2},j_{*})=sup_{i=1,...,n}N(G_{1,i},G_{2,i},j_{*,i}).$$

{\bf Remarque.} Pour $n=2$, les triplets (2) et (3) sont les m\^emes.  Les deux recettes possibles pour d\'efinir  $N(G_{1},G_{2},j_{*})$ donnent le m\^eme r\'esultat.

\ass{Lemme}{ Si un triplet endoscopique non standard $(G_{1},G_{2},j_{*})$ est issu comme ci-dessus de couples $(\eta,\epsilon)$, on a $N(G_{1},G_{2},j_{*})\leq dim(G_{SC})$.}

 Preuve. Introduisons la forme quasi-d\'eploy\'ee $G^*_{AD}$ de $G_{AD}$. Fixons-en une paire de Borel \'epingl\'ee $(B^*,T^*,(E_{\alpha})_{\alpha\in \Delta})$ d\'efinie sur $F$. On a sur $G^*_{AD}$ une action galoisienne et un automorphisme $\theta^*$ qui conservent cette paire de Borel \'epingl\'ee. Soit $t\in T^*$. On introduit la composante neutre $I(t)=G^*_{AD,t\theta^*}$ du commutant de $t\theta^*$ et son groupe dual $\hat{I}(t)$. Soit $s$ un \'el\'ement semi-simple de $\hat{I}(t)$. On introduit la composante neutre   $\hat{I}(t)_{s}$  du commutant de $s$ et son groupe dual. On note $R(s,t)$ le syst\`eme de racines de ce groupe dual. Il r\'esulte de l'hypoth\`ese et de [W1] 3.3 qu'il existe $t$ et $s$ de sorte que le syst\`eme de racines de $G_{1}$ co\"{\i}ncide avec   $R(s,t)$. De l'action galoisienne sur $G_{1}$ se d\'eduit  une action galoisienne sur  $R(s,t)$. Elle est compos\'ee de l'action galoisienne sur le syst\`eme de racines de $G^*_{AD}$ et d'un cocycle \`a valeurs dans $W^{\theta^*}$. On d\'eduit du syst\`eme de racines  $R(s,t)$ un autre syst\`eme de racines, notons-le $R(t,s)$ (dans [W1] 3.3, les deux syst\`emes sont not\'es respectivement $\Sigma_{2}$ et $\Sigma_{1}$). C'est celui du groupe $G_{2}$. 
 On a d\'ecompos\'e plus haut  $(G_{1},G_{2},j_{*})$ en produit sur $i=1,...,n$ de triplets \'equivalents \`a des triplets quasi-\'el\'ementaires $(G_{1,i},G_{2,i},j_{*,i})$. Posons $N=N(G_{1},G_{2},j_{*})=sup_{i=1,...,n}N(G_{1,i},G_{2,i},j_{*,i})$. Fixons $i$ tel que $N(G_{1,i},G_{2,i},j_{*,i})=N$. Il suffit de prouver que, s'il existe $t$ et $s$ comme ci-dessus de sorte que le syst\`eme de racines de $G_{1,i}$ soit un sous-syst\`eme de $R(s,t)$ et que le syst\`eme de racines de $G_{2,i}$  soit un sous-syst\`eme de $R(t,s)$, alors $N\leq dim(G^*_{AD})$. A ce point, on peut simplifier les notations en supposant que $(G_{1},G_{2},j_{*})=(G_{1,i},G_{2,i},j_{*,i})$ et en abandonnant l'indice $i$.  
 
 Supposons qu'il existe $t$ et $s$ v\'erifiant les conditions ci-dessus. Notons  $J$ la composante neutre du centralisateur du centre de $I(t)$. Il est stable par l'action galoisienne et par $\theta^*$. Il contient $T^*$ et le groupe $I(t)$ co\"{\i}ncide avec la composante neutre du centralisateur de $t\theta^*$ dans $J$. Il en r\'esulte que les syst\`emes de racines $R(s,t)$ et $R(t,s)$ se d\'eduisent aussi bien du groupe $J$, ou mieux de la forme quasi-d\'eploy\'ee $J^*_{AD}$ de son groupe adjoint. Si l'on prouve l'assertion pour ce groupe, on en d\'eduit la m\^eme assertion pour $G^*_{AD}$ puisque $dim(J^*_{AD})\leq dim(G^*_{AD})$. Cela nous ram\`ene au cas o\`u $J=G^*_{AD}$, ce que l'on suppose d\'esormais.
On peut d\'ecomposer $G^*_{AD}$ en produit de groupes $G^*_{j}$ pour $j=1,...,m$ tels que, pour tout $j$, les composantes irr\'eductibles de $G^*_{j}$ forment une seule orbite pour le groupe de permutations engendr\'e par l'action galoisienne et par $\theta^*$.  Les syst\`emes de racines $R(s,t)$ et $R(t,s)$ se d\'ecomposent conform\'em\'ent (y compris en tenant compte de l'action galoisienne). Il existe donc un indice $j$ tel que le groupe $G^*_{AD,j}$ v\'erifie les m\^emes hypoth\`eses que $G^*_{AD}$. De nouveau, si l'on d\'emontre l'assertion pour ce groupe $G^*_{AD,j}$, on en d\'eduit l'assertion pour $G^*$. Cela nous ram\`ene au cas o\`u les composantes  irr\'eductibles de $G^*_{AD}$ forment une seule orbite pour le groupe de permutations engendr\'e par l'action galoisienne et par $\theta^*$.  Fixons une telle composante irr\'eductible $\underline{G}$. Notons $c$ le plus petit entier strictement positif tel que $(\theta^*)^c(\underline{G})=\underline{G}$. Posons $\bar{G}=\underline{G}\times \theta^*(\underline{G})\times...\times (\theta^*)^{c-1}(\underline{G})$. Soit $F''$ l'extension de $F$ tel que $\Gamma_{F''}$ soit le sous-groupe des $\sigma\in \Gamma_{F}$ qui conservent $\bar{G}$. Alors $G^*_{AD}$ est d\'eduit de $\bar{G}$ par restriction des scalaires de $F''$ \`a $F$. Les syst\`emes de racines  $R(s,t)$ et $R(t,s)$ se d\'ecomposent en produits index\'es par $Gal(F''/F)$ de syst\`emes analogues relatifs \`a $\bar{G}$. 
Puisqu'on a suppos\'e notre triplet $(G_{1},G_{2},j_{*})$ \'equivalent \`a un triplet quasi-\'el\'ementaire, on peut supposer qu'il existe un triplet \'el\'ementaire $(G_{0,1},G_{0,2},j_{0,*})$ et une extension $F'$ de $F$ tels que le syst\`eme de racines de  $G_{0,1}$, resp. $G_{0,2}$, soit un sous-syst\`eme d'un syst\`eme $R(\bar{s},\bar{t})$, resp. $R(\bar{t},\bar{s})$, relatif \`a $\bar{G}$ et que $(G_{1},G_{2},j_{*})$ soit \'equivalent au syst\`eme d\'eduit de $(G_{0,1},G_{0,2},j_{0,*})$ par restriction des scalaires. D'apr\`es la propri\'et\'e rappel\'ee ci-dessus des actions galoisiennes, un \'el\'ement de $\Gamma_{F}$ qui conserve  $(G_{0,1},G_{0,2},j_{0,*})$ conserve aussi $\bar{G}$, autrement dit $F'\subset F''$. Notons  $(\bar{G}_{1},\bar{G}_{2},\bar{j}_{*})$ le triplet sur $F''$ d\'eduit de $(G_{0,1},G_{0,2},j_{0,*})$ par restriction des scalaires de $F'$ \`a $F''$.   On voit alors que $F''$, $\bar{G}$ et  $(\bar{G}_{1},\bar{G}_{2},\bar{j}_{*})$ v\'erifient les m\^emes hypoth\`eses que $F$, $G^*_{AD}$ et $(G_{1},G_{2},j_{*})$. Par d\'efinition, $N=[F'':F]N(\bar{G}_{1},\bar{G}_{2},\bar{j}_{*})$ (ce dernier entier \'etant relatif au corps de base $F''$) et $dim(G^*_{SC})=[F'':F]dim(\bar{G})$. Il suffit donc de d\'emontrer la relation $N(\bar{G}_{1},\bar{G}_{2},\bar{j}_{*})\leq dim(\bar{G})$. En oubliant  cette construction, on est ramen\'e au cas o\`u toutes les composantes irr\'eductibles de $G^*_{AD}$ sont permut\'ees par le groupe d'automorphismes engendr\'e par $\theta^*$. On fixe comme ci-dessus une composante irr\'eductible $\underline{G}$ et on note $c$  le plus petit entier strictement positif tel que $(\theta^*)^c(\underline{G})=\underline{G}$. Mais alors les syst\`emes de racines $R(s,t)$ et $R(t,s)$ sont exactement les m\^emes que des syst\`emes $R(\underline{s},\underline{t})$ et $R(\underline{t},\underline{s})$ d\'eduits de $\underline{G}$, de son automorphisme $(\theta^*)^c$ et d'\'el\'ements convenables $\underline{t}$ et $\underline{s}$. De nouveau, cela nous ram\`ene au cas o\`u $G^*_{AD}=\underline{G}$, autrement dit, on peut supposer $G^*_{AD}$ irr\'eductible.

   Il reste \`a \'etudier cas par cas chaque syst\`eme possible pour $I(t)$ (soumis \`a la condition $J=I(t)$) et chaque triplet quasi-\'el\'ementaire possible.   On  suppose donc que $(G_{1},G_{2},j_{*})$ est \'equivalent \`a un triplet d\'eduit par restriction des scalaires   de $F'$ \`a $F$ d'un triplet  $(G_{0,1},G_{0,2},j_{0,*})$ de type (1) \`a (5).  On peut exclure le triplet (1): dans ce cas $N=0$ et l'in\'egalit\'e $N\leq dim(G^*_{AD})$ \`a prouver est \'evidente. On peut aussi exclure le cas o\`u $\theta^*$ est l'identit\'e. On est alors dans le cas d'endoscopie non tordue et il r\'esulte des constructions de [W1] 3.3 qu'alors le triplet endoscopique non standard est forc\'ement de type (1).  D'apr\`es [L] proposition II.3.2, les possibilit\'es pour $G^*_{AD}$ et $I(t)$ sont les suivantes:
   
   (6) $G^*_{AD}$ de type $A_{2m}$, $dim(G^*_{AD})=(2m+1)^2-1$,  $I(t)$ de type $B_{m}$;
   
   (7) $G^*_{AD}$ de type $A_{2m-1}$, $dim(G^*_{AD})=4m^2-1$, $I(t)$ de type $C_{m}$;
   
   (8) $G^*_{AD}$ de type $A_{2m-1}$, $dim(G^*_{AD})=4m^2-1$, $I(t)$ de type $D_{m}$;

   (9) $G^*_{AD}$ de type $D_{m}$ avec $m\geq4$, $dim(G^*_{AD})=m(2m-1)$, $I(t)$ de type $B_{m^+}\cup B_{m^-}$, avec $m^++m^-=m-1$;
   
   (10) $G^*_{AD}$ de type $D_{4}$, $dim(G^*_{AD})=28$, $I(t)$ de type $G_{2}$;
   
   (11) $G^*_{AD}$ de type $D_{4}$, $dim(G^*_{AD})=28$, $I(t)$ de type $A_{2}$;
   
   (12) $G^*_{AD}$ de type $E_{6}$, $dim(G^*_{AD})=78$, $I(t)$ de type $F_{4}$, $C_{4}$ ou $B_{3}\cup A_{1}$. 
   
   Dans le cas o\`u $(G_{0,1},G_{0,2},j_{0,*})$  est de type (4), on peut exclure les cas classiques (6) \`a (9):  des op\'erations consistant \`a prendre des commutants ou \`a passer au groupe dual \`a partir d'un groupe classique ne sauraient cr\'eer un groupe de type $F_{4}$. Les cas (10) et (11) sont exclus car $I(t)$ y est trop petit pour contenir $F_{4}$.  Il ne reste que le cas (12).  La seule possibilit\'e est que $I(t)$ soit lui-m\^eme de type $F_{4}$. Il ne contient \'evidemment qu'une copie de  ce syst\`eme, donc $F'=F$. Mais alors, par d\'efinition, $N=78=dim(G^*_{AD})$.
  
   Dans le cas o\`u $(G_{0,1},G_{0,2},j_{0,*})$  est de type (5), on peut exclure les cas classiques  pour la m\^eme raison que ci-dessus. Un syst\`eme de type $G_{2}$ ne peut pas \^etre contenu dans des syst\`emes de types $A_{2}$, $F_{4}$, $C_{4}$ ou $B_{2}$. Il ne reste que le cas (10). De nouveau $F'=F$ et $N=28=dim(G^*_{AD})$. 
   
   Supposons    $(G_{0,1},G_{0,2},j_{0,*})$  de type (2). On peut exclure les cas (8), (10) et (11): un syst\`eme de type $B_{n}$ ne peut pas intervenir dans un syst\`eme de type $D_{m}$, $G_{2}$ ou $A_{2}$. Posons $d=[F':F]$. Alors $N=d(n+1)(2n+1)$. On a $d$ syst\`emes orthogonaux de type $B_{n}$ contenus dans le syst\`eme de racines de $I(t)$. Donc le rang de $I(t)$ est au moins $dn$. Il suffit alors de prouver
   
   - dans le cas (6), l'in\'egalit\'e $m\geq dn$ entra\^{\i}ne $(2m+1)^2-1\geq d(n+1)(2n+1)$;
   
   - dans le cas (7), l'in\'egalit\'e $m\geq dn$ entra\^{\i}ne $4m^2-1\geq d(n+1)(2n+1)$;
   
   - dans le cas  (9), l'in\'egalit\'e $m-1\geq dn$ entra\^{\i}ne $m(2m-1)\geq d(n+1)(2n+1)$;

   - dans le cas (12), l'in\'egalit\'e $4\geq dn$ entra\^{\i}ne $78\geq d(n+1)(2n+1)$.
   
  On laisse au lecteur la v\'erification \'el\'ementaire.
  
   Supposons    $(G_{0,1},G_{0,2},j_{0,*})$  de type (3). On peut supposer $n\geq3$: si $n=2$, le cas (3) se confond avec le cas (2) d\'ej\`a trait\'e. On peut alors exclure les cas (10) et (11) o\`u $I(t)$ est de trop petit rang. On peut aussi exclure les cas (6), (8) ou  (9): un syst\`eme $C_{n}$ avec $n\geq3$ n'est pas un sous-syst\`eme de $B_{m}$ ou $ D_{m}$. Dans les cas restants, le m\^eme argument que ci-dessus ram\`ene \`a prouver
   
   - dans le cas (7), l'in\'egalit\'e $m\geq dn$ entra\^{\i}ne $4m^2-1\geq d(4n^2-1)$;
   
  - dans le cas (12), l'in\'egalit\'e $4\geq dn$ entra\^{\i}ne $78\geq d(4n^2-1)$. 
  
  De nouveau, on laisse la v\'erification au lecteur. Cela prouve le lemme. $\square$
   
   \bigskip

 \subsection{Compl\'ements \`a propos de l'ensemble ${\cal Z}(\tilde{G})$}
  Soit $G$ un groupe r\'eductif connexe d\'efini sur $F$ et $\tilde{G}$ un espace tordu sous $G$. Pour $\eta\in \tilde{G}$, notons $E_{\eta}$ l'ensemble des paires de Borel \'epingl\'ees ${\cal E}$ de $G$ telles que $\eta\in Z(\tilde{G},{\cal E})$, c'est-\`a-dire que $ad_{\eta}$ conserve ${\cal E}$. 
 Notons $\Theta$ l'ensemble des $\eta\in \tilde{G}$ tels que $E_{\eta}\not=\emptyset$.  Rappelons les propri\'et\'es suivantes, cf. [KS] 1.1. Soient $\eta\in \Theta$ et ${\cal E}=(B,T,(E_{\alpha})_{\alpha\in \Delta})\in E_{\eta}$. On note $W$ le groupe de Weyl de $T$ dans $G$ et $\theta=ad_{\eta}$. Alors

   (1) les sous-groupes $Z_{G_{AD}}(\eta)$ et $(T_{ad})^\theta$ de $G_{AD}$ sont connexes ($Z_{G_{AD}}(\eta)$ est le groupe des points fixes de $\theta$ dans $G_{AD}$);

   (2)    $W^{\theta}$ s'identifie au groupe de Weyl de $T^{\theta,0}$ dans $G_{\eta}$. 
   
   Montrons que
   
   (3) pour $\eta\in \Theta$, l'ensemble $E_{\eta}$ est une classe de conjugaison sous $G_{\eta}$.
   
Preuve. Il est clair que la conjugaison par un \'el\'ement de $G_{\eta}$ conserve $E_{\eta}$. Inversement,  soient ${\cal E},{\cal E}'\in E_{\eta}$.  D'apr\`es [I] 1.3(2) et  (2) ci-dessus, il existe $x\in G_{\eta}$ tel que ${\cal E}$ et $ad_{x}({\cal E}')$ aient la m\^eme paire de Borel sous-jacente. Notons $(B,T)$ cette paire. Il existe alors $t\in T$ de sorte que ${\cal E}=ad_{tx}({\cal E}')$. Puisque $ad_{\eta}$ conserve ${\cal E}'$ et que $x\in G_{\eta}$, $ad_{\eta}$ conserve $ad_{x}({\cal E}')$. Puisque $ad_{\eta}$ conserve aussi ${\cal E}=ad_{tx}({\cal E}')$, on a $ad_{\eta}(t)\in t Z(G)$. Donc l'image $t_{ad}$ de $t$ dans $T_{ad}$ est fixe par $\theta$. D'apr\`es (1),  $T^{\theta,0}$ s'envoie surjectivement dans $T_{ad}^{\theta}$. On a donc $t\in Z(G)T^{\theta,0}$. Quitte \`a multiplier $t$ par un \'el\'ement de $Z(G)$, on peut supposer $t\in T^{\theta,0}\subset G_{\eta}$. En posant $y=tx$, on a $y\in G_{\eta}$ et ${\cal E}=ad_{y}({\cal E}')$. $\square$

Rappelons que, pour toute paire de Borel \'epingl\'ee ${\cal E}$, on a des applications
 $$(4) \qquad Z(\tilde{G},{\cal E})\to {\cal Z}(\tilde{G},{\cal E})=Z(\tilde{G},{\cal E})/(1-\theta)(Z(G))\simeq {\cal Z}(\tilde{G}).$$
 Pour une autre paire de Borel \'epingl\'ee ${\cal E}'$, on choisit $x\in G$ tel que $ad_{x}({\cal E}')={\cal E}$. Le diagramme suivant est commutatif
 $$(5)\qquad \begin{array}{ccc}Z(\tilde{G},{\cal E}')&&\\ &\searrow&\\ ad_{x}\downarrow\,&&{\cal Z}(\tilde{G})\\&\nearrow&\\ Z(\tilde{G},{\cal E})&&\\ \end{array}$$
 Pour $\eta\in \Theta$, choisissons ${\cal E}\in E_{\eta}$. Via (4), $\eta$ s'envoie alors sur un \'el\'ement de ${\cal Z}(\tilde{G})$. La propri\'et\'e (3) et la commutativit\'e du diagramme ci-dessus montrent que cet \'el\'ement ne d\'epend pas du choix de ${\cal E}$. On obtient ainsi une application $\Theta\to {\cal Z}(\tilde{G})$. Il est imm\'ediat qu'elle est \'equivariante pour les actions galoisiennes.
 L'ensemble $\Theta$ est invariant par conjugaison par $G$. Notons $\Theta/conj$ l'ensemble des classes de conjugaison.
 
 \ass{Lemme}{L'application pr\'ec\'edente se quotiente en une bijection de $\Theta/conj$ sur ${\cal Z}(\tilde{G})$.}

 Preuve.    Soient $\eta,\eta'\in \Theta$, supposons ces \'el\'ements conjugu\'es. Fixons $x\in G$ tel que $ad_{x}(\eta')=\eta$. Fixons aussi ${\cal E}\in E_{\eta}$ et ${\cal E}'\in E_{\eta'}$. La paire de Borel \'epingl\'ee $ad_{x}({\cal E}')$ est conserv\'ee par $ad_{\eta}$. D'apr\`es (3), quitte \`a multiplier $x$ \`a gauche par un \'el\'ement de $G_{\eta}$, on peut supposer $ad_{x}({\cal E}')={\cal E}$. Le diagramme (5) montre alors que $\eta$ et $\eta'$ ont m\^eme image dans ${\cal Z}(\tilde{G})$. Inversement, soient $\eta,\eta'\in \Theta$, supposons que ces \'el\'ements ont m\^eme image dans ${\cal Z}(\tilde{G})$. On fixe ${\cal E}\in E_{\eta}$ et ${\cal E}'\in E_{\eta'}$. On fixe $x\in G$ tel que $ad_{x}({\cal E}')={\cal E}$. Alors $ad_{x}(\eta')=z\eta$, avec $z\in Z(G)$. D'apr\`es (5), l'image de $\eta'$ dans ${\cal Z}(\tilde{G})$ est aussi celle de $z\eta$. Pour qu'elle soit \'egale \`a celle de $\eta$, il faut que $z$ appartienne \`a $(1-\theta)(Z(G))$. Ecrivons $z=(\theta-1)(z')$, avec $z'\in Z(G)$. Alors $ad_{z'x}(\eta')=\eta$ et $\eta$ et $\eta'$ sont  conjugu\'es. Cela prouve que l'application $\Theta\to {\cal Z}(\tilde{G})$ se quotiente en une injection $\Theta/conj\to {\cal Z}(\tilde{G})$. Fixons une 
  paire de Borel \'epingl\'ee ${\cal E}$. Alors $Z(\tilde{G},{\cal E})$ est contenu dans $\Theta$ et notre application co\"{\i}ncide sur cet ensemble avec l'application (4). Celle-ci est surjective par d\'efinition de ${\cal Z}(\tilde{G})$. Cela prouve la surjectivit\'e de l'application de l'\'enonc\'e. $\square$
  
  \bigskip
 
\subsection{D\'efinition de triplets $(G,\tilde{G},{\bf a})$ particuliers}

Consid\'erons un syst\`eme de racines et un automorphisme $\theta^*$ de ce syst\`eme. On introduit un groupe r\'eductif connexe $G$ sur $F$, simplement connexe et d\'eploy\'e, dont le syst\`eme de racines est celui fix\'e. Fixons une paire de Borel \'epingl\'ee ${\cal E}^*$ de $G$ d\'efinie sur $F$. A $\theta^*$ est alors associ\'e un automorphisme de $G$ qui conserve ${\cal E}^*$, et qui est d\'efini sur $F$. On le note encore $\theta^*$. On introduit   l'espace tordu $\tilde{G}=G\theta^*$. C'est un espace principal homog\`ene sous $G$ \`a gauche, muni du point marqu\'e $\theta^*$. L'action de $G$ \`a droite est d\'efinie par $(g\theta^*,x)\mapsto g\theta^*(x)\theta^*$ et l'action galoisienne fixe $\theta^*$. On voit que la classe d'isomorphisme du couple $(G,\tilde{G})$ est uniquement d\'etermin\'ee par le syst\`eme de racines et son automorphisme (plus pr\'ecis\'ement par la classe de conjugaison de ce dernier dans le groupe d'automorphismes du syst\`eme de racines). Plus g\'en\'eralement, consid\'erons 
 une extension finie $F'$ de $F$, un syst\`eme de racines et un automorphisme $\theta^*$ de ce syst\`eme. On introduit un couple $(G_{F'},\tilde{G}_{F'})$ d\'efini sur $F'$ associ\'e au syst\`eme de racines et \`a son automorphisme. On note $(G,\tilde{G})$ le couple sur $F$ d\'eduit de $(G_{F'},\tilde{G}_{F'})$ par restriction des scalaires. 
 
 Fixons un tel couple $(G,\tilde{G})$.  Parce que $G$ est simplement connexe, on a la propri\'et\'e suppl\'ementaire suivante.  Soient $\eta\in \Theta$ et ${\cal E}=(B,T,(E_{\alpha})_{\alpha\in \Delta})\in E_{\eta}$. Posons $\theta=ad_{\eta}$. Alors
   
   (1) $Z_{G}(\eta)$ et $T^{\theta}$ sont connexes.

 Pour $\eta\in \tilde{G}(F)$, notons $E_{\eta,F}$ l'ensemble des paires de Borel \'epingl\'ees ${\cal E}$ de $G$ d\'efinies sur $F$ telles que $\eta\in Z(\tilde{G},{\cal E})$. Notons $\Theta_{F}$ l'ensemble des $\eta\in \tilde{G}(F)$ tels que $E_{\eta,F}\not=\emptyset$.  Par construction, l'ensemble $\Theta_{F}$ n'est pas vide. On a $\Theta_{F}\subset \Theta^{\Gamma_{F}}$ et l'application de 6.2 se restreint en une application $\Theta_{F}\to {\cal Z}(\tilde{G})^{\Gamma_{F}}$.    L'ensemble $\Theta_{F}$ n'a pas de raison d'\^etre invariant par conjugaison stable. On peut n\'eanmoins introduire la relation d'\'equivalence dans $\Theta_{F}$: deux \'el\'ements sont \'equivalents si et seulement s'ils sont stablement conjugu\'es. On note $\Theta_{F}/st-conj$ l'ensemble des classes d'\'equivalence. Remarquons que, d'apr\`es (1), la classe de conjugaison stable d'un \'el\'ement $\eta\in \Theta_{F}$ est l'intersection de $\tilde{G}(F)$ avec la classe de conjugaison par $G(\bar{F})$ de $\eta$. On en d\'eduit une injection $\Theta_{F}/st-conj\to \Theta/conj$. Montrons que
 
 (2) l'application du lemme 6.2 se restreint en une bijection $\Theta_{F}/st-conj\simeq {\cal Z}(\tilde{G})^{\Gamma_{F}}$. 
 
 Preuve. Le lemme 6.2 nous dit que cette restriction est injective.   Soit $e\in {\cal Z}(\tilde{G})^{\Gamma_{F}}$. Puisque $\Theta_{F}$ n'est pas vide, fixons $\eta\in \Theta_{F}$ et ${\cal E}=(B,T,(E_{\alpha})_{\alpha\in \Delta})\in E_{\eta,F}$. On peut fixer un \'el\'ement $z\in Z(G)$ tel que $e$ soit l'image de $z\eta$. Puisque $e$ est fixe par $\Gamma_{F}$, on a $z\sigma(z)^{-1}\in (1-\theta)(Z(G))$ pour tout $\sigma\in \Gamma_{F}$. L'application $\sigma\mapsto z\sigma(z)^{-1}$ est un cocycle de $\Gamma_{F}$ dans $(1-\theta)(Z(G))$, que l'on pousse en un cocycle \`a valeurs dans $(1-\theta)(T)$. Il r\'esulte des constructions que $T$ et $(1-\theta)(T)$ sont d\'eduits par restriction des scalaires de tores d\'eploy\'es sur $F'$. Donc $(1-\theta)(T)$ est induit et $H^1(\Gamma_{F};(1-\theta)(T))=\{1\}$. On peut donc fixer $u\in T$ de sorte que $z\sigma(z)^{-1}=(1-\theta)(u\sigma(u)^{-1})$ pour tout $\sigma\in \Gamma_{F}$. En notant $u_{ad}$ l'image de $u$ dans $T_{ad}$, cette relation implique que $u_{ad}\sigma(u_{ad})^{-1}\in T_{ad}^{\theta}$. Donc $\sigma\mapsto u_{ad}\sigma(u_{ad})^{-1}$ est un cocycle \`a valeurs dans $T_{ad}^{\theta}$. De nouveau, ce tore est induit et $H^1(\Gamma_{F};T_{ad}^{\theta})=\{1\}$. On peut donc fixer $v_{ad}\in T_{ad}^{\theta}$ de sorte que $u_{ad}\sigma(u_{ad})^{-1}=v_{ad}\sigma(v_{ad})^{-1}$ pour tout $\sigma$. Relevons $v_{ad}$ en un \'el\'ement  $v\in T^{\theta}$. Posons $t=uv^{-1}$. On a alors $x_{ad}\in T_{ad}^{\Gamma_{F}}$. Puisque $v\in T^{\theta}
$, on a encore l'\'egalit\'e $z\sigma(z)^{-1}=(1-\theta)(x\sigma(x)^{-1})$ pour tout $\sigma$. Posons $\eta'=ad_{x^{-1}}(z\eta )$ et ${\cal E}'=ad_{x^{-1}}({\cal E})$. La relation pr\'ec\'edente implique que $\eta'\in \tilde{G}(F)$. Parce que $x_{ad}$ est fixe par $\Gamma_{F}$, ${\cal E}'$ est d\'efinie sur $F$. Alors   ${\cal E}'$ appartient \`a $E_{\eta',F}$ et $\eta'$ appartient \`a $\Theta_{F}$.   Le diagramme 6.2 (5) montre que l'image de $\eta'$ dans ${\cal Z}(\tilde{G})$ est la m\^eme que celle de $z\eta$, laquelle est $e$. Cela ach\`eve la preuve. $\square$

Puisque $G$ est simplement connexe, $\hat{G}$ est adjoint et $H^1(W_{F};Z(\hat{G}))=\{1\}$. Compl\'etons le couple $(G,\tilde{G})$ par l'unique cocycle possible ${\bf a}=1$. Remarquons que, pour toute donn\'ee endoscopique ${\bf G}'=(G',{\cal G}',\tilde{s})$ de $(G,\tilde{G},{\bf a})$, le choix d'un \'el\'ement $\eta\in \Theta_{F}$ permet d'identifier $\tilde{G}'$ \`a $G'$ de la fa\c{c}on suivante. On note $e$ l'image de $\eta$ dans ${\cal Z}(\tilde{G})$. Cet \'el\'ement est fixe par $\Gamma_{F}$. En se rappelant que $\tilde{G}'=G'\times_{{\cal Z}(G)}{\cal Z}(\tilde{G})$, l'application qui, \`a $x\in G'$, associe l'image de $ (x,e)$ dans $\tilde{G}'$ identifie $G'$ \`a $\tilde{G}'$.

On fixe comme toujours une paire de Borel \'epingl\'ee $\hat{{\cal E}}=(\hat{B},\hat{T},(\hat{E}_{\hat{\alpha}})_{\hat{\alpha}\in \Delta})$ de $\hat{G}$, conserv\'ee par l'action galoisienne et on note $\hat{\theta}$ l'automorphisme habituel qui la conserve. Introduisons la donn\'ee endoscopique ${\bf G}'=(G',{\cal G}',\tilde{s})$  "maximale" de $(G,\tilde{G},{\bf a})$ d\'efinie par $\tilde{s}=\hat{\theta}$ et ${\cal G}'=\hat{G}_{\hat{\theta}}\rtimes W_{F}$. {\bf On suppose d\'esormais que le syst\`eme de racines de d\'epart ne contient pas de composante de type $A_{2n}$}. On a

(3)  les applications naturelles ${\cal Z}(G)=Z(G)/(1-\theta)(Z(G))\to Z(G')$ et ${\cal Z}(\tilde{G})\to {\cal Z}(\tilde{G}')$ sont bijectives.

Preuve. Puisque les ensembles ${\cal Z}(\tilde{G})$ et ${\cal Z}(\tilde{G}')$ sont des espaces principaux homog\`enes sous respectivement ${\cal Z}(G)$ et $Z(G')$, la bijectivit\'e de la deuxi\`eme application r\'esulte de celle de la premi\`ere. Fixons des paires de Borel $(B,T)$ de $G$ et $(B',T')$ de $G'$. Notons $\Sigma(T)$ et $\Sigma(T')$ les ensembles de racines de $T$ dans $\mathfrak{g}$ et de $T'$ dans $\mathfrak{g}'$. On a un homomorphisme $\xi:T\to T'$ qui se quotiente en un isomorphisme $T/(1-\theta)(T)\simeq T'$. La description de [I] 1.6 se simplifie puisque, d'apr\`es l'hypoth\`ese sur le syst\`eme de racines, tous les \'el\'ements de $\Sigma(T)$ sont de type $1$. On obtient que $\Sigma(T')$ est l'ensemble des $N\alpha$ pour $\alpha\in \Sigma(T)$. Soit $t\in T$ tel que $\xi(t)\in Z(G')$. Alors $N\alpha(t)=1$ pour tout $\alpha\in \Sigma(T)$. Puisque $\theta$ permute la base de $X_{*}(T_{ad})$ form\'ee des copoids associ\'es aux racines simples, on voit que cette condition implique que l'image $t_{ad}$ de $t$ dans $T_{ad}$ appartient \`a $(1-\theta)(T_{ad})$. Il en r\'esulte que
 $t\in Z(G)(1-\theta)(T)$. Alors $\xi(t)$ est aussi l'image par $\xi$ d'un \'el\'ement de $Z(G)$. Cela prouve la surjectivit\'e de la premi\`ere application de (3). Pour prouver son injectivit\'e, il suffit de prouver celle de l'application $Z(G)/(1-\theta)(Z(G))\to T/(1-\theta)(T)$, ou encore de prouver que $Z(G)\cap (1-\theta)(T)=(1-\theta)(Z(G))$. Or, soit $t\in T$ tel que $(1-\theta)(t)\in Z(G)$. Alors $t_{ad}\in T_{ad}^{\theta}$. Puisque ce tore est connexe, l'application $T^{\theta}\to T_{ad}^{\theta}$ est surjective et on peut \'ecrire $t=zt'$, avec $z\in Z(G)$ et $t'\in T^{\theta}$. Alors $(1-\theta)(t)=(1-\theta)(z)\in (1-\theta)(Z(G))$. $\square$
     
Soit $\eta\in \Theta_{F}$. On l'envoie sur un \'el\'ement de ${\cal Z}(\tilde{G})$, puis sur un \'el\'ement $\epsilon\in {\cal Z}(\tilde{G}')$. On a $\epsilon\in \tilde{G}'(F)$. Fixons un \'el\'ement ${\cal E}\in E_{\eta,F}$ dont on note $(B,T)$ la paire de Borel sous-jacente et fixons une paire de Borel $(B',T')$ de $G'$ d\'efinie sur $F$. Le sextuplet $(\epsilon,B',T',B,T,\eta)$ est un diagramme. A l'aide de la description des ensembles de racines de [W1] 3.3, on v\'erifie ais\'ement que $G_{\eta}$ et $G'=G'_{\epsilon}$ sont simplement connexes. 
La correspondance $\eta\mapsto \epsilon$ se quotiente en une bijection entre l'ensemble de classes de conjugaison stable $\Theta_{F}/st-conj$ et l'ensemble des classes de conjugaison stable d'\'el\'ements de ${\cal Z}(\tilde{G}')^{\Gamma_{F}}$ (ces derni\`eres classes \'etant r\'eduites \`a un \'el\'ement). Il est assez clair qu'inversement, une classe de conjugaison stable dans $\tilde{G}'(F)$ qui correspond \`a une classe dans $\Theta_{F}/st-conj$ est la classe d\'efinie par cette bijection. 

Consid\'erons un  triplet endoscopique non standard quasi-\'el\'ementaire $(G_{1},G_{2},j_{*})$  d\'eduit par restriction des scalaires de $F'$ \`a $F$ d'un triplet de type (2), (3),  (4) ou (5). On lui associe un syst\`eme de racines et un automorphisme $\theta^*$ de ce syst\`eme:  

-  dans le cas (2), le syst\`eme est $D_{n+1}$ si $n\geq3$ et $A_{3}$ si $n=2$; $\theta^*$ est d'ordre $2$;

- dans le cas (3), le syst\`eme est $A_{2n-1}$ et $\theta^*$ est d'ordre $2$;

- dans le cas (4), le syst\`eme est $E_{6}$ et $\theta^*$ est d'ordre $2$;

- dans le cas (5), le syst\`eme est $D_{4}$ et $\theta^*$ est d'ordre $3$.

Par la construction ci-dessus, ce syst\`eme, cet automorphisme et l'extension $F'$ d\'eterminent un couple $(G,\tilde{G})$, que l'on compl\`ete par le cocycle ${\bf a}=1$. On introduit comme ci-dessus la donn\'ee endoscopique  maximale ${\bf G}'=(\hat{G'},\hat{G}_{\hat{\theta}}\rtimes W_{F},\hat{\theta})$. Soit $\eta\in \Theta_{F}$ et $\epsilon$ son image dans $\tilde{G}'(F)$. On v\'erifie que le couple $(\epsilon,\eta)$ donne naissance comme en 5.1 au triplet endoscopique non standard $(G_{1},G_{2},j_{*})$. Les groupes $G_{1}$ et $G_{2}$ sont respectivement \'egaux \`a $G_{\eta}$ et $G'_{\epsilon}=G'$, ces deux groupes \'etant simplement connexes.

  \ass{Lemme}{Soient ${\bf G}'=(G',{\cal G}',\tilde{s})$ une donn\'ee endoscopique de $(G,\tilde{G},{\bf a})$, $\eta$ un \'el\'ement semi-simple de $\tilde{G}(F)$ et $\epsilon$ un \'el\'ement semi-simple de $G'(F)$ qui se correspondent. Notons $(G'_{1},G'_{2},j'_{*})$ le triplet endoscopique non standard auquel ils donnent naissance. Supposons que ${\bf G}'$ ne soit pas \'equivalente \`a la donn\'ee maximale ou que $\eta$ ne soit pas stablement conjugu\'e \`a un \'el\'ement de $\Theta_{F}$. Alors on a l'in\'egalit\'e $N(G'_{1},G'_{2},j'_{*})<dim(G_{SC})$.}
  
  Preuve. En vertu du lemme 6.1, il s'agit d'exclure l'\'egalit\'e $N(G'_{1},G'_{2},j'_{*})=dim(G_{SC})$. On reprend la d\'emonstration de ce lemme en \'etudiant les cas o\`u les in\'egalit\'es de rang peuvent devenir des \'egalit\'es. On s'aper\c{c}oit que, si on a \'egalit\'e,    l'\'el\'ement $t$ de cette d\'emonstration est \'egal \`a $1$ et le triplet $(G'_{1},G'_{2},j'_{*})$ est notre triplet $(G_{1},G_{2},j_{*})$ de d\'epart. L'\'egalit\'e $t=1$ signifie que $\eta$ appartient \`a $\Theta$. Puisque $\eta$ appartient \`a $\tilde{G}(F)$, son image $e$ dans ${\cal Z}(\tilde{G})$ est fixe par $\Gamma_{F}$.   D'apr\`es (2), il existe $\eta_{0}\in \Theta_{F}$ qui a $e$ pour image. D'apr\`es le lemme 6.2, $\eta$ est stablement conjugu\'e \`a $\eta_{0}$. 
   On conna\^{\i}t le syst\`eme de racines de $G'_{1}$: c'est l'ensemble $\Sigma_{2}$ de [W1] 3.3. Comme on l'a dit, la description de cette r\'ef\'erence se simplifie car notre groupe $G$ n'a que des racines de type 1. En supposant $\tilde{s}=s\hat{\theta}$, avec $s\in \hat{T}$, on voit que l'\'egalit\'e  $G'_{1}=G_{1}=G_{ \eta_{0}}$ entra\^{\i}ne que, pour toute racine $\hat{\alpha}$ de $\hat{T}$, on a $N\hat{\alpha}(s)=1$. Par le m\^eme argument que dans la preuve de (3), et parce que le groupe $\hat{G}$ est adjoint, cette condition implique $s\in  (1-\hat{\theta})(\hat{T})$. Quitte \`a remplacer ${\bf G}'$ par une donn\'ee \'equivalente, on peut supposer $s=1$ et $\tilde{s}=\hat{\theta}$. Puisque $Z(\hat{G})=\{1\}$, on a une relation $s\hat{\theta}(g)w(s)^{-1}=g$ pour tout $(g,w)\in {\cal G}'$. Puisque $s=1$, cette relation se simplifie en $g\in \hat{G}^{\hat{\theta}}$. Mais ce groupe est connexe ([KS] 1.1) donc $g\in \hat{G}_{\hat{\theta}}$ et ${\cal G}'=\hat{G}'_{\hat{\theta}}\rtimes W_{F}$.  Alors la donn\'ee ${\bf G}'$ est \'equivalente \`a la donn\'ee "maximale". Mais alors, on est dans la situation que l'\'enonc\'e exclut. $\square$

\bigskip

\subsection{Mise en place des r\'ecurrences}

On aura a prouver des assertions concernant soit un triplet $(G,\tilde{G},{\bf a})$ comme en [II] 1.1, soit un triplet endoscopique non standard $(G_{1},G_{2},j_{*})$.  Concernant les triplets $(G,\tilde{G},{\bf a})$, on conserve les hypoth\`eses de r\'ecurrence pos\'ees en [II] 1.1. Mais il nous faut intercaler les hypoth\`eses concernant ces triplets et celles concernant les triplets endoscopiques non standard. 

Pour d\'emontrer une assertion concernant un triplet $(G,\tilde{G},{\bf a})$   quasi-d\'eploy\'e et \`a torsion int\'erieure, on ne pose aucune hypoth\`ese concernant les triplets endoscopiques non standard. Les seuls tels triplets intervenant dans ce cas sont triviaux (du cas (1) de 6.1) et leurs propri\'et\'es sont tautologiques.

Dans les autres cas, on  raisonne par r\'ecurrence sur un entier $N\geq0$.

Pour d\'emontrer une assertion concernant l'un  des triplets particuliers $(G,\tilde{G},{\bf a})$ d\'efinis en 6.3  tel que $dim(G_{SC})=N$, on suppose   connues toutes les assertions concernant des triplets endoscopiques non standard $(G_{1},G_{2},j_{*})$ tels que $N(G_{1},G_{2},j_{*})<N$ (en plus, naturellement, des hypoth\`eses pos\'ees en [II] 1.1). Pour d\'emontrer une assertion concernant un triplet $(G,\tilde{G},{\bf a})$ tel que $dim(G_{SC})=N$ et qui n'est pas l'un des triplets particuliers d\'efinis en 6.3, on suppose connues toutes les assertions concernant des triplets endoscopiques non standard $(G_{1},G_{2},j_{*})$ tels que $N(G_{1},G_{2},j_{*})\leq N$. 

Pour d\'emontrer une assertion concernant un triplet endoscopique non standard $(G_{1},G_{2},j_{*})$ tel que $N(G_{1},G_{2},j_{*})=N$, on suppose connues toutes les assertions concernant des triplets endoscopiques non standard $(G'_{1},G'_{2},j'_{*})$ tels que $N(G'_{1},G'_{2},j'_{*})<N$. On suppose connues toutes les assertions concernant des triplets  $(G',\tilde{G}',{\bf a}')$ quasi-d\'eploy\'es et \`a torsion int\'erieure tels que $dim(G'_{SC})\leq N$. On suppose connues toutes les assertions concernant les triplets $(G',\tilde{G}',{\bf a}')$ d\'efinis en 6.3  tels que $dim(G'_{SC})=N$. On suppose connues toutes les assertions concernant des triplets $(G',\tilde{G}',{\bf a}')$ quelconques tels que $dim(G'_{SC})<N$.

  En raisonnant ainsi, on a les deux propri\'et\'es suivantes:

- quand on travaille avec un triplet $(G,\tilde{G},{\bf a})$,  il y a au plus un nombre fini de classes  de conjugaison stable d'\'el\'ements $\eta\in \tilde{G}(F)$ qui peuvent cr\'eer  des triplets endoscopiques non-standard dont les propri\'et\'es ne sont pas connues;

- quand on travaille avec un triplet $(G_{1},G_{2},j_{*})$, on peut le d\'ecomposer en produit de triplets \'equivalents \`a des  triplets quasi- \'el\'ementaires qui, ou bien sont de type (1), auquel cas on d\'emontrera directement les propri\'et\'es en vue, ou bien sont issus d'un couple $(\epsilon, \eta)$ provenant d'un triplet $(G,\tilde{G},{\bf a})$ dont les propri\'et\'es sont d\'ej\`a connues.

\bigskip 
 
\bigskip
\subsection{ Quelques d\'efinitions}
On consid\`ere dans ce paragraphe un triplet endoscopique non standard $(G_{1},G_{2},j_{*})$.  
On utilise les notations de 6.1. De $j_{*}$ se d\'eduit une  correspondance bijective entre classes de conjugaison stables dans $\mathfrak{g}_{1}(F)$ et $\mathfrak{g}_{2}(F)$. Il s'en d\'eduit un isomorphisme $SI(\mathfrak{g}_{1}(F))\otimes Mes(G_{1}(F))\simeq SI(\mathfrak{g}_{2}(F))\otimes Mes(G_{2}(F))$, d'o\`u, par dualit\'e, un isomorphisme
$$(1)\qquad D_{g\acute{e}om}^{st}(\mathfrak{g}_{1}(F))\otimes Mes(G_{1}(F))^*\simeq D_{g\acute{e}om}^{st}(\mathfrak{g}_{2}(F))\otimes Mes(G_{2}(F))^*.$$
   L'application $j_{*}$ induit une bijection entre Levi standard de $G_{1}$ et Levi standard de $G_{2}$. Soient $M_{1}$ et $M_{2}$ deux tels Levi qui se correspondent.  On a comme ci-dessus un isomorphisme
$$(2)\qquad D_{g\acute{e}om}^{st}(\mathfrak{m}_{1}(F))\otimes Mes(M_{1}(F))^*\simeq D_{g\acute{e}om}^{st}(\mathfrak{m}_{2}(F))\otimes Mes(M_{2}(F))^*.$$
Ces isomorphismes se restreignent \'evidemment aux espaces de distributions \`a support nilpotent, not\'es $D_{nil}^{st}(\mathfrak{g}_{1}(F))$ etc...

Pour d\'efinir des int\'egrales pond\'er\'ees, on doit fixer des mesures sur ${\cal A}_{M_{i}}$ pour $i=1,2$. De l'application $j_{*}$ se d\'eduit un isomorphisme ${\cal A}_{M_{1}}\simeq {\cal A}_{M_{2}}$ et on suppose  que les mesures se correspondent par cet isomorphisme.

On doit encore d\'efinir une certaine constante. Fixons des paires de Borel invariantes par $\Gamma_{F}$ des groupes duaux $\hat{G}_{i}$, pour $i=1,2$, dont on note les tores $\hat{T}_{i}$. Soit $n>0$ un entier tel que $nb$ prenne ses valeurs dans ${\mathbb N}_{>0}$. Alors $nj^*$ envoie le r\'eseau engendr\'e par $\Sigma(T_{2})$ dans celui engendr\'e par $\Sigma(T_{1})$. De $nj^{*}$ se d\'eduit dualement un homomorphisme de $\hat{T}_{2}$ dans $\hat{T}_{1}$. On v\'erifie qu'il envoie $Z(\hat{M}_{2})$ dans $Z(\hat{M}_{1})$. Il est \'equivariant et on obtient un homomorphisme
$$\hat{j}_{n}:Z(\hat{M}_{2})^{\Gamma_{F}}\to Z(\hat{M}_{1})^{\Gamma_{F}}.$$
Puisque $G_{1}$ et $G_{2}$ sont simplement connexes, leurs groupes duaux sont adjoints et les groupes ci-dessus sont connexes. Cela entra\^{\i}ne que l'homomorphisme est surjectif. Son noyau est fini. On pose
$$c_{M_{1},M_{2}}^{G_{1},G_{2}}=n^{-a_{M_{2}}}\vert ker(\hat{j}_{n})\vert ,$$
o\`u, comme toujours, $a_{M_{2}}$ est la dimension de $A_{M_{2}}$. 
Cela ne d\'epend pas du choix de $n$. En effet, si l'on remplace $n$ par $nm$, pour un entier $m\geq1$, on a l'\'egalit\'e $\hat{j}_{nm}=(\hat{j}_{n})^m$, donc le nombre d'\'el\'ements du noyau est multipli\'e par $m$ \'elev\'e \`a la puissance $dim(Z(\hat{M}_{2})^{\Gamma_{F}})$. Or cette dimension est \'egale \`a $a_{M_{2}}$ et le terme d\'efini ci-dessus ne change pas.

Pour $i=1,2$ soit $B_{i}$ une fonction sur $\Sigma(T_{i})$ \`a valeurs dans l'ensemble ${\mathbb Q}_{>0}$ des rationnels strictement positifs. On suppose ces fonctions reli\'ees par la condition suivante:

- pour tout $\alpha_{2}\in \Sigma(T_{2})$, $B_{1}(\tau(\alpha_{2}))=\frac{B_{2}(\alpha_{2})}{b(\alpha_{2})}$.

Remarquons que cette condition est sym\'etrique en le sens suivant. Le triplet $(G_{2},G_{1},j_{*}^{-1})$ est encore endoscopique non standard. L'analogue de la bijection $\tau$ pour ce triplet est $\tau^{-1}$. L'analogue de la  fonction $b$  est la fonction $b'$ d\'efinie par $b'(\alpha_{1})=b(\tau^{-1}(\alpha_{1}))^{-1}$. Alors le couple $(B_{2},B_{1})$ v\'erifie encore l'hypoth\`ese ci-dessus pour ce triplet. 

Montrons que

(3) la fonction $B_{1}$ v\'erifie les hypoth\`eses de [II] 1.8 si et seulement si $B_{2}$ les v\'erifie. 

Preuve. D'apr\`es la sym\'etrie remarqu\'ee ci-dessus, on peut supposer que $B_{2}$ v\'erifie ces hypoth\`eses et on doit montrer que $B_{1}$ les v\'erifie aussi. 
Les conditions d'\'equivariance r\'esultent de celles v\'erifi\'ees par $B_{2}$ et par la correspondance entre racines. Il faut v\'erifier que, sur un sous-syst\`eme irr\'eductible de $\Sigma(T_{1})$ sur lequel on fixe une norme euclidienne poss\'edant les propri\'et\'es usuelles, $B_{1}$ est soit constante, soit proportionnelle \`a  la fonction $\alpha_{1}\mapsto (\alpha_{1},\alpha_{1})$. Fixons un tel sous-syst\`eme. L'ensemble des $\alpha_{2}\in \Sigma(T_{2})$ tels que $\tau(\alpha_{2})$ appartient \`a ce sous-syst\`eme forme un sous-syst\`eme irr\'eductible de $\Sigma(T_{2})$.  On est ramen\'e \`a ces deux sous-syst\`emes irr\'eductibles. En oubliant les actions galoisiennes qui ne jouent plus de r\^ole ici,  on peut aussi bien supposer $\Sigma(T_{2})$ et $\Sigma(T_{1})$ irr\'eductibles. On a rappel\'e en 6.1 tous les cas possibles, \`a homoth\'etie pr\`es (et il est clair que la question est insensible \`a une homoth\'etie). Si $j_{*}$ est l'identit\'e, l'assertion est claire. Dans les quatre autres cas, on constate qu'en munissant nos syst\`emes de produits euclidiens comme ci-dessus, il existe des constantes $c_{1},c_{2}>0$ telles que, pour tout $\alpha_{2}\in \Sigma(T_{2})$, on a les \'egalit\'es $b(\alpha_{2})=c_{2}(\alpha_{2},\alpha_{2})=c_{1}(\alpha_{1},\alpha_{1})^{-1}$, o\`u $\alpha_{1}=\tau(\alpha_{2})$. Si $B_{2}$ est constante, alors $B_{1}$ est proportionnelle \`a $\alpha_{1}\mapsto (\alpha_{1},\alpha_{1})$. Si $B_{2}$ est proportionnelle \`a $\alpha_{2}\mapsto (\alpha_{2},\alpha_{2})$, alors $B_{1}$ est constante. $\square$

  Fixons un tel couple de fonctions v\'erifiant les hypoth\`eses de [II] 1.8. Fixons deux Levi $M_{1}$ et $M_{2}$ qui se correspondent.  L'isomorphisme $j_{*}$ d\'efinit un isomorphisme encore not\'e $j_{*}:{\cal A}_{M_{1}}\to {\cal A}_{M_{2}}$. Soit $\alpha_{2}\in \Sigma(T_{2})$, posons $\alpha_{1}=\tau(\alpha_{2})$ et, pour $i=1,2$, notons $\alpha'_{i}$ la restriction de $B_{i}(\alpha_{i})^{-1}\alpha_{i}$ \`a ${\cal A}_{M_{i}}$. Il r\'esulte des d\'efinitions que $\alpha'_{2}\circ j_{*}=\alpha'_{1}$. Donc, par dualit\'e, $j_{*}$ d\'etermine une bijection de $\Sigma(A_{M_{2}},B_{2})$ sur $\Sigma(A_{M_{1}},B_{1})$. Cette bijection est compatible aux \'equivalences sur chacun de ces ensembles. On a donc aussi une bijection de ${\cal J}_{M_{2}}^{G_{2}}(B_{2})$ sur ${\cal J}_{M_{1}}^{G_{1}}(B_{1})$.   Pour $i=1,2$, soit $J_{i}\in {\cal J}_{M_{i}}^{G_{i}}(B_{i})$. On suppose que $J_{1}$ et $J_{2}$ se correspondent par cette bijection.     Soit $u\in U_{J_{2}}$. Via l'exponentielle, on consid\`ere que $u$ est un germe de fonctions d\'efini au voisinage de $0$ dans $\mathfrak{a}_{M_{2}}(F)$. On v\'erifie que la fonction $X_{1}\mapsto u(j_{*}(X_{1}))$ sur $\mathfrak{a}_{M_{1}}(F)$ appartient \`a $U_{J_{1}}$.    En composant avec l'exponentielle la d\'efinition de [II] 3.5, on obtient une application lin\'eaire
$$\sigma^{G_{i}}_{J_{i}}:D_{nil}^{st}(\mathfrak{m}_{i}(F))\otimes Mes(M_{i}(F))^*\to U_{J_{i}}\otimes (D_{nil}^{st}(\mathfrak{m}_{i}(F))\otimes Mes(M_{i}(F))^*)/Ann_{nil}^{G_{i},st}.$$
En reprenant la preuve du lemme [II] 3.1 , on voit que l'isomorphisme (2) 
 envoie $Ann_{nil}^{G_{1},st}$ sur $Ann_{nil}^{G_{2},st}$.   
 
 Remarquons que $J_{1}$ est l'\'el\'ement maximal de $\Sigma(A_{M_{1}},B_{1})$ si et seulement si $J_{2}$ est l'\'el\'ement maximal de $\Sigma(A_{M_{2}},B_{2})$.
 
 \bigskip
 
 \subsection{Les termes $\sigma_{J}$}
 Soient $(G_{1},G_{2},j_{*})$ un triplet endoscopique non standard, $B_{1}$ et $B_{2}$ deux fonctions comme en 6.6 v\'erifiant toutes deux les hypoth\`eses de [II] 1.8 et $M_{1}$ et $M_{2}$ deux Levi qui se correspondent.

\ass{Proposition (\`a prouver)}{ On suppose que $B_{1}$ est constante. Pour $i=1,2$, soient $J_{i}\in {\cal J}_{M_{i}}^{G_{i}}(B_{i})$ et $\boldsymbol{\delta}_{i}\in D_{nil}^{st}(\mathfrak{m}_{i}(F))\otimes Mes(M_{i}(F))^*$. On suppose que $J_{1}$ et $J_{2}$ se correspondent par la bijection entre les ensembles $\Sigma(A_{M_{i}},B_{i})$ et que $\boldsymbol{\delta}_{1}$ et $\boldsymbol{\delta}_{2}$ se correspondent par l'isomorphisme 6.5(2). Alors, pour tout $X_{1}\in \mathfrak{a}_{M_{1}}(F)$ en position g\'en\'erale et proche de $0$, on a l'\'egalit\'e
$$\sigma_{J_{1}}^{G_{1}}(\boldsymbol{\delta}_{1},X_{1})=c_{M_{1},M_{2}}^{G_{1},G_{2}}\sigma_{J_{2}}^{G_{2}}(\boldsymbol{\delta}_{2},j_{*}(X_{1})).$$}

Montrons que

(1) cette assertion est v\'erifi\'ee si, pour $i=1,2$, $J_{i}$ n'est pas l'\'el\'ement maximal de $\Sigma(A_{M_{i}},B_{i})$.

Preuve.  Pour $i=1,2$, on construit le  groupe $G_{i,J_{i}}$. D'apr\`es la proposition 1.2(i), on a l'\'egalit\'e
$$\sigma_{J_{i}}^{G_{i}}(\boldsymbol{\delta}_{i},X_{i})=i_{J_{i}}^{G_{i}}\sigma_{J_{i}}^{G_{i,J_{i}}}(\boldsymbol{\delta}_{i},X_{i}).$$
Notons $M_{i,sc}$ l'image r\'eciproque de $M_{i}$ dans $G_{i,J_{i},SC}$. On peut fixer $\boldsymbol{\delta}'_{i}\in D_{nil}^{st}(\mathfrak{m}_{i,sc}(F))\otimes Mes(M_{i,sc}(F))^*$ tel que
$$\boldsymbol{\delta}_{i}=\iota^*_{M_{i,sc},M_{i}}(\boldsymbol{\delta}'_{i}).$$
D'apr\`es le lemme 3.6, on a aussi
$$\sigma_{J_{i}}^{G_{i}}(\boldsymbol{\delta}_{i},X_{i})=i_{J_{i}}^{G_{i}}\iota^*_{M_{i,sc},M_{i}}(\sigma_{J_{i}}^{G_{i,J_{i},SC}}(\boldsymbol{\delta}'_{i},X_{i})).$$
Notons $T_{i,sc}$  l'image r\'eciproque de $T_{i}$  dans $G_{i,J_{i},SC}$. On v\'erifie que $j_{*}$ se restreint en un isomorphisme de $X_{*}(T_{1,sc})\otimes_{{\mathbb Z}}{\mathbb Q}$ sur $X_{*}(T_{2,sc})\otimes_{{\mathbb Z}}{\mathbb Q}$ et que, en notant encore $j_{*}$ cette restriction, le triplet $(G_{1,J_{1},SC},G_{1,J_{2},SC},j_{*})$ est encore endoscopique non standard. Les distributions $\boldsymbol{\delta}'_{1}$ et $\boldsymbol{\delta}'_{2}$ se correspondent. En supposant que les $J_{i}$ ne sont pas maximaux, on a $dim(G_{i,J_{i},SC})<dim(G_{i})$ et on peut appliquer le lemme ci-dessus:
$$\sigma_{J_{1}}^{G_{1,J_{1},SC}}(\boldsymbol{\delta}'_{1},X_{1})=c_{M_{1,sc},M_{2,sc}}^{G_{1,J_{1},SC},G_{2,J_{2},SC}}\sigma_{J_{2}}^{G_{2,J_{2},SC}}(\boldsymbol{\delta}'_{2},X_{2}).$$
Toutes ces \'egalit\'es conduisent \`a l'\'egalit\'e de l'\'enonc\'e pourvu que l'on ait
$$(2) \qquad i_{J_{1}}^{G_{1}}c_{M_{1,sc},M_{2,sc}}^{G_{1,J_{1},SC},G_{2,J_{2},SC}}=i_{J_{2}}^{G_{2}}c_{M_{1},M_{2}}^{G_{1},G_{2}}.$$
Remarquons que, puisque $G_{i}$ est simplement connexe pour $i=1,2$, on a simplement
$$i_{J_{i}}^{G_{i}}=\vert Z(\hat{G}_{i,J_{i}})^{\Gamma_{F}}\vert ^{-1}.$$
Fixons un entier $n$ assez grand. On a alors un diagramme commutatif
$$\begin{array}{ccccccccc}&&1&&1&&1&&\\&&\downarrow&&\downarrow&&\downarrow&&\\1&\to&A&\to&Z(\hat{G}_{2,J_{2}})^{\Gamma_{F}}&\stackrel{\hat{j}_{n}}{\to}&Z(\hat{G}_{1,J_{1}})^{\Gamma_{F}}&\to&1\\ &&\downarrow&&\downarrow&&\downarrow&&\\1&\to&B&\to&Z(\hat{M}_{2})^{\Gamma_{F}}&\stackrel{\hat{j}_{n}}{\to}&Z(\hat{M}_{1})^{\Gamma_{F}}&\to&1\\ &&\downarrow&&\downarrow&&\downarrow&&\\ 1&\to &C&\to&Z(\hat{M}_{2,ad})^{\Gamma_{F}}&\stackrel{\hat{j}_{n}}{\to}&Z(\hat{M}_{1,ad})^{\Gamma_{F}}&\to&1\\ &&\downarrow&&\downarrow&&\downarrow&&\\ &&1&&1&&1&&\\ \end{array}$$
o\`u $A$, $B$ et $C$ sont les noyaux des fl\`eches horizontales de droite. Les lignes de ce diagramme sont exactes. Les deux derni\`eres colonnes aussi. Il en r\'esulte que la premi\`ere colonne est exacte. D'o\`u l'\'egalit\'e $\vert B\vert =\vert A\vert \vert C\vert $. On a
$$\vert B\vert =n^{a_{M_{2}}}c_{M_{1},M_{2}}^{G_{1},G_{2}},$$
$$\vert C\vert =n^{a_{M_{2}}}c_{M_{1,sc},M_{2,sc}}^{G_{1,J_{1},SC},G_{2,J_{2},SC}},$$
et
$$\vert A\vert =\vert Z(\hat{G}_{2,J_{2}})^{\Gamma_{F}}\vert Z(\hat{G}_{1,J_{1}})^{\Gamma_{F}}\vert  ^{-1}=i_{J_{1}}^{G_{1}}(i_{J_{2}}^{G_{2}})^{-1}.$$
D'o\`u (2), ce qui ach\`eve la preuve. $\square$

\bigskip

\subsection{Germes de Shalika}
On conserve les m\^emes donn\'ees.
 Pour $i=1,2$, on sait d\'efinir le germe $Sg_{M_{i},unip}^{G_{i}}(B_{i})$ au voisinage de l'origine dans $M_{i}(F)$. Puisqu'il vit au voisinage de l'origine, on peut le descendre par l'exponentielle en un germe sur l'alg\`ebre de Lie que l'on note $Sg_{M_{i},nil}^{G_{i}}(B_{i})$. C'est un germe d'application lin\'eaire
$$D_{g\acute{e}om,G_{i}-\acute{e}qui}^{st}(\mathfrak{m}_{i}(F))\otimes Mes(M_{i}(F))^*\to D_{nil}^{st}(\mathfrak{g}_{i}(F))\otimes Mes(G_{i}(F))^*.$$
Remarquons que la notion d'\'el\'ement $G_{i}$-\'equisingulier se d\'efinit dans les alg\`ebres de Lie comme dans les groupes, cf. [II] 1.2.

\ass{Proposition (\`a prouver)}{On suppose que $B_{1}$ est constante. Pour $i=1,2$, soit $\boldsymbol{\delta}_{i}\in D_{g\acute{e}om,G_{i}-\acute{e}qui}^{st}(\mathfrak{m}_{i}(F))\otimes Mes(M_{i}(F))^*$.  Supposons que $\boldsymbol{\delta}_{1}$ et $\boldsymbol{\delta}_{2}$ se correspondent par l'isomorphisme 6.5(2). Alors, si $\boldsymbol{\delta}_{1}$ et $\boldsymbol{\delta}_{2}$ sont assez voisins de l'origine, les termes $Sg_{M_{1},nil}^{G_{1}}(\boldsymbol{\delta}_{1},B_{1})$ et $c_{M_{1},M_{2}}^{G_{1},G_{2}}Sg_{M_{2},nil}^{G_{2}}(\boldsymbol{\delta}_{2},B_{2})$ se correspondent par l'isomorphisme 6.5(1).}

\bigskip

\subsection{R\'eduction des propositions 6.6 et  6.7}
\ass{Lemme}{Supposons que la proposition 6.6, resp. la proposition 6.7, soit v\'erifi\'ee dans le cas o\`u $(G_{1},G_{2},j_{*})$ est quasi-\'el\'ementaire et $B_{1}$ est la fonction constante de valeur $1$. Alors la proposition 6.6, resp. la proposition 6.7, est  v\'erifi\'ee.}

 Preuve. On a besoin de quelques propri\'et\'es pr\'eliminaires. Consid\'erons un seul groupe $G$ r\'eductif connexe et simplement connexe d\'efini sur $F$, un Levi $M$ de $G$ et une fonction $B$ comme en [II] 1.8.     Soient  $\boldsymbol{\delta}\in D_{g\acute{e}om,G-\acute{e}qui}^{st}(\mathfrak{m}(F))\otimes Mes(M(F))^*$, $\boldsymbol{\delta}_{nil}\in D_{nil}^{st}(\mathfrak{m}(F))\otimes Mes(M(F))^*$, $H\in \mathfrak{a}_{M}(F)$ en position g\'en\'erale  et proche de $0$ et soit $J\in {\cal J}_{M}^G(B)$. On a d\'efini $\sigma_{J}^G(\boldsymbol{\delta}_{nil},H)$ et $Sg_{M,nil}^G(\boldsymbol{\delta},B)$. Ces termes d\'ependent de la mesure que l'on a fix\'ee en [II] 1.2 sur ${\cal A}_{M}^G$. Pour un instant, notons $m$ cette mesure et introduisons-la dans la notation. Si on remplace $m$ par $cm$, avec $c\in {\mathbb R}_{>0}$,
  il r\'esulte imm\'ediatement des d\'efinitions que
 $$(1)\qquad \left\lbrace\begin{array}{ccc}\sigma_{J}^G{\boldsymbol{\delta}}_{nil},H,cm)&=&c \sigma_{J}^G(\boldsymbol{\delta}_{nil},H,m),\\ Sg_{M,nil}^G(\boldsymbol{\delta},B,cm)&=&cSg_{M,nil}^G(\boldsymbol{\delta},B,m).\\ \end{array}\right.$$
 Oublions cette parenth\`ese, la mesure $m$ est maintenant fix\'ee.

 Soit $r\in F^{\times}$. On d\'efinit deux homomorphismes 
 $$\begin{array}{ccc}C_{c}^{\infty}(\mathfrak{g}(F))&\to& C_{c}^{\infty}(\mathfrak{g}(F))\\ f&\mapsto& f[r]\\ \end{array}$$
 $$\begin{array}{ccc}D_{g\acute{e}om}(\mathfrak{g}(F))\otimes Mes(G(F))^*&\to &D_{g\acute{e}om}(\mathfrak{g}(F))\otimes Mes(G(F))^*\\ \boldsymbol{\gamma}&\mapsto& \boldsymbol{\gamma}[r]\\ \end{array}$$
 de la fa\c{c}on suivante. Rappelons que les donn\'ees d'un \'el\'ement $X\in \mathfrak{g}(F)$  et d'une mesure de Haar sur $G_{X}(F)$ d\'efinissent un \'el\'ement $\boldsymbol{\gamma}\in D_{g\acute{e}om}(\mathfrak{g}(F))\otimes Mes(G(F))^*$. Pour $f\in C_{c}^{\infty}(\mathfrak{g}(F))$ et une mesure de Haar $dg$ sur $G(F)$, on a
 $$I^G(\boldsymbol{\gamma},f\otimes dg)=D^G(X)^{1/2}\int_{G_{X}(F)\backslash G(F)}f(ad_{x}^{-1}(X))\,dx,$$
 o\`u $dx$ est le quotient de $dg$ par la mesure fix\'ee sur $G_{X}(F)$.Consid\'erons d'abord un \'el\'ement $X\in \mathfrak{g}_{reg}(F)$, posons $T=G_{X}$. Munissons $\mathfrak{t}(F)$ d'une mesure de Haar. On en d\'eduit via l'exponentielle une mesure de Haar sur $T(F)$. Notons $\boldsymbol{\gamma}$ l'int\'egrale orbitale associ\'ee.   La multiplication par $r$ envoie $X$ sur $rX$ et transporte la mesure sur $\mathfrak{t}(F)$ sur cette mesure multipli\'ee par $\vert r\vert _{F}^{-dim(T)}$. On note $\boldsymbol{\gamma}[r]$ l'int\'egrale orbitale associ\'ee \`a $rX$ et cette nouvelle mesure. En posant $\delta(G)=dim(G)+dim(T)$ et en d\'efinissant $f[r]$ par l'\'egalit\'e
 $$f[r](Y)=\vert r\vert _{F}^{\delta(G)/2}f(rY)$$
 pour tout $Y\in \mathfrak{g}(F)$, on voit que
 $$(2) \qquad I^G(\boldsymbol{\gamma}[r],f\otimes dg)=I^G(\boldsymbol{\gamma},f[r]\otimes dg).$$
 Pour $\boldsymbol{\gamma}$ quelconque, $\boldsymbol{\gamma}[r]$ est d\'efini par cette relation. Soit $X\in \mathfrak{g}(F)$ quelconque, fixons une mesure sur $G_{X}(F)$  et notons $\boldsymbol{\gamma}$ l'int\'egrale orbitale associ\'ee. On v\'erifie que $\boldsymbol{\gamma}[r]$  est \'egale \`a $\vert r\vert _{F}^{\delta(G_{X_{ss}})/2}\boldsymbol{\gamma}'$, o\`u $X_{ss}$ est la partie semi-simple de $X$ et  $\boldsymbol{\gamma}'$ est l'int\'egrale orbitale associ\'ee \`a $rX$ et \`a la m\^eme mesure sur $G_{X}(F)$.  En tensorisant avec l'identit\'e de $Mes(G(F))$, on obtient des transformations de $I(\mathfrak{g}(F))\otimes Mes(G(F))$ ou  $SI(\mathfrak{g}(F))\otimes Mes(G(F))$. 
Dualement, l'application $\boldsymbol{\gamma}\mapsto \boldsymbol{\gamma}[r]$ pr\'eserve les distributions stables.   Soit $s\in Z(\hat{M})^{\Gamma_{F}}/Z(\hat{G})^{\Gamma_{F}}$. De la donn\'ee endoscopique maximale ${\bf M}$ de $M$ et de $s$ se d\'eduit une donn\'ee endoscopique ${\bf G}'(s)=(G'(s),{\cal G}'(s),s)$. Notons que, puisqu'on travaille avec des alg\`ebres de Lie, l'introduction de donn\'ees auxiliaires pour la donn\'ee endoscopique ${\bf G}'(s)$ est inutile. Les facteurs de transfert sont normalis\'es de sorte qu'ils vaillent $1$ sur la diagonale dans $\mathfrak{m}(F)\times \mathfrak{m}(F)$. Fixons $s$. On a donc un facteur de transfert $\Delta(s)$ sur un sous-ensemble de $\mathfrak{g}'(s;F)\times \mathfrak{g}(F)$. Le lemme 3.2.1 de [F] affirme qu'il existe un caract\`ere $\chi$ de $F^{\times}$ tel que $\Delta(s)(\lambda Y,\lambda X)=\chi(\lambda)\Delta(s)(Y,X)$ pour tout couple $(Y,X)$ et tout $\lambda\in F^{\times}$. En consid\'erant un couple $Y=X\in \mathfrak{m}(F)$, on obtient $\chi=1$.  Par un calcul simple, on en d\'eduit 
$$(3) \qquad  \left\lbrace\begin{array}{ccc}(f[r])^{G'(s)}&= &(f^{G'(s)})[r]\\ transfert(\boldsymbol{\tau}[r])&=&(transfert(\boldsymbol{\tau}))[r] \\ \end{array}\right.$$
pour tout $f\in I(\mathfrak{g}(F))\otimes Mes(G(F))$ et tout $\boldsymbol{\tau}\in D_{g\acute{e}om}^{st}(\mathfrak{g}'(s;F))\otimes Mes(G'(s;F))^*$.

Soit maintenant $r\in {\mathbb Q}_{>0}$. On peut appliquer les constructions ci-dessus \`a cet \'el\'ement. D'autre part, la fonction $rB$ v\'erifie les m\^emes hypoth\`eses que $B$. L'ensemble $\Sigma(A_{M},rB)$ est form\'e des $\alpha/r$ pour $\alpha\in \Sigma(A_{M},B)$. On en d\'eduit une biejction ${\cal J}_{M}^G(B)\simeq {\cal J}_{M}^G(rB)$. Notons $J/r$ l'image de $J$ dans ce dernier ensemble. Montrons que

$$(4) \qquad \left\lbrace\begin{array}{ccc}\sigma_{J/r}^G(\boldsymbol{\delta}_{nil}[r],rH)&=&(\sigma_{J}^G(\boldsymbol{\delta}_{nil},H))[r],\\ Sg_{M,nil}^G(\boldsymbol{\delta}[r],rB)&=&(Sg_{M,nil}^G(\boldsymbol{\delta},B))[r].\\ \end{array}\right.$$

Pour la premi\`ere \'egalit\'e, la d\'efinition [II] 3.5(1) nous ram\`ene par r\'ecurrence \`a prouver l'assertion analogue pour les termes $\rho_{J}^G$. On peut alors lever l'hypoth\`ese que $\boldsymbol{\delta}_{nil}$ est stable. On peut supposer que $\boldsymbol{\delta}_{nil}$ est l'int\'egrale orbitale associ\'ee \`a un \'el\'ement nilpotent $N\in \mathfrak{m}(F)$ et une mesure sur $M_{N}(F)$. Donc $\boldsymbol{\delta}_{nil}[r]$ est en tout cas une int\'egrale orbitale associ\'ee \`a $rN$.  Alors, d'apr\`es [II] 3.2(5), on a l'\'egalit\'e
$$\rho_{J/r}^G(\boldsymbol{\delta}_{nil}[r],rH)=\sum_{\underline{\alpha}'\in J/r}m(\underline{\alpha}',rN)sgn(\underline{\alpha}',rN)u_{\underline{\alpha}'}(rH)\boldsymbol{\delta}_{nil}[r].$$
Pour $\underline{\alpha}=\{\alpha_{1},...,\alpha_{n}\}\in J$, posons $\underline{\alpha}/r=\{\alpha_{1}/r,...,\alpha_{n}/r\}$. Alors $\underline{\alpha}/r\in J/r$ et la correspondance $\underline{\alpha}\mapsto \underline{\alpha}/r$ est bijective. Il est clair que $u_{\underline{\alpha}/r}(rH)=u_{\underline{\alpha}}(H)$. La formule ci-dessus se r\'ecrit
 $$\rho_{J/r}^G(\boldsymbol{\delta}_{nil}[r],rH)=\sum_{\underline{\alpha}\in J}m(\underline{\alpha}/r,rN)sgn(\underline{\alpha}/r,rN)u_{\underline{\alpha}}(H)\boldsymbol{\delta}_{nil}[r].$$
 Pour obtenir la premi\`ere formule de (4), il reste \`a prouver que
 $$m(\underline{\alpha}/r,rN)sgn(\underline{\alpha}/r,rN)=m(\underline{\alpha},N)sgn(\underline{\alpha},N).$$
 En revenant \`a la d\'efinition de ces termes, cf. [II] 3.2, il suffit de prouver que, pour tout $\alpha\in \Sigma(A_{M},B)$, on a $\rho(\alpha,N,B)=\rho(\alpha/r,rN,rB)$. Le groupe $G_{\alpha}$ associ\'e \`a $\alpha$ en [II] 1.8 est le m\^eme que le groupe $G_{\alpha/r}$ associ\'e \`a $\alpha/r$.  L'\'egalit\'e pr\'ec\'edente r\'esulte de la m\^eme \'egalit\'e pour ce groupe. Cela r\'esout le probl\`eme si $G_{\alpha}\not=G$. Supposons $G_{\alpha}=G$, donc $\alpha$ est indivisible. On vient de prouver la relation requise pour $n\alpha$ pour tout $n\geq2$. Alors la relation [II] 1.8(6) montre qu'il suffit de prouver l'\'egalit\'e suivante o\`u les fonctions $B$ et $rB$ ont disparu:
$$\rho(\alpha,rN)=\rho(\alpha,N)$$
pour tout $\alpha\in \Sigma(A_{L})$. De nouveau, un d\'evissage des d\'efinitions nous ram\`ene \`a prouver que les termes initiaux $\rho^{Art}(\beta,N)$ d\'efinis par Arthur sont insensibles au remplacement de $N$ par $rN$. Cela vient du fait qu'ils sont de nature g\'eom\'etrique et que, d'apr\`es la th\'eorie des $SL_{2}$-triplets,  $rN$ est conjugu\'e \`a $N$ par un \'el\'ement de $G(\bar{F})$. D'o\`u la premi\`ere assertion de (4).

Pour la deuxi\`eme assertion, la d\'efinition [II] 2.4(1) et la formule (3) ci-dessus nous ram\`ene \`a prouver l'assertion analogue pour les germes $g_{M,nil}^G(\boldsymbol{\delta},B)$. On peut de nouveau lever l'hypoth\`ese que $\boldsymbol{\delta}$ est stable. On a \'enonc\'e en [II] 2.3 la d\'efinition de ces germes pour des int\'egrales orbitales pond\'er\'ees invariantes, mais les m\^emes relations valent pour les int\'egrales orbitales pond\'er\'ees non invariantes, cf. [A2] proposition 9.1.  C'est-\`a-dire que pour ${\bf f}\in C_{c}^{\infty}(\mathfrak{g}(F))\otimes Mes(G(F))$ et pour $\boldsymbol{\delta}$ assez proche de $0$, on a l'\'egalit\'e
$$(5)\qquad J_{M}^G(\boldsymbol{\delta},{\bf f})=\sum_{L\in {\cal L}(M)}J_{L}^G(g_{M,nil}^L(\boldsymbol{\delta},B),B,{\bf f}).$$
Supposons prouv\'ee l'assertion suivante:

(6) soient ${\bf f}\in C_{c}^{\infty}(\mathfrak{g}(F))\otimes Mes(G(F)$ et $\boldsymbol{\tau}\in D_{g\acute{e}om}(\mathfrak{m}(F))\otimes Mes(M(F))^*$; supposons que les \'el\'ements du support de $\boldsymbol{\tau}$ sont $G$-\'equisinguliers ou nilpotents; alors on a l'\'egalit\'e
$$J_{M}^G(\boldsymbol{\tau}[r],rB,{\bf f})=J_{M}^G (\boldsymbol{\tau},B,f[r]).$$

En rempla\c{c}ant ${\bf f}$ par ${\bf f}[r]$ dans  (5) et en utilisant (6), on obtient
 $$J_{M}^G(\boldsymbol{\delta}[r],{\bf f})=\sum_{L\in {\cal L}(M)}J_{L}^G(g_{M,nil}^L(\boldsymbol{\delta},B)[r],rB,{\bf f}).$$
 Par r\'ecurrence, on peut utiliser la deuxi\`eme relation de (4) en y rempla\c{c}ant $G$ par tout Levi $L\not=G$. On obtient
 $$J_{M}^G(\boldsymbol{\delta}[r],{\bf f})=I^G((g_{M,nil}^G(\boldsymbol{\delta},B))[r],{\bf f})+\sum_{L\in {\cal L}(M),L\not=G}J_{L}^G(g_{M,nil}^L(\boldsymbol{\delta}[r],rB),rB,{\bf f}).$$
 En comparant avec (5) appliqu\'ee \`a $\boldsymbol{\delta}[r]$ et \`a la fonction $rB$, on obtient
 $$I^G((g_{M,nil}^G(\boldsymbol{\delta},B))[r],{\bf f})=I^G(g_{M,nil}^G(\boldsymbol{\delta}[r],rB),{\bf f}),$$
 ce qui prouve la deuxi\`eme assertion de (4). Il reste \`a prouver (6). Si $\boldsymbol{\tau}$ est \`a support $G$-\'equisingulier, les int\'egrales orbitales pond\'er\'ees sont d\'efinies par une honn\^ete int\'egrale et un simple calcul conduit \`a l'\'egalit\'e cherch\'ee (c'est essentiellement la m\^eme chose que pour (2), les fonctions poids ne perturbent pas le calcul). Supposons que $\boldsymbol{\tau}$ soit l'int\'egrale orbitale associ\'ee \`a un \'el\'ement nilpotent $N\in \mathfrak{m}(F)$ et \`a une mesure sur $M_{N}(F)$. Consid\'erons un \'el\'ement $X\in \mathfrak{a}_{M}(F)$ en position g\'en\'erale, notons $\boldsymbol{\tau}_{X}$ l'int\'egrale orbitale (dans $\mathfrak{m}(F)$) associ\'ee \`a $X+N$ et \`a la m\^eme mesure sur $M_{X+N}(F)=M_{N}(F)$. On a alors pour tout ${\bf f}\in C_{c}^{\infty}(\mathfrak{g}(F))\otimes Mes(G(F))$ une formule
 $$(7) \qquad J_{M}^G(\boldsymbol{\tau},B,{\bf f})=lim_{X\to 0}\sum_{L\in {\cal L}(M)}r_{M}^L(N,B,X)J_{L}^G(\boldsymbol{\tau}_{X}^L,{\bf f}).$$
 Rempla\c{c}ons $\boldsymbol{\tau}$ par $\boldsymbol{\tau}[r]$, $B$ par $rB$ et  $X$ par $rX$ dans les constructions. On v\'erifie sur la description explicite que l'on a donn\'ee plus haut que $\boldsymbol{\tau}_{X}$ est remplac\'e par $\boldsymbol{\tau}_{X}[r]$. On a aussi $(\boldsymbol{\tau}_{X}[r])^L=(\boldsymbol{\tau}_{X}^L)[r]$. La formule ci-dessus devient
 $$J_{M}^G(\boldsymbol{\tau}[r],rB,{\bf f})=lim_{X\to 0}\sum_{L\in {\cal L}(M)}r_{M}^L(rN,rB,rX)J_{L}^G(\boldsymbol{\tau}_{X}^L[r],{\bf f}).$$
 Par (6) appliqu\'e au cas d\'ej\`a prouv\'e des \'el\'ements \`a support $G$-\'equisingulier, c'est aussi
  $$J_{M}^G(\boldsymbol{\tau}[r],rB,{\bf f})=lim_{X\to 0}\sum_{L\in {\cal L}(M)}r_{M}^L(rN,rB,rX)J_{L}^G(\boldsymbol{\tau}_{X}^L,{\bf f}[r]).$$
  En comparant avec la formule (7) appliqu\'ee \`a ${\bf f}[r]$, on voit que, pour obtenir l'\'egalit\'e (6), il suffit de prouver l'\'egalit\'e
  $$r_{M}^L(rN,rB,rX)=r_{M}^L(N,B,X)$$
  pour tout $L$. Mais la preuve de cette \'egalit\'e est exactement la m\^eme que celle de la premi\`ere assertion de (4). Cela prouve la deuxi\`eme assertion de (4).

Venons-en \`a la preuve du lemme. Soit $(G_{1},G_{2},j_{*})$ un triplet endoscopique non standard pour lequel on veut prouver la proposition 6.6 ou 6.7.   Il est clair que, si  le triplet est produit de triplets $(G_{1,i},G_{2,i},j_{*,i})$ pour $i=1,...,m$, la proposition pour notre triplet r\'esulte de la m\^eme proposition pour chaque triplet $(G_{1,i},G_{2,i},j_{*,i})$. On peut donc supposer que $(G_{1},G_{2},j_{*})$ est \'equivalent \`a un triplet quasi-\'el\'ementaire, autrement dit on peut fixer $d\in {\mathbb Q}_{>0}$ tel que $(G_{1},G_{2},dj_{*})$ soit quasi-\'el\'ementaire. On a fix\'e une mesure $m_{1}$ sur ${\cal A}_{M_{1}}^{G_{1}}$, dont se d\'eduit via $j_{*}$ une mesure $m_{2}$ sur ${\cal A}_{M_{2}}^{G_{2}}$. La mesure $m'_{2}$ d\'eduite via $dj_{*}$ est \'egale \`a $d^{-a_{M}+a_{G}}m_{2}$, o\`u on a pos\'e $a_{M}=a_{M_{i}}$, $a_{G}=a_{G_{i}}$ pour $i=1,2$. On a fix\'e une fonction constante $B_{1}$. On note encore $B_{1}$ la valeur constante de cette fonction et on note ${\bf 1}$ la fonction constante sur $\Sigma(A_{M_{1}})$  de valeur $1$. Via $j_{*}$, on a d\'eduit de $B_{1}$ une fonction $B_{2}$. Via $dj_{*}$, la fonction d\'eduite de ${\bf 1}$ est $dB_{2}/B_{1}$. On a construit \`a l'aide de $j_{*}$ une constante $c_{M_{1},M_{2}}^{G_{1},G_{2}}$, notons-la plus pr\'ecis\'ement $c_{M_{1},M_{2}}^{G_{1},G_{2}}(j_{*})$ . On v\'erifie que la constante $c_{M_{1},M_{2}}^{G_{1},G_{2}}(dj_{*})$ construite \`a l'aide de $dj_{*}$ est $d^{a_{M}-a_{G}}c_{M_{1},M_{2}}^{G_{1},G_{2}}(j_{*})$. Soient $J_{1}\in {\cal J}_{M_{1}}^{G_{1}}(B_{1})$, $\boldsymbol{\delta}_{1}\in D_{nil}^{st}(\mathfrak{m}_{1}(F))\otimes Mes(M_{1}(F))^*$ et $X_{1}\in \mathfrak{a}_{M_{1}}(F)$ en position g\'en\'erale. Notons $J_{2}\in {\cal J}_{M_{2}}^{G_{2}}(B_{2})$ l'\'el\'ement correspondant \`a $J_{1}$. Notons plus pr\'ecis\'ement $transfert_{j_{*}}$ les isomorphismes 6.5(1) et 6.5(2) relatifs \`a $j_{*}$. Posons $\boldsymbol{\delta}_{2}=transfert_{j_{*}}(\boldsymbol{\delta}_{1})$. Pour prouver la proposition 6.6, on doit prouver que
$$(8) \qquad transfert_{j_{*}}(\sigma_{J_{1}}^{G_{1}}(\boldsymbol{\delta}_{1},X_{1},m_{1}))=c_{M_{1},M_{2}}^{G_{1},G_{2}}(j_{*})\sigma_{J_{2}}^{G_{2}}(\boldsymbol{\delta}_{2},j_{*}(X_{1}),m_{2}).$$
Les termes $B_{1}J_{1}$ et $ B_{1}J_{2}/d$ appartiennent respectivement \`a ${\cal J}_{M_{1}}^{G_{1}}({\bf 1})$ et ${\cal J}_{M_{2}}^{G_{2}}(dB_{2}/B_{1})$. Ils se correspondent. On v\'erifie sur les d\'efinitions que  $transfert_{dj_{*}}$ est le compos\'e de $transfert_{j_{*}}$ et de l'application  $\boldsymbol{\delta}\mapsto \boldsymbol{\delta}[d]$. On peut d'ailleurs composer dans l'un ou l'autre sens car le transfert commute \`a l'application $\boldsymbol{\delta}\mapsto \boldsymbol{\delta}[r]$ pour tout $r\in F^{\times}$. Donc $transfert_{dj_{*}}(\boldsymbol{\delta}_{1})= \boldsymbol{\delta}_{2}[d]$ et aussi  $transfert_{dj_{*}}(\boldsymbol{\delta}_{1}[1/B_{1}])=\boldsymbol{\delta}_{2}[d/B_{1}]$. Supposons la proposition 6.6 connue pour le triplet quasi-\'el\'ementaire  $(G_{1},G_{2},dj_{*})$ et la fonction ${\bf 1}$. Elle implique l'\'egalit\'e
$$transfert_{dj_{*}}(\sigma_{B_{1}J_{1}}^{G_{1}}(\boldsymbol{\delta}_{1}[1/B_{1}],X_{1}/B_{1},m_{1}))=c_{M_{1},M_{2}}^{G_{1},G_{2}}(dj_{*})\sigma_{B_{1}J_{2}/d}^{G_{2}}(\boldsymbol{\delta}_{2}[d/B_{1}],dj_{*}(X_{1}/B_{1}),m'_{2}).$$
D'apr\`es (4), on a l'\'egalit\'e
$$\sigma_{B_{1}J_{1}}^{G_{1}}(\boldsymbol{\delta}_{1}[1/B_{1}],X_{1}/B_{1},m_{1})=\sigma_{J_{1}}^{G_{1}}(\boldsymbol{\delta}_{1},X_{1},m_{1})[1/B_{1}].$$
 D'apr\`es (1) et (4), on a l'\'egalit\'e
$$\sigma_{B_{1}J_{2}/d}^{G_{2}}(\boldsymbol{\delta}_{2}[d/B_{1}],dj_{*}(X_{1}/B_{1}),m'_{2})=d^{-a_{M}+a_{G}}\sigma_{J_{2}}^{G_{2}}(\boldsymbol{\delta}_{2},j_{*}(X_{1}),m_{2})[d/B_{1}].$$
L'\'egalit\'e pr\'ec\'edente devient
$$transfert_{dj_{*}}(\sigma_{J_{1}}^{G_{1}}(\boldsymbol{\delta}_{1},X_{1},m_{1})[1/B_{1}])=c_{M_{1},M_{2}}^{G_{1},G_{2}}(j_{*})\sigma_{J_{2}}^{G_{2}}(\boldsymbol{\delta}_{2},j_{*}(X_{1}),m_{2})[d/B_{1}].$$
En vertu des propri\'et\'es d\'ej\`a signal\'ees reliant le transfert aux applications $\boldsymbol{\delta}\mapsto \boldsymbol{\delta}[r]$, cette \'egalit\'e est \'equivalente \`a (8) que l'on voulait prouver. Cela prouve la proposition 6.6 pour notre triplet $(G_{1},G_{2},j_{*})$ et nos fonctions $B_{1}$ et $B_{2}$.

L'assertion du lemme concernant la proposition 6.7 se prouve de fa\c{c}on analogue. $\square$

\bigskip

\section{Preuves des th\'eor\`emes [II] 1.10 et [II] 1.16(ii) et preuve conditionnelle du th\'eor\`eme [II] 1.16(i)}

\bigskip

\subsection{Les termes $\rho_{J}^{\tilde{G},{\cal E}}$}
On consid\`ere un triplet  $(G,\tilde{G},{\bf a})$  quelconque, un espace de Levi $\tilde{M}$  de $\tilde{G}$  et une donn\'ee endoscopique elliptique et relevante ${\bf M}'=(M',{\cal M}',\tilde{\zeta})$   de $(M,\tilde{M},{\bf a})$. 

On suppose donn\'e un diagramme $(\epsilon,B^{M'},T',B^M,T,\eta)$ joignant un \'el\'ement $\epsilon\in \tilde{M}'_{ss}(F)$ \`a un \'el\'ement $\eta\in \tilde{M}_{ss}(F)$. On suppose que $M'_{\epsilon}$ est quasi-d\'eploy\'e et que  $A_{M'_{\epsilon}}=A_{M'}$.
 On note ${\cal O}'$ la classe de conjugaison stable de $\epsilon$ dans $\tilde{M}'(F)$ et ${\cal O}$ la classe de conjugaison stable de $\eta$ dans $\tilde{M}(F)$. Comme on l'a vu, tout \'el\'ement $\tilde{s}\in \tilde{\zeta}Z(\hat{M})^{\Gamma_{F},\hat{\theta}}/Z(\hat{G})^{\Gamma_{F},\hat{\theta}}$ donne naissance \`a un triplet endoscopique non standard $(\bar{G}'(\bar{s})_{SC},G'(\tilde{s})_{\epsilon,SC},j_{*})$. On a d\'efini un syst\`eme de fonctions $B^{\tilde{G}}$ sur $\tilde{G}'(\tilde{s})$, dont on d\'eduit une fonction $B^{\tilde{G}}_{{\cal O}'}$ sur le syst\`eme de racines de $G'(\tilde{s})_{\epsilon,SC}$. Par la construction de 6.5, on en d\'eduit une fonction sur le syst\`eme de racines de $\bar{G}'(\bar{s})_{SC}$. On a
 
 (1) cette foncton est constante de valeur $1$.
 
  Le syst\`eme de racines $B^{\tilde{G}}$ a pr\'ecis\'ement \'et\'e d\'efini pour qu'il en soit ainsi.
 Pour v\'erifier cette propri\'et\'e, il suffit de reprendre la d\'efinition de [II] 1.11 du syst\`eme de fonctions $B^{\tilde{G}}$,  celle de 6.5 de la fonction associ\'ee sur le syst\`eme de racines de $\bar{G}'(\bar{s})_{SC}$ et d'utiliser les descriptions des syst\`emes de racines des groupes $G'(\tilde{s})_{\epsilon,SC}$ et $\bar{G}'(\bar{s})_{SC}$ donn\'ee en [W1] 3.3.

Suppsons $A_{G'(\tilde{s})_{\epsilon}}=A_{\tilde{G}}$.  On a  un diagramme
 $$\begin{array}{ccccc}{\cal J}_{M'_{\epsilon}}^{\tilde{G}'(\tilde{s})_{\epsilon}}(B_{{\cal O}'})&\to&{\cal J}_{\tilde{M}'}^{\tilde{G}'(\tilde{s})}(B_{{\cal O}'})&&\\ &&&\searrow&\\ \parallel&&&&{\cal J}_{\tilde{M}}^{\tilde{G}}\\&&&\nearrow&\\ {\cal J}_{\bar{M}'}^{\bar{G}'(\bar{s})}&\to&{\cal J}_{\bar{M}}^{\bar{G}}&&\\ \end{array}$$
 Les fl\`eches se d\'efinissent par endoscopie ou descente, compte tenu du fait que les ensembles en question sont insensibles au remplacement d'un groupe par le rev\^etement simplement connexe de son groupe d\'eriv\'e. On v\'erifie que ce diagramme est commutatif.  Toutes les fl\`eches sont injectives et, pour simplifier, on consid\`ere chacun des ensembles comme un sous-ensemble de ${\cal J}_{\tilde{M}}^{\tilde{G}}$.  Fixons $J\in {\cal J}_{\tilde{M}}^{\tilde{G}}$. On consid\`ere l'hypoth\`ese suivante

(2)  pour chaque triplet $(\bar{G}'(\bar{s})_{SC},G'_{\epsilon,SC},j_{*})$ ci-dessus pour lequel $A_{G'(\tilde{s})_{\epsilon}}=A_{\tilde{G}}$ et $J\in {\cal J}_{\bar{M}'}^{\bar{G}'(\bar{s})}$, la proposition 6.6 est v\'erifi\'ee pour cet \'el\'ement $J$.

  Rappelons la proposition [II] 3.8 que nous allons prouver  sous ces hypoth\`eses.

\ass{Proposition}{On suppose  $A_{M'_{\epsilon}}=A_{M'}$. Pour tout $J\in {\cal J}_{\tilde{M}}^{\tilde{G}}$ tel que (2) soit v\'erifi\'ee, tout $\boldsymbol{\delta}\in D^{st}_{G\acute{e}om}({\bf M}',{\cal O}')\otimes Mes(M'(F))^*$ et tout $a\in A_{\tilde{M}}(F)$ en position g\'en\'erale et proche de $1$, on a l'\'egalit\'e
$$\rho_{J}^{\tilde{G},{\cal E}}({\bf M}',\boldsymbol{\delta},a)=\rho_{J}^{\tilde{G}}(transfert(\boldsymbol{\delta}),a).$$}

Preuve.   Si $(G,\tilde{G},{\bf a})$ est quasi-d\'eploy\'e et \`a torsion int\'erieure et si ${\bf M}'={\bf M}$, l'\'enonc\'e est tautologique: le terme $\sigma_{J}^{\tilde{G}}(\boldsymbol{\delta},a)$ est d\'efini pour qu'il en soit ainsi. On exclut ce cas.

Rappelons la d\'efinition
$$(3) \qquad \rho_{J}^{\tilde{G},{\cal E}}({\bf M}',\boldsymbol{\delta},a)=\sum_{\tilde{s}\in \tilde{ \zeta} Z(\hat{M})^{\Gamma_{F},\hat{\theta}}/Z(\hat{G})^{\Gamma_{F},\hat{\theta}}}i_{\tilde{M}'}(\tilde{G},\tilde{G}'(\tilde{s}))\sum_{J'\in {\cal J}_{\tilde{M}'}^{\tilde{G}'(s)}(B^{\tilde{G}}_{{\cal O}'}); J'\mapsto J}$$
$$transfert(\sigma_{J'}^{{\bf G}'(s)}(\boldsymbol{\delta},a)).$$
Expliquons la notation $a$ du membre de droite. L'\'el\'ement initial $a$ appartient \`a $A_{\tilde{M}}(F)$. Il est proche de $1$, on peut l'\'ecrire $a=exp(H)$, o\`u $H\in \mathfrak{a}_{\tilde{M}}(F)$ est proche de $0$. Seules comptent les valeurs $u(a)$ pour $u\in U_{J}$, a fortiori seules comptent les valeurs $\alpha(H)$ pour $\alpha\in \Sigma(A_{\tilde{M}})$. Autrement dit, seule compte l'image de $H$ dans $\mathfrak{a}_{\tilde{M}}(F)/\mathfrak{a}_{\tilde{G}}(F)$. Pour $\tilde{s}$ apparaissant ci-dessus, avec $i_{\tilde{M}'}(\tilde{G},\tilde{G}'(\tilde{s}))\not=0$ donc ${\bf G}'(\tilde{s})$ elliptique, on a un isomorphisme naturel
$$\mathfrak{a}_{\tilde{M}}(F)/\mathfrak{a}_{\tilde{G}}(F)\simeq \mathfrak{a}_{M'}(F)/\mathfrak{a}_{G'(\tilde{s})}(F).$$
On note encore $H$ un \'el\'ement de $\mathfrak{a}_{M'}(F)$ qui a m\^eme image que le $H$ initial dans le quotient commun ci-dessus et on note encore $a$ l'\'el\'ement $exp(H)\in A_{M'}(F)$. Une convention analogue sera utilis\'ee diverses fois dans la suite du calcul.

On reprend les constructions et notations de la section 5. Apr\`es avoir fix\'e des donn\'ees auxiliaires $M'_{1},...,\Delta_{1}$, on identifie $\boldsymbol{\delta}$ \`a un \'el\'ement $\boldsymbol{\delta}_{1}\in D^{st}_{g\acute{e}om}(\tilde{M}'_{1}(F))$ auquel on applique les consid\'erations de 5.5. Ici, les parties semi-simples des \'el\'ements du support de notre \'el\'ement $\boldsymbol{\delta}$ appartiennent \`a ${\cal O}'$. Donc les termes $Z$ apparaissant en 5.5 sont nuls. Il existe donc $\boldsymbol{\delta}_{\epsilon,SC}\in D^{st}_{unip}(M'_{\epsilon,SC}(F))$ tel que 
$$\boldsymbol{\delta}_{1}=desc_{\epsilon_{1}}^{st,\tilde{M}'_{1},*}\circ\iota^*_{M'_{\epsilon,SC},M'_{1,\epsilon_{1}}}(\boldsymbol{\delta}_{\epsilon,SC}).$$
 
Soit $\tilde{s}\in \tilde{\zeta} Z(\hat{M})^{\Gamma_{F},\hat{\theta}}/Z(\hat{G})^{\Gamma_{F},\hat{\theta}}$. L'\'el\'ement $\boldsymbol{\delta}$ s'identifie aussi \`a un \'el\'ement $\boldsymbol{\delta}_{1}(\tilde{s})\in D^{st}_{g\acute{e}om}(\tilde{M}'_{1}(\tilde{s})(F))$. D'apr\`es 5.5, on peut supposer
$$\boldsymbol{\delta}_{1}(\tilde{s})=d(\tilde{s})desc_{\epsilon_{1}(\tilde{s}}^{st,\tilde{M}'_{1}(\tilde{s}),*}\circ\iota_{M'_{\epsilon,SC},M'_{1}(\tilde{s})_{\epsilon_{1}(\tilde{s})}}(\boldsymbol{\delta}_{\epsilon,SC} ).$$
Introduisons le groupe interm\'ediaire $M'(\tilde{s})_{\epsilon,sc}$, image r\'eciproque de $M'_{\epsilon}$ dans $G'(\tilde{s})_{\epsilon,SC}$. L'\'egalit\'e ci-dessus entra\^{\i}ne
$$\boldsymbol{\delta}_{1}(\tilde{s})=d(\tilde{s})desc_{\epsilon_{1}}^{st,\tilde{M}'_{1}(\tilde{s}),*}\circ\iota_{M'(\tilde{s})_{\epsilon,sc},M'_{1}(\tilde{s})_{\epsilon_{1}(\tilde{s})}}^*(\boldsymbol{\delta}(\tilde{s})_{\epsilon,sc} ),$$
o\`u
$$\boldsymbol{\delta}(\tilde{s})_{\epsilon,sc}=\iota^*_{M'_{\epsilon,SC},M'(\tilde{s})_{\epsilon,sc}}(\boldsymbol{\delta}_{\epsilon,SC}).$$
Supposons que $J$ provienne d'un \'el\'ement de ${\cal J}_{\tilde{M}'}^{\tilde{G}'(\tilde{s})}(B^{\tilde{G}}_{{\cal O}'})$. Celui-ci est alors unique et,  conform\'ement \`a ce que l'on a dit avant l'\'enonc\'e, on le note encore $J$. On a
$$\sigma_{J}^{{\bf G}'(\tilde{s})}(\boldsymbol{\delta},a)=\sigma_{J}^{\tilde{G}'_{1}(\tilde{s})}(\boldsymbol{\delta}_{1}(\tilde{s}),a).$$
Appliquons la proposition 4.3. C'est loisible  car $dim(G'(s)_{SC})<dim(G_{SC})$. En effet, le seul cas o\`u cette in\'egalit\'e n'est pas v\'erifi\'ee est celui o\`u  $(G,\tilde{G},{\bf a})$ est quasi-d\'eploy\'e et \`a torsion int\'erieure et o\`u ${\bf M}'={\bf M}$. Or on a exclu ce cas. Rappelons que l'hypoth\`ese que $J\in {\cal J}_{\tilde{M}'}^{\tilde{G}'(\tilde{s})}(B^{\tilde{G}}_{{\cal O}'})$ implique que  $A_{G'_{1}(\tilde{s})}=A_{G'_{1}(\tilde{s})_{\epsilon_{1}(\tilde{s})}}$, ce qui \'equivaut \`a  $A_{G'(\tilde{s})}=A_{G'(\tilde{s})_{\epsilon}}$, et  que $J$  provient  d'un \'el\'ement de ${\cal J}_{M'_{1}(\tilde{s})_{\epsilon_{1}(\tilde{s})}}^{G'_{1}(\tilde{s})_{\epsilon_{1}(\tilde{s})}}(B^{\tilde{G}}_{{\cal O}'})$ que l'on note encore $J$.
 Alors, d'apr\`es la proposition 4.3, on a
$$ \sigma_{J}^{\tilde{G}'_{1}(\tilde{s})}(\boldsymbol{\delta}_{1}(\tilde{s}),a)=e_{\tilde{M}'_{1}(\tilde{s})}^{\tilde{G}'_{1}(\tilde{s})}(\epsilon_{1}(\tilde{s}))d(\tilde{s})desc_{\epsilon_{1}(\tilde{s})}^{st,\tilde{M}'_{1}(\tilde{s}),*}(\sigma_{J}^{G'_{1}(\tilde{s})_{\epsilon_{1}(\tilde{s})}}(\iota^*_{M'(\tilde{s})_{\epsilon,sc},M'_{1}(\tilde{s})_{\epsilon_{1}(\tilde{s})}}(\boldsymbol{\delta}(\tilde{s})_{\epsilon,sc} ),a)).$$
On a

(4) $e_{\tilde{M}'_{1}(\tilde{s})}^{\tilde{G}'_{1}(\tilde{s})}(\epsilon_{1}(\tilde{s}))=e_{\tilde{M}'}^{\tilde{G}'(\tilde{s})}(\epsilon)$.

On a un diagramme commutatif
$$\begin{array}{ccccccccc}&&1&&1&&&&\\ &&\downarrow&&\downarrow&&&&\\ 1&\to&Z(\hat{G}'(\tilde{s}))^{\Gamma_{F}}&\to&Z(\hat{G}'_{1}(\tilde{s}))^{\Gamma_{F}}&\to&\hat{C}_{1}(\tilde{s})^{\Gamma_{F}}&\to&1\\ &&\downarrow&&\downarrow&&\parallel&&\\ 1&\to&Z(\hat{M}')^{\Gamma_{F}}&\to&Z(\hat{M}'_{1}(\tilde{s}))^{\Gamma_{F}}&\to&\hat{C}_{1}(\tilde{s})^{\Gamma_{F}}&\to&1\\ \end{array}$$
La tore $C_{1}(\tilde{s})$ est induit donc $\hat{C}_{1}(\tilde{s})^{\Gamma_{F}}$ est connexe. Les derni\`eres fl\`eches horizontales sont donc surjectives. Donc les lignes  sont exactes. Les colonnes aussi, \'evidemment. Il en r\'esulte l'\'egalit\'e
$$Z(\hat{M}')^{\Gamma_{F}}/Z(\hat{G}'(\tilde{s}))^{\Gamma_{F}}=Z(\hat{M}'_{1}(\tilde{s}))^{\Gamma_{F}}/Z(\hat{G}'_{1}(\tilde{s}))^{\Gamma_{F}}.$$
De m\^eme
$$Z(\hat{M}'_{\epsilon})^{\Gamma_{F}}/Z(\hat{G}'(\tilde{s})_{\epsilon})^{\Gamma_{F}}=Z(\hat{M}'_{1}(\tilde{s})_{\epsilon_{1}(\tilde{s})})^{\Gamma_{F}}/Z(\hat{G}'_{1}(\tilde{s})_{\epsilon_{1}(\tilde{s})})^{\Gamma_{F}}.$$
L'\'egalit\'e (4) r\'esulte alors de la d\'efinition de 4.3. 

Appliquons (4) et le lemme 3.6. On obtient
$$\sigma_{J}^{\tilde{G}'_{1}(\tilde{s})}(\boldsymbol{\delta}_{1}(\tilde{s}),a)=e_{\tilde{M}'}^{\tilde{G}'(\tilde{s})}(\epsilon)d(\tilde{s})desc_{\epsilon_{1}(\tilde{s})}^{st,\tilde{M}'_{1}(\tilde{s}),*}\circ\iota^*_{M'(\tilde{s})_{\epsilon,sc},M'_{1}(\tilde{s})_{\epsilon_{1}(\tilde{s})}}(\sigma_{J}^{G'(\tilde{s})_{\epsilon,SC}}(\boldsymbol{\delta}(\tilde{s})_{\epsilon,sc},a)).$$
 Posons $\boldsymbol{\tau}(\tilde{s})_{sc}=\sigma_{J}^{G'(\tilde{s})_{\epsilon,SC}}(\boldsymbol{\delta}(\tilde{s})_{\epsilon,sc},a)$ appartient \`a $D^{st}_{unip}(M'(\tilde{s})_{\epsilon,sc}(F))$.   Avec les notations de 5.4(3), l'\'egalit\'e pr\'ec\'edente devient
 $$(5) \qquad \sigma_{J}^{\tilde{G}'_{1}(\tilde{s})}(\boldsymbol{\delta}_{1}(\tilde{s}),a)=e_{\tilde{M}'}^{\tilde{G}'(\tilde{s})}(\epsilon)d(\tilde{s})\boldsymbol{\tau}(\tilde{s})^{\tilde{M}'_{1}(\tilde{s})}.$$
 Appliquons  5.4(3).     On obtient
$$(6) \qquad transfert(\boldsymbol{\tau}(\tilde{s})^{\tilde{M}'_{1}(\tilde{s})})=\sum_{y\in \dot{{\cal Y}}^M}c^M[y]d(\tilde{s},y)\boldsymbol{\tau}[y]^{\tilde{M}}.$$
Reprenons la construction de $\boldsymbol{\tau}[y]^{\tilde{M}}$ pour $y\in \dot{{\cal Y}}^M$. L'\'el\'ement $\bar{\boldsymbol{\tau}}(\bar{s})_{sc}$ est le transfert non standard de $\boldsymbol{\tau}(\tilde{s})_{sc}$, c'est-\`a-dire de $\sigma_{J}^{G'(\tilde{s})_{\epsilon,SC}}(\boldsymbol{\delta}(\tilde{s})_{\epsilon,sc},a)$.    
De l'\'el\'ement $J$ de ${\cal J}_{M'(\tilde{s})_{\epsilon,sc}}^{G'(\tilde{s})_{\epsilon,SC}}(B^{\tilde{G}}_{{\cal O}'})$  se d\'eduit un \'el\'ement de ${\cal J}_{\bar{M}'(\bar{s})_{sc}}^{\bar{G}(\bar{s})_{SC}}$ que l'on note encore $J$. Notons $\bar{\boldsymbol{\delta}}_{SC}$ l'image par transfert non standard de $\boldsymbol{\delta}_{\epsilon,SC}$. C'est un \'el\'ement de $D_{unip}^{st}(\bar{M}'_{SC}(F))$. En utilisant l'analogue de 3.7(4) pour le transfert non standard, on obtient que le transfert non standard de $\boldsymbol{\delta}(\tilde{s})_{\epsilon,sc}$ est $\iota^*_{\bar{M}'_{SC},\bar{M}'(\bar{s})_{sc}}(\bar{\boldsymbol{\delta}}_{SC})$. Notons $\bar{\boldsymbol{\delta}}(\bar{s})_{sc}$ cet \'el\'ement. Utilisons l'hypoth\`ese (2). Elle nous dit que le transfert non standard de $\sigma_{J}^{G'(\tilde{s})_{\epsilon,SC}}(\boldsymbol{\delta}(\tilde{s})_{\epsilon,sc},a)$ est \'egal \`a $ c\sigma_{J}^{\bar{G}'(\bar{s})_{SC}}(\bar{\boldsymbol{\delta}}(\bar{s})_{sc},a)$, o\`u
$$c=(c_{\bar{M}'(\bar{s})_{sc},M'(\tilde{s})_{\epsilon,sc}}^{\bar{G}'(\bar{s})_{SC},G'(\tilde{s})_{\epsilon,SC}})^{-1}.$$
 Autrement dit
$$\bar{\boldsymbol{\tau}}(\bar{s})_{sc}= c\sigma_{J}^{\bar{G}'(\bar{s})_{SC}}(\bar{\boldsymbol{\delta}}(\bar{s})_{sc},a).$$
On a $\bar{\boldsymbol{\tau}}(\bar{s})=\iota^*_{\bar{M}'(\bar{s})_{sc},\bar{M}'}(\bar{\boldsymbol{\tau}}(\bar{s})_{sc})$.
En utilisant le lemme 3.6, on obtient
$$\bar{\boldsymbol{\tau}}(\bar{s})=c\sigma_{J}^{\bar{G}'(\bar{s})}(\bar{\boldsymbol{\delta}},a),$$
o\`u
$$\bar{\boldsymbol{\delta}}=\iota^*_{\bar{M}'(\bar{s})_{sc},\bar{M}'}(\bar{\boldsymbol{\delta}}(\bar{s})_{sc}).$$
Remarquons que l'on a aussi
$$\bar{\boldsymbol{\delta}}=\iota^*_{\bar{M}'_{SC},\bar{M}'}(\bar{\boldsymbol{\delta}}_{SC}).$$
Cette distribution ne d\'epend pas de $\bar{s}$. Ensuite
$$\boldsymbol{\tau}[y]^{\tilde{M}}=desc_{\eta[y]}^{\tilde{M},*}\circ\iota^*_{M_{\eta[y],sc},M_{\eta[y]}}\circ transfert_{y}(\bar{\boldsymbol{\tau}}(\bar{s}))$$
$$=c \,desc_{\eta[y]}^{\tilde{M},*}\circ\iota^*_{M_{\eta[y],sc},M_{\eta[y]}}\circ transfert_{y}(\sigma_{J}^{\bar{M}'}(\bar{\boldsymbol{\delta}},a)).$$
On a ajout\'e un indice $y$ pour rappeler qu'il s'agit du transfert de $\bar{M}'$ vers $M_{\eta[y]}$.
En utilisant (5) et (6), on obtient 
 $$(7) \qquad transfert(\sigma_{J}^{\tilde{G}'_{1}(\tilde{s})}(\boldsymbol{\delta}_{1}(\tilde{s}),a))=\sum_{y\in \dot{{\cal Y}}^M}c^M[y] d(\tilde{s},y)e_{\tilde{M}'}^{\tilde{G}'(\tilde{s})}(\epsilon)d(\tilde{s})(c_{\bar{M}'(\bar{s})_{sc},M'(\tilde{s})_{\epsilon,sc}}^{\bar{G}'(\bar{s})_{SC},G'(\tilde{s})_{\epsilon,SC}})^{-1}$$
 $$desc_{\eta[y]}^{\tilde{M},*}\circ\iota^*_{M_{\eta[y],sc},M_{\eta[y]}}\circ transfert_{y}( \sigma_{J}^{\bar{G}'(\bar{s})}(\bar{\boldsymbol{\delta}},a)).$$

Soit $Y\in \mathfrak{m}'_{\epsilon}(F)$ en position g\'en\'erale et elliptique. Il lui correspond par la construction de 5.1 des \'el\'ements $\bar{Y}$, $X[y]_{sc}$ et $X[y]$. Dans le cas $y=1$, on note ce dernier terme $X$. Normalisons le facteur $\Delta(\bar{s},y)$ par l'\'egalit\'e
$$(8) \qquad \Delta(\bar{s},y)(exp(\bar{Y}),exp(X_{sc}[y]))=\Delta_{1}(exp(Y)\epsilon_{1},exp(X[y])\eta[y]).$$
On a

(9) $d(\tilde{s})d(\tilde{s},y)=1$ pour tout $y\in \dot{{\cal Y}}^M$.

Par 5.4(1), on a
$$d(\tilde{s},y)=\frac{\Delta(\tilde{s})_{1}(exp(Y)\epsilon_{1}(\tilde{s}),exp(X[y])\eta[y])}{\Delta(\bar{s},y)(exp(\bar{Y}),exp(X_{sc}[y]))}=\frac{\Delta(\tilde{s})_{1}(exp(Y)\epsilon_{1}(\tilde{s}),exp(X[y])\eta[y])}{\Delta_{1}(exp(Y)\epsilon_{1},exp(X[y])\eta[y])}.$$
C'est un rapport de facteurs de transfert pour deux s\'eries de donn\'ees auxiliaires relatives \`a la m\^eme donn\'ee ${\bf M}'$ de $(M,\tilde{M},{\bf a})$. Ils se transforment de la m\^eme fa\c{c}on par conjugaison stable en la deuxi\`eme variable. On peut donc remplacer dans le dernier terme l'\'el\'ement $exp(X[y])\eta[y]$ par l'\'el\'ement stablement conjugu\'e $exp(X)\eta$. L'\'egalit\'e (9) r\'esulte alors de 5.5(2).

On se rappelle que l'on a suppos\'e  que  $J$ provenait d'un \'el\'ement de ${\cal J}_{\tilde{M}'}^{\tilde{G}'(\tilde{s})}(B^{\tilde{G}}_{{\cal O}'})$, ce qui entra\^{\i}ne qu'il provient d'un \'el\'ement de   ${\cal J}_{M'_{1}(\tilde{s})_{\epsilon_{1}(\tilde{s})}}^{G'_{1}(\tilde{s})_{\epsilon_{1}(\tilde{s})}}(B^{\tilde{G}}_{{\cal O}'})$. Il revient au m\^eme de dire qu'il provient d'un \'el\'ement de ${\cal J}_{M'_{\epsilon}}^{G'(\tilde{s})_{\epsilon}}(B^{\tilde{G}}_{{\cal O}'})$, ou encore, d'apr\`es ce que l'on a dit avant l'\'enonc\'e, qu'il provient d'un \'el\'ement de ${\cal J}_{\bar{M}'}^{\bar{G}'(\bar{s})}$. Notons ${\cal S}_{J}$ l'ensemble des \'el\'ements $\tilde{s}\in \tilde{\zeta} Z(\hat{M})^{\Gamma_{F},\hat{\theta}}/Z(\hat{G})^{\Gamma_{F},\hat{\theta}}$ tels que ${\bf G}'(\tilde{s})$ soit elliptique, que $A_{G'(\tilde{s})}=A_{G'(\tilde{s})_{\epsilon}}$ et  que $J$ provienne d'un \'el\'ement de ${\cal J}_{M'_{\epsilon}}^{G'(\tilde{s})_{\epsilon}}(B^{\tilde{G}}_{{\cal O}'})$. A l'aide de (7) et (9), l'expression (3) se transforme en
$$(10)\qquad \rho_{J}^{\tilde{G},{\cal E}}({\bf M}',\boldsymbol{\delta},a)=\sum_{y\in \dot{{\cal Y}}^M}c^M[y]desc_{\eta[y]}^{\tilde{M},*}\circ\iota^*_{M_{\eta[y],sc},M_{\eta[y]}}\circ transfert_{y} (\boldsymbol{\xi}),$$
o\`u
$$(11)\qquad \boldsymbol{\xi}=\sum_{\tilde{s}\in  {\cal S}_{J}}i_{\tilde{M}'}(\tilde{G},\tilde{G}'(\tilde{s}))e_{\tilde{M}'}^{\tilde{G}'(\tilde{s})}(\epsilon)(c_{\bar{M}'(\bar{s})_{sc},M'(\tilde{s})_{\epsilon,sc}}^{\bar{G}'(\bar{s})_{SC},G'(\tilde{s})_{\epsilon,SC}})^{-1} \sigma_{J}^{\bar{G}'(\bar{s})}(\bar{\boldsymbol{\delta}},a).$$
L'ensemble ${\cal S}_{J}$ est un sous-ensemble de l'ensemble ${\cal S}$ d\'efini en 5.3. D'apr\`es 
 5.3(5), il est vide si $A_{\tilde{G}}\not=A_{G_{\eta}}$.  Supposons $A_{\tilde{G}}=A_{G_{\eta}}$. Notons $\bar{{\cal S}}_{J}$ l'ensemble des $\bar{s}\in \bar{\zeta}Z(\hat{\bar{M}}_{ad})^{\Gamma_{F}}$ tels que  la donn\'ee  associ\'ee ${\bf \bar{G}}'(\bar{s})$ soit elliptique et que  $J$ provienne d'un \'el\'ement de ${\cal J}_{\bar{M}'}^{\bar{G}'(\bar{s})}$. D'apr\`es 5.3(6), ${\cal S}_{J}$ est l'ensemble des $\tilde{s}\in \tilde{\zeta} Z(\hat{M},\hat{\theta})^{\Gamma_{F}}/Z(\hat{G})^{\Gamma_{F},\hat{\theta}}$ tels que l'\'el\'ement associ\'e $\bar{s}$ appartienne \`a $\bar{{\cal S}}_{J}$. Il est clair que, si  $J$ ne provient pas d'un \'el\'ement de ${ \cal J}_{\bar{M}}^{\bar{G}}$, les deux ensembles ${\cal S}_{J}$ et $\bar{{\cal S}}_{J}$ sont vides. On obtient
 
 (12) $\rho_{J}^{\tilde{G},{\cal E}}({\bf M}',\boldsymbol{\delta},a)=0$ si $A_{\tilde{G}}\not=A_{G_{\eta}}$ ou si $A_{\tilde{G}}=A_{G_{\eta}}$ et $J$ ne provient pas d'un \'el\'ement de ${ \cal J}_{\bar{M}}^{\bar{G}}$.
 
 Supposons  que $A_{\tilde{G}}=A_{G_{\eta}}$ et  que $J$  provient  d'un \'el\'ement de ${ \cal J}_{\bar{M}}^{\bar{G}}$. Alors  la  d\'efinition (11) se transforme en
$$(13) \qquad \boldsymbol{\xi}=\sum_{\bar{s}\in \bar{{\cal S}}_{J}}x(\bar{s})\sigma_{J}^{\bar{G}'(\bar{s})}(\bar{\boldsymbol{\delta}},a),$$
o\`u
$$x(\bar{s})=\sum_{\tilde{s}\in\tilde{\zeta} Z(\hat{M})^{\Gamma_{F},\hat{\theta}}/Z(\hat{G})^{\Gamma_{F},\hat{\theta}} ; \tilde{s}\mapsto \bar{s}}i_{\tilde{M}'}(\tilde{G},\tilde{G}'(\tilde{s}))e_{\tilde{M}'}^{\tilde{G}'(\tilde{s})}(\epsilon)(c_{\bar{M}'(\bar{s})_{sc},M'(\tilde{s})_{\epsilon,sc}}^{\bar{G}'(\bar{s})_{SC},G'(\tilde{s})_{\epsilon,SC}})^{-1} .$$
Soit $\bar{s}\in \bar{{\cal S}}_{J}$. Montrons que
$$(14) \qquad  x(\bar{s})=i_{\bar{M}'}( \bar{G}_{SC},\bar{G}'(\bar{s})).$$
Pour $\tilde{s}$ se projetant sur $\bar{s}$, on a le diagramme commutatif
$$\begin{array}{ccccc}&&Z(\hat{M})^{\Gamma_{F},\hat{\theta}}/Z(\hat{G})^{\Gamma_{F},\hat{\theta}}&&\\ &\swarrow&&\searrow&\\ Z(\hat{M}')^{\Gamma_{F}}/Z(\hat{G}'(\tilde{s}))^{\Gamma_{F}}&&&&Z(\hat{\bar{M}})^{\Gamma_{F}}/Z(\hat{\bar{G}})^{\Gamma_{F}}\\ \parallel&&&&\parallel\\ Z(\hat{M}'(\tilde{s})_{ad})^{\Gamma_{F}}&&&&Z(\hat{\bar{M}}_{ad})^{\Gamma_{F}}\\ \downarrow&&&&\downarrow\\ Z(\hat{M}'(\tilde{s})_{\epsilon,ad})^{\Gamma_{F}}&&&&Z(\hat{\bar{M}}')^{\Gamma_{F}}/Z(\hat{\bar{G}}'(\bar{s}))^{\Gamma_{F}}\\ \downarrow&&&&\parallel\\ Z(\hat{\bar{M}}'(\bar{s})_{ad})^{\Gamma_{F}}&&=&&Z(\hat{\bar{M}}'(\bar{s})_{ad})^{\Gamma_{F}}\\ \end{array}$$
La fl\`eche du bas \`a gauche se d\'eduit de l'homomorphisme $\hat{T}^{\hat{\theta},0}\to \hat{T}/(1-\hat{\theta})(\hat{T})$. En utilisant la description des syst\`emes de racines de [W1] 3.3, on voit que l'image r\'eciproque par cet homomorphisme d'une racine de $\hat{\bar{G}}'(\bar{s})$ ou de $\hat{\bar{M}}'$ est un multiple entier d'une racine de $\hat{G}'(\tilde{s})_{\epsilon}$ ou de $\hat{M}'_{\epsilon}$. Il en r\'esulte que cet homomorphisme envoie $Z(\hat{G}'(\tilde{s})_{\epsilon})^{\Gamma_{F}}$ dans $Z(\hat{\bar{G}}(\bar{s}))^{\Gamma_{F}}$ et $Z(\hat{M}'_{\epsilon})^{\Gamma_{F}}$ dans $Z(\hat{\bar{M}}')^{\Gamma_{F}}$. Dans le diagramme ci-dessus, tous les quotients sont connexes, les fl\`eches sont donc surjectives. Calculons le nombre d'\'el\'ements du noyau de la fl\`eche compos\'ee
$$Z(\hat{M})^{\Gamma_{F},\hat{\theta}}/Z(\hat{G})^{\Gamma_{F},\hat{\theta}}\to Z(\hat{\bar{M}}'(\bar{s})_{ad})^{\Gamma_{F}}.$$
Si on utilise le chemin de gauche, on obtient le produit des nombres d'\'el\'ements des noyaux des  trois fl\`eches descendantes. Ces nombres sont respectivement \'egaux \`a $i_{\tilde{M}'}(\tilde{G},\tilde{G}'(\tilde{s}))^{-1}$, $e_{\tilde{M}'}^{\tilde{G}'(\tilde{s})}(\epsilon)^{-1}$ et $c_{\bar{M}'(\bar{s})_{sc},M'(\tilde{s})_{\epsilon,sc}}^{\bar{G}'(\bar{s})_{SC},G'(\tilde{s})_{\epsilon,SC}}$.  Si on utilise le chemin de droite, on obtient le produit des nombres d'\'el\'ements des noyaux des deux fl\`eches descendantes. Ces nombres sont respectivement \'egaux au nombre d'\'el\'ements $d$ de toute fibre de la projection $\tilde{s}\mapsto \bar{s}$ et \`a $i_{\bar{M}'}(\bar{G}_{SC},\bar{G}'(\bar{s}))^{-1}$. On en d\'eduit l'\'egalit\'e
$$i_{\tilde{M}'}(\tilde{G},\tilde{G}'(\tilde{s}))e_{\tilde{M}'}^{\tilde{G}'(\tilde{s})}(\epsilon)(c_{\bar{M}'(\bar{s})_{sc},M'(\tilde{s})_{\epsilon,sc}}^{\bar{G}'(\bar{s})_{SC},G'(\tilde{s})_{\epsilon,SC}})^{-1}=d^{-1}i_{\bar{M}'}(\bar{G}_{SC},\bar{G}'(\bar{s})),$$
puis
$$x(\bar{s})=d^{-1}\sum_{\tilde{s}\in\tilde{\zeta} Z(\hat{M})^{\Gamma_{F},\hat{\theta}}/Z(\hat{G})^{\Gamma_{F},\hat{\theta}} ; \tilde{s}\mapsto \bar{s}}i_{\bar{M}'}(\bar{G}_{SC},\bar{G}'(\bar{s})).$$
Puisque l'ensemble de sommation a $d$ \'el\'ements, on obtient (14).

 Gr\^ace \`a (14),  l'\'egalit\'e (13) se transforme en
$$\boldsymbol{\xi}=\sum_{\bar{s}\in \bar{\zeta}Z(\hat{\bar{M}}_{ad})^{\Gamma_{F}}, J\in {\cal J}_{\bar{M}'}^{\bar{G}'(\bar{s})}} i_{\bar{M}'}( \bar{G}_{SC},\bar{G}'(\bar{s})) \sigma_{J}^{\bar{G}'(\bar{s})}(\bar{\boldsymbol{\delta}},a).$$
Soit $y\in \dot{{\cal Y}}^M$.   En se rappelant la d\'efinition de [II] 3.8, on obtient
$$transfert_{y}(\boldsymbol{\xi})=\rho_{J}^{G_{\eta[y],SC},{\cal E}}({\bf \bar{M}}',\bar{\boldsymbol{\delta}},a).$$
Ici, les groupes ne sont pas tordus et on peut appliquer la proposition 1.4(ii).  Le terme ci-dessus n'est autre que $\rho_{J}^{G_{\eta[y],SC}}(transfert_{y}(\bar{\boldsymbol{\delta}}),a)$. On se rappelle que l'on a suppos\'e $A_{\tilde{G}}=A_{G_{\eta}}$. C'est \'equivalent \`a $A_{\tilde{G}}=A_{G_{\eta[y]}}$ puisque les deux groupes $G_{\eta}$ et $G_{\eta[y]}$ sont formes int\'erieures l'un de l'autre. De m\^eme, on a suppos\'e que $J$ provenait d'un \'el\'ement de ${\cal J}_{\bar{M}}^{\bar{G}}$, ce qui \'equivaut \`a ce qu'il provienne d'un \'el\'ement de ${\cal J}_{M_{\eta[y]}}^{G_{\eta[y]}}$. On utilise 3.2(2) et le lemme 4.1.  On obtient
$$desc_{\eta[y]}^{\tilde{M},*}\circ \iota_{M_{\eta[y],sc},M_{\eta[y]}}^*(\rho_{J}^{G_{\eta[y],SC}}(transfert_{y}(\bar{\boldsymbol{\delta}}),a))=\rho_{J}^{\tilde{G}}(\boldsymbol{\tau}[y],a),$$
o\`u
$$\boldsymbol{\tau}[y]=desc_{\eta[y]}^{\tilde{M},*}\circ \iota_{M_{\eta[y],sc},M_{\eta[y]}}^*\circ transfert_{y}(\bar{\boldsymbol{\delta}}).$$
L'\'egalit\'e (10) devient
$$(15 )\qquad \rho_{J}^{\tilde{G},{\cal E}}({\bf M}',\boldsymbol{\delta},a)=\rho_{J}^{\tilde{G}}(\boldsymbol{\tau},a),$$
o\`u
$$\boldsymbol{\tau}=\sum_{y\in \dot{{\cal Y}}^M}c^M[y]\boldsymbol{\tau}[y].$$
A ce point, on peut lever l'hypoth\`ese $A_{\tilde{G}}=A_{G_{\eta}}$ et que $J$ provient d'un \'el\'ement de ${\cal J}_{\bar{M}}^{\bar{G}}$. Si elle n'est pas v\'erifi\'ee, le membre de gauche de (15) est nul d'apr\`es (12). Celui de droite l'est aussi d'apr\`es le lemme 4.1. 

Des calculs analogues \`a ceux effectu\'es ci-dessus permettent de d\'eduire de 5.4(3) l'\'egalit\'e
$ transfert(\boldsymbol{\delta})=\boldsymbol{\tau}$. On peut aussi plus simplement appliquer la relation (15) au cas $\tilde{G}=\tilde{M}$ et \`a $J=\emptyset$. Cette relation devient dans ce cas l'\'egalit\'e pr\'ec\'edente. Gr\^ace \`a celle-ci, la relation (15) est l'\'egalit\'e de l'\'enonc\'e. $\square$

\bigskip

\subsection{Les termes $\rho_{J}^{\tilde{G},{\cal E}}$, variante}
On consid\`ere un triplet  $(G,\tilde{G},{\bf a})$  quasi-d\'eploy\'e et \`a torsion int\'erieure, un syst\`eme de fonctions $B$, un espace de Levi $\tilde{M}$  de $\tilde{G}$  et une donn\'ee endoscopique elliptique et relevante ${\bf M}'=(M',{\cal M}',\tilde{\zeta})$   de $(M,\tilde{M})$. 

On suppose donn\'e un diagramme $(\epsilon,B^{M'},T',B^M,T,\eta)$ joignant un \'el\'ement $\epsilon\in \tilde{M}'_{ss}(F)$ \`a un \'el\'ement $\eta\in \tilde{M}_{ss}(F)$. On suppose $M'_{\epsilon}$ quasi-d\'eploy\'e.
 On note ${\cal O}'$ la classe de conjugaison stable de $\epsilon$ dans $\tilde{M}'(F)$ et ${\cal O}$ la classe de conjugaison stable de $\eta$ dans $\tilde{M}(F)$.
 
 \ass{Proposition}{On suppose  que $A_{M'_{\epsilon}}=A_{M'}$. Pour tout $J\in {\cal J}_{\tilde{M}}^{\tilde{G}}(B)$, tout $\boldsymbol{\delta}\in D^{st}_{g\acute{e}om}({\bf M}',{\cal O}')\otimes Mes(M'(F))^*$ et tout $a\in A_{M}(F)$ en position g\'en\'erale et proche de $1$, on a l'\'egalit\'e
$$\rho_{J}^{\tilde{G},{\cal E}}({\bf M}',\boldsymbol{\delta},a)=\rho_{J}^{\tilde{G}}(transfert(\boldsymbol{\delta}),a).$$}

La preuve est identique. On n'a plus besoin de l'hypoth\`ese (2) du paragraphe pr\'ec\'edent car, dans la situation  quasi-d\'eploy\'ee et \`a torsion int\'erieure, les triplets endoscopiques non standard qui apparaissent sont triviaux. Ils v\'erifient \'evidemment la proposition 6.6.
 
\bigskip

\subsection{ Les termes $\sigma_{J}$}
On consid\`ere un triplet  $(G,\tilde{G},{\bf a})$  quasi-d\'eploy\'e et \`a torsion int\'erieure, un syst\`eme de fonctions $B$, un espace de Levi $\tilde{M}$  de $\tilde{G}$. On consid\`ere un \'el\'ement  $\eta\in \tilde{M}(F)$, semi-simple et tel que $M_{\eta}$ soit quasi-d\'eploy\'e.    On note ${\cal O}$ la classe de conjugaison stable  de $\eta$ dans $\tilde{M}(F)$. 
\ass{Proposition}{On suppose $A_{M}=A_{M_{\eta}}$. 

(i) Soient $J\in {\cal J}_{\tilde{M}}^{\tilde{G}}(B_{{\cal O}})$, $\boldsymbol{\delta}'\in D_{unip}^{st}(M_{\eta}(F),\omega)\otimes Mes(M_{\eta}(F))^*$ et $a\in A_{\tilde{M}}(F)$ en position g\'en\'erale et assez proche de $1$. Posons $\boldsymbol{\delta}=desc_{\eta}^{st,\tilde{M},*}(\boldsymbol{\delta}')$.  On a l'\'egalit\'e
$$\sigma_{J}^{\tilde{G}}(\boldsymbol{\delta},a)=\left\lbrace\begin{array}{cc}e_{\tilde{M}}^{\tilde{G}}(\eta)desc_{\eta}^{st,\tilde{M},*}(\sigma_{J}^{G_{\eta}}(\boldsymbol{\delta}',a)),&\text{ si }A_{G}=A_{G_{\eta}} ,\\ 0,&\text{ sinon. }\\ \end{array}\right.$$

(ii) Pour tout $J\in {\cal J}_{\tilde{M}}^{\tilde{G}}(B_{{\cal O}})$, $\sigma^{\tilde{G}}_{J}$ prend ses valeurs dans 
 $$U_{J}\otimes (D_{g\acute{e}om}^{st}({\cal O})\otimes Mes(M(F))^*)/Ann_{{\cal O}}^{\tilde{G},st}.$$}
 
 Preuve. On reprend la preuve de 7.1 dans le cas que l'on avait exclu, \`a savoir ${\bf M}'={\bf M}$. On  prend pour diagramme un diagramme "trivial" $(\eta,B^M,T,B^M,T,\eta)$. La relation 7.1(3) devient
 $$(1) \qquad \rho_{J}^{\tilde{G},{\cal E}}({\bf M},\boldsymbol{\delta},a)=\sum_{s\in Z(\hat{M})^{\Gamma_{F}}/Z(\hat{G})^{\Gamma_{F}}, J\in {\cal J}_{\tilde{M}}^{\tilde{G}'(s)}(B)}i_{\tilde{M}}(\tilde{G},\tilde{G}'(s))\sigma_{J}^{\tilde{G}'(s)}(\boldsymbol{\delta},a).$$
 Le seul point qui diff\`ere de la situation de 7.1 est que l'on ne peut plus utiliser de relation   de descente pour le terme $\sigma_{J}^{\tilde{G}}(\boldsymbol{\delta},a)$ correspondant \`a $s=1$. 
Mais on peut n\'eanmoins appliquer cette relation de descente, \`a condition d'ajouter \`a l'expression obtenue la diff\'erence entre $\sigma_{J}^{\tilde{G}}(\boldsymbol{\delta},a)$ et le terme obtenu par descente. C'est-\`a-dire, posons
$$x=\left\lbrace\begin{array}{cc}\sigma_{J}^{\tilde{G}}(\boldsymbol{\delta},a)-e_{\tilde{M}}^{\tilde{G}}(\eta)desc_{\eta}^{st,\tilde{M},*}(\sigma_{J}^{G_{\eta}}(\boldsymbol{\delta},a)),&\text{ si }A_{G}=A_{G_{\eta}} ,\\ \sigma_{J}^{\tilde{G}}(\boldsymbol{\delta},a),&\text{ sinon. }\\ \end{array}\right.$$
Alors le membre de droite de (1) est la somme de $x$ et d'une  expression qui se calcule comme en 7.1. Rappelons que l'hypoth\`ese (2) de ce paragraphe est automatiquement v\'erifi\'ee dans notre situation quasi-d\'eploy\'ee et \`a torsion int\'erieure. On obtient finalement l'\'egalit\'e
$$\rho_{J}^{\tilde{G},{\cal E}}({\bf M},\boldsymbol{\delta},a)=x+\rho_{J}^{\tilde{G}}({\bf M},\boldsymbol{\delta},a).$$
Mais, pour la donn\'ee endoscopique maximale ${\bf M}$, on a tautologiquement l'\'egalit\'e 
$$\rho_{J}^{\tilde{G},{\cal E}}({\bf M},\boldsymbol{\delta},a)=\rho_{J}^{\tilde{G}}({\bf M},\boldsymbol{\delta},a)$$
car le terme $\sigma_{J}^{\tilde{G}}(\boldsymbol{\delta},a)$ est d\'efini pour qu'il en soit ainsi. D'o\`u $x=0$, ce qu'affirme le (i) de l'\'enonc\'e. 

Le membre de droite de l'\'egalit\'e du  (i) est par d\'efinition une distribution stable. Le (ii) en r\'esulte. $\square$

\bigskip

\subsection{Preuve conditionnelle des propositions [II] 2.7, [II] 3.8 et  du th\'eor\`eme [II] 1.16(i)}

On consid\`ere un triplet  $(G,\tilde{G},{\bf a})$  quelconque et un espace de Levi $\tilde{M}$  de $\tilde{G}$. On consid\`ere l'hypoth\`ese

(1) pour tout $\boldsymbol{\gamma}\in D_{g\acute{e}om}(\tilde{M}(F),\omega)\otimes Mes(M(F))^*$ \`a support form\'e d'\'el\'ements $\tilde{G}$-fortement r\'eguliers et pour tout ${\bf f}\in I(\tilde{G}(F),\omega)\otimes Mes(G(F))$, on a l'\'egalit\'e
$$I^{\tilde{G},{\cal E}}_{\tilde{M}}(\boldsymbol{\gamma},{\bf f})=I^{\tilde{G}}_{\tilde{M}}(\boldsymbol{\gamma},{\bf a}).$$

\ass{Proposition}{On suppose cette hypoth\`ese v\'erifi\'ee.

(i) Soit ${\bf M}'$ une donn\'ee endoscopique elliptique et relevante de $(M,\tilde{M},{\bf a})$ et soit ${\cal O}'$ une classe de conjugaison stable semi-simple dans $\tilde{M}'(F)$ se transf\'erant en une classe de conjugaison stable ${\cal O}$ de $\tilde{M}(F)$. Pour tout $J\in {\cal J}_{\tilde{M}}^{\tilde{G}}$, tout $\boldsymbol{\delta}\in D^{st}_{g\acute{e}om}({\bf M}',{\cal O}')\otimes Mes(M'(F))^*$ et tout $a\in A_{\tilde{M}}(F)$ en position g\'en\'erale et proche de $1$, on a l'\'egalit\'e
$$\rho_{J}^{\tilde{G},{\cal E}}({\bf M}',\boldsymbol{\delta},a)=\rho_{J}^{\tilde{G}}(transfert(\boldsymbol{\delta}),a).$$

(ii) Pour tout $\boldsymbol{\gamma}\in D_{g\acute{e}om}(\tilde{M}(F),\omega)\otimes Mes(M(F))^*$  et pour tout ${\bf f}\in I(\tilde{G}(F),\omega)\otimes Mes(G(F))$, on a l'\'egalit\'e
$$I^{\tilde{G},{\cal E}}_{\tilde{M}}(\boldsymbol{\gamma},{\bf f})=I^{\tilde{G}}_{\tilde{M}}(\boldsymbol{\gamma},{\bf a}).$$

(iii) Pour toute classe de conjugaison stable semi-simple ${\cal O}$ dans $\tilde{M}(F)$, on a l'\'egalit\'e $g_{\tilde{M},{\cal O}}^{\tilde{G},{\cal E}}=g_{\tilde{M},{\cal O}}^{\tilde{G}}$.}

Preuve. Consid\'erons la situation de (i). Fixons un diagramme $(\epsilon,B^{M'},T',B^M,T,\eta)$ reliant un \'el\'ement $\epsilon\in {\cal O}'$ tel que $M'_{\epsilon}$ soit quasi-d\'eploy\'e \`a un \'el\'ement $\eta\in {\cal O}$. Supposons d'abord que $A_{M'_{\epsilon}}\not=A_{M'}$. Introduisons le commutant $\tilde{R}'$ de $A_{M'_{\epsilon}}$ dans $\tilde{M}'$. C'est un espace de Levi propre. Du diagramme se d\'eduit un homomorphisme $\xi:T^{\theta,0}\to T'$ qui est une isog\'enie et est \'equivariant pour les actions galoisiennes. On a $A_{M'_{\epsilon}}\subset T'$. La composante neutre de $\xi^{-1}(A_{M'_{\epsilon}})$ est un tore d\'eploy\'e. On note $\tilde{R}$ son commutant dans $\tilde{M}$. C'est un espace de Levi propre  qui correspond \`a $\tilde{R}'$. En posant $B^{R'}=B^{M'}\cap R'$ et $B^R=B^M\cap R$,  le sextuplet $(\epsilon,B^{R'},T',B^R,T,\eta)$ est encore un diagramme, avec pour espaces ambiants $\tilde{R}$ et $\tilde{R}'$. Le Levi $R'$ se compl\`ete en une donn\'ee endoscopique ${\bf R}'$ de $(R,\tilde{R},{\bf a})$ qui est elliptique et relevante.  On a $M'_{\epsilon}=R'_{\epsilon}$ par construction. L'application $desc^{st,\tilde{M}',*}_{\epsilon}$ est la compos\'ee de $desc_{\epsilon}^{st,\tilde{R}',*}$ et de l'induction de $\tilde{R}'$ \`a $\tilde{M}'$.  Puisque tout \'el\'ement de $D^{st}_{g\acute{e}om}({\cal O}')\otimes Mes(M'(F))^*$ appartient \`a l'image de $desc^{st,\tilde{M}',*}_{\epsilon}$, tout tel \'el\'ement est l'induit d'un \'el\'ement de $D^{st}_{g\acute{e}om}({\cal O}_{\tilde{R}'} )\otimes Mes(R'(F))^*$, o\`u ${\cal O}_{\tilde{R}'}={\cal O}'\cap \tilde{R}'(F)$. Ceci s'adapte formellement aux donn\'ees endoscopiques. Donc $\boldsymbol{\delta}=\boldsymbol{\tau}^{{\bf M}'}$ pour un \'el\'ement $\boldsymbol{\tau}\in D^{st}_{g\acute{e}om}({\bf R}',{\cal O}_{\tilde{R}'})\otimes Mes(R'(F))^*$. On a alors
$$transfert(\boldsymbol{\delta})=(transfert(\boldsymbol{\tau}))^{\tilde{M}}.$$
Les deux membres de l'\'egalit\'e du (i) v\'erifient les formules de descente parall\`eles [II] 3.10 et [II] 3.12. On voit que cette \'egalit\'e du (i) r\'esulte d'\'egalit\'es similaires o\`u $\tilde{G}$ est remplac\'e par des espaces de Levi propres. En vertu de nos hypoth\`eses de r\'ecurrence, ces \'egalit\'es sont v\'erifi\'ees, d'o\`u (i) dans ce cas.

Supposons maintenant $A_{M'_{\epsilon}}=A_{M'}$. Comme on l'a expliqu\'e en 6.4, si $(G,\tilde{G},{\bf a})$ n'est pas l'un des triplets d\'efinis en 6.3, les donn\'ees endoscopiques non standard qui apparaissent en 7.1 v\'erifient toutes les propri\'et\'es requises. Donc l'hypoth\`ese (2) de 7.1 est v\'erifi\'ee. L'assertion (i) r\'esulte alors de cette proposition 7.1. Supposons que $(G,\tilde{G},{\bf a})$ soit l'un des triplets d\'efinis en 6.3. Rappelons que $G$ est simplement connexe. Consid\'erons un triplet endoscopique non standard $(\bar{G}'(\bar{s})_{SC},G'(\tilde{s})_{\epsilon,SC},j_{*})$ comme en 7.1. 
D'apr\`es le lemme 6.3,  on sait dans quels cas les propri\'et\'es de ce triplet ne sont pas connues. Supposons que l'on soit dans un tel  cas. Alors  la donn\'ee ${\bf G}'(\tilde{s})$ est \'equivalente \`a la donn\'ee maximale ${\bf G}'=(G',\hat{G}_{\hat{\theta}}\rtimes W_{F},\hat{\theta})$ d\'efinie en 6.3, $ad_{\eta}$ conserve une paire de Borel \'epingl\'ee de $G$ et $\epsilon$ est l'\'el\'ement de ${\cal Z}(\tilde{G}'(\tilde{s}))^{\Gamma_{F}}$ qui correspond \`a $\eta$ par l'application du lemme 6.2. 
Montrons que cela entra\^{\i}ne

(2) ${\bf M}'$ est \'equivalente \`a la donn\'ee endoscopique maximale de $(M,\tilde{M},{\bf a})$, $ad_{\eta}$ conserve une paire de Borel \'epingl\'ee de $M$ et $\epsilon$ est l'\'element de ${\cal Z}(\tilde{M}')^{\Gamma_{F}}$ qui correspond \`a $\eta$.

Avec les notations habituelles, on peut supposer $\tilde{\zeta}=\zeta\hat{\theta}$, avec $\zeta\in \hat{T}$ et on \'ecrit $\tilde{s}=z\zeta\hat{\theta}$, avec $z\in Z(\hat{M})^{\Gamma_{F},\hat{\theta}}$. Puisque ${\bf G}'(\tilde{s})$ est \'equivalent \`a la donn\'ee maximale ${\bf G}'$, on peut fixer $x\in \hat{G}$ tel que $x{\cal G}'(\tilde{s})x^{-1}=\hat{G}_{\hat{\theta}}\rtimes W_{F}$ et $xz\zeta\hat{\theta}x^{-1}=\hat{\theta}$. En particulier, $ad_{x}$ envoie $\hat{G}'(\tilde{s})$ sur $\hat{G}_{\hat{\theta}}$. Quitte \`a multiplier $x$ \`a gauche par un \'el\'ement de $\hat{G}_{\hat{\theta}}$, on peut supposer que $ad_{x}$ envoie la paire de Borel $(\hat{B}\cap \hat{G}'(\tilde{s}),\hat{T}^{\hat{\theta}})$ de $\hat{G}'(\tilde{s})$ sur la paire de Borel $(\hat{B}\cap \hat{G}_{\hat{\theta}},\hat{T}^{\hat{\theta}})$ de $\hat{G}_{\hat{\theta}}$. Alors $x$ normalise $\hat{T}$ et son  image dans $W$ est fixe par $\hat{\theta}$. Or le groupe $W^{\hat{\theta}}$ est le groupe de Weyl de $\hat{G}_{\hat{\theta}}$. Quitte \`a multiplier encore $x$ \`a gauche par un \'el\'ement de ce groupe, on peut supposer que l'image de $x$ dans $W$ est $1$, autrement dit que $x\in \hat{T}$. Mais alors $x\in \hat{M}$ et la conjugaison par $x$ conserve $^LM$. En prenant les intersections avec ce groupe, la relation $x{\cal G}'(\tilde{s})x^{-1}=\hat{G}_{\hat{\theta}}\rtimes W_{F}$ entra\^{\i}ne que $x{\cal M}'x^{-1}=\hat{M}_{\hat{\theta}}\rtimes W_{F}$. On a aussi $x\zeta\hat{\theta}x^{-1}=s^{-1}\hat{\theta}\in Z(\hat{M})\hat{\theta}$. Donc ${\bf M}'$ est \'equivalente \`a la donn\'ee maximale de $(M,\tilde{M},{\bf a})$. Soit ${\cal E}_{0}=(B_{0},T_{0},(E_{\alpha})_{\alpha\in \Delta})$ une paire de Borel \'epingl\'ee de $G$ qui est conserv\'ee par $ad_{\eta}$. Puisque $\eta$ est un \'el\'ement semi-simple de $\tilde{M}$, on peut aussi fixer une paire de Borel $(B,T)$ de $G$, pour laquelle $M$ est standard et qui est conserv\'ee par $ad_{\eta}$. D'apr\`es la preuve de 6.2(3), il existe $x\in G_{\eta}$ tel que $ad_{x}(B_{0},T_{0})=(B,T)$. Alors ${\cal E}=ad_{x}({\cal E}_{0})$ est une paire de Borel \'epingl\'ee de $G$ qui est conserv\'ee par $ad_{\eta}$. Or $M$ est standard pour cette paire. On peut donc "restreindre" celle-ci \`a $M$ et on obtient  une paire de Borel \'epingl\'ee   ${\cal E}^M$ de $M$ qui est  conserv\'ee par $ad_{\eta}$. De plus, il r\'esulte des constructions que le diagramme suivant est commutatif
$$\begin{array}{ccc}Z(\tilde{G},{\cal E})&\to&{\cal Z}(\tilde{G}')\\ \downarrow&&\downarrow\\ Z(\tilde{M},{\cal E}^M)&\to&{\cal Z}(\tilde{M}')\\ \end{array}$$
La derni\`ere assertion de (2) en r\'esulte. Cela prouve cette assertion.

 On a alors $G'(\tilde{s})_{\epsilon}\simeq G'$ et $M'_{\epsilon}=M'$ et la classe de conjugaison stable de $\epsilon$ est r\'eduite \`a ce point. 
  On  v\'erifie facilement que  l'\'el\'ement maximal de ${\cal J}_{\tilde{M}'}^{\tilde{G}'(\tilde{s})}(B^{\tilde{G}}_{ \{\epsilon\}})$ s'envoie sur l'\'el\'ement maximal de ${\cal J}_{\tilde{M}}^{\tilde{G}}$. D'apr\`es la propri\'et\'e 6.6(2), l'assertion de la proposition 6.6 est connue pour un \'el\'ement non maximal de ${\cal J}_{\tilde{M}'}^{\tilde{G}'(\tilde{s})}(B^{\tilde{G}}_{ \{\epsilon\}})$. Il en r\'esulte que l'hypoth\`ese 7.1(2) est v\'erifi\'ee si $J$ n'est pas l'\'el\'ement maximal de ${\cal J}_{\tilde{M}}^{\tilde{G}}$. La proposition 7.1 nous dit donc que l'\'egalit\'e du (i) est v\'erifi\'ee sauf pour cet \'el\'ement maximal. On abandonne notre $J$ initial et on note $J_{max}$ l'\'el\'ement maximal. 
Soit ${\bf f}\in I(\tilde{G}(F),\omega)\otimes Mes(G(F))$ et $a\in A_{\tilde{M}}(F)$ en position g\'en\'erale. Posons $\boldsymbol{\gamma}=transfert(\boldsymbol{\delta})$. Les propositions [II] 3.2 et [II] 3.9 entra\^{\i}nent que le germe en $1$ de la fonction
$$(3) \qquad a\mapsto I_{\tilde{M}}^{\tilde{G}}(a\boldsymbol{\gamma},{\bf f})- I_{\tilde{M}}^{\tilde{G},{\cal E}}(a\boldsymbol{\gamma},{\bf f})= I_{\tilde{M}}^{\tilde{G}}(a\boldsymbol{\gamma},{\bf f})-I_{\tilde{M}}^{\tilde{G},{\cal E}}({\bf M}',\xi(a)\boldsymbol{\delta},{\bf f})$$
est \'equivalent \`a
$$\sum_{\tilde{L}\in {\cal L}(\tilde{M})}\sum_{J\in {\cal J}_{\tilde{M}}^{\tilde{L}}}I_{\tilde{L}}^{\tilde{G}}(\rho_{J}^{\tilde{L}}(\boldsymbol{\gamma},a)^{\tilde{L}},{\bf f})-I_{\tilde{L}}^{\tilde{G},{\cal E}}(\rho_{J}^{\tilde{L},{\cal E}}({\bf M}',\boldsymbol{\delta},a)^{\tilde{L}},{\bf f}).$$
En vertu de nos hypoth\`eses de r\'ecurrence, les termes index\'es par $\tilde{L}\not=\tilde{M}$ et $\tilde{L}\not=\tilde{G}$ s'annulent. Ceux index\'es par $\tilde{L}=\tilde{G}$ et un $J\not=J_{max}$ s'annulent aux-aussi. L'expression ci-dessus se r\'eduit \`a
$$(4) \qquad I^{\tilde{G}}(\rho_{J_{max}}^{\tilde{G}}(\boldsymbol{\gamma},a)^{\tilde{G}}-\rho_{J_{max}}^{\tilde{G},{\cal E}}({\bf M}',\boldsymbol{\delta},a)^{\tilde{G}},{\bf f})$$
$$+I_{\tilde{M}}^{\tilde{G}}(\boldsymbol{\gamma},{\bf f})-I_{\tilde{M}}^{\tilde{G},{\cal E}}({\bf M}',\boldsymbol{\delta},{\bf f}).$$
D'apr\`es la proposition [II] 2.10, l'hypoth\`ese (1) entra\^{\i}ne la m\^eme \'egalit\'e que dans cette hypoth\`ese pour les distributions \`a support $\tilde{G}$-\'equisingulier. C'est le cas de la distribution $a\boldsymbol{\gamma}$. Donc l'expression (3) est nulle. Consid\'erons (4). Comme fonction de $a$, le premier terme appartient \`a $U_{J_{max}}$ et le second est constant. Leur somme est \'equivalente \`a $0$. En utilisant [II] 3.1(3), les deux termes sont nuls (en supposant $\tilde{M}\not=\tilde{G}$; si $\tilde{M}=\tilde{G}$, l'assertion \`a prouver est tautologique). La nullit\'e du premier pour tout ${\bf f}$ signifie que 
$$\rho_{J_{max}}^{\tilde{G}}(\boldsymbol{\gamma},a)-\rho_{J_{max}}^{\tilde{G},{\cal E}}({\bf M}',\boldsymbol{\delta},a)$$
modulo $Ann^{\tilde{G}}$. Cela ach\`eve la preuve de (i). 

Prouvons (ii). Par lin\'earit\'e, on peut supposer qu'il existe   une donn\'ee endoscopique  ${\bf M}'$ de $(M,\tilde{M},{\bf a})$ elliptique et relevante et un \'el\'ement $\boldsymbol{\delta}\in D_{g\acute{e}om}^{st}({\bf M}')\otimes Mes(M'(F))^*$ de sorte que $\boldsymbol{\gamma}=transfert(\boldsymbol{\delta})$. Toujours par lin\'earit\'e, on peut fixer une classe de conjugaison stable semi-simple ${\cal O}'$ dans  $\tilde{M}'(F)$ telle que $\boldsymbol{\delta}\in D^{st}_{g\acute{e}om}({\bf M}',{\cal O}')\otimes Mes(M'(F))^*$. On peut supposer que cette classe se transf\`ere en une classe de $\tilde{M}(F)$, sinon $\boldsymbol{\gamma}=0$ et l'\'egalit\'e \`a prouver est triviale. On reprend alors le raisonnement ci-dessus. Pour ${\bf f}\in I(\tilde{G}(F),\omega)\otimes Mes(G(F))$ et $a\in A_{\tilde{M}}(F)$ en position g\'en\'erale, on calcule le d\'eveloppement de (3). Maintenant que l'on a prouv\'e (i), ce d\'eveloppement se r\'eduit \`a
$$I_{\tilde{M}}^{\tilde{G}}(\boldsymbol{\gamma},{\bf f})-I_{\tilde{M}}^{\tilde{G},{\cal E}}({\bf M}',\boldsymbol{\delta},{\bf f}),$$
ou encore \`a
$$I_{\tilde{M}}^{\tilde{G}}(\boldsymbol{\gamma},{\bf f})-I_{\tilde{M}}^{\tilde{G},{\cal E}}(\boldsymbol{\gamma},{\bf f}).$$
Comme ci-dessus, l'hypoth\`ese (1) entra\^{\i}ne que (3) est nul donc aussi cette diff\'erence. C'est la conclusion de (ii). 

Comme on l'a vu ci-dessus, l'hypoth\`ese (1) entraine que l'hypoth\`ese du lemme [II] 2.8 est v\'erifi\'ee pour l'ensemble ${\cal D}=D_{g\acute{e}om,\tilde{G}-\acute{e}qui}(\tilde{M}(F),\omega)\otimes Mes(M(F))^*$. Le (i) de ce lemme aussi d'apr\`es  l'assertion (ii) du pr\'esent \'enonc\'e. Donc aussi le (ii) de ce lemme, qui n'est autre que la pr\'esente assertion (iii).  $\square$

\bigskip

\subsection{Preuve du th\'eor\`eme [II] 1.16(ii)}
On consid\`ere un triplet  $(G,\tilde{G},{\bf a})$  quasi-d\'eploy\'e et \`a torsion int\'erieure, un syst\`eme de fonctions $B$, un espace de Levi $\tilde{M}$  de $\tilde{G}$.

\ass{Proposition}{(i) Soit ${\bf M}'$ une donn\'ee endoscopique elliptique et relevante de $(M,\tilde{M})$ et soit ${\cal O}'$ une classe de conjugaison stable semi-simple dans $\tilde{M}'(F)$ se transf\'erant en une classe de conjugaison stable ${\cal O}$ de $\tilde{M}(F)$. Pour tout $J\in {\cal J}_{\tilde{M}}^{\tilde{G}}(B_{{\cal O}})$, tout $\boldsymbol{\delta}\in D^{st}_{g\acute{e}om}({\bf M}',{\cal O}')\otimes Mes(M'(F))^*$ et tout $a\in A_{M}(F)$ en position g\'en\'erale et proche de $1$, on a l'\'egalit\'e
$$\rho_{J}^{\tilde{G},{\cal E}}({\bf M}',\boldsymbol{\delta},a)=\rho_{J}^{\tilde{G}}(transfert(\boldsymbol{\delta}),a).$$

(ii) Pour tout $\boldsymbol{\gamma}\in D_{g\acute{e}om}(\tilde{M}(F),\omega)\otimes Mes(M(F))^*$  et pour tout ${\bf f}\in I(\tilde{G}(F),\omega)\otimes Mes(G(F))$, on a l'\'egalit\'e
$$I^{\tilde{G},{\cal E}}_{\tilde{M}}(\boldsymbol{\gamma},B,{\bf f})=I^{\tilde{G}}_{\tilde{M}}(\boldsymbol{\gamma},B,{\bf a}).$$}

La preuve est identique \`a celle de la proposition pr\'ec\'edente, en plus simple puisque notre triplet ne saurait \^etre l'un de ceux d\'efinis en 6.3. On n'a plus besoin de l'hypoth\`ese (1) de ce paragraphe: elle est v\'erifi\'ee d'apr\`es la proposition 2.9.

\bigskip

\subsection{Preuve des propositions [II] 2.4,  [II] 3.5 et du th\'eor\`eme [II] 1.10}
On consid\`ere un triplet  $(G,\tilde{G},{\bf a})$  quasi-d\'eploy\'e et \`a torsion int\'erieure, un syst\`eme de fonctions $B$, un espace de Levi $\tilde{M}$  de $\tilde{G}$.
\ass{Proposition}{(i) Soit ${\cal O}$ une classe de conjugaison stable semi-simple dans $\tilde{M}(F)$. Pour tout $J\in {\cal J}_{\tilde{M}}^{\tilde{G}}(B_{{\cal O}})$, le terme $\sigma_{J}^{\tilde{G}}$ prend ses valeurs dans
$$U_{J}\otimes (D^{st}_{g\acute{e}om}({\cal O})\otimes Mes(M(F))^*)/Ann_{{\cal O}}^{st,\tilde{G}}.$$

(ii) Pour tout $\boldsymbol{\delta}\in D^{st}_{g\acute{e}om}(\tilde{M}(F))\otimes Mes(M(F))^*$, la distribution
$${\bf f}\mapsto S_{\tilde{M}}^{\tilde{G}}(\boldsymbol{\delta},B,{\bf f})$$
est stable.

(iii) Soit ${\cal O}$ une classe de conjugaison stable semi-simple dans $\tilde{M}(F)$, notons ${\cal O}^{\tilde{G}}$ la classe de conjugaison stable dans $\tilde{G}(F)$ qui contient ${\cal O}$. Alors le germe $Sg_{\tilde{M},{\cal O}}^{\tilde{G}}(.,B)$ prend ses valeurs dans l'espace $D_{g\acute{e}om}^{st}({\cal O}^{\tilde{G}})\otimes Mes(G(F))^*$.}

Preuve. Pour la preuve de (i), on utilise le m\^eme argument qu'en 7.4. On sait comment se comportent nos termes par induction, gr\^ace \`a la proposition [II] 3.11. On se ram\`ene alors au cas o\`u la proposition 7.3 s'applique.

Prouvons (ii). Par lin\'earit\'e, on peut fixer une classe de conjugaison stable semi-simple  ${\cal O}$ dans $\tilde{M}(F)$ et supposer que $\boldsymbol{\delta}\in D^{st}_{g\acute{e}om}({\cal O})\otimes Mes(M(F))^*$. Soit ${\bf f}\in I(\tilde{G}(F))\otimes Mes(G(F))$. Supposons que les int\'egrales orbitales stables fortement  r\'eguli\`eres de ${\bf f}$ sont nulles, autrement dit que l'image de ${\bf f}$ dans $SI(\tilde{G}(F))\otimes Mes(G(F))$ est nulle. Soit $a\in A_{M}(F)$ en position g\'en\'erale. Remarquons que le (i) d\'ej\`a prouv\'e assure la validit\'e de la proposition [II] 3.5. La proposition 3.7 (ii) calcule le germe en $1$ de la fonction
$$(1) \qquad a\mapsto S_{\tilde{M}}^{\tilde{G}}(a\boldsymbol{\delta},{\bf f}).$$
Il est \'equivalent \`a
$$\sum_{\tilde{L}\in {\cal L}(\tilde{M})}\sum_{J\in {\cal J}_{\tilde{M}}^{\tilde{L}}(B_{{\cal O}})}S_{\tilde{L}}^{\tilde{G}}(\sigma_{J}^{\tilde{L}}(\boldsymbol{\delta},a)^{\tilde{L}},B,{\bf f}).$$
L'hypoth\`ese sur ${\bf f}$ et nos hypoth\`eses de recurrence assurent que tous les termes sont nuls sauf celui index\'e par $\tilde{M}$. L'expression ci-dessus se r\'eduit \`a $S_{\tilde{M}}^{\tilde{G}}(\boldsymbol{\delta},B,{\bf f})$.  Or (1) est nul d'apr\`es la proposition 2.8. Donc $S_{\tilde{M}}^{\tilde{G}}(\boldsymbol{\delta},B,{\bf f})=0$. Cette \'egalit\'e pour tout ${\bf f}$ d'image nulle dans $SI(\tilde{G}(F))\otimes Mes(G(F))$ est \'equivalente \`a l'assertion (ii). 

Le (iii) s'en d\'eduit comme en 7.4  en utilisant le lemme [II] 2.9.
$\square$

\bigskip

\subsection{Preuve de la proposition 6.6}
On se place sous les hypoth\`eses de cette proposition, dont on utilise les notations.  D'apr\`es le lemme 6.8, on peut supposer $(G_{1},G_{2},j_{*})$ quasi-\'el\'ementaire et $B_{1}$ constante de valeur $1$. On introduit le triplet $(G,\tilde{G},{\bf a})$ associ\'e \`a $(G_{1},G_{2},j_{*})$ comme en 6.3.  On fixe un \'el\'ement $\eta\in \tilde{G}(F)$ qui conserve une paire de Borel \'epingl\'ee ${\cal E}=(B,T,(E_{\alpha})_{\alpha\in \Delta})$ de $G$ d\'efinie sur $F$. De cette paire se d\'eduit une paire de Borel \'epingl\'ee de $G_{\eta}$ d\'efinie sur $F$. On peut identifier $G_{1}$ \`a $G_{\eta}$ de sorte que le Levi $M_{1}$ devienne un Levi de $G_{\eta}$ standard pour cette paire de Borel \'epingl\'ee. On note $\tilde{M}$ le commutant dans $\tilde{G}$ du tore $A_{M_{1}}$. On a $\eta\in \tilde{M}(F)$ et $M$ est standard pour  ${\cal E}$. 
On introduit les donn\'ees endoscopiques maximales   ${\bf G}'=(G',\hat{G}_{\hat{\theta}}\rtimes W_{F},\hat{\theta})$ de $(G,\tilde{G},{\bf a})$ et  ${\bf M}'=(M',\hat{M}_{\hat{\theta}}\rtimes W_{F},\hat{\theta})$ de $(M,\tilde{M},{\bf a})$.   Remarquons que ${\bf G}'$ est aussi  la donn\'ee ${\bf G}'(\hat{\theta})$  d\'eduite de ${\bf M}'$ et de l'\'el\'ement $\tilde{s}=\hat{\theta}$. 
Comme en 6.3, l'\'el\'ement $\eta\in \tilde{G}(F)$ d\'etermine un \'el\'ement $\epsilon\in {\cal Z}(\tilde{G}')^{\Gamma_{F}}$. Si l'on remplace les espaces ambiants $\tilde{G}$ et $\tilde{G}'$ par $\tilde{M}$ et $\tilde{M}'$, on obtient \'evidemment le m\^eme \'el\'ement $\epsilon$. 
On  fixe un diagramme $(\epsilon,B^{M'},T',B^M,T,\eta)$, o\`u $B^M=B\cap M$. Reprenons les constructions et notations de 7.1. En particulier, on fixe des mesures pour simplifier. Remarquons que $\bar{G}=G_{\eta}=G_{1}$, $\bar{M}=M_{\eta}=M_{1}$, $G_{2}=G'_{\epsilon}=G'$ et $M_{2}=M'_{\epsilon}=M'$. 

 On dispose d'\'el\'ements $\boldsymbol{\delta}_{1}$ et $\boldsymbol{\delta}_{2}$. On peut identifier $\boldsymbol{\delta}_{2}$ \`a un \'el\'ement de $D_{unip}^{st}(M'_{\epsilon}(F))$. On pose
$$\boldsymbol{\delta}=desc_{\epsilon}^{st,M',*} (\boldsymbol{\delta}_{2}).$$
La fonction $B_{2}$ s'identifie \`a $B^{\tilde{G}}_{{\cal O}'}$. 
L'\'el\'ement $J$ de l'\'enonc\'e de la proposition 6.6, vu comme un \'el\'ement de ${\cal J}_{M_{2}}^{G_{2}}(B_{2})$, s'identifie \`a un \'el\'ement de ${\cal J}_{M'_{\epsilon}}^{G'_{\epsilon}}(B^{\tilde{G}}_{{\cal O}'})$, qui s'envoie sur un \'el\'ement de ${\cal J}_{\tilde{M}}^{\tilde{G}}$ que l'on note encore $J$. On reprend la preuve de 7.1 pour ces \'el\'ements $J$ et $\boldsymbol{\delta}$. Remarquons que  l'homomorphisme
    $$ Z(\hat{M})^{\Gamma_{F},\hat{\theta}}/Z(\hat{G})^{\Gamma_{F},\hat{\theta}}\to Z(\hat{\bar{M}})^{\Gamma_{F}}/ Z(\hat{\bar{G}})^{\Gamma_{F}},$$
  se simplifie en
    $$(1) \qquad Z(\hat{M})^{\Gamma_{F},\hat{\theta}}\to Z(\hat{\bar{M}})^{\Gamma_{F}}$$
 puisque $G$ et $\bar{G}$ sont simplement connexes, donc leurs duaux sont adjoints. 
La preuve  marche jusqu'au point o\`u on avait utilis\'e l'hypoth\`ese (2) du paragraphe 7.1. 
Pour un $\tilde{s}\in \hat{\theta} Z(\hat{M})^{\Gamma_{F},\hat{\theta}} $ tel que ${\bf G}'(\tilde{s})$ n'est pas \'equivalent \`a ${\bf G}'$, le lemme 6.3 et nos hypoth\`eses de r\'ecurrence assurent que cette hypoth\`ese est v\'erifi\'ee. Il reste les $\tilde{s}$ tels que ${\bf G}'(\tilde{s})$ est \'equivalent \`a ${\bf G}'$.
 Notons ${\cal Z}$ le   noyau de l'homomorphisme (1), c'est-\`a-dire ${\cal Z}=Z(\hat{M})^{\Gamma_{F},\hat{\theta}} \cap (1-\hat{\theta})(\hat{T})$.  Montrons que
    
    (2) l'ensemble des $s\in Z(\hat{M})^{\Gamma_{F},\hat{\theta}} $ tels que ${\bf G}'(s\hat{\theta})$ est \'equivalente \`a ${\bf G}'$ est \'egal \`a ${\cal Z}$.
    
    Supposons ${\bf G}'(s\hat{\theta})$ \'equivalente \`a ${\bf G}'$. Il existe alors $x\in \hat{G}$ tel que $x{\cal G}'(s\hat{\theta})x^{-1}=G_{\hat{\theta}}\rtimes W_{F}$ et $xs\hat{\theta}x^{-1}=\hat{\theta}$. Le m\^eme argument que dans la preuve de 7.4(2) montre que l'on peut supposer $x\in \hat{T}$. Alors $s=(\hat{\theta}-1)(x)\in (1-\hat{\theta})(\hat{T})$. Donc $s\in {\cal Z}$. Inversement, supposons $s\in {\cal Z}$. Ecrivons $s=(\hat{\theta}-1)(x)$, avec $x\in \hat{T}$.  On a  $xs\hat{\theta}x^{-1}=\hat{\theta}$. Cela entra\^{\i}ne $x\hat{G}'(s\hat{\theta})x^{-1}=\hat{G}_{\hat{\theta}}$. Pour $g\in \hat{G}'(s\hat{\theta})$ et $w\in W_{F}$, on a $xgw(x)^{-1}=xgx^{-1}xw(x)^{-1}$. Le premier terme $xgx^{-1}$ appartient \`a $\hat{G}_{\hat{\theta}}$. L'\'egalit\'e $s=(\hat{\theta}-1)(x)$ et le fait que $s$ est invariant par $\Gamma_{F}$ entra\^{\i}ne que $xw(x)^{-1}$ appartient \`a $\hat{T}^{\hat{\theta}}$, qui est contenu dans $\hat{G}_{\hat{\theta}}$. Donc $xgw(x)^{-1}\in \hat{G}_{\hat{\theta}}$. Cela prouve que $x{\cal G}'(s\hat{\theta})x^{-1}=G_{\hat{\theta}}\rtimes W_{F}$, donc 
 ${\bf G}'(s\hat{\theta})$ est \'equivalente \`a ${\bf G}'$. D'o\`u (2).

     Consid\'erons  un $\tilde{s}$ tel que ${\bf G}'(\tilde{s})$ soit \'equivalent \`a ${\bf G}'$. D'apr\`es (2), cela \'equivaut \`a  $\bar{s}=1$.  Alors   $G'(\tilde{s})_{\epsilon}$ est isomorphe \`a $G_{2}$, $M'(\tilde{s})_{\epsilon}$ est isomorphe \`a $M_{2}$, $\bar{G}'(\bar{s})$ est isomorphe \`a $G_{1}$ et $\bar{M}'$ est isomorphe \`a $M_{1}$.   L'\'el\'ement not\'e $\boldsymbol{\delta}(\tilde{s})_{\epsilon,sc}$ en 7.1 n'est autre que $\boldsymbol{\delta}_{2}$.  L'\'el\'ement $\bar{\boldsymbol{\delta}}(\bar{s})_{sc}$ est \'egal \`a $\boldsymbol{\delta}_{1}$. L'\'el\'ement $\boldsymbol{\tau}(\tilde{s})_{sc}$ est \'egal \`a $\sigma_{J}^{G_{2}}(\boldsymbol{\delta}_{2},a)$. On ne peut plus affirmer que son transfert $\bar{\boldsymbol{\tau}}(\bar{s})_{sc}$ est \'egal \`a $(c_{M_{1},M_{2}}^{G_{1},G_{2}})^{-1}\sigma_{J}^{G_{1}}(\boldsymbol{\delta}_{1},a)$.  Mais on peut \'ecrire ce transfert sous la forme
$$(c_{M_{1},M_{2}}^{G_{1},G_{2}})^{-1}(\bar{\boldsymbol{\mu}}+\sigma_{J}^{G_{1}}(\boldsymbol{\delta}_{1},a)),$$
o\`u 
$$(3) \qquad \bar{\boldsymbol{\mu}}=c_{M_{1},M_{2}}^{G_{1},G_{2}}transfert(\sigma_{J}^{G_{2}}(\boldsymbol{\delta}_{2},a))- \sigma_{J}^{G_{1}}(\boldsymbol{\delta}_{1},a).$$
On peut alors poursuivre le calcul comme en 7.1. Il appara\^{\i}t des termes suppl\'ementaires provenant de $\boldsymbol{\mu}$. 
On obtient une \'egalit\'e similaire \`a 7.1(10):
$$(4) \qquad  \rho_{J}^{\tilde{G},{\cal E}}({\bf M}',\boldsymbol{\delta},a)=\boldsymbol{\mu}+\sum_{y\in \dot{{\cal Y}}^M}c^M[y]desc_{\eta[y]}^{\tilde{M},*} \circ transfert_{y} (\boldsymbol{\xi}),$$
o\`u
$$\boldsymbol{\mu}=x(1)\sum_{y\in \dot{{\cal Y}}^M}c^M[y]desc_{\eta[y]}^{\tilde{M},*}\circ \circ transfert_{y}  (\bar{\boldsymbol{\mu}}).$$
Les applications $\iota^*$ de 7.1(10) disparaissent ici car les groupes $G_{\eta}$ et $G'_{\epsilon}$ sont simplement connexes. 
Le terme $x(1)$ est le $x(\bar{s})$ de 7.1 pour $\bar{s}=1$. Le calcul de la somme du membre de droite de (4) se poursuit comme en 7.1. Cette somme vaut $\rho_{J}^{\tilde{G}}(transfert(\boldsymbol{\delta}),a)$. Mais, d'apr\`es nos hypoth\`eses de r\'ecurrence, toutes les propri\'et\'es sont connues pour le triplet $(G,\tilde{G},{\bf a})$. Donc 
$$\rho_{J}^{\tilde{G},{\cal E}}({\bf M}',\boldsymbol{\delta},a)=\rho_{J}^{\tilde{G}}(transfert(\boldsymbol{\delta}),a).$$
Il en r\'esulte que $\boldsymbol{\mu}=0$. 

Soit $\bar{\varphi}\in SI(\mathfrak{m}_{1}(F))$ \`a support proche de $0$. On l'identifie par l'exponentielle \`a une fonction sur  $M_{1}(F)=\bar{M}'(F)$.  On peut supposer que l'ensemble $\dot{{\cal Y}}^M$ contient l'\'el\'ement $y=1$. Pour celui-ci, $M_{\eta[1]}=M_{\eta}$ est quasi-d\'eploy\'e et on a $\bar{M}'=M_{\eta[1]}$. Modulo cette identification, le transfert $transfert_{1}$  est l'identit\'e.  On peut donc consid\'erer $\bar{\varphi}$ comme un \'el\'ement de $SI(M_{\eta[1]}(F))$, que l'on rel\`eve  en un \'el\'ement $\varphi_{1}\in I(M_{\eta[1]}(F))$. On peut \'evidemment supposer que $\varphi_{1}$ est \`a support proche de l'origine. L'application $desc_{\eta[1]}^{\tilde{M}}$ a pour image le sous-espace des \'el\'ements de $I(M_{\eta[1]}(F))$ qui sont invariants par l'action de $Z_{M}(\eta[1];F)$. Mais, parce que  $\eta[1]=\eta$ conserve une paire de Borel \'epingl\'ee de $M$ et que $T^{\theta}$ est connexe, on a $Z_{M}(\eta[1])=M_{\eta[1]}$. L'application de descente est donc surjective et on peut relever $\varphi_{1}$ en un \'el\'ement $\varphi\in I(\tilde{M}(F))$. 
Toujours parce que  $Z_{M}(\eta)$ est connexe, l'ensemble $\dot{{\cal Y}}^M$ est  un ensemble de repr\'esentants des classes de conjugaison par $M(F)$ dans la classe de conjugaison stable de $\eta$.   On peut modifier $\varphi$   de sorte que $desc_{\eta[y]}^{\tilde{M}}(\varphi)=0$ pour tout $y\in \dot{{\cal Y}}^M$ tel que $y\not=1$. Il r\'esulte alors de la d\'efinition de $\varphi$ que l'on a l'\'egalit\'e
$$I^{\tilde{M}}(\boldsymbol{\mu},\varphi)=S^{M_{1}}(\bar{\boldsymbol{\mu}},\bar{\varphi}).$$
Puisque $\boldsymbol{\mu}=0$, ceci est nul. Puisque cela est vrai pour tout $\bar{\varphi}$, on conclut $\bar{\boldsymbol{\mu}}=0$. D'apr\`es la d\'efinition (3), cette nullit\'e est l'assertion de la proposition 6.6. $\square$

    {\bf Attention.} On ne doit pas s'abuser: la d\'emonstration ci-dessus s'appuie sur des hypoth\`eses de r\'ecurrence. Elle ne deviendra une v\'eritable d\'emonstration que quand toutes les \'etapes de la r\'ecurrence auront \'et\'e \'etablies. 
    
    \bigskip

\section{Descente des germes de Shalika endoscopiques}

\bigskip

\subsection{La proposition [II] 2.7 dans un cas particulier}
On consid\`ere un triplet  $(G,\tilde{G},{\bf a})$  quelconque, un espace de Levi $\tilde{M}$  de $\tilde{G}$  et une donn\'ee endoscopique elliptique et relevante ${\bf M}'=(M',{\cal M}',\tilde{\zeta})$   de $(M,\tilde{M},{\bf a})$. 

On suppose donn\'e un diagramme $(\epsilon,B^{M'},T',B^M,T,\eta)$ joignant un \'el\'ement $\epsilon\in \tilde{M}'_{ss}(F)$ \`a un \'el\'ement $\eta\in \tilde{M}_{ss}(F)$. On suppose que $M'_{\epsilon}$ est quasi-d\'eploy\'e et que  $A_{M'_{\epsilon}}=A_{M'}$.
 On note ${\cal O}'$ la classe de conjugaison stable de $\epsilon$ dans $\tilde{M}'(F)$ et ${\cal O}$ la classe de conjugaison stable de $\eta$ dans $\tilde{M}(F)$. Comme on l'a vu, tout \'el\'ement $\tilde{s}\in \tilde{\zeta}Z(\hat{M})^{\Gamma_{F},\hat{\theta}}/Z(\hat{G})^{\Gamma_{F},\hat{\theta}}$ donne naissance \`a un triplet endoscopique non standard $(\bar{G}'(\bar{s})_{SC},G'(\tilde{s})_{\epsilon,SC},j_{*})$. Du syst\`eme de fonctions $B^{\tilde{G}}$ sur $\tilde{G}'(\tilde{s})$ se  d\'eduit une fonction $B^{\tilde{G}}_{{\cal O}'}$ sur le syst\`eme de racines de $G'(\tilde{s})_{\epsilon,SC}$, puis, par la construction de 6.5,  une fonction sur le syst\`eme de racines de $\bar{G}'(\bar{s})_{SC}$. C'est la fonction constante de valeur $1$ d'apr\`es 7.1(1). Les groupes $G'(\tilde{s})_{\epsilon,SC}$ et $\bar{G}'(\bar{s})_{SC}$ contiennent des Levi $M'(\tilde{s})_{\epsilon,sc}$ et $\bar{M}'(\bar{s})_{sc}$. On consid\`ere l'hypoth\`ese
 
 (1) pour chaque triplet $(\bar{G}'(\bar{s})_{SC},G'(\tilde{s})_{\epsilon,SC},j_{*})$ comme ci-dessus tel que $A_{G'(\tilde{s})_{\epsilon}}=A_{\tilde{G}}$, la proposition 6.7 est v\'erifi\'ee pour ces Levi.
 
 \ass{Proposition}{On suppose que $A_{M'_{\epsilon}}=A_{M'}$ et que l'hypoth\`ese (1) est v\'erifi\'ee. Soit $\boldsymbol{\delta}\in D^{st}_{g\acute{e}om,\tilde{G}-\acute{e}qui}({\bf M}')\otimes Mes(M'(F))^*$, on a l'\'egalit\'e
 $$g_{\tilde{M},{\cal O}}^{\tilde{G},{\cal E}}({\bf M}',\boldsymbol{\delta})=g_{\tilde{M},{\cal O}}^{\tilde{G}}(transfert(\boldsymbol{\delta})),$$
 pourvu que $\boldsymbol{\delta}$ soit assez proche de ${\cal O}'$.}
 
 La preuve occupe les trois paragraphes suivants.
 
 \bigskip
 
 \subsection{D\'ebut de la preuve}
 
  Si $(G,\tilde{G},{\bf a})$ est quasi-d\'eploy\'e et \`a torsion int\'erieure et si ${\bf M}'={\bf M}$, l'\'enonc\'e est tautologique: le terme $Sg_{\tilde{M},{\cal O}}^{\tilde{G}}(\boldsymbol{\delta})$ est d\'efini pour qu'il en soit ainsi. On exclut ce cas.
 
 Rappelons la d\'efinition [II] 2.6(2):
 $$(1) \qquad g_{\tilde{M},{\cal O}}^{\tilde{G},{\cal E}}({\bf M}',\boldsymbol{\delta})=\sum_{\tilde{s}\in \tilde{\zeta}Z(\hat{M})^{\Gamma_{F},\hat{\theta}}/Z(\hat{G})^{\Gamma_{F},\hat{\theta}}}i_{\tilde{M}'}(\tilde{G},\tilde{G}'(\tilde{s}))transfert(Sg_{{\bf M}',{\cal O}'}^{{\bf G}'(\tilde{s})}(\boldsymbol{\delta},B^{\tilde{G}})).$$
 On reprend les constructions et notations de la section 5. Apr\`es avoir fix\'e des donn\'ees auxiliaires $M'_{1}$,...,$\Delta_{1}$, on identifie $\boldsymbol{\delta}$ \`a un \'el\'ement $\boldsymbol{\delta}_{1}\in D^{st}_{g\acute{e}om}(\tilde{M}'_{1}(F))$. Utilisons la description de 5.5. Par lin\'earit\'e, on peut supposer qu'il existe $Z\in \mathfrak{z}(M'_{\epsilon};F)$ et $\boldsymbol{\delta}_{\epsilon,SC}\in D^{st}_{g\acute{e}om}(M'_{\epsilon,SC}(F))$ tels que 
 $$\boldsymbol{\delta}_{1}=desc_{\epsilon_{1}}^{st,\tilde{M}'_{1},*}(exp(Z)\iota^*_{M'_{\epsilon,SC},M'_{1,\epsilon_{1}}}(\boldsymbol{\delta}_{\epsilon,SC})).$$
 Soit $\tilde{s}\in \tilde{\zeta}Z(\hat{M})^{\Gamma_{F},\hat{\theta}}/Z(\hat{G})^{\Gamma_{F},\hat{\theta}}$. D'apr\`es 5.5, l'\'el\'ement $\boldsymbol{\delta}$ s'identifie aussi \`a l'\'el\'ement
 $$\boldsymbol{\delta}_{1}(\tilde{s})=d(\tilde{s})desc_{\epsilon_{1}(\tilde{s})}^{st,\tilde{M}'_{1}(\tilde{s}),*}(exp(Z) \iota^*_{M'_{\epsilon,SC},M'_{1}(\tilde{s})_{\epsilon_{1}(\tilde{s})}}(\boldsymbol{\delta}_{\epsilon,SC}))\in D_{g\acute{e}om}^{st}(\tilde{M}'_{1}(\tilde{s};F))).$$
 A l'aide de 5.3(4), on d\'ecompose $Z$ en $Z_{1}+Z_{2}(\bar{s})+Z_{3}(\bar{s})$, o\`u $Z_{1}\in \mathfrak{z}(\bar{G};F)$, $Z_{2}(\bar{s})\in \mathfrak{z}(\bar{G}'(\bar{s});F)$, $Z_{3}(\bar{s})\in \mathfrak{z}(\bar{M}'(\bar{s})_{sc};F)\simeq \mathfrak{z}(\tilde{M}'(\tilde{s})_{sc};F)$. Notons que 
 $Z_{1}+Z_{2}(\bar{s})\in \mathfrak{z}(\tilde{G}'(\tilde{s})_{\epsilon})$. On a alors
 $$\boldsymbol{\delta}_{1}(\tilde{s})=d(\tilde{s})desc_{\epsilon_{1}(\tilde{s})}^{st,\tilde{M}'_{1}(\tilde{s}),*}(exp(Z_{1}+Z_{2}(\bar{s})) \iota^*_{M'(\tilde{s})_{\epsilon,sc},M'_{1}(\tilde{s})_{\epsilon_{1}(\tilde{s})}}(\boldsymbol{\delta}(\tilde{s})_{\epsilon,sc})),$$
 o\`u
 $$\boldsymbol{\delta}(\tilde{s})_{\epsilon,sc}=exp(Z_{3}(\bar{s})) \iota^*_{M'_{\epsilon,SC},M'(\tilde{s})_{\epsilon,sc}}(\boldsymbol{\delta}_{\epsilon,SC}).$$
 On a
 $$Sg_{{\bf M}',{\cal O}'}^{{\bf G}'(\tilde{s})}(\boldsymbol{\delta},B^{\tilde{G}})=Sg_{\tilde{M}'_{1}(\tilde{s}),{\cal O}'}^{\tilde{G}'_{1}(\tilde{s})}(\boldsymbol{\delta}_{1}(\tilde{s}),B^{\tilde{G}}).$$
 Nos hypoth\`eses de r\'ecurrence autorisent \`a utiliser la proposition 4.4. Le terme ci-dessus est nul si $A_{G'_{1}(\tilde{s})}\not=A_{G'_{1}(\tilde{s})_{\epsilon_{1}(\tilde{s})}}$. Cette condition \'equivaut \`a $A_{G'(\tilde{s})}\not=A_{G'(\tilde{s})_{\epsilon}}$. Supposons que $A_{G'(\tilde{s})}=A_{G'(\tilde{s})_{\epsilon}}$. Alors le terme pr\'ec\'edent vaut
 $$e_{\tilde{M}'_{1}(\tilde{s})}^{\tilde{G}'_{1}(\tilde{s})}(\epsilon_{1}(\tilde{s}))d(\tilde{s})desc_{\epsilon}^{st,\tilde{G}'_{1}(\tilde{s}),*}( Sg_{M'_1(\tilde{s})_{\epsilon_{1}(\tilde{s})},unip}^{\tilde{G}'_{1}(\tilde{s})_{\epsilon_{1}(\tilde{s})}}(exp(Z_{1}+Z_{2}(\bar{s})) \iota^*_{M'(\tilde{s})_{\epsilon,sc},M'_{1}(\tilde{s})_{\epsilon_{1}(\tilde{s})}}(\boldsymbol{\delta}(\tilde{s})_{\epsilon,sc}),B^{\tilde{G}}_{{\cal O}'})).$$
 Comme en 7.1(4), on peut simplifier $e_{\tilde{M}'_{1}(\tilde{s})}^{\tilde{G}'_{1}(\tilde{s})}(\epsilon_{1}(\tilde{s}))$ en $e_{\tilde{M}'}^{\tilde{G}'(\tilde{s})}(\epsilon)$. 
 On applique la proposition 3.7: on a
  $$Sg_{M'_1(\tilde{s})_{\epsilon_{1}(\tilde{s})},unip}^{\tilde{G}'_{1}(\tilde{s})_{\epsilon_{1}(\tilde{s})}}(exp(Z_{1}+Z_{2}(\tilde{s})) \iota^*_{M'(\tilde{s})_{\epsilon,sc},M'_{1}(\tilde{s})_{\epsilon_{1}(\tilde{s})}}(\boldsymbol{\delta}(\tilde{s})_{\epsilon,sc}),B^{\tilde{G}}_{{\cal O}'})=\iota_{G'(\tilde{s})_{\epsilon,SC},G'_{1}(\tilde{s})_{\epsilon_{1}(\tilde{s})}}^*(\boldsymbol{\tau}(\tilde{s})_{sc} ),$$
  o\`u
  $$\boldsymbol{\tau}(\tilde{s})_{sc}=Sg_{M'(\tilde{s})_{\epsilon,sc},unip}^{G'(\tilde{s})_{\epsilon,SC}}(\boldsymbol{\delta}(\tilde{s})_{\epsilon,sc},B^{\tilde{G}}_{{\cal O}'}).$$
  Avec les notations de 5.4(2), on obtient
 $$Sg_{\tilde{M}'_{1}(\tilde{s}),{\cal O}'}^{\tilde{G}'_{1}(\tilde{s})}(\boldsymbol{\delta}_{1}(\tilde{s}),B^{\tilde{G}})= e_{\tilde{M}'}^{\tilde{G}'(\tilde{s})}(\epsilon)d(\tilde{s})\boldsymbol{\tau}(\tilde{s})^{\tilde{G}'_{1}(\tilde{s})}.$$
 Gr\^ace \`a 5.4(2), on a
 $$transfert(\boldsymbol{\tau}(\tilde{s})^{\tilde{G}'_{1}(\tilde{s})})=\sum_{y\in \dot{{\cal Y}}}c[y]d(\tilde{s},y)\boldsymbol{\tau}[y]^{\tilde{G}}.$$
Reprenons la construction des \'el\'ements $\boldsymbol{\tau}[y]$. Notons $\bar{\boldsymbol{\delta}}_{SC}$ l'image par transfert non standard de $\boldsymbol{\delta}_{\epsilon,SC}$. C'est un \'el\'ement de $D_{g\acute{e}om}^{st}(\bar{M}'_{SC}(F))$.   En utilisant l'analogue de 3.7(4) pour le transfert non standard, on obtient que le transfert non standard de $\boldsymbol{\delta}(\tilde{s})_{\epsilon,sc}$ est 
 $$\bar{\boldsymbol{\delta}}(\bar{s})_{sc}=exp(Z_{3}(\bar{s}))\iota^*_{\bar{M}'_{SC},\bar{M}'(\bar{s})_{sc}}(\bar{\boldsymbol{\delta}}_{SC}).$$ 
 Utilisons l'hypoth\`ese (1). Elle nous dit que le transfert non standard 
 $\bar{\boldsymbol{\tau}}(\bar{s})_{sc}$ de $\boldsymbol{\tau}(\tilde{s})_{sc}$ est $cSg_{\bar{M}'(\bar{s})_{sc},unip}^{\bar{G}'(\bar{s})_{SC}}(\bar{\boldsymbol{\delta}}(\bar{s})_{sc})$, o\`u
 $$c=(c_{\bar{M}'(\bar{s})_{sc},\tilde{M}'(\tilde{s})_{\epsilon,sc}}^{\bar{G}'(\bar{s})_{SC},\tilde{G}'(\tilde{s})_{\epsilon,SC}})^{-1}.$$
 En utilisant le lemme 3.7, l'\'el\'ement $\bar{\boldsymbol{\tau}}(\bar{s})=\iota^*_{\bar{G}'(\bar{s})_{SC},\bar{G}'(\bar{s})}(\bar{\boldsymbol{\tau}}(\bar{s})_{sc})$ est \'egal \`a 
 $cSg_{\bar{M}',unip}^{\bar{G}'(\bar{s})}( \bar{\boldsymbol{\delta}})$,
 o\`u
 $$\bar{\boldsymbol{\delta}}=exp(Z_{2}(\bar{s}))\iota^*_{\bar{M}'(\bar{s})_{sc},\bar{M}'}(\bar{\boldsymbol{\delta}}(\bar{s})_{sc}).$$
 Remarquons que l'on a aussi
 $\bar{\boldsymbol{\delta}}=exp(Z_{2})\iota^*_{\bar{M}'_{SC},\bar{M}'}(\bar{\boldsymbol{\delta}}_{SC})$,
 o\`u 
 $Z_{2}=Z_{2}(\bar{s})+Z_{3}(\bar{s})\in \mathfrak{z}(\bar{M}';F)$.
 Ces termes $Z_{2}$ et $\bar{\boldsymbol{\delta}}$ sont ind\'ependants de $\bar{s}$.  Ensuite
 $$\boldsymbol{\tau}[y]^{\tilde{G}}=desc_{\eta[y]}^{\tilde{G},*}\circ \iota^*_{G_{\eta[y],SC},G_{\eta[y]}}\circ transfert_{y}(\bar{\tau}(\bar{s}))$$
 $$=c\,desc_{\eta[y]}^{\tilde{G},*}\circ \iota^*_{G_{\eta[y],SC},G_{\eta[y]}}\circ transfert_{y}(Sg_{\bar{M}',unip}^{\bar{G}'(\bar{s})}( \bar{\boldsymbol{\delta}})).$$
 On obtient
$$transfert( Sg_{\tilde{M}'_{1}(\tilde{s}),{\cal O}'}^{\tilde{G}'_{1}(\tilde{s})}(\boldsymbol{\delta}_{1}(\tilde{s}),B^{\tilde{G}}))=\sum_{y\in \dot{{\cal Y}}}c[y]d(\tilde{s},y)e_{\tilde{M}'}^{\tilde{G}'(\tilde{s})}(\epsilon) d(\tilde{s})(c_{\bar{M}'(\bar{s})_{sc},\tilde{M}'(\tilde{s})_{\epsilon,sc}}^{\bar{G}'(\bar{s})_{SC},\tilde{G}'(\tilde{s})_{\epsilon,SC}})^{-1}$$
$$desc_{\eta[y]}^{\tilde{G},*}\circ \iota^*_{G_{\eta[y],SC},G_{\eta[y]}}\circ transfert_{y}(Sg_{\bar{M}',unip}^{\bar{G}'(\bar{s})}( \bar{\boldsymbol{\delta}})).$$
On se rappelle que l'on a suppos\'e $A_{G'(\tilde{s})}=A_{G'(\tilde{s})_{\epsilon}}$. Comme en 5.3, notons ${\cal S}$ l'ensemble des $\tilde{s}$ tels que ${\bf G}'(\tilde{s})$ soit elliptique et que cette \'egalit\'e soit v\'erifi\'ee. Pour tout $y\in \dot{{\cal Y}}$ et tout $\bar{s}\in Z(\hat{\bar{M}}_{ad})^{\Gamma_{F}}$, posons
$$x(\bar{s},y)=\sum_{\tilde{s}\in {\cal S}, \tilde{s}\mapsto \bar{s}}i_{\tilde{M}'}(\tilde{G},\tilde{G}'(\tilde{s}))d(\tilde{s},y)e_{\tilde{M}'}^{\tilde{G}'(\tilde{s})}(\epsilon) d(\tilde{s})(c_{\bar{M}'(\bar{s})_{sc},\tilde{M}'(\tilde{s})_{\epsilon,sc}}^{\bar{G}'(\bar{s})_{SC},\tilde{G}'(\tilde{s})_{\epsilon,SC}})^{-1}.$$
Posons
$$(2) \qquad \boldsymbol{\xi}[y]=\sum_{\bar{s}\in Z(\hat{\bar{M}}_{ad})^{\Gamma_{F}}}x(\bar{s},y)Sg_{\bar{M}',unip}^{\bar{G}'(\bar{s})}( \bar{\boldsymbol{\delta}}).$$
Alors les calculs ci-dessus transforment l'expression (1) en
$$(3) \qquad g_{\tilde{M},{\cal O}}^{\tilde{G},{\cal E}}({\bf M}',\boldsymbol{\delta})=\sum_{y\in \dot{{\cal Y}}}c[y] 
desc_{\eta[y]}^{\tilde{G},*}\circ \iota^*_{G_{\eta[y],SC},G_{\eta[y]}}\circ transfert_{y}(\boldsymbol{\xi}[y]).$$
D'apr\`es 5.3(5), l'ensemble ${\cal S}$ est vide si $A_{G_{\eta}}\not=A_{\tilde{G}}$. Cela entra\^{\i}ne

(4) $g_{\tilde{M},{\cal O}}^{\tilde{G},{\cal E}}({\bf M}',\boldsymbol{\delta})=0$ si $A_{G_{\eta}}\not=A_{\tilde{G}}$.

Dans la suite, on suppose $A_{G_{\eta}}=A_{\tilde{G}}$. Alors, d'apr\`es 5.3(6),  ${\cal S}$ est l'image r\'eciproque de l'ensemble des $\bar{s}$ tels que ${\bf G}'(\bar{s})$ soit elliptique. 
Notons ${\cal Z}$ le noyau de l'homomorphisme
$$Z(\hat{M})^{\Gamma_{F},\hat{\theta}}/Z(\hat{G})^{\Gamma_{F},\hat{\theta}}\to Z(\hat{\bar{M}})^{\Gamma_{F}}/Z(\hat{\bar{G}})^{\Gamma_{F}}.$$
 Pour $\tilde{s}\in {\cal S}$,  on a l'\'egalit\'e
$$ e_{\tilde{M}'}^{\tilde{G}'(\tilde{s})}(\epsilon)i_{\tilde{M}'}(\tilde{G},\tilde{G}'(\tilde{s}))=\vert {\cal Z}\vert ^{-1}c_{\bar{M}'(\bar{s})_{sc},\tilde{M}'(\tilde{s})_{\epsilon,sc}}^{\bar{G}'(\bar{s})_{SC},\tilde{G}'(\tilde{s})_{\epsilon,SC}}i_{\bar{M}'}(\bar{G}_{SC},\bar{G}'(\bar{s})).$$
On a vu cette \'egalit\'e dans la preuve de  7.1(14) (o\`u $\vert {\cal Z}\vert $ \'etait not\'e $d$). Gr\^ace \`a elle et \`a la description de ${\cal S}$, on transforme la d\'efinition de $x(\bar{s},y)$ en
$$(5) \qquad x(\bar{s},y)=\vert {\cal Z}\vert ^{-1} i_{\bar{M}'}(\bar{G}_{SC},\bar{G}'(\bar{s})) \sum_{\tilde{s}\in \tilde{\zeta}Z(\hat{M})^{\Gamma_{F},\hat{\theta}}/Z(\hat{G})^{\Gamma_{F},\hat{\theta}}, \tilde{s}\mapsto \bar{s}}d(\tilde{s})d(y,\tilde{s}).$$

  \bigskip
  
  \subsection{Calcul de $x(y,\bar{s})$} 
   Rappelons que l'on note ${\cal Y}^{M}$ l'analogue de ${\cal Y}$ quand on remplace $\tilde{G}$ par $\tilde{M}$, c'est-\`a-dire l'ensemble des $y\in M$ tels que $y\sigma(y)^{-1}\in I^M_{\eta}$, o\`u $I^M_{\eta}=Z(M)^{\theta}M_{\eta}$. On a fix\'e un ensemble $\dot{{\cal Y}}$ de repr\'esentants du quotient $I_{\eta}\backslash {\cal Y}/G(F)$.  On fixe de m\^eme un ensemble $\dot{{\cal Y}}^M$ du quotient $I^M_{\eta}\backslash {\cal Y}^M/M(F)$. Le lemme [I] 5.11 nous autorise \`a supposer que $\dot{{\cal Y}}^M$ est un sous-ensemble de $\dot{{\cal Y}}$.
  
 \ass{Proposition}{Soient $\bar{s}\in \bar{\zeta}Z(\hat{\bar{M}}_{ad})^{\Gamma_{F}}$ et $y\in \dot{{\cal Y}}$. Alors on peut normaliser le facteur de transfert $\Delta(\bar{s},y)$ de sorte que l'on ait  l'\'egalit\'e
 $$x(y,\bar{s})=\left\lbrace\begin{array}{cc}i_{\bar{M}'}(\bar{G}_{SC},\bar{G}'(\bar{s})) &\text{ si }y\in \dot{{\cal Y}}^M,\\ 0,&\text{ sinon.}\\ \end{array}\right.$$}
  
  Preuve. On peut suppposer ${\bf G}'(\bar{s})$ elliptique, sinon les deux membres sont nuls. Supposons d'abord $y\in \dot{{\cal Y}}^M$. On normalise le facteur $\Delta(\tilde{s},y)$ par l'\'egalit\'e 7.1(8). Cela entra\^{\i}ne $d(\tilde{s})d(\tilde{s},y)=1$ d'apr\`es 7.1(9). L'\'egalit\'e de l'\'enonc\'e r\'esulte alors directement de la formule 8.2(5).
  
    Supposons maintenant $y\not\in \dot{{\cal Y}}^M$. Fixons $\tilde{s}\in \tilde{\zeta}Z(\hat{M})^{\Gamma_{F},\hat{\theta}}/Z(\hat{G})^{\Gamma_{F},\hat{\theta}}$ se projetant sur $\bar{s}$.  L'ensemble des \'el\'ements de $\tilde{\zeta}Z(\hat{M})^{\Gamma_{F},\hat{\theta}}/Z(\hat{G})^{\Gamma_{F},\hat{\theta}}$ qui se projettent sur $\bar{s}$ est alors l'ensemble des $z\tilde{s}$ pour $z\in {\cal Z}$. On cherche \`a d\'emontrer l'\'egalit\'e
  $$\sum_{z\in {\cal Z}}d(z\tilde{s})d(z\tilde{s},y)=0.$$
  On peut supposer que $\dot{{\cal Y}}^M$ contient l'\'el\'ement $1$. Comme on vient de le voir, la fonction $z\mapsto d(z\tilde{s})d(z\tilde{s},1)$ est constante de valeur $1$. On peut donc aussi bien d\'emontrer l'\'egalit\'e
 $$(1) \qquad \sum_{z\in {\cal Z}}d(z\tilde{s},y)d(z\tilde{s},1)^{-1}=0.$$ 
 Effectuons les constructions de   5.2 dans un sens diff\'erent. On fixe un sous-tore maximal elliptique $\bar{R}'$ de $\bar{G}'(\bar{s})$. Parce qu'il est elliptique, il se transf\`ere en un tore $R^{\natural}_{sc}$ de $G_{\eta,SC}$ et en un tore $R^{\natural}[y]_{sc}$ de $G_{\eta[y],SC}$. On note $R^{\natural}$, resp. $R^{\natural}[y]$, leurs tores associ\'es dans $G_{\eta}$, resp. $G_{\eta[y]}$,  et $R$, resp. $R[y]$, les commutants de ces tores dans $G$. Pour $z\in {\cal Z}$, le tore $\bar{R}'_{sc}$ se transf\`ere par endoscopie non standard en un tore $R'(z)_{sc}$ de $G'(z\tilde{s})_{\epsilon,SC}$. On note $R'(z)$ le tore associ\'e dans $G'(z\tilde{s})_{\epsilon}$, qui est aussi un sous-tore maximal de $G'(z\tilde{s})$. On fixe $X\in \mathfrak{r}^{\theta}(F)=\mathfrak{r}^{\natural}(F)$ en position g\'en\'erale et proche de $0$, que l'on \'ecrit $X=X_{sc}+Z_{1}$, avec $X_{sc}\in \mathfrak{r}^{\natural}_{sc}(F)$ et $Z_{1}\in \mathfrak{z}(G_{\eta};F)\simeq \mathfrak{z}(\bar{G} ;F)$ (on oublie le temps de cette d\'emonstration les termes $Z_{1}$ etc... de 8.2). On tranf\`ere $X_{sc}$ en un \'el\'ement $\bar{Y}\in \mathfrak{\bar{r}}'(F)$ que l'on \'ecrit $\bar{Y}=\bar{Y}_{sc}+Z_{2}$, avec $\bar{Y}_{sc}\in \mathfrak{\bar{r}}'_{sc}(F)$ et $Z_{2}\in \mathfrak{z}(\bar{G}'(\bar{s});F)$. On transf\`ere $\bar{Y}_{sc}$ en un \'el\'ement $Y_{sc}(z)\in \mathfrak{r}'(z)_{sc}(F)$. Modulo les m\^emes identifications qu'en 5.2, on pose $Y(z)=Y_{sc}(z)+Z_{1}+Z_{2}$. C'est un \'el\'ement de $R'(z)(F)$.  On transf\`ere aussi $\bar{Y}$ en un \'el\'ement $X[y]_{sc}$ de $\mathfrak{r}^{\natural}[y]_{sc}(F)$ et on pose $X[y]=X[y]_{sc}+Z_{1}$, modulo l'isomorphisme $\mathfrak{z}(\bar{G} ;F)\simeq \mathfrak{z}(G_{\eta[y]};F)$. La d\'efinition 5.4(1) donne
$$d(z\tilde{s},y)\Delta(\bar{s},y)(exp(\bar{Y}),exp(X[y]_{sc}))=\Delta_{1}(z\tilde{s})(exp(Y(z))\epsilon_{1}(\tilde{s}),exp(X'[y])\eta[y]),$$
$$d(z\tilde{s},1)\Delta(\bar{s},1)(exp(\bar{Y}),exp(X_{sc}))=\Delta_{1}(z\tilde{s})(exp(Y(z))\epsilon_{1}(\tilde{s}),exp(X)\eta).$$
Donc
$$d(z\tilde{s},y)d(z\tilde{s},1)^{-1}=c\chi(z),$$
o\`u
$$c=\Delta(\bar{s},1)(exp(\bar{Y}),exp(X_{sc}))\Delta(\bar{s},y)(exp(\bar{Y}'),exp(X_{sc}[y]))^{-1},$$
$$\chi(z)=\Delta(z\tilde{s})_{1}(exp(Y(z))\epsilon(\tilde{s})_{1},exp(X[y])\eta[y])\Delta(z\tilde{s})_{1}(exp(Y(z))\epsilon(\tilde{s})_{1},exp(X)\eta)^{-1}$$
$$=\boldsymbol{\Delta}(z\tilde{s})_{1}(exp(Y(z))\epsilon(\tilde{s})_{1},exp(X[y])\eta[y];exp(Y(z))\epsilon(\tilde{s})_{1},exp(X)\eta).$$
L'\'egalit\'e (1) est \'equivalente \`a
$$(2) \qquad \sum_{z\in {\cal Z}}\chi(z)=0.$$
Fixons $z\in {\cal Z}$ et calculons $\chi(z)$. Le tore $R^{\natural}$ est un transfert de $R^{\natural}[y]$ par l'automorphisme int\'erieur $ad_{y}$. Quitte \`a multiplier $y$ \`a gauche par un \'el\'ement de $I_{\eta}$, ce qui ne change rien au probl\`eme, on peut supposer que $ad_{y}(R^{\natural}[y])=R^{\natural}$ et que $ad_{y}$ se restreint en un isomorphisme d\'efini sur $F$ de $R^{\natural}[y]$ sur $R^{\natural}$. Ces propri\'et\'es se prolongent automatiquement: $ad_{y}$ se restreint en un isomorphisme d\'efini sur $F$ de $R[y]$ sur $R$. On calcule $\chi(z)$ en utilisant les formules de [I] 2.2. Les tores $T$ et $\underline{T}$ sont remplac\'es par $R[y]$ et $R$. Du c\^ot\'e dual, on peut identifier les tores $\hat{R}[y]$ et $\hat{R}$. Les constructions sont les m\^emes pour les deux tores. Le cocycle $\hat{V}_{1}$ de [I] 2.2 est donc de la forme $\hat{V}_{1}(w)=(\hat{V}_{\mathfrak{R}_{1}}(w),\hat{V}_{\mathfrak{R}_{1}}(w),1)$ (on a remplac\'e la lettre $\mathfrak{T}$ de [I] 2.2 par $\mathfrak{R}$ par souci de coh\'erence) et l'\'el\'ement de $H^{1,0}(W_{F};\hat{S}_{1}\stackrel{1-\hat{\theta}}{\to}\hat{U})$ est $(\hat{V}_{1}, {\bf zs})$, o\`u ${\bf zs}=(z_{sc}s_{sc},z_{sc}s_{sc})$. Du c\^ot\'e des groupes sur $F$, on doit faire un peu attention. On  peut    identifier les tores $R$ et $R[y]$ par l'automorphisme $ad_{y}$. Mais les cocycles ne s'identifient pas exactement. On a des cocycles $V_{R[y]} $ et $V_{R}$ d\'efinis par des formules
$$V_{R[y]}(\sigma)=r_{R[y]}(\sigma)n_{{\cal E}}(\omega_{R[y]}(\sigma))u_{{\cal E}}(\sigma),$$
$$V_{R}(\sigma)=r_{R}(\sigma)n_{\underline{{\cal E}}}(\omega_{R}(\sigma))u_{\underline{{\cal E}}}(\sigma),$$
en adaptant les notations de [I] 2.2 \`a la pr\'esente situation. Fixons une d\'ecomposition $y=y_{sc}d$, avec $y_{sc}\in G_{SC}$ et $d\in Z(G)$. On a d\'efini $u_{\underline{{\cal E}}}(\sigma)$ par $u_{\underline{{\cal E}}}(\sigma)=y_{sc}u_{{\cal E}}(\sigma)\sigma(y_{sc})^{-1}$. On v\'erifie alors que l'on a l'\'egalit\'e
$$y^{-1}V_{R}(\sigma)y=V_{R[y]}(\sigma)\sigma(y_{sc})^{-1}y_{sc}.$$
De m\^eme, on a pos\'e $exp(X'[y])\eta[y]=\nu e$ et $exp(X')\eta=\underline{\nu}\underline{e}$.  On a $exp(X')\eta=ad_{y}(exp(X'[y])\eta[y])$, mais $\underline{e}=ad_{y_{sc}}(e)$. On en d\'eduit $\underline{\nu}=d\theta(d)^{-1}ad_{y}(\nu)$. En identifiant maintenant les deux tores via $ad_{y}$, on obtient que le cocycle $V$ est de la forme $V(\sigma)=(V_{R[y]}(\sigma),V_{R[y]}(\sigma)^{-1}y_{sc}^{-1}\sigma(y_{sc})) $ et que l'\'el\'ement $\boldsymbol{\nu}_{1}$ est de la forme $(\nu_{1},\nu_{1}^{-1}\theta(d)d^{-1})$. Pour $\sigma\in \Gamma_{F}$, posons $\tau(\sigma)=\sigma(y_{sc})^{-1}y_{sc}$. Alors $\tau$ est un cocycle \`a valeurs dans $R[y]_{sc}$. On v\'erifie que le couple $(\tau, \theta(d)^{-1}d)$ est un cocycle qui d\'efinit un \'el\'ement de $H^{1,0}(\Gamma_{F};R[y]_{sc}\stackrel{1-\theta}{\to}(1-\theta)(R[y]))$. On a un homomorphisme naturel 
 $$j: H^{1,0}(\Gamma_{F};R[y]_{sc}\stackrel{1-\theta}{\to}(1-\theta)(R[y]))\to H^{1,0}(\Gamma_{F};U\stackrel{1-\theta}{\to}S_{1})$$
  (via les secondes composantes, cf. les formules ci-dessus).
On  peut d\'ecomposer le cocycle $(V,\boldsymbol{\nu}_{1})\in Z^{1,0}(\Gamma_{F};U\stackrel{1-\theta}{\to}S_{1})$ en le produit de l'image naturelle de l'inverse du  cocycle pr\'ec\'edent et du cocycle $(V_{0},\boldsymbol{ \nu}_{0})$ d\'efini par $V_{0}(\sigma)=(V_{R[y]}(\sigma),V_{R[y]}(\sigma)^{-1})$ et $\boldsymbol{\nu}_{0}=(\nu_{1},\nu_{1}^{-1})$.  On a alors
$$\chi(z)=<(V,\boldsymbol{\nu}_{1}),(\hat{V}_{1},{\bf zs})>^{-1}=<(V_{0},\boldsymbol{\nu}_{0}),(\hat{V}_{1},{\bf zs})>^{-1}<(\tau,\theta(d)^{-1}d),j^*(\hat{V}_{1},{\bf zs})>,$$
o\`u $j^*$ est l'homomorphisme dual de $j$. On reconna\^{\i}t le premier terme du membre de droite: c'est 
$$\boldsymbol{\Delta}_{1}(z\tilde{s})(exp(Y'(z))\epsilon_{1}(\tilde{s}),exp(X'[y])\eta[y];exp(Y'(z))\epsilon_{1}(\tilde{s}),exp(X'[y])\eta[y]),$$
et on sait que ce facteur vaut $1$. L'\'el\'ement $j^*(\hat{V}_{1},{\bf zs})$ appartient \`a $H^{1,0}(W_{F};\hat{R}[y]/\hat{R}[y]^{\hat{\theta},0}\stackrel{1-\hat{\theta}}{\to}\hat{R}[y]_{ad})$. On voit que c'est le couple $(Y(z),z_{ad}s_{ad})$, o\`u $Y(z)$ est l'image dans $\hat{R}[y]/\hat{R}[y]^{\hat{\theta},0}$ du cocycle $t_{R[y]}$ construit en [I] 2.2 et $z_{ad}$ et $s_{ad}$ sont les images de $z$ et $s$ dans $\hat{G}_{AD}$. On a ajout\'e un $z$ dans la notation $Y(z)$ parce qu'il d\'epend en effet de $z$ et parce que cela va nous \^etre utile. On obtient
$$\chi(z)=<(\tau,\theta(d)^{-1}d),(Y(z),z_{ad}s_{ad})>,$$
d'o\`u
$$(3)\qquad \chi(z)=\chi(1)<(\tau,\theta(d)^{-1}d),( Y(1)^{-1}Y(z),z_{ad})>.$$
Calculons $Y(1)^{-1}Y(z)$. On ajoute des indices $z$ ou $1$ dans les notations pour distinguer les termes relatifs \`a ${\bf G}'(z\tilde{s})$ de ceux relatifs \`a ${\bf G}'(\tilde{s})$. La d\'efinition de [I] 2.2 donne, pour $w\in W_{F}$,
$$t_{R[y],z}(w)=\hat{r}_{R[y]}(w)\hat{n}(\omega_{R[y]}(w))g_{z}(w)^{-1}\hat{n}_{G'(z\tilde{s})}(\omega_{R[y],G'(z\tilde{s})}(w))^{-1}\hat{r}_{R[y],G'(z\tilde{s})}(w)^{-1}.$$
Pour simplifier les notations, on pose
$$\hat{u}_{z}(w)=\hat{r}_{R[y],G'(z\tilde{s})}(w)\hat{n}_{G'(z\tilde{s})}(\omega_{R[y],G'(z\tilde{s})}(w)).$$
 Donc $Y(1)^{-1}(w)Y(z)(w)$ est la projection de
$$\hat{u}_{ 1}(w)g_{1}(w)g_{z}(w)^{-1}\hat{u}_{ z}(w)^{-1}.$$
On se rappelle que $g_{z,w}=(g_{z}(w),w)$ est un \'el\'ement de ${\cal G}'(z\tilde{s})$ tel que $ad_{g_{z,w}}\circ w_{G}$ agisse comme $w_{G'(z\tilde{s})}$ sur $\hat{G}'(z\tilde{s})$. On introduit de m\^eme un \'el\'ement $m_{w}=(m(w),w)\in {\cal M}'$ tel que $ad_{m_{w}}\circ w_{M}$ agisse comme $w_{M'}$ sur $\hat{M}'$. Puisque ${\cal G}'(z\tilde{s})=\hat{G}'(z\tilde{s}){\cal M}'$ par d\'efinition, les \'el\'ements $g_{z,w}$ et $m_{w}$ appartiennent tous deux \`a ${\cal G}'(z\tilde{s})$ et conservent la m\^eme paire de Borel (celle que l'on a fix\'ee pour laquelle $\hat{M}'$ est un Levi standard). Il en r\'esulte que $g_{z}(w)\in \hat{T}^{\hat{\theta},0}m(w)$. Donc $g_{1}(w)g_{z}(w)^{-1}\in \hat{T}^{\hat{\theta},0}$. Les \'el\'ements $\hat{u}_{ z}(w)$ et $\hat{u}_{ 1}(w)$ normalisent ce tore. Puisqu'on projette  dans $\hat{T}/\hat{T}^{\hat{\theta},0}$, on peut aussi bien supprimer le terme $g_{1}(w)g_{z}(w)^{-1}$ et on obtient que  $Y(1)^{-1}(w)Y(z)(w)$ est la projection de
 $ \hat{u}_{ 1}(w) \hat{u}_{ z}(w)^{-1}$.
Les deux \'el\'ements $ \hat{u}_{ z}(w) $ et $\hat{u}_{ 1}(w)$  se rel\`event naturellement dans $\hat{G}_{SC}$. D\'efinissons une cocha\^{\i}ne $\underline{Y}:W_{F}\to \hat{R}[y]_{sc}/\hat{R}[y]_{sc}^{\hat{\theta}}$ ainsi: $\underline{Y}(w)$ est la projection dans $\hat{R}[y]_{sc}/\hat{R}[y]_{sc}^{\hat{\theta}}$ de $ \hat{u}_{ 1}(w) \hat{u}_{z}(w)^{-1}$, vu comme un \'el\'ement de $\hat{R}[y]_{sc}$. L'\'el\'ement $z$ appartient par d\'efinition \`a $Z(\hat{M})^{\Gamma_{F},\hat{\theta}}/Z(\hat{G})^{\Gamma_{F},\hat{\theta}}$. On a d\'ej\`a dit plusieurs fois que ce tore n'est autre que $Z(\hat{M}_{ad})^{\Gamma_{F},\hat{\theta}}$. Il est connexe. On peut donc relever $z_{ad}$ en un \'el\'ement $z_{sc}\in Z(\hat{M}_{sc})^{\Gamma_{F},\hat{\theta},0}$. Ce groupe est un sous-groupe de $\hat{R}[y]_{sc}$. Montrons que

(4) le couple $(\underline{Y},z_{sc})$ appartient \`a $Z^{1,0}(W_{F}; \hat{R}[y]_{sc}/\hat{R}[y]_{sc}^{\hat{\theta}}\stackrel{1-\hat{\theta}}{\to}\hat{R}[y]_{sc})$. 

On note $w\mapsto w_{R}$ l'action galoisienne sur $\hat{R}[y]$ (ou les tores reli\'es tels que $\hat{R}[y]_{sc}$ etc...). On a les \'egalit\'es $w_{R} =ad_{\hat{u}_{z}(w)}\circ w_{G'(z\tilde{s})}=ad_{\hat{u}_{z}(w)}\circ ad_{g_{z}(w)}\circ w_{G}$. Remarquons  que $w_{R}$ n'agit que sur $\hat{R}$ mais le dernier op\'erateur s'etend \`a tout $\hat{G}$.  Pour $w,w'\in W_{F}$, on a
$$\underline{Y}(w)w_{R}(\underline{Y}(w'))=\underline{Y}(w)\hat{u}_{z}(w)g_{z}(w)w_{G}(\underline{Y}(w'))g_{z}(w)^{-1}\hat{u}_{z}(w)^{-1}.$$
C'est la projection dans $\hat{R}[y]_{sc}/\hat{R}[y]_{sc}^{\hat{\theta}}$ de
$$\hat{u}_{1}(w)g_{z}(w)w_{G}(\hat{u}_{1}(w')) w_{G}(\hat{u}_{z}(w'))^{-1}g_{z}(w)^{-1}\hat{u}_{z}(w)^{-1}$$
$$=\hat{u}_{1}(w)g_{z}(w)w_{G}(\hat{u}_{1}(w')) g_{z}(w)^{-1}g_{z}(w)w_{G}(\hat{u}_{z}(w'))^{-1}g_{z}(w)^{-1}\hat{u}_{z}(w)^{-1}$$
$$=\hat{u}_{1}(w)g_{z}(w)w_{G}(\hat{u}_{1}(w')) g_{z}(w)^{-1}w_{G'(z\tilde{s})}(\hat{u}_{z}(w'))^{-1}\hat{u}_{z}(w)^{-1}.$$
 On a vu ci-dessus que $g_{z}(w)g_{1}(w)^{-1}$ appartenait \`a $\hat{T}^{\hat{\theta},0}$. En relevant l'image de cet \'el\'ement dans $\hat{G}_{AD}$ en un \'el\'ement de $\hat{T}_{sc}^{\hat{\theta}}$, on obtient qu'il existe $t\in \hat{T}_{sc}^{\hat{\theta}}$ tel que $ad_{g_{z}(w)}=ad_{t}\circ ad_{g_{1}(w)}$. Donc
$$g_{z}(w)w_{G}(\hat{u}_{1}(w')) g_{z}(w)^{-1}=tg_{1}(w)w_{G}(\hat{u}_{1}(w'))g_{1}(w)^{-1}t^{-1}=tw_{G'(\tilde{s})}(\hat{u}_{1}(w'))t^{-1}.$$
On obtient que $\underline{Y}(w)w_{R}(\underline{Y}(w'))$ est la projection dans $\hat{R}[y]_{sc}/\hat{R}[y]_{sc}^{\hat{\theta}}$ de
$$ \hat{u}_{1}(w)tw_{G'(\tilde{s})}(\hat{u}_{1}(w'))t^{-1}w_{G'(z\tilde{s})}(\hat{u}_{z}(w'))^{-1}\hat{u}_{z}(w)^{-1}.$$
  Comme plus haut, les \'el\'ements $\hat{u}_{z}(w)$ etc... normalisent le tore $\hat{T}_{sc}^{\hat{\theta}}$, donc les \'el\'ements $t$ de la formule ci-dessus disparaissent par projection. Il reste
 $$ \hat{u}_{1}(w)w_{G'(\tilde{s})}(\hat{u}_{1}(w'))w_{G'(z\tilde{s})}(\hat{u}_{z}(w'))^{-1}\hat{u}_{z}(w)^{-1}.$$ 
 Or, d'apr\`es la construction de Langlands et Shelstad, les applications $\hat{u}_{1}$, resp. $\hat{u}_{z}$, sont des cocycles (\`a valeurs dans $\hat{G}'(\tilde{s})$, resp. $\hat{G}'(z\tilde{s})$). Le terme ci-dessus est donc \'egal \`a $\hat{u}_{1}(ww')\hat{u}_{z}(ww')^{-1}$. Sa projection est $\underline{Y}(ww')$. Cela prouve que $\underline{Y}$ est un cocycle.  Pour $w\in W_{F}$, on a l'\'egalit\'e
 $$(1-\hat{\theta})(\underline{Y}(w))=\hat{u}_{ 1}(w) \hat{u}_{ z}(w)^{-1}\hat{\theta}(\hat{u}_{z}(w)\hat{u}_{1}(w)^{-1}).$$
Fixons un rel\`evement $s_{sc}$ de $s_{ad}$ dans $\hat{G}_{SC}$. On peut remplacer $\hat{\theta}$ par $ad_{z_{sc}s_{sc}}\circ\hat{\theta}$ puisque le terme auquel on applique cet op\'erateur commute \`a $z_{sc}s_{sc}$: il appartient \`a $\hat{T}_{sc}$. Parce que $\hat{u}_{z}(w)\in \hat{G}'(z\tilde{s})_{sc}$, ce terme est fixe par $ad_{z_{sc}s_{sc}}\circ \hat{\theta}$. La formule se simplifie en
 $$(1-\hat{\theta})(\underline{Y}(w))=\hat{u}_{ 1}(w)z_{sc}s_{sc}\hat{\theta}(\hat{u}_{1}(w)^{-1})s_{sc}^{-1}z_{sc}^{-1}.$$
 Pour la m\^eme raison que ci-dessus, $\hat{u}_{1}(w)$ est fixe par $ad_{s_{sc}}\circ \hat{\theta}$. On obtient
 $$(1-\hat{\theta})(\underline{Y}(w))=\hat{u}_{ 1}(w)z_{sc}\hat{u}_{1}(w)^{-1}z_{sc}^{-1}.$$
 L'\'el\'ement $z_{sc}$ appartient \`a $Z(\hat{M}_{sc})^{\Gamma_{F},\hat{\theta},0}$. A fortiori, $w_{G}(z_{sc})=z_{sc}$. On a vu ci-dessus que $g_{1}(w)\in \hat{T}^{\hat{\theta},0}m(w)$, donc $g_{1}(w)\in \hat{M}$. Puisque $z_{sc}\in Z(\hat{M}_{sc})$, ces deux \'el\'ements commutent, d'o\`u $ad_{g_{1}(w)}w_{G}(z_{sc})=z_{sc}$, c'est-\`a-dire $w_{G'(\tilde{s})}(z_{sc})=z_{sc}$. Mais alors
 $$ad_{\hat{u}_{1}(w)}(z_{sc})=w_{R}(z_{sc}).$$
 On obtient
 $$(1-\hat{\theta})(\underline{Y}(w))=w_{R}(z_{sc})z_{sc}^{-1},$$
ce qui prouve que $(\underline{Y},z_{sc}^{-1})$ est un cocycle. Cela d\'emontre (4).

Il y a un homomorphisme naturel
$$H^{1,0}(W_{F};\hat{R}[y]_{sc}/\hat{R}[y]_{sc}^{\hat{\theta}}\stackrel{1-\hat{\theta}}{\to}\hat{R}[y]_{sc})\to H^{1,0}(W_{F};\hat{R}[y]/\hat{R}[y]^{\hat{\theta},0}\stackrel{1-\hat{\theta}}{\to}\hat{R}[y]_{ad}).$$
Le cocycle  $(Y(1)^{-1}Y(z),z_{ad})$ est l'image par cet homomorphisme de $(\underline{Y},z_{sc})$.  Notons $\tau_{ad}$ l'image de $\tau$ dans $R[y]_{ad}$. Puisque $d\in Z(G)$, l'image de  $(\tau,\theta(d)^{-1}d)$ dans $H^{1,0}(\Gamma_{F};R[y]_{ad}\stackrel{1-\theta}{\to}(1-\theta)(R[y]_{ad}))$ par l'homomorphisme dual du pr\'ec\'edent est $(\tau_{ad},1)$. Gr\^ace \`a (3), on obtient
$$\chi(z)=\chi(1)<(\tau_{ad},1),(\underline{Y},z_{sc})>.$$
Rappelons que $R[y]_{ad}^{\theta}$ est connexe. Le fait que $(\tau_{ad},1)$ soit un cocycle implique que $\tau_{ad}$ prend ses valeurs dans ce tore. On a un homomorphisme naturel
$$H^{1,0}(\Gamma_{F};R[y]_{ad}^{\theta}\to \{1\})\to H^{1,0}(\Gamma_{F};R[y]_{ad}\stackrel{1-\theta}{\to}(1-\theta)(R[y]_{ad}))$$
et $(\tau_{ad},1)$ est l'image du m\^eme cocycle, vu comme un \'el\'ement du premier groupe. L'image de $(\underline{Y},z_{sc}^{-1})$ dans $H^{1,0}(W_{F};\{1\}\to \hat{R}[y]_{sc}/(1-\hat{\theta})(\hat{R}[y]_{sc}))$ par l'homomorphisme dual du pr\'ec\'edent est $(1,\bar{z}_{sc})$, o\`u $\bar{z}_{sc}$ est l'image de $z_{sc}$ dans $\hat{R}[y]_{sc}/(1-\hat{\theta})(\hat{R}[y]_{sc})$. D'o\`u
$$\chi(z)=\chi(1)<(\tau_{ad},1),(1,\bar{z}_{sc})>,$$
le produit \'etant celui sur
$$H^{1,0}(\Gamma_{F};R[y]_{ad}^{\theta}\to \{1\})\times H^{1,0}(W_{F};\{1\}\to \hat{R}[y]_{sc}/(1-\hat{\theta})(\hat{R}[y]_{sc})).$$
En appliquant [KS] A.3.14, cela se simplifie en
$$\chi(z)=\chi(1)<\tau_{ad},\bar{z}_{sc}>,$$
le produit \'etant celui sur
$$H^1(\Gamma_{F};R[y]_{ad}^{\theta})\times H^0(W_{F};\hat{R}[y]_{sc}/(1-\hat{\theta})(\hat{R}[y]_{sc})).$$
Transf\'erons  le tore $R[y]^{\theta,0}$ de $G_{\eta[y]}$ en le  tore $R^{\theta,0}$ de $G_{\eta}$. Cela remplace $\tau_{ad}$ par $\tau'_{ad}$ d\'efini par $\tau'_{ad}(\sigma)=ad_{y}(\tau_{ad}(\sigma))$ pour tout $\sigma\in \Gamma_{F}$.  On note $\hat{\bar{R}}$ le tore dual de $R^{\theta,0}$. On peut consid\'erer que c'est un sous-tore de $\hat{\bar{G}}$. On a une suite
$$G_{\eta}\to G_{\eta,ad}\to G_{\eta,AD},$$
o\`u $G_{\eta,ad}$ est l'image de $G_{\eta}$ dans $G_{AD}$. On a une suite similaire pour les formes quasi-d\'eploy\'ees
$$\bar{G}\to \bar{G}_{ad}\to \bar{G}_{AD}$$
et une suite duale
$$\hat{\bar{G}}_{SC}\to \hat{\bar{G}}_{sc}\to \hat{\bar{G}}.$$
  En notant $\hat{\bar{R}}_{sc}$ l'image r\'eciproque de $\hat{\bar{R}}$ dans $\hat{\bar{G}}_{sc}$, le produit ci-dessus devient celui sur
$$H^1(\Gamma_{F};R_{ad}^{\theta})\times H^0(W_{F};\hat{\bar{R}}_{sc}).$$
Puisque $z\in {\cal Z}$, l'image de $z$ dans $\hat{\bar{R}}$ appartient \`a $Z(\hat{\bar{G}})$. Il en r\'esulte que $\bar{z}_{sc}\in Z(\hat{\bar{G}}_{sc})$. On a une dualit\'e sur
$$H^1(\Gamma_{F};G_{\eta,ad})\times Z(\hat{\bar{G}}_{sc})^{\Gamma_{F}}/Z(\hat{\bar{G}}_{sc})^{\Gamma_{F},0}$$
qui est compatible au produit pr\'ec\'edent. Pour simplifier, nous ne changerons pas les notations: on a encore
$$(5) \qquad \chi(z)=\chi(1)<\tau'_{ad},\bar{z}_{sc}>,$$
o\`u cette fois, $\tau'_{ad}$ est vu comme un \'el\'ement du premier groupe ci-dessus et $\bar{z}_{sc}$ comme un \'el\'ement du second. On a choisi le rel\`evement $z_{sc}$ mais la formule obtenue montre que le membre de droite ci-dessus ne d\'epend pas de ce choix. Pour deux \'el\'ements $z,z'\in {\cal Z}$, on peut choisir $z_{sc}z'_{sc}$ comme rel\`evement de $zz'$. On voit alors que,
\`a la constante $\chi(1)$ pr\`es, $\chi(z)$ est  la valeur en $z$ d'un caract\`ere de ${\cal Z}$. Pour obtenir la relation de nullit\'e (2), il reste \`a prouver que ce caract\`ere n'est pas trivial. On a un homomorphisme
$$(6) \qquad Z(\hat{\bar{G}}_{sc})^{\Gamma_{F}}/Z(\hat{\bar{G}}_{sc})^{\Gamma_{F},0}\to Z(\hat{\bar{M}}_{sc})^{\Gamma_{F}}/Z(\hat{\bar{M}}_{sc})^{\Gamma_{F},0},$$
o\`u $\hat{\bar{M}}_{sc}$ est ici l'image r\'eciproque de $\hat{\bar{M}}$ dans $\hat{\bar{G}}_{sc}$.
Montrons que

(7) tout \'el\'ement du noyau de (6) est de la forme $\bar{z}_{sc}$ pour un choix convenable de $z\in {\cal Z}$ et de rel\`evement $z_{sc}\in Z(\hat{M}_{sc})^{\Gamma_{F},\hat{\theta},0}$. 

Soit $x\in Z(\hat{\bar{G}}_{sc})^{\Gamma_{F}}$ relevant un \'el\'ement du noyau. Alors $x\in Z(\hat{\bar{M}}_{sc})^{\Gamma_{F},0}$. On a la m\^eme relation que 5.3(2) au niveau du groupe $\hat{G}_{SC}$, c'est-\`a-dire que l'homomorphisme
$$Z(\hat{M}_{sc})^{\Gamma_{F},\hat{\theta},0}\to Z(\hat{\bar{M}}_{sc})^{\Gamma_{F},0}$$
est surjectif. On peut donc relever $x$ en un \'el\'ement $z_{sc}\in Z(\hat{M}_{sc})^{\Gamma_{F},\hat{\theta},0}$. Notons $z$ son image dans $Z(\hat{M})^{\Gamma_{F},\hat{\theta}}$. L'image naturelle de $z$ dans $Z(\hat{\bar{M}})^{\Gamma_{F}}$ est \'egale \`a celle de $x$. Or $x$ appartient \`a $Z(\hat{\bar{G}}_{sc})^{\Gamma_{F}}$. Donc cette image appartient \`a $Z(\hat{\bar{G}})^{\Gamma_{F}}$. Par d\'efinition de ${\cal Z}$, cela entra\^{\i}ne que $z\in {\cal Z}$. Bien s\^ur, l'image de $x$ dans $Z(\hat{\bar{G}}_{sc})^{\Gamma_{F}}/Z(\hat{\bar{G}}_{sc})^{\Gamma_{F},0}$ est $\bar{z}_{sc}$. Cela d\'emontre (7). 

D'apr\`es (5) et (7), $\chi$ ne peut \^etre constant que si le caract\`ere de $Z(\hat{\bar{G}}_{sc})^{\Gamma_{F}}/Z(\hat{\bar{G}}_{sc})^{\Gamma_{F},0}$ d\'efini par $\tau'_{ad}$ annule le noyau de (6).  On a un diagramme
$$\begin{array}{ccc}H^1(\Gamma_{F};G_{\eta,ad})&\times& Z(\hat{\bar{G}}_{sc})^{\Gamma_{F}}/Z(\hat{\bar{G}}_{sc})^{\Gamma_{F},0}\\ \uparrow&&\downarrow\\ H^1(\Gamma_{F};M_{\eta,ad})&\times &Z(\hat{\bar{M}}_{sc})^{\Gamma_{F}}/Z(\hat{\bar{M}}_{sc})^{\Gamma_{F},0}\\ \end{array}$$
qui est compatible aux dualit\'es. Alors $\tau_{ad}$ annule le noyau de (6) si et seulement s'il provient d'un \'el\'ement de $H^1(\Gamma_{F};M_{\eta,ad})$. Rappelons qu'au d\'epart, on avait $\tau(\sigma)=\sigma(y_{sc})^{-1}y_{sc}$ pour tout $\sigma\in \Gamma_{F}$. Puisque $y_{sc}\in yZ(G)$, on a simplement $\tau_{ad}(\sigma)=\sigma(y_{ad})^{-1}y_{ad}$. D'o\`u $\tau'_{ad}(\sigma)=y_{ad}\sigma(y_{ad})^{-1}$. Autrement dit, c'est le cocycle provenant naturellement du torseur int\'erieur $ad_{y^{-1}}:G_{\eta}\to G_{\eta[y]}$. Si ce cocycle provient d'un \'el\'ement de $H^1(\Gamma_{F};M_{\eta,ad})$, alors le groupe de Levi $M_{\eta}$ de $G_{\eta}$ se transf\`ere \`a $G_{\eta[y]}$. C'est interdit par notre hypoth\`ese $y\not\in \dot{{\cal Y}}^M$ et le lemme [I] 5.11. Cela ach\`eve la preuve de (2) et de la proposition. $\square$

\bigskip

\subsection{Fin de la preuve de la proposition 8.1}
On normalise les facteurs $\Delta(\bar{s},y)$ de sorte que la proposition pr\'ec\'edente soit v\'erifi\'ee. En utilisant cette proposition et la formule 8.2(3), on a $\boldsymbol{\xi}[y]=0$ si $y\not\in \dot{{\cal Y}}^M$. Supposons $y\in \dot{{\cal Y}}^M$. Alors
$$\boldsymbol{\xi}[y]=\sum_{\bar{s}\in Z(\hat{\bar{M}}_{ad}^{\Gamma_{F}}} i_{\bar{M}'}(\bar{G}_{SC},\bar{G}'(\bar{s}))Sg_{\bar{M}',unip}^{\bar{G}'(\bar{s})}(\bar{\boldsymbol{\delta}}).$$
Par d\'efinition
$$transfert_{y}(\boldsymbol{\xi}[y])=g_{M_{\eta[y],sc},unip}^{G_{\eta[y],SC},{\cal E}}({\bf \bar{M}}',\bar{\boldsymbol{\delta}}).$$
ici, les groupes ne sont pas tordus. On peut utiliser le corollaire 1.5: le terme ci-dessus vaut $g_{M_{\eta[y],sc},unip}^{G_{\eta[y],SC}}( \boldsymbol{\delta}[y]_{sc})$, o\`u
$$\boldsymbol{\delta}[y]_{sc}=transfert_{y}(\bar{\boldsymbol{\delta}}).$$
 En utilisant la proposition 3.4, on a
$$\iota^*_{G_{\eta[y],SC},G_{\eta[y]}}(g_{M_{\eta[y],sc},unip}^{G_{\eta[y],SC}}( \boldsymbol{\delta}[y]_{sc}))=g_{M_{\eta[y]},unip}^{G_{\eta[y]}}(\boldsymbol{\delta}[y]),$$
o\`u
$$\boldsymbol{\delta}[y]=exp(Z_{1})\iota^*_{M_{\eta[y],sc},M_{\eta[y]}}(\boldsymbol{\delta}[y]_{sc}).$$
On a identifi\'e ici $Z_{1}$ \`a un \'el\'ement de $\mathfrak{z}(G_{\eta[y]};F)$. Enfin, on se rappelle que l'on a suppos\'e $A_{G_{\eta}}=A_{\tilde{G}}$, ce qui \'equivaut \`a $A_{G_{\eta[y]}}=A_{\tilde{G}}$ puisque $G_{\eta[y]}$ est une forme int\'erieure de $G_{\eta}$. D'apr\`es la proposition 4.2, on a donc
$$desc_{\eta[y]}^{\tilde{G},*}(g_{M_{\eta[y]},unip}^{G_{\eta[y]}}(\boldsymbol{\delta}[y]))=g_{\tilde{M},{\cal O}}^{\tilde{G}}(desc_{\eta[y]}^{\tilde{M},*}(\boldsymbol{\delta}[y])).$$
 La formule 8.1(3) devient
$$(1) \qquad g_{\tilde{M},{\cal O}}^{\tilde{G},{\cal E}}({\bf M}',\boldsymbol{\delta})=g_{\tilde{M},{\cal O}}^{\tilde{G}}(\boldsymbol{\tau}),$$
o\`u
$$\boldsymbol{\tau}=
\sum_{y\in \dot{{\cal Y}}^M}c[y]desc_{\eta[y]}^{\tilde{M},*}(\boldsymbol{\delta}[y]).$$
A ce point, on peut lever l'hypoth\`ese que $A_{G_{\eta}}=A_{\tilde{G}}$. Si elle n'est pas v\'erifi\'ee, le membre de gauche de  (1) est nul d'apr\`es 8.2(4). Celui de droite l'est d'apr\`es la proposition 4.2. 

Il est facile de reprendre tous ces calculs en rempla\c{c}ant l'espace $\tilde{G}$ par $\tilde{M}$ pour calculer $transfert(\boldsymbol{\delta})$. C'est d'ailleurs le m\^eme calcul qu'en 7.1, aux translations pr\`es par l'\'el\'ement central $Z$ qui est \'evidemmment inoffensif. On n'obtient pas tout-\`a-fait l'\'egalit\'e $transfert(\boldsymbol{\delta})=\boldsymbol{\tau}$. Le remplacement de $\tilde{G}$ par $\tilde{M}$ conduit \`a l'\'egalit\'e
$$transfert(\boldsymbol{\delta})=\sum_{y\in \dot{{\cal Y}}^M}c^M[y]desc_{\eta[y]}^{\tilde{M},*}(\boldsymbol{\delta}[y]).$$
Mais on a

(2) $c^M[y]=c[y]$ pour tout $y\in \dot{{\cal Y}}^M$. 

Pour simplifier la notation, on peut supposer $y=1$. On doit prouver que l'application naturelle
$$M_{\eta}(F)\backslash I^M_{\eta}(F)\to G_{\eta}(F)\backslash I_{\eta}(F)$$
est bijective. L'injectivit\'e r\'esulte de l'\'egalit\'e $G_{\eta}\cap M=M_{\eta}$ (un \'el\'ement de $G_{\eta}\cap M$ commute \`a $A_{\tilde{M}}=A_{M_{\eta}}$ donc appartient \`a $M_{\eta}$). Soit $u\in I_{\eta}(F)$. On \'ecrit $u=gz$ avec $g\in G_{\eta}$ et $z\in Z(G)^{\theta}$. On d\'efinit un cocycle $\xi$ sur $\Gamma_{F}$ par $\xi(\sigma)=z^{-1}\sigma(z)=g\sigma(g)^{-1}$. Cette formule montre qu'il  prend ses valeurs dans $Z(G_{\eta})$ et qu'il est cohomologiquement trivial dans $G_{\eta}$. C'est donc un \'el\'ement du noyau de $H^1(\Gamma_{F};Z(G_{\eta}))\to H^1(\Gamma_{F};G_{\eta})$. Cette application se factorise en 
$$H^1(\Gamma_{F};Z(G_{\eta}))\to H^1(\Gamma_{F};M_{\eta})\to H^1(\Gamma_{F};G_{\eta}).$$
La deuxi\`eme application est injective car $M_{\eta}$ est un Levi de $G_{\eta}$. Donc $\xi$ appartient au noyau de la premi\`ere application. On peut donc trouver $m\in M_{\eta}$ tel que $\xi(\sigma)=m\sigma(m)^{-1}$ pour tout $\sigma\in \Gamma_{F}$. Posons $v=mz$. Alors $v\in I^M_{\eta[y]}$. On a $u=gm^{-1}v$ et les relations $\xi(\sigma)=m\sigma(m)^{-1}=g\sigma(g)^{-1}$ entra\^{\i}nent que $gm^{-1}$ appartient \`a $G_{\eta}(F)$. Donc l'image de $u$ dans $ G_{\eta}(F)\backslash I_{\eta}(F)$ est \'egale \`a celle de $v$, ce qui d\'emontre la surjectivit\'e cherch\'ee et (2).  

Donc $transfert(\boldsymbol{\delta})=\boldsymbol{\tau}$ et la formule (1) devient  
$$g_{\tilde{M},{\cal O}}^{\tilde{G},{\cal E}}({\bf M}',\boldsymbol{\delta})=g_{\tilde{M},{\cal O}}^{\tilde{G}}(transfert(\boldsymbol{\delta})).$$
Cela prouve la proposition 8.1. $\square$

 \bigskip
 
 \subsection{Egalit\'e de germes et de germes endoscopiques}
 On consid\`ere un triplet quelconque $(G,\tilde{G},{\bf a})$, un espace de Levi $\tilde{M}$ et une classe de conjugaison stable semi-simple ${\cal O}\subset \tilde{M}(F)$. 
 
 Il y un cas particulier que nous devons exclure. C'est  celui o\`u $(G,\tilde{G},{\bf a})$ est l'un des triplets d\'efinis en 6.3 et o\`u ${\cal O}$ est la classe de conjugaison stable d'un \'el\'ement $\eta\in \tilde{M}(F)$ qui conserve une paire de Borel \'epingl\'ee de $G$ d\'efinie sur $F$.  
  
 \ass{Proposition}{On suppose que l'on n'est pas dans le cas particulier ci-dessus. Soit $\boldsymbol{\gamma}\in D_{g\acute{e}om,\tilde{G}-\acute{e}qui}(\tilde{M}(F),\omega)\otimes Mes(M(F))^*$. On suppose que les \'el\'ements du support de $\boldsymbol{\gamma}$ sont $\tilde{G}$-\'equisinguliers et proches de ${\cal O}$. Alors on a l'\'egalit\'e
 $$g_{\tilde{M},{\cal O}}^{\tilde{G}}(\boldsymbol{\gamma})=g_{\tilde{M},{\cal O}}^{\tilde{G},{\cal E}}(\boldsymbol{\gamma}).$$}
 
 Preuve. Par lin\'earit\'e, on peut fixer une donn\'ee endoscopique ${\bf M}'=(M',{\cal M}',\tilde{\zeta})$ de $(M,\tilde{M},{\bf a})$, elliptique et relevante, et un \'el\'ement $\boldsymbol{\delta}\in D^{st}_{g\acute{e}om}({\bf M}')\otimes Mes(M'(F))^*$ de sorte que $\boldsymbol{\gamma}=transfert(\boldsymbol{\delta})$. On peut aussi fixer une classe de conjugaison stable semi-simple dans $\tilde{M}'(F)$ se transf\'erant sur ${\cal O}$ et supposer $\boldsymbol{\delta}$ proche de ${\cal O}'$. Soit $\epsilon\in {\cal O}'$. Supposons d'abord $A_{M'_{\epsilon}}\not=A_{M'}$. Alors, comme dans la preuve de 7.4, les deux membres de l'\'enonc\'e sont calcul\'es par des formules de descente, \`a savoir celles des propositions [II] 2.11 et [II] 2.12. Par r\'ecurrence, on en d\'eduit l'\'egalit\'e de l'\'enonc\'e. Supposons $A_{M'_{\epsilon}}=A_{M'}$. D'apr\`es nos hypoth\`eses de r\'ecurrence et le lemme 6.3, l'hypoth\`ese (1) de 8.1 est v\'erifi\'ee sauf dans le cas particulier que l'on a exclu. On peut donc appliquer la proposition 8.1 qui conclut. $\square$
 
 \bigskip
 
 \subsection{Preuve de la proposition 4.4}
 On renvoie \`a 4.4 pour l'\'enonc\'e de cette proposition. Ici, le triplet $(G,\tilde{G},{\bf a})$ est quasi-d\'eploy\'e et \`a torsion int\'erieure, muni d'un syst\`eme de fonctions $B$. La preuve de la proposition est similaire \`a celle de la proposition 7.3. On reprend la preuve de 8.1 dans le cas   ${\bf M}'={\bf M}$ que l'on avait exclu. Elle conduit \`a une \'egalit\'e 
 $$g_{\tilde{M},{\cal O}}^{\tilde{G},{\cal E}}({\bf M},\boldsymbol{\delta},B)=x+g_{\tilde{M},{\cal O}}^{\tilde{G}}(transfert(\boldsymbol{\delta}),B),$$
 o\`u $x$ est la diff\'erence entre les deux membres de l'\'egalit\'e que l'on veut prouver. Mais 
 $g_{\tilde{M},{\cal O}}^{\tilde{G},{\cal E}}({\bf M},\boldsymbol{\delta},B)$ est ici tautologiquement \'egal \`a $g_{\tilde{M},{\cal O}}^{\tilde{G}}(transfert(\boldsymbol{\delta}),B)$, ce qui entra\^{\i}ne $x=0$. 
 
 \bigskip
 
 \subsection{Preuve de la proposition 6.7}
 L'argument est le m\^eme qu'en 7.7. Gr\^ace au lemme 6.8, on peut supposer $(G_{1},G_{2},j_{*})$ quasi-\'el\'ementaire. On introduit le triplet $(G,\tilde{G},{\bf a})$ qui lui est associ\'e comme en 6.3. On introduit les m\^emes donn\'ees $\eta$, $\tilde{M}$, ${\bf M}'$ et $\epsilon$ qu'en 7.7. 
 On reprend la preuve de 8.1 pour ces donn\'ees. On obtient une \'egalit\'e
 $$g_{\tilde{M},{\cal O}}^{\tilde{G},{\cal E}}({\bf M}',\boldsymbol{\delta})=\boldsymbol{\mu}+g_{\tilde{M},{\cal O}}^{\tilde{G}}(transfert(\boldsymbol{\delta})),$$
 o\`u $\boldsymbol{\mu}$ est un cetain terme compl\'ementaire. Par nos hypoth\`eses de r\'ecurrence, on conna\^{\i}t l'\'egalit\'e
 $$g_{\tilde{M},{\cal O}}^{\tilde{G},{\cal E}}({\bf M}',\boldsymbol{\delta})= g_{\tilde{M},{\cal O}}^{\tilde{G}}(transfert(\boldsymbol{\delta})).$$
 D'o\`u $\boldsymbol{\mu}=0$. En choisissant convenablement $\boldsymbol{\delta}$, on en d\'eduit l'assertion cherch\'ee concernant notre triplet $(G_{1},G_{2},j_{*})$. On laisse les d\'etails au lecteur.

\bigskip

{\bf Bibliographie}

[A1] J. Arthur: {\it A stable trace formula I. General expansions} Journal of the Inst. of Math. Jussieu 1 (2002), p. 175-277

[A2] -----------: {\it The local behaviour of weighted orbital integrals}, Duke Math. Journal 56 (1988), p. 223-293

[F] A. Ferrari: {\it Th\'eor\`eme de l'indice et formule des traces}, manuscripta math. 124 (2007), p. 363-390

[KS] R. Kottwitz, D. Shelstad: {\it Foundations of twisted endoscopy}, Ast\'erisque 255 (1999)

[L] J.-P. Labesse: {\it Stable twisted trace formula: elliptic terms}, Journal of the Inst. of Math. Jussieu 3 (2004), p. 473-530

[S] J.-P. Serre: {\it Corps locaux}, Hermann 1968

 [W1] J.-L. Waldspurger: {\it L'endoscopie tordue n'est pas si tordue}, Memoirs AMS 908 (2008)

[W2] -----------------------: {\it La formule des traces locale tordue}, pr\'epublication 2012

[I] ---------------------:{\it  Stabilisation de la formule des traces tordue I: endoscopie tordue sur un corps local}, pr\'epublication 2014

[II] --------------------: {\it  Stabilisation de la formule des traces tordue II: int\'egrales orbitales et endoscopie sur un corps local non-archim\'edien; d\'efinitions et \'enonc\'es des r\'esultats}, pr\'epublication 2014

\bigskip

Institut de Math\'ematiques de Jussieu, CNRS

2 place Jussieu 75005 Paris

e-mail: waldspur@math.jussieu.fr

\bigskip
\end{document}